\documentclass[11pt,a4paper]{article}
\usepackage[a4paper,margin=1in]{geometry}
\usepackage[T1]{fontenc}
\usepackage[utf8]{inputenc}
\usepackage{lmodern,microtype}
\usepackage{amsmath,amssymb,amsthm,mathtools}
\usepackage{tikz-cd}
\usepackage{booktabs,longtable,array,tabularx}
\usepackage{enumitem}
\usepackage[hidelinks]{hyperref}
\usepackage{bookmark}
\usepackage[nameinlink,noabbrev]{cleveref}
\setlist{itemsep=3pt,topsep=4pt}
\newtheorem{theorem}{Theorem}[section]
\newtheorem{proposition}[theorem]{Proposition}
\newtheorem{lemma}[theorem]{Lemma}
\newtheorem{corollary}[theorem]{Corollary}
\theoremstyle{definition}
\newtheorem{definition}[theorem]{Definition}
\newtheorem{example}[theorem]{Example}
\theoremstyle{remark}
\newtheorem{remark}[theorem]{Remark}
\newcommand{\Set}{\mathbf{Set}}
\newcommand{\TS}{\mathsf{TS}}
\newcommand{\SX}{\mathsf{SX}}
\newcommand{\id}{\mathrm{id}}
\newcommand{\Hom}{\operatorname{Hom}}
\newcommand{\base}{\operatorname{base}}
\newcommand{\heap}{\operatorname{heap}}
\newcommand{\Seq}{\mathsf{Seq}}
\newcommand{\Cau}{\mathsf{Cau}}
\newcommand{\Rep}{\mathsf{Rep}}
\newcommand{\Tok}{\mathsf{Tok}}

\newcommand{\ATC}{\mathbf{ATC}}
\newcommand{\Tkn}{\mathbf{Tkn}}
\newcommand{\TTS}{\mathsf{TTS}}
\newcommand{\tmorph}{\mathbf{tmorph}}
\newcommand{\corr}{\mathbf{corr}}
\newcommand{\ary}{\mathbf{ary}}
\newcommand{\tary}{\mathbf{tary}}
\newcommand{\universal}{\mathbf{universal}}
\newcommand{\Tr}{\mathcal{T}}
\newcommand{\Gt}{\mathcal{G}}
\newcommand{\Softmax}{\operatorname{softmax}}
\newcommand{\Attn}{\operatorname{Attn}}
\newcommand{\MHA}{\operatorname{MHA}}
\newcommand{\FFN}{\operatorname{FFN}}
\newcommand{\RMSNorm}{\operatorname{RMSNorm}}

\DeclareMathOperator{\im}{im}
\DeclareMathOperator{\supp}{supp}
\title{Token Space: A Category Theory Framework\\for AI Computations}
\author{Wuming Pan\\College of Computer Science, Sichuan University\\Chengdu, P.R. China, 610065\\\texttt{panwuming@scu.edu.cn}}
\date{}
\hypersetup{pdftitle={Token Space: A Category Theory Framework for AI Computations},pdfauthor={Wuming Pan},pdfkeywords={Token Space, algebraic tokenization, category theory, Token-level structure, finite mappings, execution gluing, frontier functor, concurrency, elastic computation, unbounded computing cores, Transformers, structural distillation}}
\begin{document}
\maketitle

\begin{abstract}
Token Space is a categorical framework and mathematical language for AI computations. It connects internal relations, program descriptions, execution states and observations, locating questions about structure, behavior and cost at their appropriate levels. Five guiding positions concern structural interiors, categorical self-description, interfaces, occurrence identity and extensible computation. Tokens are finite records of carrier elements and fixed symbols; Token maps preserve selected heaps. The elementary category is a quasitopos, hence locally cartesian closed, but not a topos. Algebraic tokenization is fully faithful for fixed finitary signatures. Represented finite mappings admit concurrent graph execution, gluing, functorial frontiers and exact state migration characterized by kernel inclusion. Transformers are one implementation family. Coherent occurrence prefixes yield natural numerical maps, and a cache invariant proves agreement with full-prefix evaluation. A prescribed access policy determines the least retained index set under a no-reconstruction discipline. Parameterised state expresses changing interfaces. For a specified future-observation heap, structural indiscernibility is behavioral equivalence. An encoding admits exact incremental execution precisely when its kernel is a right congruence contained in that equivalence; every reachable exact realization maps uniquely onto the behavioral quotient. Teacher-induced heaps and declared readouts connect these constructions to distillation and explicit knowledge. The definitions, examples and theorems demonstrate how the language joins lines of reasoning while distinguishing representation, execution and observation. Effective implementations and quantitative performance remain further questions.
\end{abstract}

\textbf{Keywords:} Token Space; Algebraic Tokenization; Category Theory; Token-Level Structure; Finite Mappings; Concurrent Computation; Computational Elasticity; Unbounded Computing Cores; Transformer Models; Structural Distillation.

\section{Introduction}\label{introduction}

During the last sixty years category theory has become one of the essential building blocks of modern mathematics \cite{lawvere1997,awodey2006,maclane1998}: it shifts attention from objects to the structure-preserving maps between them and thereby offers one descriptive framework for very different structures. Its use in machine learning is by now well established. Optics \cite{riley2018} describe modular access to and modification of parts of data structures and the forward and backward passes of layered networks; functorial accounts of backpropagation, gradient-based learning and architectures \cite{fong2019,cruttwell2022,shiebler2021,gavranovic2024} treat a network as a morphism in a suitable category. These approaches already address compositional structure and architectural constraints. Our emphasis is a concrete record language inside the objects: which elements of a layer's input are grouped, ordered, or held fixed, and how those records determine the maps that are admitted. This emphasis complements the existing categorical accounts.

This paper introduces the Token Space to develop that record-level perspective. Its initiative is to imbue AI computations with meanings at the Token level, such as token grouping, order and parameter types, without giving up the ability to represent every function as a morphism. A \emph{Token} is a finite tuple whose entries are elements of a carrier set or symbols of a fixed \emph{core}; for example $(+,x,y,z)$ records the fact $x+y=z$, and $(c_k,s,t)$ records that two sequences agree on their first $k$ entries. A \emph{Token class} is a carrier together with a \emph{heap} of such Tokens, and a \emph{Token map} is a function on carriers that preserves every Token of the heap. Which structure a computation respects is therefore decided by the choice of heap, and the morphisms acquire an internal categorical structure: layers of a network become Token maps between classes of sequences, architectures become Token maps built from products, exponentials, evaluation and composition, and parameters become points of bases.

\paragraph{A mathematical language for connecting questions.}
Our aim is to provide a common language in which different lines of
reasoning about AI computation can meet. A useful language makes it
possible to express a structure, move between justified descriptions,
and locate a question at the level where its assumptions and answers
belong. Here the carrier identifies what is being studied, the heap
selects the relations of interest, and a Token map states what a
transformation preserves. Full faithfulness, evaluation, gluing and
state factorization then establish specific connections between these
choices. The examples show how those connections work; counterexamples
explain when a proposed identification or transfer fails. The framework
builds on established relational category theory and develops this
record language for structures, programs and AI computations. Its
usefulness can be assessed through the questions it clarifies and the
constructions it connects, alongside subsequent algorithmic studies.

\subsection*{Position: five theses}

The construction is driven by five convictions about category theory and AI computation. They determine what the framework is meant to make visible and operable. The results below give these convictions a mathematical form; their stated hypotheses identify the part of the programme established here.

\begin{enumerate}[label=\textbf{T\arabic*.},ref=T\arabic*,leftmargin=2.5em]
\item\label{thesis:interior} \textbf{AI computation needs the interior of objects as data.} The categorical study of an object through its maps is powerful, but a computing machine also needs a presentation of the internal relationships it is to manipulate. We want one record language in which grouping, order, operation roles and fixed parameters can be inspected, transported and changed across the chosen kinds of structure. Token Space makes this requirement explicit: a Token is a finite record and a heap is the selected interior available as data. For a fixed finitary signature and all its homomorphisms, Theorem~\ref{thm:tokenization} proves that this presentation preserves exactly the specified structure and maps. The operator calculus of Section~\ref{sec:operators} then acts on the records themselves. The choice of interior is part of the model, and a finite record language does not require the whole heap to be finite.

\item\label{thesis:executable} \textbf{Category theory should compute with its own objects.} Our ambition is that categories, maps and constructions should enter a machine as structured data, as directly as other mathematical objects do. Proposition~\ref{prop:cat-tokenization} represents small categories and functors fully faithfully and represents natural transformations through tokenized arrow categories. Proposition~\ref{prop:finite-cat-constructions} supplies terminating constructions on finite category records, including a construction of the functor category with its natural transformations. This gives a concrete executable part of the programme, while the general small-category statement concerns representation. We continue the computational direction of Rydeheard--Burstall and functorial data migration \cite{rydeheard1988,spivak2012}, with the aim of placing categorical data and finite-mapping implementations in the same record framework.

\item\label{thesis:interfaces} \textbf{Functions with arbitrary interfaces have a categorical semantics, and their structural obligations belong in the interfaces.} A computation should be described by its domain and codomain together with the relations it is required to preserve. Cartesian closed structure supplies higher-order evaluation and composition for the pure, total function fragment, in correspondence with typed $\lambda$-calculus \cite{lambek1986}. Token Space enriches that account by making the selected relations explicit as heaps. We use \emph{structural refinement} for this enrichment: the carrier remains available, while the heap determines admissible maps. This is related to the motivations of refinement types and relational parametricity \cite{freeman1991,reynolds1983}; the condition proved here is preservation of specified records. The denotation map $\delta$ of Section~\ref{sec:representation-system} assigns an implementation its admissible interface map, so one interface can support several representations of a computation.

\item\label{thesis:occurrence} \textbf{The same vector does not determine the same Token occurrence.} Numerical equality alone does not identify states located at different layers, positions or inputs of a language model. The model must retain that location and must specify how states are compared. We make this discipline explicit through the located state of Table~\ref{tab:token-levels}, tagged layer carriers and comparison maps in Section~\ref{sec:pipeline}. Proposition~\ref{prop:observation-factorization} states exactly when an observation can forget location and depend on the vector alone. Here an occurrence is distinguished from the \emph{structural Token} that records a fact about it. Set-theoretic tags supply the identity data; the categorical organization makes the interfaces and permitted comparisons visible. This combination is why we use Token Space to describe language-model computations.

\item\label{thesis:cores} \textbf{Computation should admit an unbounded, dynamically organized population of Token computing cores.} We envisage a machine whose logical computing units can be generated, connected, expanded, fused and reclaimed as a task unfolds, without a preset bound on their population across all admitted workloads. Each unit has a typed interface, an executable mapping description and an invocation identity. Independent ready units may proceed concurrently, while physical resources determine their actual execution. Section~\ref{sec:unbounded-cores} makes this position precise through finite configurations, effective growth, a uniform witness of unbounded ready width and a colimit of compatible structural stages. The finite graph theorems provide local correctness; certified transformations connect changes of granularity. The broader objective is a machine that can reorganize its own represented computations while preserving their obligations and making progress under finite resources.
\end{enumerate}

These theses have a common direction: expose structure, make that structure available to computation, specify what its transformations must preserve, retain the identities needed to interpret those transformations, and organize an extensible population of computational units. The represented-mapping machine of Sections~\ref{sec:finite-representations}--\ref{sec:elastic-execution} is the computational expression of this direction. Its representations and execution boundaries can change under explicit correctness conditions, while the relations that matter remain available for inspection.

The construction rests on two elementary rules---a category with one object and its identity is a bicartesian closed category and a topos, and products of bicartesian closed categories are bicartesian closed---followed by a \emph{subsets extension} that selects the heaps. We establish that the fixed-core elementary Token category is complete, cocomplete and a quasitopos, hence locally cartesian closed, while it is not a topos. The classification applies established relational category theory, with a separate factor for nullary facts; dependent products describe varying structural interfaces, while the strong-subobject classifier distinguishes induced structure from arbitrary subobjects (Section~\ref{elementary-Token-space-and-Token-topoi}). A bounded variant keeps every carrier and heap finite, and the algebra $\ATC$ of concrete Token categories organizes iterated constructions; categorical closure properties of nested extensions require separate proofs.

The central theme of the paper is \emph{algebraic tokenization}. Recording all constants, relations and graphs of operations of a structured set as Tokens headed by a core symbol yields a Token class, and we prove full faithfulness for a fixed finitary signature with all its homomorphisms as arrows. Binary products and equalizers are preserved when they are the componentwise and induced constructions in the source category (Theorem~\ref{thm:tokenization}, Proposition~\ref{prop:tokenization-limits}). Algebraic tokenization thus places algebras, ordered sets, graphs and vector spaces (and, contravariantly and faithfully, topological spaces) in one ambient category together with data sets, finite mappings and neural network layers, where they can be compared and combined by a common calculus of operators (Section~\ref{sec:operators}). On objects and morphisms, replacing operations by their graphs is the classical passage from algebraic to relational signatures \cite{hodges1993,lawvere1963}; what the present construction makes explicit is a common core for the chosen collection of signatures, so that data may be moved between carrier and signature, exponential classes whose global elements are the homomorphisms, and tree Tokens for nested terms (Remark~\ref{rem:relational}). Section~\ref{sec:operation-graphs-closure} explains the distinction
between faithful operation-graph encoding and ambient cartesian closure:
a two-element group example shows why a relational function object may
leave the encoded algebraic category. The word ``tokenization'' is chosen deliberately: the tokens of natural language processing \cite{sennrich2016,kudo2018} are elements of the carrier of a sequence class, and the Tokens of this paper are the structural facts about them (Section~\ref{sec:token-levels}).

Three applications test the framework. First, a machine built from \emph{represented finite mappings}: a mapping on finite interfaces may be stored as a table, a circuit, a decision diagram \cite{bryant1986} or a parameterized model, while its implementations remain distinct objects with certified transformations between them. We prove evaluation and sharing results for finite computation graphs and relate cross-language correctness to Token maps between evaluation classes (Section~\ref{sec:finite-representations}). Second, Transformers \cite{vaswani2017}: permutation heaps characterize equivariant layers, prefix-agreement heaps characterize causal decoder blocks, and we verify the latter for a complete decoder block with masked attention, rotary positions, normalization, residual paths and feed-forward maps (Section~\ref{attention-in-token-space}). Third, knowledge distillation as structure-preserving compression: a student is faithful to a teacher iff it is a Token map from the teacher-induced class, symmetries of the teacher are inherited by the student, and the smallest relation-preserving and relation-reflecting state quotient is determined by an indiscernibility congruence (Section~\ref{distillation}). This quotient concerns state identification, separately from implementation cost.

\paragraph{Causal computation and exact retained state.}
The language-model application is developed in this framework manuscript,
while the mathematical companion focuses on dependency presentations.
Sections~\ref{llm:prefix}--\ref{sec:parameterised-cache-state} connect
coherent occurrences to numerical prefixes, exact cache simulation,
safe eviction and parameterised state. These constructions reuse the
frontier and migration results, rather than introducing a second
execution theory. Theorem~\ref{obs:observation-compression} then
identifies an explicit future-observation heap for which structural
indiscernibility and behavioral equivalence coincide. Correct observation
factorization and correct state update are separate obligations. The
result puts exact state compression within the record framework while
leaving efficient implementations and approximate compression as further
questions. The state and coalgebra constructions are standard; their
roles and connections to the selected heaps are stated explicitly.

A further connection concerns explicit knowledge and learned internal computation (Section~\ref{sec:explicit-hidden}). Knowledge graph facts and selected ontology rules enter the record language, while declared readouts make correspondences with contextual model states precise. A commuting square specifies when an internal computation agrees with an explicit one; its validity and the admissibility of its arrows are additional conditions. This connection extends the application programme without requiring every hidden state to have a named concept or treating structural preservation as logical soundness.

Thesis~\ref{thesis:cores} states the kind of computing machine that this construction is intended to support. We envisage finite mappings with several compressed or expanded representations, whose internal components can execute independently when their dependencies permit. Its computational scale should be adjustable: worker allocation, active tasks, representation granularity and input size need not be fixed together. Section~\ref{sec:elastic-execution} formalizes a first part of that programme. Pure acyclic graph executions have a schedule-independent result. Their completed sets form a finite distributive lattice realized by a fully faithful diagram in a slice of Token Space. We prove universal gluing of completed parts, a completed-set adjunction for sharing, and functorial composition of frontier transitions. Explicit module expansion preserves the computed mapping; a completed boundary determines a residual mapping that may be replaced by a certified equivalent implementation. A block-prefix family demonstrates exact changes of granularity with explicit work and span bounds. These claims concern specified descriptions and execution rules, and do not follow from the finiteness of an extensional mapping alone. They give concurrency and elasticity a place in the framework alongside representation compression. Section~\ref{sec:unbounded-cores} adds the population-level foundation: every reached configuration is finite, but no common finite bound need constrain all admitted workloads, and the available concurrency can grow without such a bound.

The paper is organized as follows. Section~\ref{sec:construction-roadmap} states the two construction rules and the purpose of each step. Section~\ref{sec:token-levels} separates the levels at which the word ``token'' is used in language models. Sections~\ref{construction-operators-and-initial-categories} and~\ref{elementary-Token-space-and-Token-topoi} construct the elementary Token space, prove its categorical properties, introduce parameter objects and the bounded variant, and define the algebra $\ATC$ of Token categories. Section~\ref{representing-categories-of-structured-objects-in-Token-space} develops algebraic tokenization. Section~\ref{sec:operators} introduces the operator calculus and applies it to data sets. Section~\ref{interior-structure-mapping-and-tree-Token-classes} treats tree Tokens and reification. Section~\ref{sec:finite-representations} specifies represented finite-mapping machines; Section~\ref{sec:elastic-execution} develops their concurrent and elastic execution, including information and concurrency criteria for migrations. Section~\ref{v13:causal:section} summarizes causal gluing, designated-root function spaces and finite-stage criteria from the self-contained mathematical companion \cite{pancausal}, and explains their role in the machine programme. Sections~\ref{attention-in-token-space} and~\ref{distillation} connect the framework to Transformers and structural distillation. Section~\ref{sec:explicit-hidden} relates explicit knowledge structures to hidden-state computation, Section~\ref{sec:limits} states the scope, and Section~\ref{conclusion} lists open problems.

\subsection*{Main results at a glance}

\begin{itemize}[leftmargin=1.5em]
\item \emph{The category.} $\TS_E$ is complete, cocomplete and a quasitopos, with dependent products and a classifier of induced subobjects (Theorem~\ref{v13:thm-quasitopos}); its monomorphisms are the injective Token maps and it is not a topos (Proposition~\ref{prop:mono}); a bounded variant has finite carriers and heaps (Proposition~\ref{prop:bounded}); parameter objects classify admissible families (Proposition~\ref{prop:parameters}); the algebra $\ATC$ closes the construction under nesting (Theorem~\ref{thm:atc}).
\item \emph{Algebraic tokenization.} Categories of structures for a fixed finitary signature with all homomorphisms embed fully faithfully; componentwise binary products and induced equalizers are preserved under the stated closure hypotheses (Theorem~\ref{thm:tokenization}, Proposition~\ref{prop:tokenization-limits}). Topological spaces have a faithful contravariant representation (Section~\ref{sec:top}). Small categories, functors and natural transformations have explicit encodings, with terminating finite constructions (Propositions~\ref{prop:cat-tokenization}, \ref{prop:finite-cat-constructions}); the choice between core and carrier decides the admissible symmetries (Remark~\ref{rem:signature}).
\item \emph{Operators and trees.} A calculus of operators compares structures across signatures and turns a data set into its covariant Token classes (Section~\ref{sec:operators}); tree Tokens are generated by three operators, and reification embeds the supported tree substitution category fully faithfully into an extended-core Token Space (Theorem~\ref{thm:reification}).
\item \emph{Represented mappings.} Evaluation of a valid computation graph is a unique Token map and graph Token maps certify sharing (Theorem~\ref{thm:graph-evaluation}, Proposition~\ref{prop:sharing}); module expansion preserves denotation (Theorem~\ref{thm:module-expansion}); correct re-encoding is a Token map between evaluation classes (Proposition~\ref{prop:reencoding}).
\item \emph{Concurrency and elasticity.} Partial executions are subclasses, completed parts glue along squares that are both pullbacks and pushouts, and every maximal progressing legal execution ends at the same Token map (Theorems~\ref{thm:schedule-independence}, \ref{thm:execution-gluing}); frontier interfaces form a functor with natural residual maps (Theorem~\ref{thm:frontier-functor}); sharing induces an adjunction and a closure on completion lattices (Theorem~\ref{thm:sharing-adjunction}); unfinished work may be replaced at a completed boundary (Theorem~\ref{thm:residual-replacement}); a block-prefix family realizes selectable granularity (Proposition~\ref{prop:block-prefix}).
\item \emph{Unbounded core organization.} Effective finite transitions preserve finite configurations; a uniform product family exposes arbitrarily many ready invocations; compatible growing record classes have their union as a colimit in $\TS_E$ (Propositions~\ref{prop:finite-core-stages}, \ref{prop:unbounded-core-width}, \ref{prop:core-union}). Finite generated requests complete within a continuing population under an explicit progress condition (Proposition~\ref{prop:finite-request-progress}). These are structural and operational foundations for T5, not a claim of unlimited physical capacity.
\item \emph{Transformers and distillation.} Located states distinguish occurrence identity from numerical equality; an observation factors through the vector exactly when it is constant on its coordinate fibers (Proposition~\ref{prop:observation-factorization}). Permutation heaps characterize equivariant layers and prefix-agreement heaps characterize causal computation, verified for a full decoder block (Proposition~\ref{prop:equivariance}, Theorem~\ref{thm:causal}, Proposition~\ref{prop:decoder-causal}); a faithful student is a Token map from the teacher-induced class, symmetries are inherited, and the smallest exactly faithful compression is the quotient by the indiscernibility congruence (Propositions~\ref{prop:distill-equivariance}, \ref{prop:congruence}).

\item \emph{Exact retained state.} Coherent occurrence prefixes induce
natural numerical maps (Theorem~\ref{llm:prefix-naturality}); the
cache invariant and deletion square compare execution organizations
(Theorems~\ref{llm:cache-simulation} and~\ref{llm:eviction-theorem}).
Changing interfaces carry standard parameterised state
(Proposition~\ref{prop:parameterised-cache-state}). A specifically
chosen future-observation heap identifies the structural quotient with
the minimal reachable behavioral realization, with a right-stability
condition for exact updates (Theorem~\ref{obs:observation-compression}).
\end{itemize}

\noindent Two structural results from the mathematical companion connect record design to computation: relative rootedness exactly characterizes gluing of all induced completed fragments, and a specified result position per record shape gives a cartesian closed category with explicit exponentials (summarized in Section~\ref{v13:causal:section}; complete proofs in \cite{pancausal}). This is a sufficient construction, not a necessary boundary or an assertion that executable graphs are closed under exponentials. For elastic execution, exact reconstruction of a target reachable state is governed by a kernel inclusion, while a graph map has an exact ready-cube image formula (Section~\ref{v13:mig:criteria}). These results keep the retained information and exposed concurrency visible when representations change.

\paragraph{Companion papers and scope.}
This paper concentrates on the Token Space construction, the five theses and the represented-mapping machine. The structural theory is developed in the self-contained mathematical companion \emph{Causal Token Structures} \cite{pancausal}; its results are summarized here at their points of computational use. The elementary classification is calibrated against established relational theory, with the full dependent-product construction in that companion. The graph execution, state migration, sharing and concurrency proofs remain in this paper. The separate working manuscript \emph{Types as Token Classes} \cite{pantypes} develops data, codata and constraint satisfaction, providing a third focus on types and programming. Each manuscript defines the structures needed for its own principal arguments; the shared record notation is not a circular proof dependency.

\paragraph{Four uses for studying language models.}
The framework offers a language for connecting questions that arise at
different levels of LLM computation. Table~\ref{tab:llm-uses} identifies
four uses of that language and points to the constructions that make
each question precise.

\begin{table}[htbp]
\centering\small\renewcommand{\arraystretch}{1.16}
\begin{tabularx}{\linewidth}{@{}>{\raggedright\arraybackslash}p{0.26\linewidth}>{\raggedright\arraybackslash}X>{\raggedright\arraybackslash}p{0.19\linewidth}@{}}
\toprule
Research question & What the framework offers & Where developed \\
\midrule
Which execution changes preserve specified outputs and allow an inference process to continue?
& Cache simulation and reachable-state reconstruction separate equality of complete computations from correctness of changing a running state. The relevant square specifies what must continue to agree.
& Theorems~\ref{llm:cache-simulation}, \ref{v13:mig:kernel}. \\
\addlinespace
Which distinctions between states are necessary for future outputs, and can a compressed state still be updated?
& Future-observation heaps identify the required behavioral distinctions. The right-congruence condition adds the update obligation; safe eviction separately concerns addresses retained by a prescribed implementation.
& Theorems~\ref{obs:observation-compression}, \ref{llm:eviction-theorem}. \\
\addlinespace
How can alternative representations and changing resource arrangements be studied within one computational description?
& Representation systems retain distinct implementations of one mapping. Frontiers and residual interfaces specify where execution can continue after a certified change, while dependencies expose the available independent work.
& Sections~\ref{sec:representation-system}, \ref{sec:frontier-functor}, \ref{sec:online-replacement}. \\
\addlinespace
Along which structural dimensions can LLM descriptions, inference protocols, and state realizations be compared?
& Prefix coherence, dependency structure, state simulation, and observation quotients provide distinct tests. They organize comparisons by the property being investigated, with model parameters and protocol assumptions stated explicitly.
& Sections~\ref{llm:prefix}, \ref{sec:llm-elasticity}, \ref{obs:future-observation}. \\
\bottomrule
\end{tabularx}
\caption{Four uses of Token Space for locating and connecting questions about LLM computation.}
\label{tab:llm-uses}
\end{table}

These uses range from proved correctness criteria to ways of organizing
further research. Practical gains depend on the representations,
algorithms, and costs chosen for a particular implementation.

\subsection*{Reading routes}

\begin{itemize}[leftmargin=1.5em]
\item \emph{For the categorical construction:} Sections~\ref{sec:construction-roadmap}, \ref{construction-operators-and-initial-categories}--\ref{representing-categories-of-structured-objects-in-Token-space}, then the overview in Section~\ref{v13:causal:section}, its mathematical companion \cite{pancausal}, and the open problems of Section~\ref{conclusion}.
\item \emph{For machine learning:} Section~\ref{sec:token-levels}, Remark~\ref{rem:signature}, Section~\ref{sec:dataset-examples}, then Sections~\ref{attention-in-token-space}, \ref{distillation} and~\ref{sec:explicit-hidden}.
\item \emph{For the computing-machine proposal:} Section~\ref{sec:finite-representations}, then Section~\ref{sec:elastic-execution} with Appendix~\ref{app:scheduling}, Section~\ref{v13:causal:section} for the structural criteria, and the finite checks of Appendix~\ref{app:finite-checks}.
\item \emph{For the five theses:} Thesis~\ref{thesis:interior} is realized in Sections~\ref{sec:construction-roadmap} and~\ref{representing-categories-of-structured-objects-in-Token-space}; Thesis~\ref{thesis:executable} in Propositions~\ref{prop:cat-tokenization}--\ref{prop:finite-cat-constructions} and Section~\ref{sec:finite-representations}; Thesis~\ref{thesis:interfaces} in Sections~\ref{sec:representation-system} and~\ref{elementary-Token-space-and-Token-topoi}; Thesis~\ref{thesis:occurrence} in Sections~\ref{sec:token-levels} and~\ref{sec:pipeline}; Thesis~\ref{thesis:cores} in Section~\ref{sec:unbounded-cores}, building on the finite execution and replacement theorems.
\end{itemize}

\section{How and Why Is Token Space Constructed?}\label{sec:construction-roadmap}

The initiative behind the construction of the Token Space is to imbue AI computations with meanings at the Token level, such as token grouping, order, and parameter types. Here ``meaning'' refers to stipulated computational and mathematical roles, not a uniquely determined external referent. An abstract morphism by itself does not display the chosen records. Generalized elements $S\to X$, the Yoneda embedding and concrete underlying-set functors offer ways of describing structure; they do not by themselves choose a common record format and the operations a machine will use to edit it. Global elements $1\to X$ are the special case with terminal domain, and need not recover all carrier elements (Remark~\ref{rem:points}). Our requirement is that the selected interior be available as data in a shared language (Thesis~\ref{thesis:interior}). Heaps implement this requirement for the structures treated here, so that their records can be compared, transported and transformed. Empty-heap classes keep every carrier function available, while chosen nonempty heaps express preservation obligations. To establish this categorical structure at the Token level we employ two fundamental rules:
\begin{enumerate}
\item A category with a single object and only its identity morphism (a copy of the terminal category, whose object may be any set $E$) is a bicartesian closed category and a topos.
\item A product of bicartesian closed categories is bicartesian closed; a product of topoi is a topos.
\end{enumerate}
The first rule fails for a discrete category with two or more objects, which has no terminal object. Token Space uses it to attach a fixed set $E$ of \emph{core} symbols to every object without changing the categorical structure of $\Set$ (the product $\Set\times\mathbf E$ is isomorphic to $\Set$), while remaining able to represent any function as a morphism. The genuinely new structure is added in a third step, the \emph{subsets extension} of Section~\ref{subsets-extension-of-subcategories-of-mathbfset}, which turns the elements of $\Set\times\mathbf E$ into labelled sets of Tokens and the functions into Token maps.

A second motivation is computational. A finite mapping $f:A\to B$ is specified by its table, but a circuit, a decision diagram, a factorized calculation, or a neural architecture with fixed parameters may represent the same mapping. Comparing such implementations requires both an account of their behavior and explicit records of their internal organization; the heaps supply the latter. Exposing this organization also makes dependencies, independent work and alternative granularities available to a machine. The intention is to let the execution structure adapt while the specified behavior remains stable, under the correctness conditions developed in Sections~\ref{sec:finite-representations} and~\ref{sec:elastic-execution}. Table~\ref{tab:construction-roadmap} summarizes the route; the definitions and proofs follow in Sections~\ref{construction-operators-and-initial-categories} and~\ref{elementary-Token-space-and-Token-topoi}.

\begin{table}[htbp]
\centering\small\renewcommand{\arraystretch}{1.13}
\begin{tabularx}{\linewidth}{@{}p{0.19\linewidth}p{0.29\linewidth}X@{}}
\toprule
Construction & Mathematical contribution & Computational purpose \\
\midrule
Sets and functions & A carrier $B$ and transformations $f:B\to C$. & Specify inputs, states, outputs, and their extensional behavior. \\
Fixed core & Add $E$; use $f\sqcup\id_E$ on $B\sqcup E$. & Separate transformable data from fixed relation names, operators, and protocol labels. \\
Finite words & Form $F_E(B)=(B\sqcup E)^*$. & Supply ordered records with repeated entries and explicit labels. \\
Selected heaps & Choose $H\subseteq F_E(B)$. & State which connections, dependencies, operation graphs, or constraints belong to this description. \\
Subsets extension & Require $F_E(f)(H)\subseteq K$. & Select the computations and translations that preserve the specified records. \\
Products, exponentials and dependent products & Combine interfaces and obtain map objects and structured families over a base. & Organize parallel composition, parameterized computation, and context-dependent structural obligations. \\
Trees and reification & Record nesting and replace it by explicit node relations. & Describe composite expressions and connect their syntax to graph-based representations. \\
\bottomrule
\end{tabularx}
\caption{Construction of Token Space and the role of each step. The last two rows are constructions within, or extensions of, the elementary framework; their properties require their own proofs.}\label{tab:construction-roadmap}
\end{table}

A mapping enters the construction in two ways. An execution is a morphism $f:T_A\to T_B$, expressing preservation of the chosen input and output relations; the graph of the same mapping can also be stored as data, through records $(\mathrm{apply},\iota_A(a),\iota_B(b))$ with $b=f(a)$ on the tagged carrier $A\sqcup B$. The first view supports composition of computations; the second lets a computation, its representation, or an observation of it become an object of another computation. Exponentials connect the two views categorically. Finite descriptions of Token classes can in turn be encoded as carrier elements, with evaluation and comparison records specified separately; iterating $\ATC$ does not automatically provide such an encoding.

This is why merely retaining all functions in $\Set$ is insufficient. A function space already contains every input--output behavior; what it does not specify is which relations a computation must preserve, how a description encodes it, or which changes to a description preserve its behavior. Heaps make the first question explicit; algebraic tokenization (Section~\ref{representing-categories-of-structured-objects-in-Token-space}) covers the fixed-signature homomorphisms specified in Definition~\ref{def:structure}; the operator calculus (Section~\ref{sec:operators}) and the representation systems of Section~\ref{sec:finite-representations} address the second and third.

\paragraph{Choosing the level of a question.}
The same calculation can be studied through its values, its description,
or the state retained during execution. These choices need different
carriers and heaps. Equality of computed functions concerns behavior;
a map between descriptions can certify a particular transformation;
and replacing a running state additionally requires enough information
to continue. The common record language connects these questions while
keeping their objects and preservation conditions explicit. The prefix
computation in Section~\ref{sec:parity} provides a concrete thread through
these levels.

\paragraph{How to read the diagrams.}
An arrow inside a computation graph is a dependency between operation nodes. An arrow in a categorical diagram is instead a map between the displayed objects: whole Token classes, description sets, or value interfaces. The ambient category is stated in each caption. A commuting diagram asserts equality of the composites along its routes; a pullback or pushout additionally asserts a universal factorization property. These distinctions make the diagrams useful for checking the construction, rather than only illustrating its vocabulary.

\begin{table}[htbp]
\centering\small\renewcommand{\arraystretch}{1.13}
\begin{tabularx}{\linewidth}{@{}>{\raggedright\arraybackslash}p{0.24\linewidth}>{\raggedright\arraybackslash}p{0.28\linewidth}>{\raggedright\arraybackslash}X@{}}
\toprule
Question & Diagram to consult & What its arrows express \\
\midrule
What does a Token map preserve? & Figure~\ref{fig:heap-square} & An ambient carrier function transports each selected record. \\
How is a natural transformation data? & Figure~\ref{fig:category-lift} & A lift between encoded categories has prescribed endpoint functors. \\
Why can independent work be joined? & Figures~\ref{fig:execution-gluing} and~\ref{fig:ready-diamond} & Completed classes embed and glue; their values agree. \\
Why does expansion retain behavior? & Figure~\ref{fig:expansion-evaluation} & Changing a description commutes with evaluation. \\
How do live interfaces change? & Figure~\ref{fig:frontier-coherence} & Frontier advancement composes and preserves the remaining output. \\
What work can sharing reuse? & Figure~\ref{fig:sharing-adjunction} & Images of completed sets and consistent values transfer together. \\
Why can an invocation be adapted? & Figure~\ref{fig:online-replacement} & The retained boundary feeds an equivalent residual implementation. \\
When can a retained state be reconstructed? & Figure~\ref{v13:mig:state-output-diagram} & Boundary maps expose which original-input distinctions remain available. \\
How can finite stages have an unbounded horizon? & Figure~\ref{fig:core-union} & Compatible inclusions assemble into a union colimit; values extend uniquely. \\
When do causal fragments glue? & Figure~\ref{v13:causal:gluing-figure} & Induced completed pieces cover each record under the relative-root criterion. \\
What does a designated-root function object do? & Figure~\ref{v13:causal:exponential-figure} & Currying and evaluation preserve the stated order and record obligations. \\
When can location be forgotten? & Figures~\ref{fig:occurrence-forgetting} and~\ref{fig:layer-comparison} & Observations must factor through coordinates; layer comparisons are specified maps. \\
How is LLM causality expressed? & Figure~\ref{fig:causal-factorization} & The first output coordinates factor through the corresponding input prefix. \\
When does internal computation agree with explicit computation? & Figure~\ref{fig:explicit-hidden-square} & Declared readouts connect two computations; commutation is a condition to establish. \\
\bottomrule
\end{tabularx}
\caption{A reading route through the categorical diagrams. Objects and maps are specified in each case. The associated propositions establish the proven diagrams; Figure~\ref{fig:explicit-hidden-square} instead states a correspondence condition to investigate.}\label{tab:diagram-guide}
\end{table}

\section{What Is a Token in This Framework?}\label{sec:token-levels}

\subsection{Vocabulary items, occurrences, states, and structural records}

Fix a deterministic tokenizer $\tau:\mathcal X\to V^*$, where $\mathcal X$ is the set of admissible texts and $V$ is its vocabulary. Its output is a sequence $s=(v_1,\ldots,v_n)$, usually represented by integer IDs. Tokenization may include normalization and special symbols; a vocabulary item need not be a word or an independently meaningful unit \cite{sennrich2016,kudo2018}. We do not assume that tokenization preserves concatenation, or that decoding recovers an unnormalized source exactly.

\begin{table}[htbp]
\centering\small\renewcommand{\arraystretch}{1.13}
\begin{tabularx}{\linewidth}{@{}p{0.20\linewidth}p{0.25\linewidth}X@{}}
\toprule
Object & Mathematical representation & Role in a language model or in Token Space \\
\midrule
Tokenizer token type & $v\in V$ & A vocabulary entry; its ID specifies an embedding-table lookup. \\
Token occurrence & $o=(s,i)$ with value $v_i$ & One occurrence at a particular input position; repeated IDs give different occurrences. \\
Input embedding & $e(v_i)\in\mathbb R^d$ & A learned vector associated with the vocabulary item, before contextual mixing. \\
Located hidden state & $\omega=(\theta,\ell,s,i,h_i^{(\ell)}(s))$ & A state occurrence with its model, layer, input, and position; its vector alone omits this location. \\
Sequence state & $H^{(\ell)}\in(\mathbb R^d)^n$ & The input or output of a Transformer block. This is one carrier element in the sequence classes used below. \\
Structural Token & $r\in\Upsilon\subseteq(B\sqcup E)^*$ & A finite record in the mathematical description, such as $(c_k,H,H')$ or $(\mathrm{next},o_i,o_{i+1})$. \\
\bottomrule
\end{tabularx}
\caption{The objects called ``tokens'' in informal discussion occupy distinct levels. A Token class is an object of $\TS_E$; an individual structural Token is an element of its heap.}\label{tab:token-levels}
\end{table}

For example, suppose a particular tokenizer produces the three IDs $(a,b,a)$. There is one vocabulary type $a$ at two positions, two corresponding occurrences $(s,1)$ and $(s,3)$, and a shared lookup vector $e(a)$. Positional processing and attention can produce different hidden states at the two positions. A class with carrier $V$ describes vocabulary types; a class with carrier $O_s=\{(s,i):1\le i\le n\}$ describes occurrences; a class with carrier $V^n$ describes complete sequences. None of these choices should be silently substituted for another.

A vocabulary type can be encoded as the one-entry formal tuple $(v)$. This is an encoding of the tokenizer token into the tuple language, not an identification of all structural Tokens with tokenizer outputs. Relations among occurrences require additional records, and contextual computation requires an appropriately enlarged carrier.

\subsection{Structure, relative atomicity, and semantic commitments}

A tokenizer token is atomic relative to a chosen input interface: the model receives its ID as one discrete input unit. It is not thereby a smallest semantic unit. Likewise, a structural Token can serve as one elementary record without being internally indivisible: a tuple has entries, and a tree Token has subtrees. This relative atomicity also applies to executable descriptions: a mapping used as one module at an interface may expose a finer implementation graph when expanded. The module and the expanded graph must have the same specified denotation, while their physical evaluation costs remain explicit (Section~\ref{sec:granularity}). When a headed tuple represents a single atomic formula $R(b_1,\ldots,b_k)$ in a fixed relational signature, we may call it an \emph{atomic relational record}. We do not use ``atomic'' as the general definition of Token.

The basic data $(B,E,\Upsilon)$ do not include a map assigning each carrier element or tuple a privileged external referent or natural-language content. They can nevertheless carry precise mathematical interpretations and empirically grounded relations. In an algebraic example, $(+,x,y,z)$ records a stipulated operation graph. In an interpretability application, $(\mathrm{readout},\bar P,\omega,\bar c)$ records an output of a specified readout $P$. It does not, solely by membership in a heap, establish that the intrinsic content of $\omega$ is $c$. The absence of a semantic assignment from a formal definition must be distinguished from the claim that its realizations have no content.

\section{Construction Operators}\label{construction-operators-and-initial-categories}

\subsection{Fixed cores, tagged unions, and finite words}\label{identity-set-categories}\label{products-of-categories}

All signatures and carriers are sets in a sufficiently large ambient universe. The notation $A\sqcup E$ means a tagged disjoint union; an element with the same printed name in the carrier and the core remains distinguishable. We write
\[
A^*=\coprod_{m\ge0}A^m
\]
for finite words, including the empty word. Writing a word as a tuple gives the equivalent notation $A^{\&}$ used for a \emph{Tokenoid}. A function $f:A\to B$ extends entrywise to $f^*:A^*\to B^*$. These extensions preserve word length. They form the image of the finite-word functor; they are not all homomorphisms between free monoids, since a general monoid homomorphism can send one generator to a longer word.

\begin{remark}[The free-monoid route]\label{rem:fmon}
Equivalently, let $\mathbf{FMon}_{\Set}$ be the category of free monoids $\mathcal H(A)=\langle A^*,\cdot,\varepsilon\rangle$ and of the homomorphisms $f^*$ extended from functions $f$ between generator sets, let $\mathcal F:\mathbf{FMon}_{\Set}\to\Set$ be the forgetful functor, and let $(\cdot)$ write a word $x_1\cdot x_2\cdots x_n$ as the tuple $(x_1,\ldots,x_n)$. The composite $\mathcal M=(\cdot)\circ\mathcal F\circ\mathcal H$ is an isomorphism of $\Set$ onto the \emph{pre-Token category} $\Tok_{\Set}$, whose objects are the Tokenoids $A^{\&}$; hence $\Tok_{\Set}$, and $\mathcal M(U_E(\Set\times\mathbf E))$ for the functor $U_E$ below, are bicartesian closed topoi whose structure components are those of $\Set$ transported along $\mathcal M$. The finite-word functor $W$ is $\mathcal F\circ\mathcal H$. The algebraic description is what makes the tree Tokens of Section~\ref{interior-structure-mapping-and-tree-Token-classes}, obtained by adding one unary operator to the monoid signature, a natural continuation.
\end{remark}

For a fixed set $E$, the identity set category $\mathbf E=\mathbf{idset}(E)$ has the single object $E$ and only $\id_E$. The full subcategory $\mathbf C_*$ of $\Set$ on the singleton sets is equivalent to it, with $\mathbf C_0=\{\{*\}\}$ as a skeleton. Abstractly $\mathbf E$ is a terminal category: $E$ is simultaneously its initial object, terminal object, product, coproduct, exponential $E^E$ and subobject classifier $\Omega_E$, so $\mathbf E$ is a degenerate topos and a bicartesian closed category \cite{lambek1986,goldblatt1984,maclane1998,johnstone2002}. The product $\Set\times\mathbf E$ is isomorphic to $\Set$. Its concrete realization sends $(B,E)$ to $B\sqcup E$ and $(f,\id_E)$ to $f\sqcup\id_E$. It attaches fixed labels, but by itself imposes no new condition on $f$.

For products of concrete categories, tagged disjoint union of their underlying sets gives a faithful realization of the product. If the original category objects carry structure, their labels are retained: a faithful underlying-set functor need not be injective on objects. Abstract products of bicartesian closed categories are bicartesian closed, but adding selected subsets to a concrete realization requires a separate argument.

\subsection{The construction as a sequence of functors}\label{sec:functorial-route}

Let $U_E:\Set\times\mathbf E\to\Set$ send $(B,E)$ to $B\sqcup E$ and $(f,\id_E)$ to $f\sqcup\id_E$. Let $W:\Set\to\Set$ be the finite-word functor. The concrete route is
\[
\Set\times\mathbf E\xrightarrow{\ U_E\ }\Set
\xrightarrow{\ W\ }\Set,
\qquad
(B,E)\longmapsto B\sqcup E\longmapsto(B\sqcup E)^*.
\]
It is useful to separate the indexing category $\Set\times\mathbf E$ from the concrete sets used to display its objects. Both functors are faithful. Through the isomorphism $\Set\times\mathbf E\cong\Set$, their composite is $F_E(B)=(B\sqcup E)^*$.

At this stage no heap has been chosen and no carrier function has been excluded. The fixed core forces the same symbols to recur unchanged; word formation supplies ordered positions and repeated occurrences. Selecting a subset of each word set is the step that specifies a structure and restricts its admissible maps. With the subsets extension defined below, the full construction can be written
\[
\TS_E\cong\SX_{W\circ U_E}(\Set\times\mathbf E)
=\SX_{F_E}(\Set)
\]
where the last equality uses the identification of the indexing categories.

Products of several identity set categories similarly separate several fixed symbol sets before their tagged union is taken, and products with structured categories distinguish several variable sorts. Closing the construction under these operations yields the algebra $\ATC$ of Token categories of Section~\ref{Token-categories}.

\subsection{Subsets extension with ambient maps retained}\label{subsets-extension-of-subcategories-of-mathbfset}

\begin{definition}[Subsets extension]\label{def:sx}
Let $J:\mathbf C\to\Set$ be a faithful functor. The category $\SX_J(\mathbf C)$ has objects $(A,S)$, where $A\in\mathbf C$ and $S\subseteq J(A)$. A morphism
\[
(A,S)\longrightarrow(B,T)
\]
is a morphism $f:A\to B$ in $\mathbf C$ such that $J(f)(S)\subseteq T$. Its restriction $J(f)|_S:S\to T$ is determined by $f$, but is not substituted for $f$ in the definition of equality of morphisms. Composition is composition in $\mathbf C$.
\end{definition}

Identity and composition preserve the required inclusions, so this is a category. The projection $\mathrm{Or}(A,S)=A$, $\mathrm{Or}(f)=f$, is a faithful functor, and $A\mapsto(A,J(A))$ is fully faithful. One may regard a morphism as the pair $(f,J(f)|_S)$, with the second component redundant. The notation $\mathbf{on}(f,S,T)$, when used, denotes this morphism with its ambient map retained.

\begin{example}[Why restrictions alone are insufficient]\label{ex:empty-heap}
For $A=\{a\}$, $B=\{b,c\}$, and $S=T=\varnothing$, there are two morphisms $(A,S)\to(B,T)$, given by the two functions $A\to B$. Their restrictions to $S$ are the same empty function. Identifying those restrictions would produce a different category and would lose the carrier computation at unrecorded elements.
\end{example}

When $\mathbf C=\Set$ and $J=\id$, objects are sets with one distinguished subset. Their product is $(A\times B,S\times T)$, their coproduct is $(A\sqcup B,S\sqcup T)$, their terminal object is $(\{*\},\{*\})$, and their initial object is $(\varnothing,\varnothing)$. The exponential is
\[
(B,T)^{(A,S)}=
\bigl(B^A,\{f:A\to B\mid f(S)\subseteq T\}\bigr).
\]
Indeed, $g:C\times A\to B$ preserves the distinguished subsets exactly when its curry sends the distinguished subset of $C$ to the displayed subset of $B^A$. Equalizers have the induced distinguished subset. Thus this extension has finite limits and is bicartesian closed. It is not a topos: the identity map $(\{*\},\varnothing)\to(\{*\},\{*\})$ is both monic and epic but not invertible. Every topos is balanced \cite{johnstone2002}.

The extension depends on the concrete functor $J$, not only on the abstract category $\mathbf C$. In particular, Cartesian closedness of an arbitrary $\mathbf C$ does not by itself imply Cartesian closedness of $\SX_J(\mathbf C)$. We prove the required result for Token Space directly.

\section{Elementary Token Space}\label{elementary-Token-space-and-Token-topoi}

\subsection{Objects, morphisms, and relational interpretation}

Fix a core $E$ and put $F_E(B)=(B\sqcup E)^*$ and $F_E(f)=(f\sqcup\id_E)^*$. The functor $F_E:\Set\to\Set$ is faithful, since its restriction to one-entry carrier words recovers $f$.

\begin{definition}[Token class and Token map]\label{def:token-class}
A Token class over $E$ is
\[
T=\langle B,E,\Upsilon\rangle,\qquad \Upsilon\subseteq F_E(B).
\]
The set $B$ is its \emph{base} or carrier, $E$ its \emph{core}, and $\Upsilon$ its \emph{heap}. Elements of $\Upsilon$ are structural Tokens. A Token map $T_1\to T_2$ is a function $f:B_1\to B_2$ satisfying
\begin{equation}\label{eq:hom}
F_E(f)(\Upsilon_1)\subseteq\Upsilon_2.
\end{equation}
It can be displayed as $\langle f,\id_E,p\rangle$, where $p=F_E(f)|_{\Upsilon_1}:\Upsilon_1\to\Upsilon_2$. Equality of Token maps means equality of the carrier functions. The category of these objects and maps is $\TS_E=\SX_{F_E}(\Set)$.
\end{definition}

\begin{figure}[htbp]
\centering
\begin{tikzcd}[column sep=large,row sep=large]
\Upsilon_A \arrow[r,"p"] \arrow[d,hook] & \Upsilon_B \arrow[d,hook] \\
F_E(A) \arrow[r,"F_E(f)"'] & F_E(B)
\end{tikzcd}
\caption{The record-preservation square in $\Set$ for a carrier function $f:A\to B$. A Token map exists exactly when $F_E(f)$ restricts to the displayed $p$. The carrier function $f$ remains part of the map; the restriction $p$ alone can lose information, as Example~\ref{ex:empty-heap} shows.}\label{fig:heap-square}
\end{figure}

The map $p$ is a total function on the heap, although it is a restriction of a function on all words. A Token map is therefore not a partially defined carrier computation. A heap need not contain every word, need not be finite, and need not be closed under concatenation. Finiteness is a property of each individual Token.

For each word length, choice of carrier positions, and assignment of core symbols to the remaining positions, a heap determines a relation on the carrier positions. All-core words correspond to nullary relations. Conversely, this family of relations reconstructs the heap. Hence $\TS_E$ is a category of relational structures for a fixed, possibly infinite, finitary signature \cite{hodges1993,adamek1990}. This observation explains the constructions below; it does not supply a natural-language interpretation of the relations.

\subsection{Finite limits, coproducts, and exponentials}

\begin{theorem}\label{thm:tse}
The category $\TS_E$ has all finite limits, finite coproducts, and exponentials. It is bicartesian closed but is not a topos.
\end{theorem}

\begin{proof}
The initial and terminal objects are respectively
\[
0_E=\langle\varnothing,E,\varnothing\rangle,
\qquad
1_E=\langle\{*\},E,F_E(\{*\})\rangle.
\]
The unique empty carrier function preserves an empty heap, and the unique function to $\{*\}$ preserves a full target heap.

Let $T_A=\langle A,E,H_A\rangle$ and $T_B=\langle B,E,H_B\rangle$. With projections $\pi_A,\pi_B$ from $A\times B$, the product, also denoted $T_A\otimes T_B$, has carrier $A\times B$ and heap
\begin{equation}\label{eq:product-heap}
F_E(\pi_A)^{-1}(H_A)\cap F_E(\pi_B)^{-1}(H_B).
\end{equation}
The pairing of two carrier maps preserves this intersection precisely when both component maps preserve their heaps. In a product Token the component words must have the same length and the same core symbols at the same positions; the other entries are pairs of carrier elements.

The equalizer of $f,g:T_A\to T_B$ has carrier $A_0=\{a:f(a)=g(a)\}$ and heap $H_A\cap F_E(A_0)$, identifying $F_E(A_0)$ with its image under inclusion. Any map equalizing $f$ and $g$ factors uniquely through this carrier, and heap preservation passes to the factor. Together with the terminal object and products, this gives finite limits.

The coproduct $T_A\uplus T_B$ has carrier $A\sqcup B$ and heap
\begin{equation}\label{eq:coproduct-heap}
F_E(\iota_A)(H_A)\cup F_E(\iota_B)(H_B).
\end{equation}
A copairing preserves this heap exactly when both maps preserve their source heaps. The same formula, with a union over the indices, constructs set-indexed coproducts.

For the exponential $T_B^{T_A}$, take carrier $B^A$. Call a word over $B^A\sqcup E$ and a word over $A\sqcup E$ \emph{compatible} when they have the same length, the same core positions, and equal core symbols at those positions. For compatible words $u=(u_j)$ and $a=(a_j)$ define $\operatorname{ev}_*(u,a)$ entrywise: use $u_j(a_j)$ at carrier positions and the common symbol at core positions. Put
\begin{equation}\label{eq:exp-heap}
H_B\vartriangleleft H_A
=\{u\in F_E(B^A):
 \operatorname{ev}_*(u,a)\in H_B\text{ for every compatible }a\in H_A\}.
\end{equation}
This condition includes the case of all-core words and is vacuous when no compatible source Token exists. Evaluation is a Token map by construction. A function $g:C\times A\to B$ preserves the product heap iff its curry $\widehat g:C\to B^A$ sends every Token of $T_C$ into \eqref{eq:exp-heap}. The usual natural currying bijection therefore restricts to
\[
\Hom_{\TS_E}(T_C\otimes T_A,T_B)
\cong\Hom_{\TS_E}(T_C,T_B^{T_A}).
\]
The failure of the topos property follows from Proposition~\ref{prop:mono} below and the subsequent counterexample.
\end{proof}

\begin{table}[htbp]
\centering\small\renewcommand{\arraystretch}{1.13}
\begin{tabularx}{\linewidth}{@{}p{0.16\linewidth}p{0.20\linewidth}X@{}}
\toprule
Construction & Carrier & Heap and computational role \\
\midrule
Initial object $0_E$ & $\varnothing$ & Empty heap. The empty interface has a unique map into any class. \\
Terminal object $1_E$ & $\{*\}$ & $F_E(\{*\})$. A target for discarding an input; maps out of it select global elements. \\
Product & $A\times B$ & $F_E(\pi_A)^{-1}(H_A)\cap F_E(\pi_B)^{-1}(H_B)$. Combines two interfaces with matched record shapes. \\
Coproduct & $A\sqcup B$ & $F_E(\iota_A)(H_A)\cup F_E(\iota_B)(H_B)$. Keeps alternatives tagged, with casewise maps. \\
Equalizer of $f,g$ & $A_0\subseteq A$ & $H_A\cap F_E(A_0)$, where $A_0=\{a:f(a)=g(a)\}$. Restricts the interface to agreement states. \\
Exponential $T_B^{T_A}$ & $B^A$ & $H_B\vartriangleleft H_A$, defined in~\eqref{eq:exp-heap}. Supports evaluation and currying; its global elements are the admissible maps. \\
\bottomrule
\end{tabularx}
\caption{Elementary Token Space: carriers, heaps, and uses of the universal constructions. A carrier element of the exponential need not be a global element. There is no classifier for all subobjects; induced subobjects do have a classifier (Theorem~\ref{v13:thm-quasitopos}).}\label{tab:elementary-constructions}
\end{table}

Table~\ref{tab:elementary-constructions} separates the underlying sets from the relations that make the constructions universal. In particular, the exponential is not obtained by simply declaring its carrier to be the set of Token maps. Its carrier contains all functions, while its heap controls admissible families and its global elements select the Token maps. This distinction will allow an implementation description to name a computation without identifying that description with the function it computes.
\begin{proposition}\label{prop:mono}
A Token map is monic iff its carrier function is injective, and epic iff that function is surjective.
\end{proposition}
\begin{proof}
The forward implications are tested by maps from $\langle\{*\},E,\varnothing\rangle$ and maps into $\langle\{0,1\},E,F_E(\{0,1\})\rangle$, respectively. These allow every carrier function of the indicated types. The reverse implications follow because Token maps are determined by their carrier functions.
\end{proof}

The identity carrier function
\[
\langle\{b\},E,\varnothing\rangle
\longrightarrow\langle\{b\},E,\{(b)\}\rangle
\]
is monic and epic but has no inverse Token map. Since every topos is balanced, $\TS_E$ is not a topos. In view of its finite limits and exponentials, it cannot have a subobject classifier. We call $T_1=\langle B_1,E,\Upsilon_1\rangle$ a \emph{subclass} of $T_2=\langle B_2,E,\Upsilon_2\rangle$, written $T_1\Subset T_2$, when $B_1\subseteq B_2$ and $\Upsilon_1\subseteq\Upsilon_2$; the inclusion is then a monomorphism, and every monomorphism $f:T_1\to T_2$ is isomorphic, as a subobject of $T_2$, to the subclass $\langle f(B_1),E,F_E(f)(\Upsilon_1)\rangle\Subset T_2$. A subclass need not carry the entire induced heap.

\begin{remark}[Carrier elements are not always global elements]\label{rem:points}
The carrier of $T_B^{T_A}$ is the set of \emph{all} functions $A\to B$. Its global elements $1_E\to T_B^{T_A}$ correspond exactly to Token maps $T_A\to T_B$. In general a carrier element of a Token class does not define a global element: the required map from $1_E$ must preserve its full heap. This distinction is essential for parameter objects in Section~\ref{sec:parameter-objects}.
\end{remark}

\subsection{Relational classification and varying interfaces}\label{v13:sec-quasitopos}

The fixed-core category admits a stronger classification than the finite constructions above. This is an application of established relational category theory, with nullary facts retained explicitly. A \emph{shape} is a finite word in $E\sqcup\{\square\}$; its arity is its number of blanks. A heap $H_A$ specifies, for each shape $s$ of arity $n$, the relation
\[
R_s^A=\{(a_1,\ldots,a_n):s[a_1,\ldots,a_n]\in H_A\}.
\]
Pure-core words, including the empty word, are nullary facts. Letting $\Sigma_E^+$ denote the positive-arity shapes gives
\begin{equation}\label{v13:eq-shape-decomposition}
\TS_E\cong\operatorname{Str}(\Sigma_E^+)\times(\mathcal P(E^*),\subseteq).
\end{equation}
A carrier function preserves the heap exactly when it preserves these relations and the source nullary facts are included in the target ones.

\begin{theorem}[Elementary classification]\label{v13:thm-quasitopos}
For each fixed set $E$, with arbitrary set carriers and arbitrary heaps of finite words, $\TS_E$ is complete, cocomplete and a quasitopos: it has finite limits and colimits, is locally cartesian closed, and has a classifier for strong subobjects. Its regular and strong monomorphisms coincide and are precisely the induced embeddings, namely injective Token maps $m:T_A\to T_B$ with $H_A=F_E(m)^{-1}(H_B)$, including agreement on nullary facts.
\end{theorem}
\begin{proof}
Parker's relational classification \cite[Theorem~4.5]{parker2022} applies to the first factor of~\eqref{v13:eq-shape-decomposition}. The powerset factor and all its slices are complete Boolean algebras; its only strong subobjects are identities. Products therefore give the asserted classification. Limits use carrier limits with records preserved by every projection; colimits use carrier colimits with the union of the image heaps. Equalizers are induced. Conversely, an induced subset is the equalizer of its characteristic function and the constant-$1$ map into the full-heap class on $\{0,1\}$. A strong mono factors through its induced carrier image; its bijective first factor is epic, and the strong lifting property makes that factor invertible. Thus the strong and regular monos are exactly the induced ones. The map from the terminal class selecting $1$ in the full-heap class on $\{0,1\}$ classifies them.
\end{proof}

Local cartesian closure supplies a useful interface construction. For a Token map $p:T_A\to T_B$, the pullback functor has a right adjoint $p^*\dashv\Pi_p$. If $q:T_X\to T_A$ is a structured family, the carrier of its dependent product consists of pairs $(b,\sigma)$ with $\sigma:p^{-1}(b)\to X$ and $q(\sigma(a))=a$ for every $a\in p^{-1}(b)$. Its records require a base record in $T_B$ and preservation by the selected sections of every matching record in $T_A$; a nullary shape has truth value $R_s^B\wedge(R_s^A\Rightarrow R_s^X)$. These are set-theoretic sections, selected into admissible families by the heap. The full construction and adjunction proof appear in the mathematical companion \cite[Appendix~A]{pancausal}. For T2 and T3 this supplies structured choices over changing contexts; a finite description, a dependent-type syntax and efficient evaluation are additional tasks.

The strong classifier does not classify all injective Token maps: such a map can omit records without changing the carrier. The category is accordingly not a topos. Nor does the classification transfer automatically to nested, tree, finite-heap or executable-graph subcategories.

\subsection{The boundary of ambient gluing}\label{v13:sec-adhesive-boundary}

By the quasitopos result, pushouts along regular monos are stable under pullback and are themselves pullbacks: the category is \emph{rm-quasiadhesive}. It is neither rm-adhesive nor adhesive in the terminology of Garner and Lack \cite{garner2012}. To see the first boundary, consider the induced singleton subobjects of the class on $\{u,v\}$ with heap $\{(u,v)\}$. Their pushout over the empty intersection has empty heap, so its comparison mono into the original class is not induced. Regular subobjects are thus not closed under binary unions, contradicting the additional condition for rm-adhesivity \cite[Corollary~3.4]{garner2012}. For the second, the noninvertible bijective Token map from the singleton empty-heap class to the singleton class with heap $\{(b)\}$ is both mono and epi. Its self-pushout has two identity legs and is not a pullback, which cannot happen in an adhesive category.

The computational issue is concrete: a relation can cross two pieces and disappear from both restrictions. The execution gluing proofs below therefore use the dependency structure to ensure that each completed record lies in one piece. Section~\ref{v13:causal:section} gives the exact support criterion; the mathematical companion supplies its general proof and further categorical boundaries \cite{pancausal}.

\subsection{Maps as data and parameter families}\label{sec:parameter-objects}

By Cartesian closedness,
\[
\Hom_{\TS_E}(T_A,T_B)\cong\Hom_{\TS_E}(1_E,T_B^{T_A}).
\]
Thus each admissible computation is a global element of an exponential. To describe an entire parameter family, merely using the empty-heap object $D_E(\Theta)$ is insufficient: every function $\Theta\to\base(T_B^{T_A})$ is then a Token map, including families whose members fail to preserve the desired structure.

\begin{definition}[Parameter object]\label{def:parameter}
For a set $\Theta$, define
\[
P_E(\Theta)=\coprod_{\theta\in\Theta}1_E.
\]
Its carrier is $\Theta$. Its heap is the union, over $\theta$, of all words whose carrier entries are all that same $\theta$, with arbitrary core entries. Equivalently,
\[
\heap(P_E(\Theta))=\bigcup_{\theta\in\Theta}F_E(j_\theta)(F_E(\{*\})),
\qquad j_\theta(*)=\theta.
\]
For $\Theta=\varnothing$ the heap is empty. Each parameter value determines a global element $j_\theta:1_E\to P_E(\Theta)$.
\end{definition}

\begin{proposition}[Parameter families]\label{prop:parameters}
For any $T_A,T_B$ there are natural bijections
\begin{align*}
\Hom_{\TS_E}(P_E(\Theta),T_B^{T_A})
&\cong\Set(\Theta,\Hom_{\TS_E}(T_A,T_B))\\
&\cong\Hom_{\TS_E}(P_E(\Theta)\otimes T_A,T_B).
\end{align*}
Thus a joint function $\Phi(\theta,a)$ is a Token map from the displayed product iff every fixed-parameter function $a\mapsto\Phi(\theta,a)$ is a Token map.
\end{proposition}
\begin{proof}
The coproduct property identifies the first hom set with the product over $\theta$ of $\Hom(1_E,T_B^{T_A})$. The exponential property identifies each factor with $\Hom(T_A,T_B)$ and identifies the outer hom set with the joint evaluation map. These bijections use the usual carrier functions.
\end{proof}

The parameter object encodes admissibility of a family, whether its elements are numerical parameters, circuits, or other computational descriptions. It does not supply a metric, an optimization rule, or identifiability of different descriptions. Composition and evaluation organize forward computation; differential accounts of learning add further structure \cite{cruttwell2022}.

\subsection{Finite carriers and bounded record shapes}\label{sec:bounded}

The finiteness of a carrier does not make its Tokenoid finite: words may have arbitrary length. A finite heap is one possible representation, but the finite-heap objects alone do not inherit the displayed terminal object, whose heap contains words of every length. There is a direct bounded variant suitable for finite specifications.

Fix a finite core $E$ and an integer $m\ge1$. Replace $F_E$ by
\[
F_{E,\le m}(B)=\coprod_{j=0}^{m}(B\sqcup E)^j,
\]
and restrict carriers to finite sets. Denote the resulting category by $\TS^{\mathrm{fin}}_{E,\le m}$. Its arrows are the same entrywise heap-preserving carrier functions, restricted to these objects.

\begin{proposition}[A finite version]\label{prop:bounded}
The category $\TS^{\mathrm{fin}}_{E,\le m}$ has finite limits, finite coproducts, and exponentials. All of their carriers and heaps are finite. It is not a topos.
\end{proposition}
\begin{proof}
Use the constructions of Theorem~\ref{thm:tse}, with every word restricted to length at most $m$. Projection, injection, and entrywise evaluation preserve word length, and the currying argument uses only words of matching shape, so all the same universal properties hold. Products, disjoint unions, subsets, and function sets of finite sets are finite. The noninvertible bijective map from the empty heap on a singleton to the heap containing its one-entry word remains available because $m\ge1$.
\end{proof}

The bounded category has its own terminal and exponential heaps; they must not be replaced by the unbounded ambient heaps. For a finite collection of named operations and relations, $m$ can cover all record lengths in use. This gives a finite mathematical specification, not an efficiency guarantee: a function carrier $B^A$ can still be very large. Alternatively, an unbounded heap may be described by a finite rule, with membership and evaluation implemented symbolically. Either choice requires distinguishing a represented set from the data structure used to store it.

\subsection{Forgetting, filling, and a common core}\label{Token-space}

Define $D_E(B)=\langle B,E,\varnothing\rangle$ and $C_E(B)=\langle B,E,F_E(B)\rangle$. Both are fully faithful embeddings of $\Set$, and
\[
D_E\dashv\base\dashv C_E.
\]
The two adjunctions follow because any carrier function out of an empty heap, or into a full heap, is a Token map. The empty-heap and full-heap classes are thus two different copies of $\Set$ in the ambient category. The functors $\varnothing=D_E\base$ (\emph{forgetting}) and $\circledast=C_E\base$ (\emph{stuffing}) empty and fill the heap, and the \emph{abstracting} functor $\circledcirc(T)=D_E(\heap(T))$, $\circledcirc(f)=p$, turns the heap into a carrier. The image of $C_E$, the full subcategory of \emph{full} classes, is isomorphic to $\Set$ and therefore a topos. By Proposition~\ref{prop:mono} a Token map $f$ is a subclass inclusion up to isomorphism iff $\circledast f$ is monic, and then $\circledcirc f$ is monic as well; inside the copy of $\Set$, the subobjects of $\circledast T$ and of $\circledast\circledcirc T$ \emph{are} classified by pullback squares along $\mathrm{true}:1_E\to C_E(\{0,1\})$, one for the base and one for the heap. These two set-theoretic classifications describe the base inclusion and the heap inclusion, subject to compatibility of the latter with the former. They are not a proof that a subobject classifier is impossible; that conclusion follows from the noninvertible bijective Token map above. Abstracting can forget distinctions between carrier functions, as Example~\ref{ex:empty-heap} demonstrates.

For size bookkeeping choose Grothendieck universes $\mathcal V\in\mathcal W$ \cite[I.6]{maclane1998}. Use $\Set=\Set_{\mathcal W}$ throughout the construction and take $U=\mathcal V$, a $\mathcal W$-small set, as a common core for $\mathcal V$-small signatures. Then $F_U$ really is an endofunctor of $\Set_{\mathcal W}$; taking the same universe for both $\Set$ and its own entire core would not justify that assertion. A signature $E\in\mathcal V$ has its symbols in $U$ and is extended by $I:E\hookrightarrow U$:
\[
\langle B,E,H\rangle\longmapsto\langle B,U,(\id_B\sqcup I)^*(H)\rangle.
\]
Carrier and core tags remain distinct even if the underlying sets overlap. This extension is fully faithful for each fixed $E$. In $\TS_U$, a map between heaps with different used label sets is a carrier function preserving the recorded tuples while fixing all core labels. Renaming labels or adding fresh metavocabulary uses an explicit tagged core extension; the labels can be recoded into a sufficiently large common core. We write $\TS$ for the chosen common-core category. No set of literally all sets or all possible signatures is required.

\subsection{Token categories and the algebra \texorpdfstring{$\ATC$}{ATC}}\label{Token-categories}

Treat a concrete category as a pair $(\mathbf C,J)$ with a specified faithful $J:\mathbf C\to\Set$. A finite product of such pairs has the faithful tagged-union realization $U((A_i))=\bigsqcup_iJ_i(A_i)$. Write
\[
\SX^{\&}(\mathbf C,J)=\SX_{W\circ J}(\mathbf C).
\]
The underlying-set functor of this extension sends $(A,H)$ to $J(A)$ and retains the ambient map. These operations may be iterated. The algebra $\ATC$ of \emph{Token categories} is generated from $(\Set,\id)$ and identity set categories by finite products, such concrete realizations and $\SX^{\&}$. Objects retain their structural labels at every stage. Examples include
\[
\mathbf A,\quad\Set\times\mathbf B,\quad
\SX^{\&}(\Set\times\mathbf A),\quad
\SX^{\&}\bigl(\Set\times\mathbf B\times\SX^{\&}(\Set\times\mathbf A)\bigr),
\]
with the concrete functors understood. In the notation $\Tkn^{\&}_{\mathbf{TC}}=U(\mathbf{TC})$, the image is \emph{labelled}: its objects retain the original object together with its realized set. The notation $\Tkn_{\mathbf{TC}}=\SX^{\&}(\mathbf{TC})$ denotes the corresponding subsets extension.

\begin{theorem}\label{thm:atc}
Every category so constructed is concrete. A labelled realization is isomorphic to its indexing category and transports whatever limits, coproducts and exponentials that category has. For $\mathbf{TC}=\Set\times\mathbf E_1\times\cdots\times\mathbf E_n$,
\[
\Tkn_{\mathbf{TC}}\cong\TS_{E_1\sqcup\cdots\sqcup E_n};
\]
hence this elementary case has finite limits and is bicartesian closed, but is not a topos.
\end{theorem}
\begin{proof}
Tagged realization of a finite product is faithful, since restricting to each tag recovers the component arrow. The projection from a subsets extension is faithful, and its composite with the previous concrete functor is faithful. Induction proves concreteness. Retaining object labels and the represented arrow set gives the labelled realization its isomorphism with the indexing category. In the displayed elementary case the concrete functor is exactly word formation on a variable carrier and the fixed tagged core; Theorem~\ref{thm:tse} applies.
\end{proof}

A faithful functor to raw sets need not distinguish objects: different heaps can have the same underlying carrier. Removing their labels would invalidate the isomorphism statement. Likewise, a product with a category of Token classes combines underlying sorts; it does not automatically make those classes themselves data elements. To manipulate classes as data, choose encodings or identifiers and record their structure explicitly, as the machine descriptions in Section~\ref{sec:finite-representations} do. Whether a nested subsets extension has the desired closure properties depends on its concrete functor and remains a separate question.

\section{Algebraic Tokenization: Representing Categories of Structured Objects in Token Space}\label{representing-categories-of-structured-objects-in-Token-space}

\subsection{Structures, structure-preserving maps, and what the axioms do not do}

How are structures and structure-preserving maps defined in mathematics? The general pattern is the following.
\begin{itemize}
\item There are some sets $A,B,C,\ldots$ under consideration.
\item In each of them there are some special elements, and a structure-preserving map must map these elements correspondingly. In $\Set_{\&}$ the only special element of $A$ is its base point $\&_A$, and a map $A\to B$ must send $\&_A$ to $\&_B$; in algebras the special elements are the nullary operations.
\item On each of them there are some relations or operations. An operation is a relation connecting its arguments to its result: $+$ connects $x,y$ to $z$ when $x+y=z$. A structure-preserving map must preserve these relations, so that the same connections hold between the images.
\item There are some properties that the sets, special elements, relations and operations must satisfy. \textbf{For a fixed signature and already admissible objects, the homomorphism equations need not repeat the object axioms.} A semigroup $\langle S,\cdot\rangle$ satisfies $x\cdot(y\cdot z)=(x\cdot y)\cdot z$, but the definition of a semigroup homomorphism makes no reference to associativity.
\end{itemize}
The last point motivates the construction for the categories specified below. Their arrows preserve named constants, relations and operations, each of which can be recorded as a Token. Axioms still determine which objects belong to the category, and additional requirements on arrows must be accounted for separately.

\begin{definition}[Structure, signature]\label{def:structure}
A \emph{signature} is a pair $\langle\Theta,\Gamma\rangle$ of disjoint tagged symbol sets: constants in $\Theta$, and relation or operation symbols in $\Gamma$, with their kinds and finite arities specified. Symbols with the same printed name but different roles receive different tags. A \emph{structure} for it on a set $A$ interprets each $\theta\in\Theta$ as an element $\theta_A\in A$, each $m$-ary relation symbol $R$ as $R_A\subseteq A^m$, and each $m$-ary operation symbol $\omega$ as its graph
\[
\omega_A=\{(a_1,\ldots,a_m,\omega(a_1,\ldots,a_m)):a_i\in A\}\subseteq A^{m+1}.
\]
A \emph{structure-preserving map} $f:A\to B$ satisfies $f(\theta_A)=\theta_B$ for all $\theta\in\Theta$ and $(f(a_1),\ldots,f(a_k))\in\gamma_B$ whenever $(a_1,\ldots,a_k)\in\gamma_A$, for all $\gamma\in\Gamma$. A \emph{category of structured objects} under $\langle\Theta,\Gamma\rangle$ is a category of structures for that signature with all structure-preserving maps as morphisms. The signature is \emph{finitary} when every arity is finite; the number of symbols may be infinite, and nullary relations are allowed.
\end{definition}

\begin{definition}[Algebraic tokenization]\label{def:tokenization}
Let $\mathbf C$ be a category of structured objects under $\langle\Theta,\Gamma\rangle$ and let $E=\Theta\cup\Gamma$, one core symbol for each symbol of the signature. The \emph{algebraic tokenization} of $\mathbf C$ is the functor $\Rep_{\mathbf C}:\mathbf C\to\TS_E$ with
\begin{align*}
\Rep_{\mathbf C}(A)&=\langle A,E,\Upsilon_A\rangle,\\
\Upsilon_A&=\{(\theta,\theta_A):\theta\in\Theta\}\ \cup\
\{(\gamma,a_1,\ldots,a_k):(a_1,\ldots,a_k)\in\gamma_A,\ \gamma\in\Gamma\},
\end{align*}
and $\Rep_{\mathbf C}(f)=\langle f,\id_E,F_E(f)|_{\Upsilon_A}\rangle$. The class $\Rep_{\mathbf C}(A)$ is the \emph{tokenization} of the structured set $A$, the core $E$ is its \emph{signature}, and a Token $(\gamma,a_1,\ldots,a_k)\in\Upsilon_A$ is a \emph{structural Token}: it records the single fact $(a_1,\ldots,a_k)\in\gamma_A$, e.g.\ $x+y=z$ for the Token $(+,x,y,z)$. A nullary relation gives the core-only word $(\gamma)$. When $\mathbf C$ is clear we write $\Rep$.
\end{definition}

Algebraic tokenization turns every axiom-free part of a structure---its constants, relations and the graphs of its operations---into a heap of Tokens over the carrier and forgets the axioms. The axioms are not erased as properties: the graph of an operation still determines whether it is associative. For the chosen class of admissible objects the specified homomorphism conditions are preserved exactly; the axioms remain conditions on that class.

\begin{theorem}[Algebraic tokenization is a full embedding]\label{thm:tokenization}
Let $\mathbf C$ be a category of structured objects under $\langle\Theta,\Gamma\rangle$ and $E=\Theta\cup\Gamma$. Then $\Rep_{\mathbf C}:\mathbf C\to\TS_E$ is injective on objects and fully faithful: a function $f:A\to B$ is structure preserving iff it is the base of a Token map $\Rep_{\mathbf C}(A)\to\Rep_{\mathbf C}(B)$. Consequently $\mathbf C$ is isomorphic to a full subcategory of the Token Space.
\end{theorem}
\begin{proof}
The base of $\Rep(A)$ is $A$, and each $\theta_A$ and $\gamma_A$ is recovered from the heap as the set of Tokens headed by $\theta$, resp.\ $\gamma$; hence $\Rep$ is injective on objects. By \eqref{eq:hom} a Token map $\Rep(A)\to\Rep(B)$ is a function $f$ with $F_E(f)(\Upsilon_A)\subseteq\Upsilon_B$. Now $F_E(f)$ sends $(\theta,\theta_A)$ to $(\theta,f(\theta_A))$, which lies in $\Upsilon_B$ iff $f(\theta_A)=\theta_B$, and sends $(\gamma,a_1,\ldots,a_k)$ to $(\gamma,f(a_1),\ldots,f(a_k))$, which lies in $\Upsilon_B$ iff $(f(a_1),\ldots,f(a_k))\in\gamma_B$. So the inclusion holds iff $f$ is structure preserving, and $\Rep(f)$ is the unique Token map with base $f$. For the graph of an operation $\omega$ the condition reads
\[
f(\omega_A(a_1,\ldots,a_m))=\omega_B(f(a_1),\ldots,f(a_m)),
\]
the usual homomorphism condition.
\end{proof}

The hypothesis on morphisms matters. A category in which arrows are required to reflect relations, be injective, be continuous, or satisfy an additional analytic condition is not automatically covered by this theorem. Such conditions must be encoded or treated separately; topological spaces (Section~\ref{sec:top}) are the example treated below.

For one $k$-ary operation, the record condition can be read as the familiar square below. Algebraic tokenization gives the same condition by preserving its operation-headed Tokens, while retaining a common record language for other kinds of structure.
\begin{figure}[htbp]
\centering
\begin{tikzcd}[column sep=large,row sep=large]
A^k \arrow[r,"\omega_A"] \arrow[d,"f^k"'] & A \arrow[d,"f"] \\
B^k \arrow[r,"\omega_B"'] & B
\end{tikzcd}
\caption{An operation homomorphism as a square in $\Set$. Applying the operation and then $f$ agrees with transforming its operands first. For $k=0$, the left-hand objects are singleton sets and the same square expresses preservation of a constant.}\label{fig:operation-homomorphism}
\end{figure}

\begin{proposition}[Tokenization preserves products and equalizers]\label{prop:tokenization-limits}
If $\mathbf C$ is closed under componentwise binary products ($\theta_{A\times B}=(\theta_A,\theta_B)$, and $((a_1,b_1),\ldots,(a_k,b_k))\in\gamma_{A\times B}$ iff $(a_i)\in\gamma_A$ and $(b_i)\in\gamma_B$) and under induced substructures on equalizers, then
\[
\Rep(A\times B)=\Rep(A)\otimes\Rep(B)
\]
with $\otimes$ the product of $\TS_E$, and $\Rep$ preserves equalizers, hence pullbacks. It does \emph{not} preserve the terminal object, and need not preserve coproducts.
\end{proposition}
\begin{proof}
A Token of $\Rep(A)\otimes\Rep(B)$ is a tuple $(\langle x_1,y_1\rangle,\ldots,\langle x_m,y_m\rangle)$ with $(x_i)\in\Upsilon_A$, $(y_i)\in\Upsilon_B$ and $\langle x_i,y_i\rangle=x_i$ whenever $x_i=y_i\in E$. Since all Tokens of $\Upsilon_A$ and $\Upsilon_B$ start with a core symbol and contain no other core symbol, the two Tokens must be headed by the \emph{same} symbol, and the matched Token is $(\gamma,(a_1,b_1),\ldots,(a_k,b_k))$ with $(a_i)\in\gamma_A$, $(b_i)\in\gamma_B$, i.e.\ exactly a Token of $\Rep(A\times B)$; the projections agree. For $f,g:A\to B$ the equalizer in $\mathbf C$ is the induced substructure on $A'=\{a:f(a)=g(a)\}$, whose tokenization has heap $\Upsilon_A\cap F_E(A')$, the equalizer of Theorem~\ref{thm:tse}. The heap of $\Rep(1)$ contains only Tokens headed by a symbol of the signature, whereas the terminal object of $\TS_E$ has the full heap $F_E(\{*\})$. Coproducts of algebras may add elements generated by mixed expressions, unlike the disjoint carrier union in $\TS_E$.
\end{proof}

\begin{remark}[What tokenization does not preserve]\label{rem:not-preserved}
\begin{itemize}
\item \emph{Coproducts.} A coproduct of monoids can introduce mixed products and identify units, so its carrier is not generally the tagged union used by Token Space. Computing such an algebraic coproduct requires its own construction.
\item \emph{Subobjects.} The subclasses of $\Rep(A)$ are the pairs $S\subseteq A$, $\Upsilon'\subseteq\Upsilon_A\cap F_E(S)$. Induced substructures are the subclasses with $\Upsilon'=\Upsilon_A\cap F_E(S)$, and subalgebras additionally require $S$ to be closed under the operations. The Token Space therefore has many more subobjects than $\mathbf C$.
\item \emph{Exponentials.} The ambient object $\Rep(B)^{\Rep(A)}$
need not be the tokenization of an algebra of the chosen kind, even
though its global points are precisely the homomorphisms $A\to B$.
Example~\ref{ex:c2-relational-exponential} shows an operation relation
losing totality. An internal exponential, if one exists in the source
category, requires its own universal property; full faithfulness alone
does not identify it with the ambient construction.
\item \emph{Axioms.} Associativity, commutativity, etc.\ are properties of the heap, not part of the Token structure; two tokenizations are isomorphic in $\TS$ only if the original structures are (Theorem~\ref{thm:tokenization}), and selecting the tokenizations of semigroups still requires imposing associativity on the recorded operation graphs.
\end{itemize}
\end{remark}

\begin{remark}[Relation to relational structures]\label{rem:relational}
Replacing operations by their graphs is the classical passage from algebraic to relational signatures in model theory \cite{hodges1993}, and the resulting categories of structures are the classical setting of functorial semantics \cite{lawvere1963} and of concrete categories \cite{adamek1990}; on objects and morphisms algebraic tokenization is this passage. What the Token Space adds is (a) a \emph{single} ambient category for the chosen collection of signatures: the symbols are elements of the common core $U$, so structures of different signatures are compared through the renaming $\oslash$, the lifting $\Uparrow$ and the release $\Downarrow$ of Section~\ref{sec:operators}, and data may be moved between the carrier and the signature (Remark~\ref{rem:signature}); (b) the exponentials of Theorem~\ref{thm:tse}, whose carrier is all functions $A\to B$ and whose global elements are exactly the homomorphisms, with parameter objects selecting admissible families; (c) the tree Tokens of Section~\ref{interior-structure-mapping-and-tree-Token-classes}, which tokenize terms and nested expressions and are flattened back by reification.
\end{remark}

\subsection{Operation graphs and relational function objects}
\label{sec:operation-graphs-closure}

The use of operation graphs has a mathematical purpose. It puts data
relations and program operations into one record language while retaining
the exact homomorphism conditions. For an $m$-ary operation, the relation
\[
 R_\omega(a_1,\ldots,a_m,b)
 \quad\Longleftrightarrow\quad
 b=\omega(a_1,\ldots,a_m)
\]
records both the operands and their result. Applying one carrier map to
all $m+1$ entries yields exactly the homomorphism equation of
Theorem~\ref{thm:tokenization} and Figure~\ref{fig:operation-homomorphism}.
Recording only the possible results, or an operation name without its
operand--result tuples, would not generally retain that equation. Within
the present relational language this is the standard direct encoding of
an operation; no uniqueness claim among all possible encodings is needed.

Two uses of this relation must be distinguished. In $\Rep(A)$, each
operation relation is the graph of a total function satisfying the
axioms of the chosen algebra. In the ambient category $\TS_E$, a heap is
an arbitrary selection of records. Its exponential is defined by the
universal evaluation condition~\eqref{eq:exp-heap}; the resulting
relations need not again be graphs of total operations. Thus faithful
representation of an algebra and the existence of ambient function
objects are compatible without requiring the algebraic image to be
closed under those function objects. This is an instance of the
distinction between unrestricted relational structures and models of
additional axioms \cite{parker2022}.

\begin{example}[An ambient exponential leaves the group image]
\label{ex:c2-relational-exponential}
Let $A=B=C_2=\{0,1\}$, with addition modulo two, and use the common
group core $E=\{\mathtt{plus},\mathtt{neg},\mathtt{zero}\}$. The heap
of $T=\Rep(C_2)$ records addition, inverse and the identity element in
the format of Definition~\ref{def:tokenization}. The exponential
$Z=T^T$ has the four functions $C_2\to C_2$ as its carrier. For its
addition-shaped records, formula~\eqref{eq:exp-heap} specializes to
\begin{equation}\label{eq:c2-exponential-relation}
 (\mathtt{plus},f,g,h)\in H_Z
 \quad\Longleftrightarrow\quad
 h(x+y)=f(x)+g(y)\quad\text{for all }x,y\in C_2.
\end{equation}
The two operands are evaluated at independently chosen $x$ and $y$;
this is the condition required by joint evaluation on a categorical
product. It is not the pointwise formula $h(t)=f(t)+g(t)$.

Take $f=z$, the constant-zero map, and $g=i$, the identity map.
If a function $h$ made this a record, its equations at two input pairs
would give
\[
 \begin{array}{rcll}
 (x,y)=(0,0)&:&h(0)=0+0=0,\\
 (x,y)=(1,1)&:&h(0)=0+1=1.
 \end{array}
\]
These conditions contradict one another. No such $h$ exists, although
every carrier function is available as a candidate. Consequently the
addition-shaped relation of $Z$ is not the graph of a total binary
operation. In particular $Z$ is not isomorphic in $\TS_E$ to the
tokenization of a group with this signature: totality of this named
relation is preserved by a Token isomorphism.

Both $z$ and $i$ are group homomorphisms. By full faithfulness and
cartesian closure, they are precisely the two global elements of $Z$:
\[
 \Hom_{\TS_E}(1_E,Z)
 \cong\Hom_{\TS_E}(T,T)
 \cong\Hom_{\mathbf{Grp}}(C_2,C_2).
\]
Thus the failure remains even when the two operands are restricted to
global elements. Table~\ref{tab:c2-exponential-relation} makes this
restriction explicit. The inverse and identity records remain part of
the common signature; they do not remove the contradiction in the
addition-shaped relation.
\end{example}

\begin{table}[htbp]
\centering
\begin{tabular}{@{}ccc@{}}
\toprule
$f$ & $g$ & $h$ allowed by~\eqref{eq:c2-exponential-relation}\\
\midrule
$z$ & $z$ & $z$\\
$z$ & $i$ & none\\
$i$ & $z$ & none\\
$i$ & $i$ & $i$\\
\bottomrule
\end{tabular}
\caption{The addition-shaped relation on the two global points of
$\Rep(C_2)^{\Rep(C_2)}$. This table is a relational compatibility
condition, not pointwise addition of homomorphisms. The latter would
give $z+i=i$ and $i+i=z$, which do not satisfy the corresponding rows.}
\label{tab:c2-exponential-relation}
\end{table}

The pointwise group on the homomorphism set does exist in this example,
but it fails the required joint evaluation test. Let
$M=\{z,i\}=\Hom_{\mathbf{Grp}}(C_2,C_2)$ with pointwise addition and
let $e:M\times C_2\to C_2$ be evaluation. In the product group,
\[
 e(i,0)+e(z,1)=0+0=0,
 \qquad
 e\bigl((i,0)+(z,1)\bigr)=e(i,1)=1.
\]
Thus $e$ is not a group homomorphism. Having a group of homomorphisms
does not by itself give an exponential for cartesian products: its
evaluation must also be a morphism in the same category.

Example~\ref{ex:c2-relational-exponential} proves that the \emph{ambient}
exponential leaves the algebraic image. That fact alone would not rule
out a different exponential internal to a subcategory. For groups,
nonexistence follows independently from the following short argument.

\begin{proposition}[A zero-object obstruction to cartesian closure]
\label{prop:zero-object-ccc}
If a cartesian closed category has a zero object, every hom-set has
exactly one element. In particular, the category of groups and the
category of vector spaces over a field are not cartesian closed.
\end{proposition}
\begin{proof}
A zero object is both initial and terminal. Taking it as the terminal
object $1$, the exponential adjunction gives
\[
 \Hom(A,B)\cong\Hom(1\times A,B)\cong\Hom(1,B^A).
\]
The last set is a singleton because $1$ is also initial. Groups have
the trivial group as a zero object, but $C_2\to C_2$ has the two maps
$z$ and $i$. Vector spaces have the zero space as a zero object, but a
one-dimensional space has distinct zero and identity linear maps.
Both conclusions therefore follow by contradiction.
\end{proof}

This obstruction uses cartesian products. It does not prohibit other
closed monoidal structures, and it does not say that all algebraic
categories fail to be cartesian closed. Sets with one specified unary
operation, for example, form the functor category for the one-object
category of the free monoid on one generator; that is a presheaf
category after taking the opposite indexing category, hence a topos
\cite{maclane1998}.
Additional algebraic conditions must therefore be assessed in their
particular categories.

There is no conflict with the full embedding of groups into $\TS_E$.
The image of the trivial group is terminal within the group image but
is not the ambient terminal object $1_E$, as already noted in
Proposition~\ref{prop:tokenization-limits}. The latter has every possible
record, whereas $\Rep(1_{\mathbf{Grp}})$ has only the group-signature
records. Full faithfulness neither identifies these two objects nor
requires their ambient exponentials to remain in the image.

This explains the design choice relevant to program--data unification.
Operation graphs retain the structure and maps of concrete algebras;
the broader relational category supplies objects for structural
functions and their families. These function objects can then be
treated as structured data, with global points selecting admissible
maps. A finite program representation and an executable evaluator
remain additional constructions, supplied for the fragments of
Section~\ref{sec:finite-representations}; cartesian closure alone does
not turn every function object into an executable graph. The
mathematical companion \cite{pancausal} studies which further
dependency constraints preserve function objects and which do not.

\subsection{Tokenizing familiar categories}\label{sec:tokenization-examples}

In each example $\Rep_{\mathbf C}(f)$ is the Token map with base $f$; by Theorem~\ref{thm:tokenization} it exists exactly when $f$ is a morphism of $\mathbf C$.

\textbf{Sets.} A set $A$ has the empty structure $\langle\varnothing,\varnothing\rangle$ and $\Rep_{\Set}(A)=\langle A,\varnothing,\varnothing\rangle$: every function is a Token map.

\textbf{Pointed sets.} $\Rep_{\Set_{\&}}(A)=\langle A,\{\&\},\{(\&,\&_A)\}\rangle$.

\textbf{Ordered sets.} $\Rep_{\mathbf{Ords}}\langle P,\le\rangle=\langle P,\{\le\},\{(\le,x,y):x\le y\}\rangle$; Token maps are the monotone maps.

\textbf{Graphs.} $\Rep_{\mathbf{Gph}}\langle V,Edg\rangle=\langle V,\{Edg\},\{(Edg,x,y):(x,y)\in Edg\}\rangle$; Token maps are graph homomorphisms. Multigraphs need an edge sort with incidence records, so that distinct edge occurrences are not collapsed into a set of vertex pairs.

\textbf{Rings.} $\Rep_{\mathbf{Rng}}\langle R,+,\cdot,-,0\rangle=\langle R,\{+,\cdot,-,0\},\Upsilon\rangle$ with
\[
\Upsilon=\{(+,x,y,z):x+y=z\}\cup\{(\cdot,x,y,z):x\cdot y=z\}\cup\{(-,x,y):-x=y\}\cup\{(0,0_R)\};
\]
if unital homomorphisms are intended, add the constant $1$.

\textbf{Vector spaces.} For $\mathbf V=\langle V,0,+,\mathbb R\rangle$ take the signature $\langle\{0\},\{+\}\cup\mathbb R\rangle$ with one unary operation symbol per scalar, with constant, addition and scalar-operation tags understood:
\[
\Rep_{\mathbf{Vct}}(\mathbf V)=\bigl\langle V,\ \{0,+\}\cup\mathbb R,\ \{(+,x,y,z):x+y=z\}\cup\{(a,x,y):ax=y,\ a\in\mathbb R\}\cup\{(0,0_V)\}\bigr\rangle.
\]
Token maps are exactly the linear maps. Alternatively the scalars may be placed in the carrier: regarding $\mathbb R$ as the ring $\langle\mathbb R,+_{\mathbb R},\cdot_{\mathbb R},-_{\mathbb R},0_{\mathbb R},1_{\mathbb R}\rangle$ and $\mathbf V$ as the two-sorted set $V\sqcup\mathbb R$ with unary sort predicates $V$ and $\mathbb R$,
\[
\Rep'_{\mathbf{Vct}}(\mathbf V)=\langle V\sqcup\mathbb R,E',\Upsilon'\rangle,\qquad
E'=\{0_V,0_{\mathbb R},1_{\mathbb R},+_V,\cdot_V,+_{\mathbb R},\cdot_{\mathbb R},-_{\mathbb R},V,\mathbb R\},
\]
where $\Upsilon'$ consists of $(+_V,x,y,z)$ for $x+y=z$, $(\cdot_V,a,x,y)$ for $ax=y$, the ring records $(+_{\mathbb R},a,b,c)$, $(\cdot_{\mathbb R},a,b,c)$, $(-_{\mathbb R},a,b)$, the constants $(0_V,0_{\mathbf V})$, $(0_{\mathbb R},0)$, $(1_{\mathbb R},1)$, and the sort records $(V,x)$, $(\mathbb R,a)$. A linear map $f$ is represented by the Token map with base $f\sqcup\id_{\mathbb R}$, which fixes every scalar.

\paragraph{Small categories as structured data.}
Fix a universe for small categories and a larger ambient universe when collecting them into $\mathbf{Cat}$. Let $O_{\mathbf C}$ and $M_{\mathbf C}$ be tagged copies of the objects and morphisms of a small category $\mathbf C$. The tags distinguish an object from its identity arrow even when the original notation uses related names. Use the fixed core
\[
E_{\mathrm{cat}}=\{\mathrm{Ob},\mathrm{Mor},\mathrm{dom},\mathrm{cod},\mathrm{id},\mathrm{comp}\}.
\]
Write the domain and codomain of $f$ as $d(f),c(f)$. On $O_{\mathbf C}\sqcup M_{\mathbf C}$, take the heap
\begin{align*}
\Upsilon_{\mathbf C}={}&\{(\mathrm{Ob},x):x\in O_{\mathbf C}\}
\cup\{(\mathrm{Mor},f):f\in M_{\mathbf C}\}\\
&\cup\{(\mathrm{dom},f,d(f)),(\mathrm{cod},f,c(f)):f\in M_{\mathbf C}\}\\
&\cup\{(\mathrm{id},x,1_x):x\in O_{\mathbf C}\}\\
&\cup\{(\mathrm{comp},g,f,g\circ f):c(f)=d(g)\}.
\end{align*}
All these are relations on one tagged carrier; in particular, the domain, codomain and identity records are graphs on their specified sorts, not total unary operations on the union. Category axioms are conditions on the recorded data. An arbitrary heap over $E_{\mathrm{cat}}$ need not satisfy them. Define
\[
\Rep_{\mathbf{Cat}}(\mathbf C)
=\langle O_{\mathbf C}\sqcup M_{\mathbf C},E_{\mathrm{cat}},\Upsilon_{\mathbf C}\rangle.
\]

\begin{proposition}[Categories are Token classes]\label{prop:cat-tokenization}
The functor $\Rep_{\mathbf{Cat}}:\mathbf{Cat}\to\TS_{E_{\mathrm{cat}}}$ is injective on objects and fully faithful: a function on the tagged carriers is a Token map between the encodings iff its restrictions are the object and arrow maps of a functor. For functors $F,G:\mathbf C\to\mathbf D$, natural transformations $F\Rightarrow G$ correspond bijectively to Token maps
\[
\widehat\alpha:\Rep_{\mathbf{Cat}}(\mathbf C)\longrightarrow
\Rep_{\mathbf{Cat}}(\mathbf D^{\to})
\]
whose two endpoint composites are $\Rep_{\mathbf{Cat}}(F)$ and $\Rep_{\mathbf{Cat}}(G)$.
\end{proposition}
\begin{proof}
The heap recovers both sorts and all the category data. Sort preservation sends objects to objects and arrows to arrows. The domain, codomain and identity records enforce the corresponding functor equations, and the composition records enforce $F(g\circ f)=F(g)\circ F(f)$ for every composable pair. Conversely a functor preserves every record. This is the relational instance of Theorem~\ref{thm:tokenization}; identities and composition are respected on both sorts.

An object of $\mathbf D^{\to}$ is an arrow of $\mathbf D$, and an arrow from $u:a\to b$ to $v:a'\to b'$ is a pair $(r:a\to a',s:b\to b')$ with $s\circ u=v\circ r$. A natural transformation $\alpha:F\Rightarrow G$ therefore defines a functor $\mathbf C\to\mathbf D^{\to}$ by $x\mapsto\alpha_x$ and $f\mapsto(Ff,Gf)$. Conversely a functor with these endpoints gives exactly these components and their naturality equations. Full faithfulness translates this bijection to Token maps.
\end{proof}

\begin{figure}[htbp]
\centering
\begin{tikzcd}[column sep=large,row sep=large]
\Rep(\mathbf C) \arrow[r,"\widehat\alpha"] \arrow[dr,"{\langle\Rep(F),\Rep(G)\rangle}"'] &
\Rep(\mathbf D^{\to}) \arrow[d,"{\langle\Rep(d),\Rep(c)\rangle}"]\\
& \Rep(\mathbf D)\otimes\Rep(\mathbf D)
\end{tikzcd}
\caption{A natural transformation as a lift with fixed endpoints in $\TS_{E_{\mathrm{cat}}}$. Here $d,c:\mathbf D^{\to}\to\mathbf D$ are the endpoint functors, and $\otimes$ is the Token Space product. An endpoint pair alone does not specify a lift; its existence and choice encode the natural transformation.}\label{fig:category-lift}
\end{figure}

\begin{remark}[Representing functors and constructing a functor category]\label{rem:cat-exponential}
Cartesian closedness gives the useful bijection
\[
\operatorname{Fun}(\mathbf C,\mathbf D)
\cong\Hom_{\TS_{E_{\mathrm{cat}}}}\bigl(1_{E_{\mathrm{cat}}},\Rep(\mathbf D)^{\Rep(\mathbf C)}\bigr).
\]
It identifies functors with global elements of an exponential. The exponential is generally \emph{not} the tokenization of the functor category $[\mathbf C,\mathbf D]$, whose arrows are natural transformations. For the one-object, one-arrow category $\mathbf 1$, the carrier of $\Rep(\mathbf 1)$ has two tagged elements. Hence $\Rep(\mathbf 1)^{\Rep(\mathbf 1)}$ has four carrier elements, whereas $\Rep([\mathbf 1,\mathbf 1])$ has two. They cannot be isomorphic in Token Space. Full faithfulness alone also does not identify the encoded categorical product with the ambient product: the latter can contain pairs of different sorts. To construct a category of functors as data, construct its object, arrow and composition records explicitly, as follows.
\end{remark}

\begin{proposition}[Executable finite category constructions]\label{prop:finite-cat-constructions}
Given valid finite category heaps with effective labels, there are terminating algorithms producing valid finite heaps for their categorical products, slices over a specified object, arrow categories and functor categories. In the functor-category construction, the object records enumerate exactly the functors and the arrow records exactly the natural transformations, with identity and vertical composition computed componentwise. The input and output category encodings lie in $\TS^{\mathrm{fin}}_{E_{\mathrm{cat}},\le4}$.
\end{proposition}
\begin{proof}
Decode the finite sorts and the domain, codomain, identity and partial composition tables. Validity itself is decidable for these finite tables: check typing, unique domains, codomains and identities, one composite for each composable pair, unit laws and associativity on every composable triple. Products use pairs in each sort and componentwise operations. In $\mathbf C/x$, objects are arrows $u:a\to x$ and arrows $u\to v:b\to x$ are triples $(u,v,h)$ with $h:a\to b$ and $v\circ h=u$; composition uses that of $\mathbf C$. The arrow category uses the commutative squares described in the preceding proof. These are finite enumerations, and their category laws follow from those of the inputs.

For $[\mathbf C,\mathbf D]$, enumerate all pairs of functions
\[
F_0:O_{\mathbf C}\to O_{\mathbf D},\qquad
F_1:M_{\mathbf C}\to M_{\mathbf D}
\]
and retain exactly those preserving domain, codomain, identity and composition. These are precisely the functors by Proposition~\ref{prop:cat-tokenization}. For every retained pair $F,G$, enumerate families $\alpha_x:Fx\to Gx$ and retain exactly those satisfying
\[
G(f)\circ\alpha_x=\alpha_y\circ F(f)
\qquad(f:x\to y\text{ in }\mathbf C).
\]
Record each arrow with its endpoints $(F,G,\alpha)$, including the empty-family cases. Identity families and the pointwise composites $(\beta\alpha)_x=\beta_x\circ\alpha_x$ satisfy naturality; the category laws follow componentwise from $\mathbf D$. Every functor and natural transformation is considered, so the output is exactly the functor category, up to the chosen finite identifiers. All enumerations terminate. Finally, every category record has length at most four over the six-symbol core, proving the assertion about the bounded Token Space.
\end{proof}

These are classical categorical constructions expressed as algorithms on the chosen records \cite{rydeheard1988}. Their role here is to give Thesis~\ref{thesis:executable} an actual input--output construction within the same language used for computation graphs. After constructing $[\mathbf C,\mathbf D]$, the evaluation functor $[\mathbf C,\mathbf D]\times\mathbf C\to\mathbf D$ can likewise be encoded: it sends $(F,x)$ to $Fx$ and $(\alpha:F\Rightarrow G,f:x\to y)$ to $G(f)\circ\alpha_x=\alpha_y\circ F(f)$. The product here is first formed in $\mathbf{Cat}$, then tokenized. A size bound on the encoded input categories gives a finite input domain and a finite set of possible encoded outputs, hence an instance of the represented finite mappings of Section~\ref{sec:finite-representations}. Without a fixed bound, the algorithms operate on an unbounded family of finite descriptions. Exhaustive enumeration establishes executability, not an efficient complexity bound.

Thus category theory can be presented as data on which this framework operates: categories have records, functors preserve them, and natural transformations have endpoint-constrained records of their own. The representation theorem applies to small categories, including infinite ones; the algorithms just proved apply to finite tables. Larger effectively presented categories require algorithms appropriate to those presentations. Collecting a category of structures and encoding it at a further level also requires that it be small in the universe used at that level. This is an explicit programme of categorical self-description, with the executable fragment established above; it does not assert a universe-free self-containing category or identify the encoding with a theory of internal categories in $\TS$.
\begin{remark}[Uncountable versus finite signatures]\label{rem:signature}
The two tokenizations of vector spaces illustrate a design choice that arises for every structure with a parameter set: a scalar $a$ can be a \emph{core symbol}, as in $\Rep_{\mathbf{Vct}}$, or an \emph{element of the base}, as in $\Rep'_{\mathbf{Vct}}$. The two choices differ in three respects.
\begin{itemize}
\item \emph{Size of the signature.} $\Rep_{\mathbf{Vct}}$ has the uncountable signature $\{0,+\}\cup\mathbb R$: one unary operation symbol for every scalar, exactly as in the presentation of $\mathbb R$-modules as algebras with unary operations. $\Rep'_{\mathbf{Vct}}$ has the finite signature $E'$ of ten symbols, at the price of a larger base $V\sqcup\mathbb R$. Theorem~\ref{thm:tokenization} does not care about the size of the signature (finitary refers to arities, not to the number of symbols); a finite signature can simplify comparison and permits the scalar sort to be studied as data. The second encoding records the scalar operations; their field axioms remain properties of the records, rather than individual Tokens. The bounded category $\TS^{\mathrm{fin}}_{E,\le m}$ of Section~\ref{sec:bounded} requires a finite signature.
\item \emph{Rigidity.} Core symbols are fixed by every Token map, base elements are not. Hence $\Rep_{\mathbf{Vct}}$ can only see maps that fix each scalar, whereas a Token map $\Rep'_{\mathbf{Vct}}(\mathbf V_1)\to\Rep'_{\mathbf{Vct}}(\mathbf V_2)$ is a pair $(\sigma,f)$ with $\sigma:\mathbb R\to\mathbb R$ a ring endomorphism preserving $0,1$ and $f(ax)=\sigma(a)f(x)$, i.e.\ a \emph{semilinear} map. Since the only ring endomorphism of $\mathbb R$ is the identity, $\Rep'_{\mathbf{Vct}}$ is still fully faithful; over a field $k$ with nontrivial endomorphisms, e.g.\ $k=\mathbb C$, the same construction embeds $k$-$\mathbf{Vct}$ faithfully but \emph{not} fully, the full ambient subcategory on those objects having the semilinear maps, while the linear maps form the intended smaller arrow set. Which symmetries count as morphisms is thus decided by what is put into the signature.
\item \emph{Moving between the two.} Lifting and release change whether scalar occurrences are fixed labels or variable data. They do not alone convert the action record $(\cdot_V,a,x,y)$ into $(a,x,y)$: that requires an explicit change of record format and treatment of the scalar and sort records. Releasing scalar labels without adding their ring-operation graphs does not force the scalar map to be an endomorphism, and hence does not by itself produce semilinear hom sets. These observations distinguish moving data between base and core from changing a mathematical signature.
\end{itemize}
The same trade-off appears for graphs with edge labels, for weighted relations such as the attention classes of Section~\ref{sec:attention-records}, and for data sets whose feature values are real numbers: labels and values in the core give rigid, exactly comparable Tokens over an infinite signature; labels and values in the base give a finite signature and allow Token maps to transform them.
\end{remark}

This choice is relevant to neural networks. Linear projections are maps between tokenized vector spaces. A complete attention block is not linear and is not a Token map between those same vector-space structures; to describe a nonlinear block as a Token map we use a different heap on its sequence carrier, for example the permutation heap of Section~\ref{sec:equivariance} or the causal heap of Section~\ref{sec:causal}. A computational carrier may support several classes, and the heap determines the claim made about a map.

\subsection{Topological spaces: faithful but not full}\label{sec:top}

Continuity can be expressed by preservation of open subsets under inverse image. We use this contravariant construction, rather than assert that the point carrier alone with the preceding finitary encoding captures continuity. For a space $X$ take the carrier $\mathcal P(X)$, record its Boolean operations $\varnothing,X,\cap,\cup,{}^c$ by their graphs, and add $(\mathrm{open},S)$ for each open subset $S$; call this class $\Rep_{\mathbf{Top}^{op}}(X)$. For a continuous $f:X\to Y$ define the carrier map
\[
f^{op}=f^{-1}:\mathcal P(Y)\longrightarrow\mathcal P(X).
\]
It preserves the Boolean records and the open-set records, $(\id_X)^{-1}=\id_{\mathcal P(X)}$ and $(g\circ f)^{-1}=f^{-1}\circ g^{-1}$, so $X\mapsto\Rep_{\mathbf{Top}^{op}}(X)$, $f\mapsto f^{op}$ is a functor $\mathbf{Top}^{op}\to\TS$; it is the reversing operator $\rightleftharpoons$ of Section~\ref{sec:operators} applied to the underlying sets, with the open sets added as a unary relation. It is faithful: if $f(x)\ne g(x)$, the inverse images of $\{f(x)\}$ distinguish $f^{op}$ and $g^{op}$.

\begin{remark}[Topological spaces are embedded faithfully only]\label{rem:top-nonfull}
Unlike the previous examples, $\mathbf{Top}^{op}$ is not a category of structured objects in the sense of Definition~\ref{def:structure}, and $\Rep_{\mathbf{Top}^{op}}$ is faithful but \emph{not} full. Under the axiom of choice a nonprincipal ultrafilter on an infinite discrete set $Y$ gives a Boolean homomorphism $\mathcal P(Y)\to\{0,1\}=\mathcal P(\{*\})$; all subsets of both spaces are open, so it preserves every record, but it is not evaluation at a point of $Y$ and hence not $f^{op}$ for any $f:\{*\}\to Y$. Fullness would require recording that $f^{-1}$ is a \emph{complete} Boolean algebra homomorphism (such homomorphisms are exactly the preimage maps), and arbitrary unions are not among the finitary operations recorded in this particular encoding. This does not rule out other encodings with additional sorts and records. Theorem~\ref{thm:tokenization} therefore applies to the finitary structures above, while topological spaces are embedded faithfully only.
\end{remark}

\section{Exploring Structure Relations of Token Classes: The Operator Calculus}\label{sec:operators}\label{exploring-structure-relations-of-Token-classes}

We need operators to explore the relations between structured objects inside the Token Space: to compare, combine, relabel and transport Token classes, and, in Section~\ref{sec:dataset-examples}, to turn a data set into a family of structured objects. All results are unique up to isomorphism. The operators act on explicitly chosen encodings and do not automatically preserve the axioms of orders, algebras or probabilities. Table~\ref{tab:operators} collects them; throughout, $T=\langle B,E,\Upsilon\rangle$, $T_i=\langle B_i,E,\Upsilon_i\rangle$, and $f=\langle f,\id_E,p\rangle:T_1\to T_2$ is a Token map.

\begin{table}[htbp]
\centering\small\renewcommand{\arraystretch}{1.13}
\begin{tabularx}{\linewidth}{@{}p{0.26\linewidth}p{0.13\linewidth}X@{}}
\toprule
Name & Symbol & Effect \\
\midrule
Merging / meeting / deleting & $\Cup,\ \Cap,\ \ominus$ & Union, intersection, difference of heaps (and bases). \\
Forgetting / stuffing / abstracting & $\varnothing,\ \circledast,\ \circledcirc$ & Empty the heap; fill the heap; make the heap the base. \\
Matchup (product) & $\otimes$ & Carrier $B_1\times B_2$; a Token iff \emph{both} projections are Tokens. \\
Blending (extending product) & $\oplus$ & Carrier $B_1\times B_2$; a Token iff \emph{some} projection is a Token. \\
Union / Y-union & $\uplus$ & Tagged union of bases and heaps. \\
Referring (exponential) & $\vartriangleleft$ & Carrier $B_1^{B_2}$; heap $\Upsilon_1\vartriangleleft\Upsilon_2$ of \eqref{eq:exp-heap}. \\
Pullback / direct image & $f^*,\ f_!$ & Largest heap making $f$ a map into $T$; smallest heap making $f$ a map out of $T$. \\
Renaming structure & $\oslash_\alpha$ & Rename core symbols along $\alpha:E\to E'$. \\
Lifting / absolute lifting & $\Uparrow,\ \upuparrows$ & Move $S\subseteq B$ into the core; name $S$ in the core while keeping $B$. \\
Release / obscuring & $\Downarrow,\ \leadsto$ & Move $E'\subseteq E$ into the base; release the whole core. \\
Introducing unknown & $\odot$ & Add a placeholder $e$ with the record $(\varepsilon,e)$. \\
Powering & $\nabla$ & Add the subsets of $B$ with membership records. \\
Reversing & $\rightleftharpoons$ & Send $f$ to $f^{-1}$ between powerset classes (contravariant). \\
Generalizing & $\rtimes$ & Replace the elements of a subclass by sets containing them. \\
\bottomrule
\end{tabularx}
\caption{The operator calculus on Token classes.}\label{tab:operators}
\end{table}

\subsection{Combining, removing, and transporting records}

\textbf{Basic set-theoretic operations.} For classes over a common core whose carriers lie in a common tagged set,
\begin{align*}
T_1\mathbin{\Cup}T_2&=\langle B_1\cup B_2,E,\Upsilon_1\cup\Upsilon_2\rangle &&\text{(\emph{merging})},\\
T_1\mathbin{\Cap}T_2&=\langle B_1\cap B_2,E,\Upsilon_1\cap\Upsilon_2\rangle &&\text{(\emph{meeting})},\\
T_2\mathbin{\ominus}T_1&=\langle B_2,E,\Upsilon_2\setminus\Upsilon_1\rangle &&\text{(\emph{deleting})}.
\end{align*}
For a fixed carrier these are operations on the lattice of heaps, and an inclusion of heaps gives an identity-carrier Token map in the same direction. Heap union on a shared carrier is different from the coproduct, whose carrier is a disjoint union. An unrecorded tuple is not thereby a false claim; reading absence as falsity requires a closed-world convention.

\textbf{Binary operations.} \emph{Matchup} $\otimes$ is the product, \emph{union} $\uplus$ the coproduct and \emph{referring} $\vartriangleleft$ the exponential of \eqref{eq:product-heap}--\eqref{eq:exp-heap}: a tuple $(\langle x_1,y_1\rangle,\ldots,\langle x_n,y_n\rangle)$ over $B_1\times B_2$, with $\langle x_i,y_i\rangle$ nominally denoting $x_i$ when $x_i=y_i\in E$, is a Token of $T_1\otimes T_2$ iff $(x_i)\in\Upsilon_1$ \emph{and} $(y_i)\in\Upsilon_2$. \emph{Blending} $T_1\oplus T_2$, the \emph{extending product}, has the same carrier and admits the tuple iff $(x_i)\in\Upsilon_1$ \emph{or} $(y_i)\in\Upsilon_2$, i.e.\ its heap is $F_E(\pi_1)^{-1}(\Upsilon_1)\cup F_E(\pi_2)^{-1}(\Upsilon_2)$; it is not the categorical product. Both $\otimes$ and $\oplus$ are associative and commutative up to the carrier bijections, their heap conditions being conjunction and disjunction of the component conditions.

\begin{example}\label{ex:order-product}
Given ordered sets $\mathbf P=\langle P,\ge\rangle$ and $\mathbf P'=\langle P',\ge\rangle$ with tokenizations $T_{\mathbf P}=\langle P,\{\ge\},\{(\ge,x,y):x\ge y\}\rangle$ and $T_{\mathbf P'}$, the matchup $T_{\mathbf P}\otimes T_{\mathbf P'}$ is the tokenization of the product order $\mathbf P\times\mathbf P'$: $(\ge,\langle x_1,x_2\rangle,\langle y_1,y_2\rangle)$ is a Token iff $x_1\ge y_1$ and $x_2\ge y_2$ (Proposition~\ref{prop:tokenization-limits}). The blending $T_{\mathbf P}\oplus T_{\mathbf P'}$ instead relates two pairs as soon as one coordinate does.
\end{example}

\textbf{Pullback and direct image.} For a function $f:A\to B$ and $T_B=\langle B,E,\Upsilon_B\rangle$ the \emph{pullback heap} and the \emph{induced class} are
\begin{equation}\label{eq:pullback}
f^*\Upsilon_B=F_E(f)^{-1}(\Upsilon_B),\qquad
f^*T_B=\langle A,E,f^*\Upsilon_B\rangle,
\end{equation}
and for a heap $\Upsilon_A$ on $A$ the \emph{direct-image heap} is $f_!\Upsilon_A=F_E(f)(\Upsilon_A)$. They satisfy
\begin{equation}\label{eq:heap-adjunction}
f_!\Upsilon_A\subseteq\Upsilon_B\quad\Longleftrightarrow\quad \Upsilon_A\subseteq f^*\Upsilon_B,
\end{equation}
so the induced class is the largest heap on $A$ making $f$ a Token map into $T_B$, and the direct image is the smallest heap on $B$ making $f$ a Token map out of $T_A$; $f_!(T_1)$ denotes $\langle f(B_1),E,f_!\Upsilon_1\rangle$.

\subsection{Changing labels and moving values between base and core}

\begin{itemize}
\item \textbf{Renaming structure} $\oslash$: for an injection $\alpha:E\to E'$,
\[
\oslash(T,\alpha)=\oslash_\alpha T=\langle B,E',(\id_B\sqcup\alpha)^*(\Upsilon)\rangle .
\]
It leaves carrier maps unchanged and is fully faithful; it is a change of core, not a Token map within a fixed $\TS_E$. Noninjective label maps can also push heaps forward, but may identify different constraints and need not give a full embedding.
\item \textbf{Lifting} $\Uparrow$: for $S\subseteq B$, the elements of $S$ are moved from the base to the core while the Tokens are kept,
\[
\Uparrow(T,S)=\bigl\langle B\setminus S,\ E\sqcup S,\ (\id_{B\setminus S}\sqcup\id_E\sqcup inc_{S\to E\sqcup S})^*(\Upsilon)\bigr\rangle,
\]
using the retagging bijection $(B\setminus S)\sqcup(E\sqcup S)\cong B\sqcup E$. All lifted elements are subsequently fixed by Token maps.
\item \textbf{Absolute lifting} $\upuparrows$: retain $B$ and name the elements of $S$ in the core,
\[
\upuparrows(T,S)=\bigl\langle B,\ E\sqcup S,\ \Upsilon\cup\{(inc_{S\to E\sqcup S}(x),x):x\in S\}\bigr\rangle .
\]
\item \textbf{Release} $\Downarrow$: for $E'\subseteq E$, $\Downarrow(T,E')=\langle B\sqcup E',E\setminus E',\Upsilon\rangle$ moves chosen core elements into a tagged carrier summand.
\item \textbf{Obscuring} $\leadsto$: $\leadsto(T)=\langle B\sqcup E,\varnothing,\Upsilon\rangle$, $\leadsto(f)=f\sqcup\id_E$, releases the entire core.
\end{itemize}
Moving a specified summand into the core and releasing that same summand undo each other on objects up to retagging. Absolute lifting instead adds naming records. Different operation-record formats, such as the two vector-space encodings in Remark~\ref{rem:signature}, require additional explicit transformations. Treating an attention weight as a core label requires exact preservation of that value, whereas treating it as a carrier element permits its transformation subject to the relations recorded on that sort. These operations are always defined on objects; their action on a map requires the map to respect the stated choices.

\subsection{Unknowns, powering, reversing, and generalizing}

\begin{itemize}
\item \textbf{Introducing unknown} $\odot$: for fresh $\varepsilon,e\notin B\cup E$,
\[
\odot(T)=\langle B\cup\{e\},E\cup\{\varepsilon\},\Upsilon\cup\{(\varepsilon,e)\}\rangle,\qquad
\odot(f)=f\cup\{e\mapsto e\}.
\]
We say that $e$ is an \emph{unknown element} introduced in $\odot T$; it is a distinguished placeholder, not a theory of uncertainty.
\item \textbf{Powering} $\nabla$: with a fresh core label $\gamma$ (membership) and the tagged carrier $B\sqcup\mathcal P(B)$,
\begin{align*}
\nabla(T)&=\bigl\langle B\sqcup\mathcal P(B),\ E\cup\{\gamma\},\ \Upsilon_\gamma\cup\Upsilon_\nabla\bigr\rangle,\\
\Upsilon_\gamma&=\{(\gamma,S,x):x\in S\},\\
\mathcal R(x)&=\{\iota_B(x)\}\cup\{\iota_{\mathcal P(B)}(S):x\in S\subseteq B\},\\
\Upsilon_\nabla&=\{r[z_i/x_i]:r\in\Upsilon,\ z_i\in\mathcal R(x_i)\text{ at each carrier position}\},
\end{align*}
where replacement ranges only over carrier positions of $r$ and every core position is unchanged. Thus each carrier entry may stay in the original sort or be replaced by a subset containing it; $\nabla(f)=f\sqcup\mathcal P(f)$ with $\mathcal P(f)(S)=f[S]$ preserves membership and lifted records, so $\nabla$ is a functor between the extended-core categories. Tags are essential: a carrier element that is itself a set is not a member of the powerset sort. $\odot\nabla$ adds the unknown $e$ of the earlier presentation.
\item \textbf{Reversing} $\rightleftharpoons$: let $\wp(B)=\langle\mathcal P(B),E_{\mathrm{Bool}},\Upsilon_{\mathrm{Bool}}\rangle$ be the tokenized Boolean algebra of subsets of $B$ (graphs of $\varnothing,B,\cap,\cup,{}^c$). For any function $f:B_1\to B_2$,
\[
\rightleftharpoons(f)=f^{-1}:\wp(B_2)\longrightarrow\wp(B_1)
\]
is a Token map, $\rightleftharpoons(\id)=\id$ and $\rightleftharpoons(g\circ f)=\rightleftharpoons(f)\circ\rightleftharpoons(g)$; reversing is a contravariant functor $\Set^{op}\to\TS$ which is faithful by Section~\ref{sec:top}. It is defined on the powerset sort only: sending the individual points of $B_2$ to a common placeholder would not preserve identities, which is why $\Rep_{\mathbf{Top}^{op}}$ uses $\wp(X)$ rather than $\nabla\varnothing\langle X,\varnothing,\varnothing\rangle$ as its carrier.
\item \textbf{Generalizing} $\rtimes$: let $T'\Subset T$, extend their core to $E\cup\{\gamma\}$ without adding records, and let $f:T\to\nabla T$ be a Token map whose base sends each $x\in B$ to $x$ or to a subset containing $x$. Then
\[
\rtimes(f,T')=f_!(T')=\bigl\langle f(B'),E\cup\{\gamma\},F_{E\cup\{\gamma\}}(f)(\Upsilon')\bigr\rangle
\]
is the \emph{generalization} of $T'$ along $f$: its elements are replaced by sets (bins, intervals, clusters) containing them and its Tokens are transported. If the base of $f$ does not satisfy the condition, $\rtimes(f,T')=f_!(\varnothing T')$, and only the carrier is generalized.
\end{itemize}

\subsection{Data sets as structured objects}\label{sec:dataset-examples}

The operators turn a data set into a family of Token classes that record the structure a learning algorithm is asked to respect. We consider the $n$-dimensional real-valued vector space $V=X_1\times\cdots\times X_n$.

\begin{example}[The $2^n$ coordinatewise orders]\label{ex:orders}
For each $X_i$ consider the ordered set $\mathbf P_i^0=\langle X_i,\ge\rangle$ and its dual $\mathbf P_i^1=\langle X_i,\le\rangle$. There are $2^n$ choices of coordinatewise order on $V$ (some may coincide on degenerate carriers),
\[
\mathbf M_j=\mathbf P_1^{j_1}\times\mathbf P_2^{j_2}\times\cdots\times\mathbf P_n^{j_n}=\langle V,\geqslant_j\rangle,\qquad j=\textstyle\sum_i j_i2^{i-1}\in\{0,\ldots,2^n-1\},
\]
and by Example~\ref{ex:order-product} their tokenizations are the $2^n$ matchups
\[
T_{\mathbf M_j}=\oslash(T_{\mathbf P_1^{j_1}},f_1)\otimes\oslash(T_{\mathbf P_2^{j_2}},f_2)\otimes\cdots\otimes\oslash(T_{\mathbf P_n^{j_n}},f_n),
\]
where $f_i$ renames the symbol $\ge^{j_i}$ ($\ge^0$ is $\ge$, $\ge^1$ is $\le$) to the common symbol $\geqslant_j$ so that the factors can be matched.
\end{example}

\begin{example}[Covariant Token classes of a data set]\label{ex:covariant}
Let $\mathcal D$ be a data set of $l$ instances $(k,v)$ with $v\in V$ and $k\in\{1,\ldots,l\}$ the identifier of the instance. With $T_I^j=\langle\mathbb N,\{\geqslant_j\},\varnothing\rangle$, the blending $T_I^j\oplus T_{\mathbf M_j}$ has carrier $\mathbb N\times V$ and relates $(k_1,v_1)$ to $(k_2,v_2)$ iff $v_1\geqslant_j v_2$: it is a \emph{preorder} on identified instances, generally not a total preorder. The data set with the $j$-th order is the subclass
\[
D_j=\bigl\langle\mathcal D,\{\geqslant_j\},\heap(T_I^j\oplus T_{\mathbf M_j})\cap(\{\geqslant_j\}\times\mathcal D\times\mathcal D)\bigr\rangle\Subset T_I^j\oplus T_{\mathbf M_j}.
\]
Instances with the same feature vector are equivalent in this preorder. Form the quotient poset $Q_j$ and its Hasse diagram, then pull its cover relation back to instance identifiers: $E_j$ consists of pairs whose feature-equivalence classes are related by a cover. Let $H_j$ be the subclass of $D_j$ retaining precisely these records, and put $\mathbf{GC}_j=\oslash(H_j,\geqslant_j\mapsto\succcurlyeq_j)$. We call it the \emph{covariant Token class} for this orientation, without identifying this usage with statistical covariance.

A monotone feature map preserves the full order records, but need not preserve covers. For example the monotone map on $0<1<2$ given by $0\mapsto0$, $1\mapsto2$, $2\mapsto2$ sends the cover $0\prec1$ to a non-cover and collapses $1\prec2$. Order isomorphisms do preserve covers. Thus strictly monotone coordinate recodings that induce an order isomorphism onto the transformed data, including invertible sign flips with the corresponding orientation change, give Token maps between the covariant classes after their relation labels are aligned. General normalization or binning requires checking cover preservation or retaining the full order heaps instead.
\end{example}

\begin{example}[Labels and classification]\label{ex:labels}
Let $\langle\mathcal D,f,\mathcal C\rangle$ be training data with label set $\mathcal C$ and $f:\mathcal D\to\mathcal C$. Then $T_{\mathcal C}=\langle\mathcal C,\{\succcurlyeq\},\{(\succcurlyeq,C,C):C\in\mathcal C\}\rangle$ is the tokenization of the discrete order on labels. For a test instance $d=(l+1,v)$ let $\mathbf{GC}'_j$, $j=0,\ldots,2^n-1$, be the covariant Token classes of $\mathcal D\cup\{d\}$ and let
\[
\mathbf{GC}=\oslash(\mathbf{GC}'_0,g_0)\uplus\oslash(\mathbf{GC}'_1,g_1)\uplus\cdots\uplus\oslash(\mathbf{GC}'_{2^n-1},g_{2^n-1}),
\]
where $g_j$ renames $\succcurlyeq_j$ to $\succcurlyeq$. For each candidate label $C_i$ of $d$, extend $f$ to $p_i:\mathcal D\cup\{d\}\to\mathcal C$. Its copairing $\bar p_i$ on the tagged coproduct carrier sends every copy of an instance to its $p_i$-label. Because the target heap is diagonal, $\bar p_i$ is a Token map exactly when every recorded edge joins equal labels. Restricting to the subclass $\mathbf{GC}^{(d)}$ containing only edges incident to $d$ gives the more local condition: all its order-cover neighbours must have label $C_i$. These are neighbours in the chosen order graphs, not metric nearest neighbours. Conflicting neighbouring labels yield no admissible candidate; with no neighbours every candidate passes. A prediction procedure therefore needs an explicit conflict or tie rule. Fewer graph summands, deletion by $\ominus$, or bins introduced by $\rtimes$ specify weaker variants, whose statistical merits require evaluation.

\end{example}

\section{Interior Structure Mapping: Tree Tokens and Reification}\label{interior-structure-mapping-and-tree-Token-classes}

\subsection{Generation of tree Tokens}\label{generation-of-tree-Tokens}

A tree Token on a set $B$ is a tuple whose entries may be elements of $B$ or again tree Tokens, such as $(x_1,(x_2,(x_3,x_4)),x_5)$. Tuple nodes may have any finite number of children; operation symbols acquire a fixed rank only when a particular term signature is imposed \cite{comon2007}. To keep formation algebraic, we add one unary operator to the monoid signature of Remark~\ref{rem:fmon}.

\begin{definition}[Pre-treepoid, treepoid]\label{def:treepoid}
The \emph{pre-treepoid} $\Gt(B)$ is the free operated monoid $\langle M,\cdot,\&,\varepsilon\rangle$ on $B$: concatenation is associative with unit $\varepsilon$, and $\&$ obeys no additional equations. Its elements are forests with a unique factorization into letters and factors $\&(r)$. Define length, depth and top-level arity by
\begin{align*}
|\varepsilon|&=\|\varepsilon\|=\ary(\varepsilon)=0,&
|x|&=\|x\|=\ary(x)=1\quad(x\in B),\\
|\&(r)|&=\max(|r|,1),&
\|\&(r)\|&=\|r\|+1,\quad\ary(\&(r))=1,\\
|rs|&=|r|+|s|,&
\|rs\|&=\max(\|r\|,\|s\|),\quad\ary(rs)=\ary(r)+\ary(s).
\end{align*}
The \emph{treepoid} is $\Tr(B)=B^{\pitchfork}=\{r:\ary(r)\le1\}$. It contains the empty forest $\varepsilon$, the letters, and the tuple trees $\&(r)$. These three syntactic forms are distinct. For a tuple tree put $\tary(\&(r))=\ary(r)$; put $\tary(x)=\tary(\varepsilon)=0$. The inverse notation $\&^{-1}$ is used only on the image of $\&$.
\end{definition}

Write $\&(r_1\cdots r_n)$, with each $\ary(r_i)=1$, as $(r_1,\ldots,r_n)$. In particular $()=\&(\varepsilon)$ is an empty tuple of depth one, distinct from the empty forest. A nonempty flat tuple of letters has depth two. Thus the algebraic presentation with $\cdot,\&,\varepsilon$ gives ordered tuple trees, together with an explicit empty forest; depth alone does not distinguish letters from empty tuples.

Every function $f_0:B_1\to\Gt(B_2)$ extends uniquely to an operated-monoid homomorphism. If every $f_0(x)$ has arity one, the extension preserves top-level arity and restricts to a \emph{tree morphism} $\tmorph(f_0):\Tr(B_1)\to\Tr(B_2)$. It may substitute $(y,z)$ for $x$, but may not substitute the forest $yz$ or $\varepsilon$. Depth and tuple arity need not be preserved.

\begin{definition}[Connection and flattening]\label{def:join}
For tuple trees $r,s$, their \emph{connection} is
\[
r\Join s=\&(\&^{-1}(r)\cdot\&^{-1}(s)).
\]
It is associative, with unit $()$, and extends pointwise to sets of tuple trees. Let $\lambda:\Gt(B)\to B^*$ erase all tuple brackets: $\lambda(x)=x$, $\lambda(\varepsilon)=\varepsilon$, $\lambda(rs)=\lambda(r)\lambda(s)$, and $\lambda(\&(r))=\lambda(r)$. Define \emph{flattening} by $\Join(r)=\&(\lambda(r))$, and apply it pointwise to heaps.
\end{definition}

Flattening has depth two if there is at least one letter, and returns $()$ otherwise. It condenses empty tuple nodes: $\Join(a,(),b)=(a,b)$. It commutes with letter relabelling. General substitutions can change the number of flattened letters, so they do not in general induce ordinary Token maps between the flattened heaps: $a\mapsto(b,c)$ sends $(a)$ to $(b,c)$ after flattening, while ordinary Token maps preserve word length. The relation $Eq_{\Join}(r,s)$ defined by equality of flattenings is nevertheless an equivalence relation.

\begin{definition}[Correspondence, universal Token, scattering]\label{def:scattering}
A \emph{tree correspondence} from $r$ to $s$ is a substitution on precisely the letters occurring in $r$, with nonempty trees as values, whose extension sends $r$ to $s$. Put $1{:}n=(1,\ldots,n)$, including $1{:}0=()$. A \emph{universal Token} has each letter $1,\ldots,n$ exactly once in left-to-right order, for some $n\ge0$; equivalently its flattening is $1{:}n$. For a correspondence $g:u\to r$ with universal source, define
\[
\ll_g(r)=\{g(i):1\le i\le n\},\qquad
\ll^*(r)=\bigcup_{u,g:u\to r}\ll_g(r).
\]
The \emph{flat scattering} $\ll(r)$ is the set of letters occurring in $r$. Here $n$ counts letter occurrences, not $|u|$: empty tuple nodes contribute to length but supply no substitution variable.
\end{definition}

\begin{theorem}[Uniqueness of correspondence]\label{thm:corr}
There is at most one tree correspondence $\corr(r,s)$ between fixed trees. Each $r$ has a unique universal Token $\universal(r)$ admitting a letter-to-letter correspondence to $r$. It is obtained by numbering letter occurrences from left to right and retaining every tuple node. It has the same length and shape as $r$.
\end{theorem}
\begin{proof}
The empty forest can only map to itself. A letter determines its own substitution value. For a tuple, injectivity of $\&$ and unique forest factorization determine the images of its ordered children; induction determines the substitutions on their letters. Repeated occurrences must agree or no correspondence exists. Letter relabelling preserves tuple shape, and the universal condition then fixes the numbering uniquely, including the case of no letters.
\end{proof}

For $r\ne\varepsilon$, $\ll^*(r)$ is exactly its set of nonempty subtrees, including $r$ and its letters. A substitution image occurs at a position of $r$; conversely any subtree can be replaced by one fresh numbered letter, with the remaining letters numbered in order. For $r=\varepsilon$ both scatterings are empty. Thus $\ll(r)\subseteq\ll^*(r)$ and subtree depth never exceeds the depth of $r$.

\subsection{Tree Token classes and tree Token maps}\label{Tokens-maps-between-tree-Token-classes}

All unions of base and core are tagged. No literal disjointness assumption on their underlying sets is needed.

\begin{definition}[t-class, t-map, tree Token space]\label{def:tclass}
A \emph{t-class} is $T=\langle B,E,H\rangle$ with $H\subseteq\Tr(B\sqcup E)$. A \emph{t-map} $T_1\to T_2$ is a substitution
\[
f_0:B_1\to\Tr(B_2\sqcup E)\setminus(E\cup\{\varepsilon\})
\]
whose extension $\tmorph(f_0\sqcup\id_E)$ sends $H_1$ into $H_2$. Equality retains the entire function $f_0$, including values on unused letters. Write $\TTS_E$ for this category, and $\TTS$ when using a common core. A \emph{letter map} has $f_0(B_1)\subseteq B_2$. A t-class is \emph{supported} if every element of $B$ occurs in some member of $H$; these objects form the full subcategory $\TTS_E^{\mathrm{supp}}$.
\end{definition}

\begin{proposition}\label{prop:tts-category}
Identities extend $x\mapsto x$, and composition substitutes into substitution values:
\[
(g\star f)_0(x)=\tmorph(g)(f_0(x)).
\]
These operations make $\TTS_E$ a category. Letter maps form a subcategory. Encoding a flat word $(z_1,\ldots,z_n)$ as the tuple tree $\&(z_1\cdots z_n)$, including $()$ for the empty word, defines a faithful functor $\Xi:\TS_E\to\TTS_E$; it is fully faithful when the target arrows are restricted to letter maps. Coproducts are tagged unions of bases and of the renamed heaps.
\end{proposition}
\begin{proof}
Freeness makes extension commute with substitution and gives associativity. An allowed substitution value remains either a base letter or a tuple tree after another substitution, and hence never becomes a bare core letter or the empty forest. Heap preservation is closed under composition. On flat tuples letter substitution is the ordinary entrywise action, proving the assertion about $\Xi$. Maps out of the tagged union are specified independently on its two generator sets, which proves the coproduct property.
\end{proof}

General substitutions need not preserve tree shape, and $\Xi$ is not asserted full into the unrestricted substitution category. Located neural states and structured reasoning traces can be recorded by trees when a corresponding syntax has actually been specified; the existence of a numerical hidden state does not alone provide such a syntax.

\begin{remark}[Conventions]\label{rem:tree-conventions}
The support $\base(H)$ is the set of base letters occurring in $H$. When only a heap is given we use the supported class $\langle\base(H),E,H\rangle$. Passing from an arbitrary class to its support discards unused generator information and is not an equivalence of the full substitution categories. Metavocabulary such as the reification labels below is chosen fresh, using tagged extensions of the core.
\end{remark}

\subsection{Shared-subterm reification}\label{reification-of-tree-Token-classes}

We retain shared subterms as single carrier elements. A full embedding for arbitrary substitutions requires the support hypothesis: an unused source letter could otherwise be sent to a tree absent from the target heap and its subterms.

\begin{definition}[Reification]\label{def:reification}
For a supported $T=\langle B,E,H\rangle$, let $N_H$ be the set of tuple subtrees of members of $H$, with a tagged identifier $[p]$ for each distinct tuple $p$. The carrier is $B\sqcup N_H$. Extend the core to $E^+=E\sqcup\{\Finv,\mathrm{root},\mathrm{empty}\}$. Represent a base letter by itself, a core letter by its core tag, a tuple $p$ by $[p]$, and the empty forest by the core symbol $\mathrm{empty}$; denote this assignment by $\ulcorner-\urcorner$. Define $\Finv(T)$ by the records
\begin{align*}
&(\Finv,[p],\ulcorner r_1\urcorner,\ldots,\ulcorner r_m\urcorner)
&&\text{for each }p=(r_1,\ldots,r_m)\in N_H,\\
&(\mathrm{root},\ulcorner r\urcorner)&&\text{for each }r\in H.
\end{align*}
For $p=()$ the first record is $(\Finv,[()])$. The identifier tags and the core symbol $\mathrm{empty}$ are distinct.
\end{definition}

The node records recover the ordered children of each tuple, and the root records recover the heap. Sharing identifies equal subterms, rather than equal occurrences. No unheaded copy of a tuple is needed: its entries are recovered from its node record. Bare core roots and the empty forest are represented without turning core symbols into variable nodes.

\begin{lemma}\label{lem:reif-subtrees}
For a t-map between supported classes, the image of each non-core, nonempty source subtree is a non-core, nonempty subtree of some target heap member.
\end{lemma}
\begin{proof}
Every source base letter occurs in a heap member. Substitution places its image, and the image of every surrounding tuple, at a position of the substituted heap member. Tuple roots remain tuples, and allowed base images exclude bare core letters and the empty forest.
\end{proof}

\begin{theorem}[Reification on supported classes]\label{thm:reification}
Reification is a fully faithful functor
\[
\Finv:\TTS_E^{\mathrm{supp}}\longrightarrow\TS_{E^+},\qquad
\TTS_E(T_1,T_2)\cong\TS_{E^+}(\Finv T_1,\Finv T_2).
\]
It sends a substitution to its action on represented subtrees. A finite supported heap has a finite reification carrier and heap.
\end{theorem}
\begin{proof}
The lemma supplies the carrier map. Substitution commutes with tuple formation, so it preserves node records; heap preservation gives preservation of roots. Identities and composition follow from substitution.

Conversely let $h:\Finv T_1\to\Finv T_2$ preserve the records. Decode $h(b)$ for each $b\in B_1$ as a target base letter or tuple subtree, obtaining an allowed $f_0(b)$. Each source node record forces its image to be the unique target tuple having the images of its children in the same order. Induction on the source tree therefore gives $h([p])=\ulcorner\tmorph(f_0)(p)\urcorner$. Root preservation forces the substituted roots into $H_2$, including core and empty roots, so $f_0$ is a t-map. All carrier elements are either base letters or tuple nodes, which proves uniqueness. A finite heap uses finitely many letters and tuple subtrees.
\end{proof}

\begin{remark}[Unused letters and occurrence reification]\label{rem:occurrence-reif}
The support condition cannot simply be omitted. In $T_1=\langle\{a,u\},\varnothing,\{a\}\rangle$ and $T_2=\langle\{b,c\},\varnothing,\{b\}\rangle$, two substitutions agreeing on $a\mapsto b$ but differing on $u$ become indistinguishable if reification keeps only heap subtrees. Merely adding unused letters also fails to accommodate arbitrary substitution images absent from the target subtree set. For letter maps one can retain all base letters and add unary leaf records; that is a different, restricted construction.

To retain occurrences, give each tuple-node occurrence an identifier $(r,\pi)$, where $\pi$ is its path within the heap member $r$, and use the same child and root records. Keep all base letters as a separate summand. Letter relabelling then induces a functor by $(r,\pi)\mapsto(\tmorph(f)(r),\pi)$. No fullness claim is needed for this variant. Sending occurrence identifiers to their subterm identifiers gives the canonical sharing Token map on supported classes. Arbitrary heaps over this signature need not be valid trees: acyclicity and consistent child records are additional conditions. Section~\ref{sec:graph-evaluation} supplies evaluation conditions for computation graphs and proves when sharing preserves the computed function.
\end{remark}

\section{Finite Mappings and Computing Machines}\label{sec:finite-representations}

\subsection{Behavior, description, and the research proposal}

A total mapping $f:A\to B$ on finite encoded interfaces has an exhaustive table. Its implementation can instead be a circuit, a decision diagram, a factorized calculation, or an architecture with parameters. Our proposed machine treats these implementations as structured data: it evaluates them, composes them, and transforms them while retaining an explicit account of their behavior and cost. A Transformer is one possible description family within this proposal.

We use \emph{compact evaluable representation} (CER) for an executable description that is shorter than a declared reference for the same Token map under a fixed encoding (Definition~\ref{def:cer}). The common machine interface concerns representation, evaluation and certified transformation; each representation language may choose its own internal syntax.

The basic ingredients have established precedents. Decision diagrams exploit shared structure in Boolean functions \cite{bryant1986}; knowledge compilation compares representation languages by succinctness and supported operations \cite{darwiche2002}; monoidal computers distinguish programs from computations and supply categorical evaluators \cite{pavlovic2013}; verified compilation proves preservation of behavior \cite{leroy2009}; equality saturation organizes equivalent expressions for optimization \cite{willsey2021}. We do not claim these ideas as new. The purpose here is to connect them to the paper's construction: operation graphs give a Token specification of evaluation, reification exposes repeated subterms, and heap transport describes selected transformations of implementations. The resulting theorems are foundations for this research proposal, not a claim of priority for sharing or compiler correctness.

\paragraph{Functions with structural interfaces.}
Thesis~\ref{thesis:interfaces} asks that an interface record what a computation must preserve. For pure total functions, Cartesian closed structure supports the higher-order operations corresponding to typed $\lambda$-calculus \cite{lambek1986}. Our finite-mapping proposal uses that fragment. Partiality, effects and stochastic behavior require an explicit extension of the state or semantic model.

Within Token Space we call the choice of a heap on a carrier a \emph{structural refinement of the interface}. The carrier still specifies the available values; the heap specifies relations whose images must be admitted by the target. Thus the interface contract is
\[
\delta(r)=f:A\to B,\qquad F_E(f)(\Upsilon_A)\subseteq\Upsilon_B.
\]
This use of refinement concerns whole relational structures. It can express order, grouping or equal-prefix relations between several values, rather than only a predicate selecting individual values from a base type. Refinement types provide a related way of enriching type information \cite{freeman1991}; establishing a correspondence with a particular refinement-type calculus would require its syntax, typing rules and interpretation.

The direction of the obligations is explicit. At fixed carriers, enlarging the source heap adds preservation requirements, whereas enlarging the target heap admits more images:
\begin{align*}
\Upsilon_A\subseteq\Upsilon_A'&\ \Longrightarrow\
\Hom(\langle A,E,\Upsilon_A'\rangle,T_B)\subseteq\Hom(\langle A,E,\Upsilon_A\rangle,T_B),\\
\Upsilon_B\subseteq\Upsilon_B'&\ \Longrightarrow\
\Hom(T_A,\langle B,E,\Upsilon_B\rangle)\subseteq\Hom(T_A,\langle B,E,\Upsilon_B'\rangle).
\end{align*}
The hom sets are regarded here as sets of carrier functions. These inclusions follow directly from heap preservation and show why the two interface roles must be specified, rather than inferred from the word ``refinement.''

Heap preservation is also related to the use of relations in parametricity \cite{reynolds1983}, but the result established here is preservation of the \emph{chosen} records. Reynolds-style relational parametricity additionally concerns uniformity across type interpretations, with relation liftings for the type constructors and an abstraction theorem. Those ingredients do not follow from preserving a fixed heap: an empty source heap alone admits every carrier function. Developing such a bridge is a further task.

The denotation map $\delta:R_{A,B}\to\Hom_{\TS_E}(T_A,T_B)$ defined next is a denotational semantics for represented implementations with these interface obligations. The parameter object $P_E(R_{A,B})$ of Section~\ref{sec:parameter-objects} expresses that every member of the specified description family is admissible. A circuit, a lookup table and a Transformer can therefore have the same structural interface while exposing different internal representations to the machine. The interface determines the preserved relations; a representation determines how the mapping is evaluated and how its implementation may be changed.
\begin{table}[htbp]
\centering\small\renewcommand{\arraystretch}{1.13}
\begin{tabularx}{\linewidth}{@{}p{0.21\linewidth}p{0.34\linewidth}X@{}}
\toprule
Representation & What is stored & How evaluation proceeds \\
\midrule
Explicit table & One output for each encoded input. & Index or search the table. \\
Circuit & Primitive operations and directed connections. & Evaluate nodes in dependency order. \\
Decision diagram & Tests, branches, and shared subgraphs. & Follow the tests selected by the input. \\
Factorized calculation & Factors and compatible intermediate interfaces. & Execute the factors without materializing their composite table. \\
Parameterized model & Architecture, numerical parameters, and an execution protocol. & Run the prescribed numerical operations. \\
Hybrid description & Typed components in several representations. & Dispatch their evaluators and compose the results. \\
\bottomrule
\end{tabularx}
\caption{Representation choices for finite mappings. Expressivity, correctness and cost must be established for each choice; none is assumed universally compact.}\label{tab:representations}
\end{table}

Three comparisons differ. \emph{Exact representation compression} lowers a specified cost while retaining the same mapping on the same input domain. \emph{State aggregation} identifies states relative to selected observations. \emph{Approximation} changes the mapping under an error criterion. Sharing syntax nodes below is the first kind: it does not merge distinct external inputs. Section~\ref{sec:quotient-compression} separately studies quotients of semantic state carriers.

\subsection{Typed representation systems}\label{sec:representation-system}

Choose Token classes $T_A,T_B,\ldots$ with finite, effectively encoded carriers over a common core $E$. Their heaps specify which mappings are admissible. Individual descriptions are finite strings or finite encoded structures; the set of all descriptions need not be finite.

\begin{definition}[Representation system]\label{def:representation-system}
A typed representation system has sets $R_{A,B}$ of valid descriptions, denotations
\[
\delta_{A,B}:R_{A,B}\to\Hom_{\TS_E}(T_A,T_B),
\]
and a total evaluator on valid descriptions and admitted inputs, with
\begin{equation}\label{eq:representation-eval}
\operatorname{eval}_{A,B}(r,a)=\delta_{A,B}(r)(a).
\end{equation}
It supplies identity descriptions $e_A$ and composition constructors
\[
\operatorname{comp}_{A,B,C}:R_{B,C}\times R_{A,B}\to R_{A,C}
\]
satisfying
\begin{equation}\label{eq:representation-laws}
\delta(e_A)=\id_A,\qquad
\delta(\operatorname{comp}(s,r))=\delta(s)\circ\delta(r).
\end{equation}
An effective realization specifies algorithms for evaluation and these constructors, and a finite certificate or decidable check for description validity. Where parallel and alternative interfaces are used, it also supplies constructors denoting the categorical products of maps and their coproducts. Their correctness is an additional requirement, not a consequence of their names.
\end{definition}

Surjectivity of $\delta$ is optional: a restricted library or fixed neural architecture may realize only some admissible maps. Injectivity is usually undesirable, since distinct implementations can have the same behavior. Write $r\equiv r'$ if $\delta(r)=\delta(r')$.

\begin{proposition}[The behavioral quotient]\label{prop:representation-category}
The sets $R_{A,B}/{\equiv}$ are the hom sets of a category, with composition induced by $\operatorname{comp}$. Denotation defines a faithful functor into $\TS_E$, full on the chosen objects when each $\delta_{A,B}$ is surjective.
\end{proposition}
\begin{proof}
The composition equation makes the result independent of representatives. Associativity and identities hold after applying $\delta$, hence in the quotient. Equality of denotations is precisely the equivalence used to form each hom set.
\end{proof}

This proposition checks the organization of behaviors. It deliberately forgets implementation differences, so it is not the structure in which to optimize descriptions. Raw composition may retain different bracketings as different syntax.

Parameter objects connect execution to the exponential construction:
\begin{equation}\label{eq:representation-family}
\begin{tikzcd}[column sep=large,row sep=large]
P_E(R_{A,B})\otimes T_A \arrow[r,"\widehat\delta\otimes\id"] \arrow[dr,"\operatorname{eval}_{A,B}"'] & T_B^{T_A}\otimes T_A \arrow[d,"\operatorname{ev}"] \\
& T_B.
\end{tikzcd}
\end{equation}
Proposition~\ref{prop:parameters} supplies the displayed Token maps because every $\delta(r)$ is admissible. The curry $\widehat\delta$ names the mapping realized by a description. Its existence supplies neither an algorithm for finding descriptions nor a cost bound.

\paragraph{Compact evaluable representations of Token maps.}
We give the compact-representation viewpoint a name within this effective
setting. A program retains a structure from which the specified mapping
can be evaluated on an input. Its internal notation may use tables,
connections, factorizations or parameters. Compactness compares the size
of that structure with a declared reference representation.

\begin{definition}[Compact evaluable representation (CER)]\label{def:cer}
Fix an effective realization of Definition~\ref{def:representation-system},
with a fixed interpreter, operation library and binary description
encoding. Let $L(r)$ be the encoded length of a valid description $r$.
Fix a Token map $f:T_A\to T_B$ and a reference description
$b\in R_{A,B}$ with $\delta(b)=f$. A \emph{compact evaluable representation}
(\emph{CER}) of $f$ relative to this reference and encoding is a valid
description $r\in R_{A,B}$ such that
\begin{equation}\label{eq:cer}
 \begin{gathered}
 \delta(r)=\delta(b)=f,\qquad L(r)<L(b),\\
 \operatorname{eval}_{A,B}(r,a)=f(a)\quad(a\in A).
 \end{gathered}
\end{equation}
Evaluation is performed by the fixed effective evaluator on $(r,a)$;
the reference $b$ is not an additional input to that evaluator.
All description-specific parameters, tables, module bodies and decoding
information required for execution are included in $L$, as are any
certificates that must be stored with the executable description. Only
the fixed, declared common environment is excluded from the comparison. When these
choices are understood, we say simply \emph{a CER of the Token map $f$}.
\end{definition}

The equality in~\eqref{eq:cer} supplies exact behavior, and membership in
$\Hom_{\TS_E}(T_A,T_B)$ supplies heap preservation. The strict inequality
supplies the particular compactness claim. An explicit table on the finite
input interface is a natural reference, but another declared implementation
can serve the same role. Several representation languages can be compared
through a fixed tagged description language and its evaluator, with their
format tags and nonshared support data counted. Interface specifications
are held fixed; moving description-specific information into an uncounted
library does not establish a saving.

Here ``compact'' concerns description length, not topological compactness
or a claim of a shortest program. It is a property of a representation
relative to the chosen reference. Evaluable descriptions remain valid
programs when the inequality fails. The evaluator can use the description
on the supplied input without taking an exhaustive mapping table as a
separate argument; time and workspace bounds require further analysis
of that evaluator. For an effectively generated family of valid
descriptions $(r_n)$ and valid references $(b_n)$ satisfying
$\delta(r_n)=\delta(b_n)=f_n$ in the same fixed setting, savings can also
be stated asymptotically: $L(r_n)=o(L(b_n))$, with $L(b_n)>0$ eventually,
gives CERs for all sufficiently large $n$. The prefix-XOR construction in
Section~\ref{sec:parity} provides a concrete family.

\paragraph{Interfaces without prescribed variable syntax.}
The representation system fixes inputs, outputs and evaluation while
leaving internal syntax to each description language. Input ports and
connections, table indices, decision tests and numerical parameters can
serve their respective representations; named variables are an optional
choice. When a whole Token sequence or configuration is the argument,
encode it as a carrier value and choose heaps expressing its interface
obligations. This allows computations to produce configurations with
different internal sizes. The induced action on the outer structural
records still fixes core symbols and preserves record shape, as required
by Definition~\ref{def:token-class}. Program descriptions and ordinary inputs
can share a record language once their encodings and roles are specified
(Section~\ref{sec:module-expansion}).

The compactness condition introduces no new closure assumption. Correct
composition supplies an evaluable description of the composite map, but
does not by itself prove a length advantage against the chosen composite
reference. Likewise, a machine may expand a CER to expose independent
operations and later pack it into another representation. The
implementation categories below retain these distinct descriptions and
their certified transformations, including intermediate descriptions that
are not compact. CER names the size comparison within that organization;
the resource profile of Section~\ref{sec:granularity} also records work,
span, memory and communication. This terminology connects established
representation and compilation ideas to Token-map interfaces rather than
asserting a new compression algorithm.

\subsection{Retaining implementations and their transformations}\label{sec:implementation-fibres}

For an admissible $f:T_A\to T_B$, retain the set
\[
\mathcal I(f)=\{r\in R_{A,B}:\delta(r)=f\}.
\]
Equip it with specified elementary rewrites $r\xRightarrow{c}r'$, where $c$ certifies validity of both descriptions and preservation of denotation. A certificate can invoke a proved rewrite schema; it need not enumerate all inputs. Finite paths of such rewrites, with empty paths as identities and concatenation as composition, form a category $\mathsf{Impl}(f)$. Distinct descriptions and distinct paths remain distinct. This is a category of implementations of one behavior, whereas Proposition~\ref{prop:representation-category} is a category of behaviors between interfaces.

\begin{proposition}[Replacement in computational contexts]\label{prop:context-rewrite}
If $\delta(r)=\delta(r')$, replacing $r$ by $r'$ in any well-typed expression built from correct composition, product and coproduct constructors leaves its denotation unchanged. A finite sequence of certified replacements has the same property.
\end{proposition}
\begin{proof}
Induct on the expression. Composition respects equality in each argument. Products and coproducts of maps respect equality componentwise. Induction on the number of replacements gives the second claim.
\end{proof}

This permits a local correctness certificate to be reused in a larger machine. It does not say that equal functions have isomorphic implementation graphs, or that every equality can be reached by the chosen rewrites. Semantic identities such as commutativity may require a different certificate from the graph maps studied next.

\subsection{Computation graphs and evaluation as a Token map}\label{sec:graph-evaluation}

Fix a nonempty finite value set $D$ and a finite library $\Omega$ of total operations $\omega_D:D^{k_\omega}\to D$, including nullary operations if needed. The interpreter and the operation library are fixed in every cost comparison. We first use one value sort; several finite sorts can be tagged and supplied with sort records.

\begin{definition}[A valid computation graph]\label{def:computation-graph}
A graph description $G:D^n\to D^m$ has a finite node set $N_G$ and a heap consisting of the following record forms
\[
(\mathrm{in}_i,v_i),\qquad
(\omega,v,u_1,\ldots,u_{k_\omega}),\qquad
(\mathrm{out}_j,o_j).
\]
The $n$ input nodes $v_i$ are distinct. Each node is defined exactly once, either as an input node or by one operation record. Every operand and output refers to a node in $N_G$, each output index has exactly one record, and the directed dependency relation from operands to operation results is acyclic. Output nodes may coincide and unused input nodes are retained. The labels name fixed input positions, output positions and primitive operations, using disjoint tags for these roles.
\end{definition}

These are finite, checkable conditions. In particular, arbitrary Token heaps are not automatically executable graphs. The operation record is the graph of an operation written with its result first; this is a fixed coordinate permutation of the algebraic tokenization of Section~\ref{representing-categories-of-structured-objects-in-Token-space}. Tree terms over $\Omega$ give occurrence graphs, and the shared-subterm construction gives shared graphs, with interface records added explicitly.

For an input $a=(a_1,\ldots,a_n)$, define the value class $\mathcal A_a$ on carrier $D$ by
\begin{align*}
H_{\mathcal A_a}={}&\{(\mathrm{in}_i,a_i):1\le i\le n\}\\
&\cup\{(\omega,\omega_D(d_1,\ldots,d_k),d_1,\ldots,d_k):\omega\in\Omega,\ d_i\in D\}\\
&\cup\{(\mathrm{out}_j,d):1\le j\le m,\ d\in D\}.
\end{align*}
Output records here impose no restriction on the value; they allow all possible results. The graph and value classes share their core.

\begin{theorem}[Evaluation and structural transformations]\label{thm:graph-evaluation}
For every valid $G$ and input $a$ there is a unique Token map $v_{G,a}:G\to\mathcal A_a$. Its output tuple defines a total function
\[
\delta(G)(a)=(v_{G,a}(o_1),\ldots,v_{G,a}(o_m)).
\]
If $q:G\to G'$ is a Token map between valid graphs with the same interfaces and library, then
\begin{equation}\label{eq:graph-semantics}
v_{G,a}=v_{G',a}\circ q\quad\text{for every }a,
\qquad \delta(G)=\delta(G').
\end{equation}
\end{theorem}
\begin{proof}
Input records force $v_{G,a}(v_i)=a_i$. In a topological ordering, each operation record forces the value of its result from already determined operands. Totality of each primitive supplies this value, and functionality makes it unique, including for nullary operations. This constructs the unique Token map. Since $v_{G',a}\circ q$ is also a Token map to $\mathcal A_a$, uniqueness gives the equality. Preservation of the uniquely indexed output records gives equal output tuples.
\end{proof}

The theorem is a sufficient correctness criterion using the description's local structure. Its verification does not require listing the $|D|^n$ inputs. It is not necessary for semantic equivalence: two graphs can compute the same function without a graph Token map between them. Validity of both endpoints matters, since identifying arbitrary nodes can introduce cycles or inconsistent definitions.

The operator calculus also supplies a graph composition procedure. For $G:D^n\to D^m$ and $K:D^m\to D^p$, first make their nodes disjoint. Delete $G$'s output records and $K$'s input records with $\ominus$, merge the remaining heaps, and use direct image along the map identifying each input node of $K$ with the corresponding output node of $G$. Keep $G$'s inputs and $K$'s outputs. The resulting graph is valid: the first graph has no dependencies on the second, and substituted input nodes introduce no second definitions. Its evaluation first evaluates $G$ and then $K$, so its denotation is $\delta(K)\circ\delta(G)$. Identity graphs have inputs as outputs; disjoint graphs with suitably renamed interface indices give parallel products. Thus these descriptions realize the representation laws by explicit finite operations on Token records.

\subsection{Sharing subterms: a complete finite example}\label{sec:parity}

Consider a finite ordered list of terms over $\Omega$ in variables $x_1,\ldots,x_n$. Its occurrence graph $G_{\mathrm{occ}}$ has one node for each operation occurrence, common input nodes for the variables, and a labelled output for each term. Let $G_{\mathrm{sh}}$ have one node for each distinct operation subterm, retaining all declared inputs. Define $q$ by sending an occurrence to its syntactic subterm and fixing input nodes.

\begin{proposition}[Certified sharing]\label{prop:sharing}
The map $q:G_{\mathrm{occ}}\to G_{\mathrm{sh}}$ is a surjective Token map, and
\[
H_{G_{\mathrm{sh}}}=F_E(q)(H_{G_{\mathrm{occ}}}),\qquad
\delta(G_{\mathrm{occ}})=\delta(G_{\mathrm{sh}}).
\]
If $t$ counts operation occurrences and $s$ counts distinct operation subterms, topological evaluation executes $t$ and $s$ primitive calls, respectively. For a fixed bounded-arity library and $m$ outputs, graph encodings using node indices have length
\[
O\bigl((n+t+m)\log(n+t+2)\bigr),\qquad
O\bigl((n+s+m)\log(n+s+2)\bigr),
\]
respectively, with the library supplied separately.
\end{proposition}
\begin{proof}
Every operation occurrence has the same label and ordered operand subterms as its image, so all records are preserved, and every shared record has an occurrence preimage. A shared subterm depends only on proper subterms; their smaller term size proves acyclicity. Each distinct term has one defining operation, hence the shared graph is valid. Theorem~\ref{thm:graph-evaluation} proves behavioral equality. Evaluating each operation node once gives the call counts. Bounded arity bounds the number of indices per node, and each index takes logarithmic space.
\end{proof}

Operationally, the sharing procedure traverses terms from children to parents and interns each key $(\omega,\text{child identifiers})$ in a dictionary. Identical keys reuse an identifier; distinct keys create a new node. Comparing keys exactly makes the procedure correct even when a hash implementation has collisions. It returns both the shared graph and the map $q$ as a certificate. A direct certificate check examines records and validity of the two finite graphs, rather than the full input--output table. The procedure exploits syntactic equality; additional algebraic identities require separately justified rewrites.

For $D=\{0,1\}$ with XOR as the binary primitive, put
\[
t_1=x_1,\qquad t_{i+1}=t_i\oplus x_{i+1},\qquad
P_n(x)=(t_1(x),\ldots,t_n(x)),\quad n\ge2.
\]
An output table for $P_n$ has $n2^n$ bits in a fixed input order. Writing each prefix term independently repeats earlier subterms: it has $\sum_{i=1}^n(i-1)=n(n-1)/2$ XOR occurrences. Sharing retains exactly $t_2,\ldots,t_n$, so it uses $n-1$ XOR nodes. The three representations have the same external domain $\{0,1\}^n$ and the same outputs.

\begin{table}[htbp]
\centering\small\renewcommand{\arraystretch}{1.13}
\begin{tabularx}{\linewidth}{@{}p{0.25\linewidth}p{0.30\linewidth}X@{}}
\toprule
Representation of $P_n$ & Storage under the specified encoding & Evaluation work \\
\midrule
Full output table & $n2^n$ output bits & One table access plus input indexing and $n$ output bits. \\
Separate prefix expression trees & $O(n^2\log n)$ bits using indexed variables & $n(n-1)/2$ XOR calls. \\
Shared prefix graph & $O(n\log n)$ bits using node indices & $n-1$ XOR calls plus output access. \\
\bottomrule
\end{tabularx}
\caption{A family with explicitly proved sharing and cost bounds. These are costs of particular encodings and evaluators, not representation-independent lower bounds.}\label{tab:sharing-cost}
\end{table}

In the terminology of Definition~\ref{def:cer}, this is an effectively generated CER family relative to explicit tables for all sufficiently large $n$. Fix a common tagged encoding for graphs and tables, the XOR library, and the interface convention. The graph descriptions have length $O(n\log n)$, whereas table descriptions contain $n2^n$ output bits, so their length ratio tends to zero. This asymptotic comparison does not assert a storage saving for every small $n$; finite cases require the actual code lengths, including headers and format tags.

For $n=4$, the tree description has six XOR occurrences and the shared graph has three. Its heap is explicitly
\[
\begin{gathered}
(\mathrm{in}_i,x_i)\ (1\le i\le4),\\
(\oplus,u_2,x_1,x_2),\quad(\oplus,u_3,u_2,x_3),\quad(\oplus,u_4,u_3,x_4),\\
(\mathrm{out}_1,x_1),\quad(\mathrm{out}_2,u_2),\quad
(\mathrm{out}_3,u_3),\quad(\mathrm{out}_4,u_4).
\end{gathered}
\]
The table has $64$ output bits; the asymptotic comparison does not assert that graph serialization is already shorter at this small size. The final-output-only parity function has a $2^n$-bit table and the same $n-1$ XOR nodes. Every prefix output depends only on that prefix, so $P_n$ is a Token map $\Cau_n(\{0,1\})\to\Cau_n(\{0,1\})$ and obeys the cross-length equation~\eqref{eq:prefix-natural}.

This example uses two different carriers deliberately. $G_{\mathrm{occ}}$ and $G_{\mathrm{sh}}$ have syntax nodes as carrier elements; $\Cau_n(\{0,1\})$ has complete input or output sequences as carrier elements. The syntax quotient lowers repeated computation while the sequence mapping retains every input. Reification, direct image, evaluation and the causal heap therefore participate in one verified construction, without identifying graph-node sharing with quotienting the input state space.

\paragraph{Following the example into execution and observation.}
The prefix-parity example will also distinguish the roles of retained
state. Once its past outputs have been delivered, an incremental
implementation can retain the running parity and update it with the
next input. Sharing expression nodes, choosing a live graph frontier,
and identifying histories with the same future parity outputs are
three different questions about this computation. The frontier and
future-observation constructions in Sections~\ref{sec:residual-mapping}
and~\ref{obs:future-observation} specify the latter two levels. The
decoder cache in Section~\ref{llm:llm-cache} asks the corresponding
continuation question for a richer computation with its own access
rules. This common route lets an example guide the choice of objects
and maps in a different application.

\subsection{Re-encoding across representation languages}\label{sec:reencoding}

A graph Token map verifies a particular kind of structural transformation. To compare different languages, such as a table and a circuit, use their evaluation records. Let $R,R'$ represent mappings on the same finite interfaces $A,B$. The class $\mathcal G_R$ has tagged carrier $R\sqcup A\sqcup B$, unary sort records and
\[
(\mathrm{eval},r,a,b)\in H_{\mathcal G_R}
\quad\Longleftrightarrow\quad \operatorname{eval}_R(r,a)=b.
\]
Define $\mathcal G_{R'}$ in the same way.

\begin{proposition}[Interface-fixing re-encoding]\label{prop:reencoding}
For $\kappa:R\to R'$, the function $\kappa\sqcup\id_A\sqcup\id_B$ is a Token map $\mathcal G_R\to\mathcal G_{R'}$ iff
\begin{equation}\label{eq:correct-reencoding}
\operatorname{eval}_{R'}(\kappa(r),a)=\operatorname{eval}_R(r,a)
\quad\text{for all }r,a.
\end{equation}
\end{proposition}
\begin{proof}
Sort records are preserved. The image of an evaluation record is $(\mathrm{eval},\kappa(r),a,b)$, which belongs to the target heap exactly when the new evaluator returns the same $b$. Totality and functionality give the equivalence.
\end{proof}

This proposition is a relational expression of semantic preservation. It is a specification, not an algorithm for proving arbitrary compiler correctness. The sharing theorem gives a structural certificate for a particular $\kappa$; exhaustive evaluation is another finite but potentially expensive method. A partial optimizer can be made total by returning its input unchanged when its applicability condition fails, using a common tagged description language. An arbitrary map between evaluation classes may also change inputs or outputs; interface fixing is essential to the displayed claim.

For a fixed interpreter and coding scheme, attach description length $L(r)$, worst-case execution time $T(r)$, and workspace $M(r)$ to implementations. These quantities generally do not descend to the behavioral quotient. A rewrite is an exact storage improvement on $r$ if it preserves $\delta(r)$ and $L(r')<L(r)$; time or memory improvements are separate comparisons. A sequence strictly decreasing a nonnegative integer $L$ terminates in at most $L(r)$ steps, although it need not find a global minimum. Rewrite discovery and certificate checking also have costs. No fixed lossless coding scheme makes every finite mapping substantially shorter than its table, since $|B|^{|A|}$ behaviors must be distinguished.

With a chosen distance and input distribution, one can instead require
\[
\mathbb E_{a\sim\mu}\bigl[d(\operatorname{eval}_{R'}(\kappa(r),a),
\operatorname{eval}_R(r,a))\bigr]\le\epsilon.
\]
This is an approximate transformation, not an arrow of the exact implementation category. An error bound need not survive a later component unchanged; composing such bounds requires assumptions on that component, such as a Lipschitz bound.

\subsection{A machine specification and the Transformer instance}\label{sec:machine-specification}

\begin{definition}[Machine of represented finite mappings]\label{def:represented-machine}
A machine specification consists of a typed representation system, a library of certified description transformations, a chosen resource accounting, and an effective control protocol. Its configurations contain typed data registers, stored descriptions, and a control state. Instructions may evaluate a stored description, construct a composite, replace a description by a certified rewrite, copy or route data, and branch or halt. Each instruction has an explicitly defined transition on valid configurations; error and halt states are included when needed to make the transition total.
\end{definition}

The machine can execute both CERs and other valid descriptions. It can retain a compact representation for storage, expand it to expose work, and replace it by another certified implementation. The CER condition compares representation sizes; effective evaluation gives those representations their computational role.

This separates the finite interfaces of component computations from the storage and control of a running machine. With finite alphabets, a finite control-state set and fixed bounds on all registers and description stores, the configuration set is finite. Unbounded execution with growing storage requires an explicit memory model or a uniform family of such bounded machines. Cartesian closedness alone supplies neither recursion nor unbounded memory.

The graph evaluator, graph composition and sharing procedure above give executable instances of three instructions. For a fixed external input, any finite sequence of certified replacements performed before evaluating a component preserves its output, by Proposition~\ref{prop:context-rewrite}. Replacing a component midway through an execution that retains internal state needs an additional state-migration invariant; function equality by itself does not justify that operation. Section~\ref{sec:online-replacement} provides a sufficient rule by retaining the values at a completed boundary and replacing the residual mapping. The control protocol may choose implementations according to estimated or measured costs, but correctness of replacement and quality of the cost estimate are different obligations.

The machine can mix a table for a small irregular mapping, a circuit for a regular Boolean mapping, and a numerical model for another component. A concrete Transformer with fixed vocabulary, length bound, numerical formats, parameters, and deterministic execution and decoding rules denotes a finite mapping, provided every admitted input has a defined result. Randomness must be supplied as data or modeled stochastically. The ideal real-valued treatment in Section~\ref{attention-in-token-space} uses the larger set-based Token Space.

A Transformer architecture together with its parameters is consequently one description form, and its forward algorithm is its evaluator. Calling it a compressed realization of a specified target requires a correctness or approximation claim and an actual cost comparison. The present construction proves a sharing transformation for a concrete non-neural family and gives a common specification for broader implementations; it does not claim a new compression algorithm for trained Transformers. The research task is to extend the certified transformation library and determine when Token-level structure leads to useful decompositions, representation choices and resource improvements beyond these initial examples. The next section develops concurrency and elasticity as explicit parts of that task, retaining both compressed and expanded implementations when their resource profiles make them useful.

\section{Concurrent and Elastic Execution of Represented Mappings}\label{sec:elastic-execution}

The machine envisaged here has two further design principles: \emph{intrinsic concurrency} and \emph{elasticity of computational scale}. A mapping specifies a behavior, while its representation exposes a particular organization of the work. Dependencies determine which operations must wait; they need not impose a total order on operations that are independent. The machine should also be able to adjust the amount of active work and the granularity of its representations to the resources available. These principles connect the construction of Token Space to a question about machines: how can a stable computation admit changing execution structures?

Both principles are motivated by the Token-level view of Sections~\ref{sec:token-levels} and~\ref{interior-structure-mapping-and-tree-Token-classes}. Atomicity is relative: a tree Token can be treated as one entry while retaining nested structure, and an explicitly encoded finite description can become a carrier element of a further Token class. The iteration of $\ATC$ alone does not supply that encoding. Shared-subterm reification $\Finv$ exposes syntactic structure on its supported domain; executable expansion additionally requires an interpreted operation signature, a module body and a port convention. Section~\ref{sec:module-expansion} supplies these ingredients and proves denotation preservation. The basic computational unit is thus relative to the chosen interface and representation, while the denotation supplies a stable criterion for correctness. In a valid operation graph, concurrency is exposed by nodes between which there is no directed dependency path in either direction; sharing an ancestor does not make two nodes dependent on each other. Module packing and certified expansion can hide or expose this internal work without changing the specified behavior.

This section supplies a first operational account in that language. Partial executions are subclasses of the graph, schedule independence is the unique extension of Token maps along subclass inclusions, and the replacement of unfinished work is a construction with the operators of Section~\ref{sec:operators}. The standard resource bounds of parallel scheduling are stated for completeness, with proofs in Appendix~\ref{app:scheduling}. The results concern valid acyclic graph descriptions over a fixed library of pure, deterministic, total operations; an unordered heap is not by itself an execution semantics, and stateful effects require explicit state and dependency records.

\subsection{Structural concurrency}\label{sec:structural-concurrency}

Call concurrency \emph{intrinsic to a representation} when the representation specifies the necessary dependencies without prescribing an additional execution order among independent operations. Actual parallelism is the simultaneous execution of some of those operations by a physical implementation; a serial evaluator realizes a representation with intrinsic concurrency without exploiting it.

In $\TS_E$ the matchup $\otimes$ is the categorical product, so the interchange law
\begin{equation}\label{eq:parallel-interchange}
(g_1\otimes g_2)\circ(f_1\otimes f_2)
=(g_1\circ f_1)\otimes(g_2\circ f_2)
\end{equation}
and the diagonal $T\to T\otimes T$ are available without further structure: parallel composition of independent interfaces and sharing of an input are already morphisms. What they do not supply is the dependency structure \emph{inside} one computation; that is the role of the graph records. A Cartesian product of semantic interfaces organizes a product of behaviors, whereas a graph description of their independent execution uses disjoint operation nodes with renamed ports (Section~\ref{sec:graph-evaluation}); the two are related through denotation, and a representation system must connect its product constructor to the semantic product correctly.

\subsection{Partial executions as subclasses}\label{sec:partial-valuations}

Fix a valid graph $G:D^n\to D^m$ and an input $a\in D^n$. Write $I_G$ for its input nodes, $O_G=N_G\setminus I_G$ for its operation nodes, and $\operatorname{pred}_G(v)$ for the set of operands of an operation node; a nullary operation has an empty predecessor set. The dependency relation generated by the operation records is a strict partial order $\prec_G$ on $N_G$, and $u,v$ are \emph{independent} when neither $u\prec_G v$ nor $v\prec_G u$.

\begin{definition}[Partial execution as a subclass]\label{def:partial-state}
A set $U$ with $I_G\subseteq U\subseteq N_G$ is \emph{predecessor-closed} if it is a down-set of $\prec_G$. For such $U$ let
\[
G|_U\Subset G
\]
be the subclass with carrier $U$ whose heap consists of the input records and the operation records whose result lies in $U$ (output records omitted). A \emph{partial execution state} of $G$ on $a$, written $(U,\sigma)$, is a Token map
\[
\sigma:G|_U\longrightarrow\mathcal A_a .
\]
The \emph{ready set} is
\begin{equation}\label{eq:ready-set}
\mathcal R_G(U)=\{v\in O_G\setminus U:\operatorname{pred}_G(v)\subseteq U\},
\end{equation}
the set of uncompleted operation nodes $v$ for which $U\cup\{v\}$ is again predecessor-closed. \emph{Firing} $v\in\mathcal R_G(U)$ extends $\sigma$ along the inclusion $G|_U\Subset G|_{U\cup\{v\}}$ by the value prescribed by the operation record of $v$; a \emph{batch firing} extends along $G|_U\Subset G|_{U\cup R}$ for any $R\subseteq\mathcal R_G(U)$, using only values present before the batch.
\end{definition}

Unwinding the definition, $\sigma:U\to D$ has $\sigma(v_i)=a_i$ at every input and satisfies
\[
\sigma(v)=\omega_D(\sigma(u_1),\ldots,\sigma(u_k))
\]
for every executed operation record $(\omega,v,u_1,\ldots,u_k)$. It is the usual consistent partial valuation. Write $G^-$ for $G$ with only its output records deleted. Then $G|_{N_G}=G^-$, not literally $G$. Since $\mathcal A_a$ admits every value at each output label, the final valuation extends to a Token map on $G$ with the same carrier function. Reading the ordered output ports yields the result.

\begin{proposition}[The lattice of completed sets]\label{prop:execution-lattice}
The collection
\[
\mathcal L_G=\{U:I_G\subseteq U\subseteq N_G,\ U\text{ predecessor-closed}\}
\]
is a finite distributive lattice with meet $U\cap V$, join $U\cup V$, bottom $I_G$ and top $N_G$. Its covers are exactly the inclusions $U\subset U\cup\{v\}$ with $v\in\mathcal R_G(U)$. If $R\subseteq\mathcal R_G(U)$ then
\begin{equation}\label{eq:ready-boolean-interval}
\mathcal P(R)\ \cong\ [U,U\cup R],\qquad A\longmapsto U\cup A.
\end{equation}
\end{proposition}
\begin{proof}
Intersections and unions of predecessor-closed sets containing $I_G$ have the same property. Their set operations satisfy distributivity. For $U\subsetneq V$, a minimal element of $V\setminus U$ has all predecessors in $U$, so it is ready. If the inclusion is a cover, adding that element already gives $V$; the converse is immediate. All predecessors of a ready node lie in $U$, so every $U\cup A$ with $A\subseteq R$ is predecessor-closed. Every intermediate set has this form uniquely, proving the interval isomorphism.
\end{proof}

These lattice points are possible completed sets, not arbitrary combinations of numerical values. Readiness is decidable from the finite graph and the completed set; a runtime maintains predecessor counters instead of rescanning the heap. Nodes are executed at most once in this semantics.

\begin{lemma}[Unique extension]\label{lem:partial-valuation}
For every predecessor-closed $U$ there is exactly one partial execution state on $G|_U$, namely
\begin{equation}\label{eq:partial-invariant}
\sigma=v_{G,a}|_U ,
\end{equation}
the restriction of the unique Token map of Theorem~\ref{thm:graph-evaluation}. Firing a ready node is the unique extension of $\sigma$ along $G|_U\Subset G|_{U\cup\{v\}}$. If $U\ne N_G$ then $\mathcal R_G(U)\ne\varnothing$.
\end{lemma}
\begin{proof}
The restriction of a Token map to a subclass is a Token map, so $v_{G,a}|_U$ is a partial execution state. Conversely a partial execution state is determined on $U$ by induction along a linear extension of $\prec_G$ restricted to $U$: input values are forced by the input records, and each operation in $U$ has all its operands in $U$ and its value forced by its record. Uniqueness gives both equalities. If unfinished nodes exist, a $\prec_G$-minimal one among them is not an input and has no unfinished predecessor, so it is ready.
\end{proof}

For $U\subseteq V$ in $\mathcal L_G$, let $\iota_{U,V}:G|_U\to G|_V$ be the inclusion Token map. The unique states satisfy the compatibility law
\begin{equation}\label{eq:partial-compatibility}
\sigma_V\circ\iota_{U,V}=\sigma_U.
\end{equation}
Thus the Boolean intervals identify independently available choices of progress, while this law identifies the values shared by every route through them.

\begin{theorem}[Commutation and schedule independence]\label{thm:schedule-independence}
If distinct nodes $u,v$ are ready in the same state, their firings commute: both orders, and the batch $\{u,v\}$, yield the unique state on $G|_{U\cup\{u,v\}}$. Every maximal legal execution path of nonempty batch firings, starting at $I_G$, terminates at $N_G$ with carrier valuation $v_{G,a}$ and output $\delta(G)(a)$. In an execution that also permits waiting, the same conclusion holds provided waiting does not prevent eventual progress on unfinished work.
\end{theorem}
\begin{proof}
Simultaneously ready nodes are independent, so $U\cup\{u\}$, $U\cup\{v\}$ and $U\cup\{u,v\}$ are predecessor-closed, and by Lemma~\ref{lem:partial-valuation} each carries exactly one state, the restriction of $v_{G,a}$; the three routes from $U$ to $U\cup\{u,v\}$ therefore agree, and the assertion extends to any ready batch. Each nonempty batch strictly increases $|U|$, so such a path has at most $|O_G|$ steps, and a path that stops below $N_G$ is not maximal by the lemma. Its terminal state is $v_{G,a}$ independently of the batches chosen. Waiting changes no stored value, and the progress assumption rules out an infinite unfinished wait.
\end{proof}

The theorem turns the unique Token map of Theorem~\ref{thm:graph-evaluation} into an operational invariant: an execution is a path in the lattice $\mathcal L_G$ of subclasses, and all maximal progressing paths from the initial state end with the same valuation. An asynchronous implementation can issue a ready operation when a worker is available and publish its result when it finishes; it must guarantee unique publication, retention of operands until their reads complete, and correct readiness tracking. Hardware barriers, locks and communication costs are not supplied by the categorical proof. For several simultaneous invocations, node identity must include an invocation identifier, e.g.\ through runtime records $(\mathrm{value},j,v,d)$ with $j$ the invocation; since heaps are sets, distinct occurrences are distinguished by identifiers.

\subsection{Execution diagrams, gluing and consistent values}\label{sec:execution-diagrams}

A dependency graph has operations as nodes and data dependencies as edges. The diagrams below have entire Token classes as objects and Token maps as arrows. This change of level allows a local execution, including its records, to be compared with another local execution. We abbreviate $G|_U$ to $G_U$ throughout this subsection; all heaps still omit output records.

\begin{figure}[htbp]
\centering
\begin{tikzcd}[column sep=large,row sep=large]
G_U \arrow[r,hook,"\iota_{U,V}"] \arrow[dr,"\sigma_U"'] & G_V \arrow[d,"\sigma_V"] \\
& \mathcal A_a
\end{tikzcd}
\caption{A partial execution extends without changing completed values. The objects and arrows lie in $\TS_E$. The arrow across the top enlarges the completed structure; the arrows into $\mathcal A_a$ assign its values. Commutativity is the invariant~\eqref{eq:partial-compatibility}.}\label{fig:partial-cone}
\end{figure}

\begin{definition}[Execution diagrams]\label{def:execution-diagram}
Regard $\mathcal L_G$ as a category with a unique arrow $U\to V$ exactly when $U\subseteq V$. Define
\[
J_G:\mathcal L_G\longrightarrow\TS_E/G^-,\qquad
J_G(U)=(G_U\hookrightarrow G^-).
\]
In this slice category an arrow from $i:X\to G^-$ to $j:Y\to G^-$ is a Token map $h:X\to Y$ such that $j h=i$. Composition with the total valuation $v_{G,a}:G^-\to\mathcal A_a$ gives a second diagram
\begin{equation}\label{eq:valued-execution-functor}
\mathcal E_{G,a}:\mathcal L_G\longrightarrow\TS_E/\mathcal A_a,
\qquad U\longmapsto(G_U\xrightarrow{\sigma_U}\mathcal A_a).
\end{equation}
Both functors send $U\subseteq V$ to the inclusion $\iota_{U,V}$.
\end{definition}

The first diagram retains where every node occurs in the fixed graph; the second retains its value. These are different forms of information. Functoriality follows from composition of inclusions and the valuation invariant. In particular, maps between objects of $\TS_E/\mathcal A_a$ need not preserve node identity merely because they preserve values.

\begin{theorem}[Execution gluing]\label{thm:execution-gluing}
The functor $J_G$ is fully faithful. For any $U,V\in\mathcal L_G$, the inclusion square in Figure~\ref{fig:execution-gluing} is both a pullback and a pushout in $\TS_E$, and also in $\TS_E/G^-$. In particular,
\begin{equation}\label{eq:execution-amalgamation}
G_{U\cup V}\ \cong\ G_U\sqcup_{G_{U\cap V}}G_V.
\end{equation}
The valuations on $G_U$ and $G_V$ agree on $G_{U\cap V}$ and glue uniquely to $\sigma_{U\cup V}$. The square with these valuations is also a pushout in $\TS_E/\mathcal A_a$.
\end{theorem}

\begin{figure}[htbp]
\centering
\begin{tikzcd}[column sep=large,row sep=large]
G_{U\cap V} \arrow[r,hook] \arrow[d,hook] & G_U \arrow[d,hook] \\
G_V \arrow[r,hook] & G_{U\cup V}
\end{tikzcd}
\caption{Two completed parts are joined along their common completed part. This square is bicartesian: it has both the pullback and pushout universal properties. These properties refer to entire Token classes and their maps; the square does not assert a hardware synchronization mechanism.}\label{fig:execution-gluing}
\end{figure}

\begin{proof}
An arrow $J_G(U)\to J_G(V)$ must commute with the inclusions into $G^-$. Its carrier function is consequently the inclusion of $U$ in $V$, if that inclusion exists, and no such function exists otherwise. When $U\subseteq V$ it preserves all records. Thus the hom sets agree with those of $\mathcal L_G$.

Write $H_U$ for the heap of $G_U$. Every retained non-input record is identified by its result node, so
\[
H_{U\cap V}=H_U\cap H_V,\qquad
H_{U\cup V}=H_U\cup H_V.
\]
For the pullback, carrier maps into $U$ and $V$ with equal composites into $U\cup V$ factor uniquely through $U\cap V$. The image of any source record lies in both heaps, so the factor is a Token map. For the pushout, two Token maps from $G_U,G_V$ to an arbitrary $T\in\TS_E$ which agree on $G_{U\cap V}$ give a unique function $U\cup V\to\base(T)$. Each source record lies in one of the two heaps and is preserved there; hence the glued function is a Token map. This proves both universal properties directly.

The same factorizations commute with the maps to $G^-$, proving the slice assertions. By~\eqref{eq:partial-compatibility}, the two partial valuations agree on the intersection. Their unique gluing is the restriction of $v_{G,a}$ to the union, hence equals $\sigma_{U\cup V}$. For maps in the value slice, the glued map also commutes with the maps into $\mathcal A_a$, giving its pushout property.
\end{proof}

This result identifies meets of completed sets with pullbacks over $G^-$, and joins with unions amalgamated over the meet. It does not identify a union with the ordinary coproduct $G_U\sqcup G_V$, which would keep separate copies of shared nodes. Nor does it assert that $\mathcal E_{G,a}$ is fully faithful: an additional map may identify distinct nodes that have equal values.

\begin{figure}[htbp]
\centering
\begin{tikzcd}[column sep=large,row sep=large]
G_U \arrow[r,hook,"+u"] \arrow[d,hook,"+v"'] & G_{U\cup\{u\}} \arrow[d,hook,"+v"] \\
G_{U\cup\{v\}} \arrow[r,hook,"+u"'] & G_{U\cup\{u,v\}}
\end{tikzcd}
\caption{The execution diamond for two simultaneously ready nodes. The labels $+u,+v$ name extensions by a completed operation; the arrows themselves are inclusions. Together with Figure~\ref{fig:partial-cone}, the square says that either order, or their ready batch, gives the same structure and the same values.}\label{fig:ready-diamond}
\end{figure}

\begin{corollary}[Ready cubes]\label{cor:ready-cubes}
For a finite ready set $R\subseteq\mathcal R_G(U)$, the objects $G_{U\cup A}$, indexed by $A\subseteq R$, form a $|R|$-dimensional cube of inclusions. Every two-dimensional face is an execution-gluing square, and all its partial valuations are compatible. Three ready nodes give an ordinary cube; the Boolean interval~\eqref{eq:ready-boolean-interval} supplies the general case.
\end{corollary}
\begin{proof}
For $A,B\subseteq R$, intersections and unions of $U\cup A,U\cup B$ are $U\cup(A\cap B)$ and $U\cup(A\cup B)$. Apply Theorem~\ref{thm:execution-gluing} and the unique partial valuation.
\end{proof}

\begin{example}[Two independent XORs]\label{ex:gluing-xor}
Let $u=x_1\oplus x_2$, $v=x_3\oplus x_4$, and $y=u\oplus v$, with all four inputs initially available. The diamond starts at the input set $I_G$ and ends at $I_G\cup\{u,v\}$. The final operation $y$ becomes ready only at that endpoint. For input $(1,0,1,1)$, the two intermediate results are $u=1,v=0$ in either order and the final result is $1$. The common inputs occur once in the glued class, and the valuation is the unique map to $\mathcal A_a$. With one worker the diamond is traversed along two edges; with two workers its two ready operations may finish as one batch.
\end{example}

Order ideals and configuration-based accounts of concurrency are established ideas in event-structure semantics \cite{winskel1987}. The present finite pure-DAG fragment has no conflict relation or alternative-event semantics. Its additional construction is an explicit realization by Token classes, their universal gluing maps, and their compatible valuations. Pushouts also have an established role in categorical graph rewriting \cite{lack2003}. The theorem proves the indicated squares; it does not establish all the stability and Van Kampen conditions needed for an adhesive-category account of arbitrary rewrites.

\subsection{Work, span and a changing number of workers}\label{sec:elastic-scheduling}

For the following bounds, all input values are initially available, each operation takes one time slot on one worker, data reads and result routing are free, and there is enough memory for the described schedule. Define $W(G)=|O_G|$ and let $S(G)$ be the maximum number of operation nodes on a dependency path, with $S(G)=0$ when there are none: the \emph{work} and \emph{span} of this representation under this evaluator. A coarse node cannot be assigned unit cost merely because it denotes one finite mapping; its evaluator must justify that cost or expose its internal work. The work--span model and greedy scheduling estimates are standard tools of parallel algorithm analysis \cite{blelloch1999}; we state them for the present graph language, with proofs in Appendix~\ref{app:scheduling}, and separate variable resource supply from a change of representation.

\begin{proposition}[Fixed-resource bound]\label{prop:work-span}
For $p\ge1$ workers, every schedule completing $G$ takes at least $\max\{W/p,S\}$ slots. A greedy schedule, which executes $\min\{p,|\mathcal R_G(U)|\}$ nodes in each unfinished slot, satisfies
\begin{equation}\label{eq:work-span}
T_p\le \frac{W}{p}+\left(1-\frac1p\right)S
\le \frac{W}{p}+S.
\end{equation}
\end{proposition}

\begin{proposition}[Variable resource supply]\label{prop:variable-capacity}
Suppose slot $t$ has $p_t$ workers, where $0\le p_t\le P$ for a fixed integer $P\ge1$, and a greedy scheduler executes up to $p_t$ ready nodes. Define the supplied capacity through a horizon $h$ by $C(h)=\sum_{t=0}^{h-1}p_t$. Execution has finished by any horizon satisfying
\begin{equation}\label{eq:capacity-bound}
C(h)\ge W+(P-1)S.
\end{equation}
In particular, unbounded cumulative capacity implies completion, and the output is independent of the resource-supply sequence.
\end{proposition}

Both are bounds on the ideal schedule for a fixed graph, not on a distributed machine's measured latency; unequal operation times, fan-out and communication require a weighted model, and peak memory depends on the schedule rather than on $W$ and $S$ alone. The second proposition allows workers to be added or withdrawn between slots without changing the denotation, provided completed values remain available; capacity that arrives only after a long outage cannot shorten that outage. Migration of an operation in flight, and failure recovery, are separate protocols.

Categorically, a schedule through $T$ slots is a monotone path $\rho:\{0<1<\cdots<T\}\to\mathcal L_G$, subject to additional operational conditions: each increment is a subset of the ready set at its source and has size at most the supplied worker count. Waiting is an identity step. Changing worker supply selects a different such path through the same execution diagram $\mathcal E_{G,a}$. Monotonicity alone does not imply readiness or a resource bound; the operational conditions and progress assumptions remain essential. Figure~\ref{fig:ready-diamond} already displays the smallest case of choosing different paths or a ready batch.

\subsection{Four meanings of computational elasticity}\label{sec:elasticity-levels}

\begin{table}[htbp]
\centering\small\renewcommand{\arraystretch}{1.13}
\begin{tabularx}{\linewidth}{@{}>{\raggedright\arraybackslash}p{0.18\linewidth}>{\raggedright\arraybackslash}p{0.29\linewidth}>{\raggedright\arraybackslash}X@{}}
\toprule
What changes & What stays fixed & Required mechanism or invariant \\
\midrule
Worker allocation & The graph, its input and output behavior & Readiness-based scheduling, retained values and continued progress. \\
Active work & The requested result of one invocation & Demand propagation, branch selection, node lifetimes and valid residual computations. \\
Representation and granularity & The denoted mapping at the chosen interface & Certified implementations in $\mathsf{Impl}(f)$; a migration condition for a running invocation. \\
Input or interface size & A specified uniform family or protocol & Effective construction of the family and any required cross-size compatibility. \\
\bottomrule
\end{tabularx}
\caption{Computational elasticity has several distinct objects. More workers, more active nodes, a larger description and a larger input domain are different changes.}\label{tab:elasticity-levels}
\end{table}

Table~\ref{tab:elasticity-levels} separates changes that can otherwise be hidden by the phrase ``dynamic scale.'' The propositions just proved concern its first row. A graph may expose many ready operations while running on one worker, or expose only a serial chain while running on many workers. Similarly, a large stored description may have a small active part for a particular request.

For a prescribed set of output ports $J$, let $A_J$ contain their output nodes and all their ancestors. Keeping $A_J\cup I_G$, all defining records on those nodes, every declared input, and only the requested output records yields a valid graph for those outputs, after renumbering output ports if necessary. Topological evaluation agrees with the corresponding coordinates of $\delta(G)$. Thus irrelevant work can be omitted for that request. This restriction changes the requested output interface unless all original outputs are retained; it is not automatically an element of $\mathcal I(\delta(G))$ with the original interface.

For input-dependent activity, a representation language may contain a conditional description $\operatorname{if}(b,r_0,r_1)$, with $b:A\to\{0,1\}$ and $r_0,r_1:A\to B$. Its evaluator first evaluates the guard and then the chosen branch, with denotation $a\mapsto\delta(r_{b(a)})(a)$. A tagged sum can carry the selected alternative into a copairing. This lazy evaluator must be specified: an ordinary total selector node in a strict graph does not itself prevent both of its arguments from being computed. General branching, loops and recursion are therefore extensions of the acyclic graph fragment, with separate progress obligations.

Finite interfaces impose no common bound on all descriptions or all invocations. A machine can generate successively larger finite structures through effective allocation rules. If both storage and all interface sizes are globally bounded, its configurations instead form the finite state space discussed in Section~\ref{sec:machine-specification}. Increasing an input size from $n$ to $n+1$ concerns a family $f_n$, not just two implementations of one fixed $f$. Prefix compatibility such as~\eqref{eq:prefix-natural} is appropriate for sequence computations when earlier outputs must remain stable; other applications require their own compatibility condition. Merely collecting the family in a coproduct does not prove such a condition.

\subsection{Granularity, expansion and certified representation changes}\label{sec:granularity}

A stored finite mapping can be treated as a single callable module or represented by a network of smaller mappings. Atomicity is consequently relative to an interface and an evaluator, just as a tree Token is atomic as an entry of an enclosing Token and composite as a tree. To use a module description $r$ as a carrier element we choose an explicit code; to expose its operations we choose a valid body and an expansion procedure. Section~\ref{sec:module-expansion} defines these choices, including a restricted compilation route from $\Finv$ to primitive graphs. Fusion hides a suitable region behind its boundary interface; its cost includes the evaluator of that module. Token Space supplies a common language for the records involved, while the denotation-preservation theorem states when a change of level is computationally correct.

A rewrite in $\mathsf{Impl}(f)$ may decrease description length, increase it to expose useful structure, or exchange work for locality. Reassociation under a valid algebraic law, block decomposition, table compilation and a change between two exact evaluators can all be certified at an interface. Arbitrary maps between syntax carriers do not certify all these transformations. Graph Token maps give the particular sufficient condition of Theorem~\ref{thm:graph-evaluation}; more general changes use other certificates in $\mathsf{Impl}(f)$.

For a fixed execution platform and accounting protocol, it is useful to retain the resource profile
\begin{equation}\label{eq:resource-profile}
\mathcal C(r,\pi,a)=
\bigl(L(r),W(r,a),S(r,a),M(r,\pi,a),Q(r,\pi,a)\bigr),
\end{equation}
where $\pi$ is a schedule and placement, $M$ is peak live storage and $Q$ is communication under the specified model. For straight-line graphs the operation counts and structural span above are input-independent; conditional evaluators need not have that property. There is usually no single implementation minimizing all components. Policies can choose among a finite catalog of certified implementations using predicted latency, measured costs or a memory constraint. Correctness of their choices follows from the certificates; quality of the predictions requires measurement.

Sharing syntactically identical pure subterms reduces repeated work, but does not by itself imply a shorter span. In the occurrence-to-subterm construction, a subterm has the same recursively determined depth at every occurrence, so shared and occurrence graphs have the same maximum operation depth. Physical consequences can differ: a shared result has more consumers, can need longer retention, and may require communication to remote consumers. Selective duplication can improve locality in some deployments while increasing work and storage. No such benefit follows without a placement and communication model.

Description construction, rewrite search, certificate checking and migration all consume resources. For $K$ identical invocations under a fixed resource regime, a proposed change with additional preparation cost $C_{\mathrm{change}}$ and per-invocation times $T_{\mathrm{old}},T_{\mathrm{new}}$ has a positive time saving only if
\[
C_{\mathrm{change}}<K(T_{\mathrm{old}}-T_{\mathrm{new}}).
\]
This elementary accounting also explains why a larger table can be worthwhile for a frequently reused small mapping while being unsuitable for a large input interface. Allowing expansion removes the strict description-length descent used to prove termination of compression-only rewriting in Section~\ref{sec:reencoding}. An elastic policy therefore needs another progress condition or a finite transformation budget.

\subsection{Explicit module codes and denotation-preserving expansion}\label{sec:module-expansion}

To make relative atomicity a construction, we now specify how a finite graph becomes a callable component, how its description becomes data, and when that component can be expanded. This gives a first precise connection between the representation systems of Section~\ref{sec:representation-system} and the change of descriptive level envisaged in $\ATC$. We use single-output modules over the same nonempty finite value set $D$. An invocation graph may still have several output ports. A module with several outputs would require a product-valued interface or a record defining several result nodes; that extension is not silently identified with the single-result graph syntax used here.

\begin{definition}[A ranked module library]\label{def:module-library}
Fix the primitive library $\Omega$ of Definition~\ref{def:computation-graph}. A finite \emph{ranked module library} is a list
\[
\mathcal M=((m_1,k_1,B_1),\ldots,(m_q,k_q,B_q))
\]
with fresh operation symbols $m_j$, arities $k_j\ge0$, and bodies $B_j$ satisfying the following conditions. The body $B_j$ is a valid finite acyclic graph with $k_j$ ordered input ports and one output port; its operation labels belong to $\Omega\sqcup\{m_1,\ldots,m_{j-1}\}$. All declared input nodes are retained, including unused ones. Recursively interpret
\begin{equation}\label{eq:module-meaning}
(m_j)_D=\delta_{\mathcal M_{<j}}(B_j):D^{k_j}\longrightarrow D,
\end{equation}
where the singleton output tuple is identified with its component. A \emph{module graph} is a valid graph over $\Omega\sqcup\{m_1,\ldots,m_q\}$, evaluated with these interpretations. Its denotation is written $\delta_{\mathcal M}$.
\end{definition}

Each body is finite and calls only lower-ranked modules, so the recursive interpretation is well defined and total. The unique-valuation and schedule-independence theorems apply to a module graph with these symbols treated as interpreted operations. That application establishes correctness at the chosen granularity; it does not make each module call a unit of physical time.

Choose an injective, effectively decodable encoding of finite records, ordered port lists and operation labels. A description includes its module library, or an immutable identifier resolving that library. Let $\mathcal R_{n,m}$ be the set of well-formed encoded module graphs with the indicated interface. These codes are ordinary carrier elements: for example they may form the carrier of a class with an empty heap, and their evaluation records may be added on the tagged carrier
\[
\mathcal R_{n,m}\sqcup D^n\sqcup D^m.
\]
Use a fresh core symbol $\mathrm{eval}$ and the heap
\begin{equation}\label{eq:module-evaluation-records}
\{(\mathrm{eval},\iota_R(r),\iota_A(a),\iota_B(b)):
 b=\delta_{\mathcal M}(r)(a)\}.
\end{equation}
The tags distinguish descriptions, inputs and outputs even if their underlying encodings coincide. A finite catalog has a finite evaluation heap; the set of all finite codes need not be finite. These are specification records, not a requirement to materialize an exhaustive table in an implementation. Explicit coding, rather than the iteration of $\ATC$ alone, permits a graph or Token class to occur as data in a further Token class.

\paragraph{Expansion.}
For an occurrence $(m_j,v,u_1,\ldots,u_{k_j})$, read the body $B_j$ from the library and make a fresh copy of every internal operation node. Wire its ordered formal inputs to $u_1,\ldots,u_{k_j}$; these actual operands need not be distinct. Replace each use of $v$, including its uses at output ports, by the copied body's output node. Remove the call record and its former result node. If the body's output is one of its inputs, this replacement is an alias to the corresponding actual operand. Retain the surrounding graph's declared inputs. Internal input and output records are used for wiring and are then removed; the surrounding graph keeps its ordered external ports. Recursively expand lower-ranked calls to obtain a primitive graph, denoted $\operatorname{Expand}_{\mathcal M}(G)$.

An equivalent constructive implementation traverses the surrounding graph in topological order. It maintains a map from old nodes to newly built nodes, emitting a primitive node directly or recursively copying a module body at a call. The image of a call node is its body's output. This avoids assigning two definitions to a node when a body merely returns an argument. It also preserves repeated uses and output aliases.

\begin{theorem}[Correct expansion]\label{thm:module-expansion}
For every ranked module library and every valid module graph $G:D^n\to D^m$, expansion terminates, yields a valid primitive graph with the same ordered external interface, and satisfies
\begin{equation}\label{eq:expansion-correctness}
\delta_{\Omega}\bigl(\operatorname{Expand}_{\mathcal M}(G)\bigr)
=\delta_{\mathcal M}(G).
\end{equation}
Expanding any selected collection of call occurrences, leaving the others interpreted by their library, also preserves the denotation.
\end{theorem}
\begin{proof}
Induct on the largest module rank. A primitive graph requires only fresh naming. In a body of rank $j$, every recursive expansion has smaller rank; each body and surrounding graph has finitely many nodes. Thus the recursive copying procedure terminates and creates finitely many nodes.

For validity, process the surrounding graph in a topological order and each copied body in its own topological order. Its actual inputs are already available in the node map. Every new primitive node has fresh identity, one defining record and operands that precede it. Argument aliases create no new definitions. Every external input remains declared, and each external output is mapped to an existing node. This gives an acyclic graph with the required unique definitions and ports.

For correctness, fix an external input and induct through this traversal. Primitive nodes preserve their values by the fixed interpretation of $\Omega$. At a call of $m_j$, the actual operands have their original values by the outer induction. By the induction on rank, the copied body returns $(m_j)_D$ of exactly those values, including when its result aliases an argument. Replacing all uses of the call result therefore preserves their values. Reading the ordered output ports proves~\eqref{eq:expansion-correctness}. The same argument at each selected occurrence proves the last assertion; the unexpanded calls retain their stipulated interpretations.
\end{proof}

\begin{figure}[htbp]
\centering
\begin{tikzcd}[column sep=large,row sep=large]
\mathcal R_{\mathcal M}^{n,m}\times D^n
 \arrow[r,"\operatorname{Expand}_{\mathcal M}\times\id"] \arrow[d,"\operatorname{eval}_{\mathcal M}"']
 & \mathcal R_{\Omega}^{n,m}\times D^n \arrow[d,"\operatorname{eval}_{\Omega}"] \\
D^m \arrow[r,equal] & D^m
\end{tikzcd}
\caption{Module expansion commutes with evaluation in $\Set$. Here $\mathcal R_{\mathcal M}^{n,m}$ and $\mathcal R_{\Omega}^{n,m}$ are valid coded descriptions at the fixed interface, for the fixed library and for primitive graphs respectively. Theorem~\ref{thm:module-expansion} proves this square. These description sets need not be finite even though every description and value interface is finite.}\label{fig:expansion-evaluation}
\end{figure}

The horizontal arrow transforms descriptions, and the vertical arrows evaluate them on inputs. These are maps between different kinds of objects. This diagram is not automatically a square of graph Token maps: its promotion to an enriched description category would require chosen heaps and preservation proofs.

\paragraph{Packing and fusion.}
For a primitive graph $B:D^k\to D$, add a fresh module with body $B$. Its \emph{packed invocation} has $k$ declared inputs, one call node and that node as output. Expansion recovers a graph isomorphic to $B$ by an input- and output-preserving renaming of internal nodes. Conversely, packing an expanded graph gives an equivalent callable description, not the original library, module boundaries or names. Packing and expansion are therefore not inverse functions on all descriptions. Nor is packing an inverse of $\Finv$, whose domain is a supported tree class and whose result is a syntactic record structure.

Fusion of a region in a larger graph requires an interface condition. One sufficient single-output case is a finite set $S$ of operation nodes with a designated node $z\in S$, such that every node of $S$ reaches $z$ along a path contained in $S$ (allowing the length-zero path), and every use of an $S$-node outside $S$, including an external output use, is a use of $z$. List the distinct operands entering $S$ from outside as boundary inputs. The induced body on $S$ and these boundary inputs is then a valid single-output module. Replacing $S$ by its call preserves behavior, and expanding the call restores the region up to fresh naming. The quotient graph is acyclic: a newly introduced cycle would give a path from $z$ through the outside back into $S$; the reentered node has a path within $S$ to $z$, so this would already be a cycle in the original graph. For several exit values one must first specify a multi-output interface. An arbitrary selection of nodes cannot be fused while ignoring its boundary. In particular, a unique exit alone is insufficient when dead internal work is included: in a path $z\to u\to w$ with output $z$, fusing $\{z,w\}$ would make the new call depend on $u$, while $u$ depends on that call. The reachability condition excludes this case.

\paragraph{From tree reification to executable modules.}
There is also a precise restricted route from $\Finv$ to the executable graph language. Let $t$ be a finite term over the primitive signature, with distinct declared variables $x_1,\ldots,x_k$. Encode a variable as its base letter and an operation application as the tuple
\[
\operatorname{enc}(\omega(t_1,\ldots,t_s))
=(\omega,\operatorname{enc}(t_1),\ldots,\operatorname{enc}(t_s)).
\]
Here operation names are core symbols; a nullary operation is the singleton tuple $(\omega)$. Form the supported tree class $T_t$ whose heap is $\{\operatorname{enc}(t)\}$ and whose base contains precisely the variables occurring in $t$. In $\Finv(T_t)$, each tuple node has a primitive operation as its first child and exactly the required number of operand children. A compiler $C_{\Omega,k}$ checks this shape, replaces the corresponding syntactic node record by the operation record $(\omega,[p],\ulcorner\operatorname{enc}(t_1)\urcorner,\ldots,\ulcorner\operatorname{enc}(t_s)\urcorner)$, adds all $k$ declared input nodes and records, and makes the unique root the output port. In particular, unused declared variables are restored from the external interface, not inferred from the supported class.

\begin{proposition}[Evaluation after term reification]\label{prop:term-compilation}
For every such term and declared interface, $C_{\Omega,k}(\Finv(T_t))$ is a valid primitive graph and
\begin{equation}\label{eq:term-compilation}
\delta_\Omega\bigl(C_{\Omega,k}(\Finv(T_t))\bigr)(a)
=\operatorname{val}_D(t,a),\qquad a\in D^k.
\end{equation}
Here $\operatorname{val}_D(t,a)$ is the usual recursive interpretation of the term. Equal syntactic subterms share a single result node.
\end{proposition}
\begin{proof}
A tuple node depends only on its proper term subtrees, so increasing term size provides a topological order. Distinct tuple subterms have distinct tagged identifiers and exactly one defining operation record. Variables have input records, including the separately retained unused inputs. Induction on subterms identifies each node value with its term value. A variable root yields the corresponding input value, and a nullary root yields its stipulated constant; hence the output assertion includes these boundary cases. Sharing is exactly the identification of equal tuple subterms in Definition~\ref{def:reification}.
\end{proof}

The two results supply a chain of constructions with explicit boundaries:
\[
\text{term syntax}\ \xrightarrow{\ \Finv\ }\
\text{syntactic records}\ \xrightarrow{\ C_{\Omega,k}\ }\
\text{valid graph}\ \xrightarrow{\ \mathrm{pack}\ }\
\text{callable module}.
\]
The first arrow exposes syntax; the next uses an interpreted signature and an ordered interface; the last adds a library entry. General tree Tokens can express data unrelated to executable operations, so the compiler is defined only on the stated valid fragment. Likewise, a table or a parameterized model may be a callable description in Section~\ref{sec:representation-system} without an expansion into this particular $\Omega$. Such an expansion requires a supplied body and a correctness certificate; neither finiteness nor reification proves that every primitive library can express every mapping.

\begin{example}[A module boundary can hide available concurrency]\label{ex:module-concurrency}
Over $D=\{0,1\}$ with XOR, define
\[
m(a,b,c,d)=(a\mathbin\oplus b)\mathbin\oplus(c\mathbin\oplus d).
\]
The packed invocation has one call node. Its expanded body has two independent pair operations followed by a final XOR, so primitive work is $3$ and span is $2$. A left-associated body computes the same mapping with work $3$ and span $3$, by associativity. These are different certified descriptions of one behavior. Expanding the balanced body exposes its existing independence; changing the left-associated body into the balanced one additionally requires an algebraic certificate. Counting the packed call as one node says nothing about its evaluator's running time. A scheduler may invoke the module as a unit, or expose its body when finer dispatch is useful, while the same interface specifies the result.
\end{example}

\begin{table}[htbp]
\centering\small\renewcommand{\arraystretch}{1.13}
\begin{tabularx}{\linewidth}{@{}p{0.22\linewidth}p{0.32\linewidth}X@{}}
\toprule
Change of description & Required data & What is established \\
\midrule
Tree reification & A supported tree Token class. & Ordered child records and shared subterms; no evaluator is assumed. \\
Term compilation & Primitive term shape, interpretation and declared ports. & A valid graph with the term's value. \\
Module expansion & A ranked library, fresh nodes and ordered port wiring. & A valid graph with the call's denotation. \\
Region fusion & A valid region and its complete boundary interface. & A callable component with the same boundary behavior. \\
Online replacement & A consistent completed frontier and an equivalent residual implementation. & The original invocation's remaining result is preserved. \\
\bottomrule
\end{tabularx}
\caption{Relative atomicity requires an explicit construction at each level. The last row is proved in Section~\ref{sec:online-replacement}; it is not a consequence of packing alone.}\label{tab:module-bridges}
\end{table}

A hierarchy can share module bodies in storage even when full expansion copies them at many call sites. Repeated nested calls can make the expanded graph much larger than its library description. Expansion is therefore a choice of execution granularity, not a universal compression method. Lazy expansion before dispatch can instantiate only the selected call occurrences. Changing a component already in execution additionally needs the completed-boundary rule below. In both cases, the runtime must account for expansion, scheduling and retained values; denotational equality alone supplies no performance gain.

\subsection{A live frontier and the residual mapping}\label{sec:residual-mapping}

Equality of whole input--output functions is sufficient for replacing an unstarted component. To replace part of a running invocation, the machine must identify the values that connect completed work to the remaining work. The graph language makes this boundary explicit.

\begin{definition}[Live frontier and residual graph]\label{def:residual}
For a partial execution state $(U,\sigma)$, define
\begin{equation}\label{eq:live-frontier}
\partial_G U=\{u\in U:
u\text{ is an output node, or }u\in\operatorname{pred}_G(w)
\text{ for some }w\notin U\}.
\end{equation}
Choose an ordering $(b_1,\ldots,b_k)$ of these distinct nodes and put
$\beta_U=(\sigma(b_1),\ldots,\sigma(b_k))$.
The residual graph $G/U:D^k\to D^m$ has carrier
$(N_G\setminus U)\cup\partial_G U$.
Its input records are $(\mathrm{in}_i,b_i)$, its operation records are exactly those with result outside $U$, and its output records are the original ones. Old input records and operation definitions of boundary nodes are removed. Interface labels for all finite arities are taken from a common core.
\end{definition}

An output already computed remains an input of the residual graph, so that its value is still available for final reporting. The boundary may be empty; $D^0$ then denotes a singleton input interface. The definition retains unfinished dead operations if present; they may subsequently be removed by an independently correct output-demand restriction.

\begin{proposition}[Residual correctness and reclamation]\label{prop:residual}
The graph $G/U$ is valid, and
\begin{equation}\label{eq:residual-correctness}
\delta(G/U)(\beta_U)=\delta(G)(a).
\end{equation}
The values at nodes of $U\setminus\partial_G U$ are unnecessary for this remaining evaluation and final output. They may be discarded in this pure graph model once no in-flight read or external observer still requires them.
\end{proposition}
\begin{proof}
Every unfinished node is an operation node. Each operand either is unfinished or belongs to the frontier by definition. Every original output either is unfinished or belongs to the frontier. Nodes in the frontier are uniquely defined as new inputs, and all other retained nodes keep their unique operation definition. Removing completed nodes and edges into the new inputs cannot introduce a cycle. Restrict $v_{G,a}$ to the retained carrier. By Lemma~\ref{lem:partial-valuation} it has input tuple $\beta_U$, and it satisfies every retained operation equation. Uniqueness in Theorem~\ref{thm:graph-evaluation} proves the displayed equation. A discarded completed node is neither an output nor an operand of any unfinished operation, proving the reclamation claim under the stated runtime conditions.
\end{proof}

This is an explicit use of the construction operators: restrict the node carrier, delete obsolete definition and input records, add boundary input records, and keep the remaining operation and output records. Heap deletion alone would not establish correctness; the partial valuation invariant and the selected boundary do so. The boundary size is also a first logical memory measure. Actual peak memory includes in-flight operands, temporary values, stored descriptions and runtime metadata.

\subsection{Frontier interfaces form a functor}\label{sec:frontier-functor}

The frontier of one completed set supplies a replacement interface. Comparing several completed sets supplies more structure: the interfaces themselves have coherent transition maps. Choose once an ordering of each $\partial_G U$, and write
\[
C_U=D^{|\partial_G U|},\qquad f_U=\delta(G/U):C_U\to D^m.
\]
These are finite sets and functions. Their dimensions may increase or decrease as execution progresses; a frontier is not generally included in the next frontier.

For $U\subseteq V$, define an \emph{advancement graph} $K_{U,V}:C_U\to C_V$. Its nodes are $\partial_G U\cup(V\setminus U)$, its inputs are the ordered $\partial_G U$, and its operation records have exactly the result nodes in $V\setminus U$. Its outputs are the ordered nodes of $\partial_G V$. Define
\begin{equation}\label{eq:frontier-transition}
b_{U,V}=\delta(K_{U,V}):C_U\longrightarrow C_V.
\end{equation}
Old frontier values that cease to be needed can be discarded by this map; newly computed live values enter its output. Input and output records are renamed to the indicated interfaces, as in the residual construction.

\begin{theorem}[Frontier composition and residual factorization]\label{thm:frontier-functor}
Every $K_{U,V}$ is a valid graph. The assignments $U\mapsto C_U$ and $(U\subseteq V)\mapsto b_{U,V}$ define a functor
\[
\mathcal B_G:\mathcal L_G\longrightarrow\mathbf{FinSet}.
\]
For $U\subseteq V\subseteq W$,
\begin{equation}\label{eq:frontier-composition}
b_{U,U}=\id_{C_U},\qquad
b_{U,W}=b_{V,W}\circ b_{U,V},\qquad
f_U=f_V\circ b_{U,V}.
\end{equation}
Consequently the residual output maps form a natural transformation
$f:\mathcal B_G\Rightarrow\Delta(D^m)$ to the constant functor.
\end{theorem}

\begin{figure}[htbp]
\centering
\begin{tikzcd}[column sep=large,row sep=large]
C_U \arrow[r,"b_{U,V}"] \arrow[dr,"b_{U,W}"'] & C_V \arrow[d,"b_{V,W}"] \\
& C_W
\end{tikzcd}
\qquad
\begin{tikzcd}[column sep=large,row sep=large]
C_U \arrow[r,"b_{U,V}"] \arrow[d,"f_U"'] & C_V \arrow[d,"f_V"] \\
D^m \arrow[r,equal] & D^m
\end{tikzcd}
\caption{Coherence of changing frontier interfaces, in $\mathbf{FinSet}$. Left: advancing in two stages equals advancing directly. Right: advancing the boundary and then finishing preserves the residual output. Commutativity holds for every input of $C_U$, including boundary tuples not reachable from the original external input.}\label{fig:frontier-coherence}
\end{figure}

\begin{proof}
If an operand of a node in $V\setminus U$ lies in $U$, it belongs to $\partial_G U$; otherwise predecessor closure puts it in $V\setminus U$. A node of $\partial_G V$ that already lies in $U$ is an original output or has a consumer outside $V$, hence outside $U$, so it belongs to $\partial_G U$. All required operands and output nodes are therefore present in $K_{U,V}$. Giving boundary nodes input records and retaining the other defining records establishes unique definitions and acyclicity.

Fix any tuple in $C_U$ and evaluate $G/U$. Restricting this valuation to $\partial_G U\cup(V\setminus U)$ gives the unique valuation of $K_{U,V}$, so its output is the boundary tuple at $V$. Subsequent evaluation of $G/V$ agrees with the original valuation on every remaining node by topological induction. This proves $f_U=f_V b_{U,V}$ without assuming that the tuple came from an external input of $G$.

The same argument restricted to the nodes through $W$ proves that $K_{U,W}$ and the composite of $K_{U,V}$ with $K_{V,W}$ have the same ordered output values. For $U=V$, the graph consists of the frontier inputs wired to the same ordered outputs and has no operation nodes. These observations prove identity and composition, and the last equality is exactly naturality of $f$.
\end{proof}

This proof also gives a structural residual law. In $G/U$, the completed set corresponding to $V$ is
\[
W_{U,V}=\partial_G U\cup(V\setminus U).
\]
Its frontier is exactly $\partial_G V$, and, using the fixed ordering at $V$,
\begin{equation}\label{eq:iterated-residual}
(G/U)/W_{U,V}=G/V
\end{equation}
as labelled graphs with the same original node identifiers. Indeed, both retain precisely the operation nodes outside $V$, make $\partial_G V$ the inputs and preserve the original output ports. This law connects the functorial equation to the actual record construction.

\begin{example}[An interface can shrink]\label{ex:frontier-shrinking}
For the graph in Example~\ref{ex:gluing-xor}, choose the boundary order $(x_1,x_2,x_3,x_4)$ at $I_G$, $(u,x_3,x_4)$ after $u$, $(u,v)$ after both pair operations, and $(y)$ at completion. The frontier diagram includes
\[
D^4\longrightarrow D^3\longrightarrow D^2\longrightarrow D,
\]
with maps $(a,b,c,d)\mapsto(a\oplus b,c,d)$, $(e,c,d)\mapsto(e,c\oplus d)$ and $(e,h)\mapsto e\oplus h$. The other order of pair evaluation gives the same composite. Forgetting the consumed input values is justified by the frontier, while the computation's output interface remains fixed. Other graphs can increase the number of live values before it decreases.
\end{example}

The execution diagram $J_G$ grows completed structure by inclusion. The frontier functor $\mathcal B_G$ retains only the changing interface needed to continue. It need not recover a completed subgraph from its interface: different completed sets can have the same frontier dimension, and the function $b_{U,V}$ does not record the work, communication or memory needed to realize it. Keeping its graph description $K_{U,V}$ alongside the function supplies that implementation information.

\begin{remark}[Typed numerical frontiers]\label{rem:typed-frontier}
The evaluation and frontier proofs above also apply to a finite acyclic
presentation whose occurrence $v$ has a specified nonempty value set
$D_v$ and whose operations are total maps of their declared sorts.
Replace $D^{|U|}$ by the product of those value sets. Topological
induction gives the same residual factorization and composition laws
in $\Set$; finiteness of a value set is not used in these proofs.
This extension is used for numerical decoder states below. It does not
transfer the finite-table algorithms or resource bounds to real-valued
state spaces. Already completed outputs required by the final interface
remain in the frontier under the same definition.
\end{remark}

\subsection{Safe replacement during an invocation}\label{sec:online-replacement}

A graph Token map can support a more direct transfer of already computed values before any residual representation is rebuilt.

\begin{proposition}[State transfer along a graph Token map]\label{prop:graph-state-transfer}
Let $q:G\to G'$ be a Token map between valid graphs on the same input and output interfaces. For a partial execution state $(U,\sigma)$ of $G$ on $a$, set $U'=q(U)$. Then
\begin{equation}\label{eq:graph-state-transfer}
\sigma'(q(u))=\sigma(u)\qquad(u\in U)
\end{equation}
defines a partial execution state $(U',\sigma')$ for $G'$ on $a$. Continuing that state to completion returns $\delta(G)(a)$.
\end{proposition}
\begin{proof}
Theorem~\ref{thm:graph-evaluation} gives $v_{G,a}=v_{G',a}\circ q$, and Lemma~\ref{lem:partial-valuation} identifies $\sigma$ with the restriction of $v_{G,a}$. Thus values agree whenever two nodes of $U$ have the same $q$-image, making $\sigma'$ well defined. Input records ensure that $q(U)$ contains every target input with the correct value. An operation record at $u\in U$ maps to the defining operation record at $q(u)$; validity of $G'$ ensures that it is that node's unique definition. Its operand nodes are the images of operands in $U$, hence belong to $q(U)$. Input nodes cannot also be operation nodes in either valid graph. It follows that $U'$ is predecessor-closed and its values satisfy all completed equations. The transferred state is consequently a partial execution state, and schedule independence and equality of graph denotations prove the output assertion.
\end{proof}

\paragraph{The adjunction behind sharing.}
The state-transfer construction has an order-theoretic form. It distinguishes the completed work visible in a shared graph from the source occurrences whose values become available through that work.

\begin{theorem}[Completed-set adjunction and saturation]\label{thm:sharing-adjunction}
For a graph Token map $q:G\to G'$ as in Proposition~\ref{prop:graph-state-transfer}, direct and inverse image restrict to monotone maps
\[
q_!:\mathcal L_G\longrightarrow\mathcal L_{G'},\quad U\longmapsto q(U),
\qquad
q^*:\mathcal L_{G'}\longrightarrow\mathcal L_G,\quad V\longmapsto q^{-1}(V).
\]
They are adjoint as functors between poset categories:
\begin{equation}\label{eq:completed-adjunction}
q_!U\subseteq V\quad\Longleftrightarrow\quad U\subseteq q^*V.
\end{equation}
The operator $c_q=q^*q_!$ is an extensive, monotone and idempotent closure on $\mathcal L_G$. If $q$ is surjective on nodes, then $q_!q^*=\id$ and
\begin{equation}\label{eq:saturated-completions}
\operatorname{Fix}(c_q)\cong\mathcal L_{G'}.
\end{equation}
For any partial valuation $\sigma_U$, the values on $c_q(U)$ are determined without further primitive evaluation by
\begin{equation}\label{eq:saturation-values}
\overline\sigma(v)=\sigma_U(u)\quad
\text{for any }u\in U\text{ with }q(u)=q(v).
\end{equation}
This is the consistent partial valuation on $G_{c_q(U)}$.
\end{theorem}

\begin{figure}[htbp]
\centering
\begin{tikzcd}[column sep=huge]
\mathcal L_G \arrow[r,shift left=1.2ex,"q_!"] & \mathcal L_{G'} \arrow[l,shift left=1.2ex,"q^*"]
\end{tikzcd}
\qquad $q_!\dashv q^*$
\par\medskip
\begin{tikzcd}[column sep=large,row sep=large]
G_U \arrow[r,"q|_U"] \arrow[d,"\sigma_U"'] & G'_{q(U)} \arrow[d,"\sigma'_{q(U)}"] \\
\mathcal A_a \arrow[r,equal] & \mathcal A_a
\end{tikzcd}
\caption{Two levels of a certified sharing transformation. Above: an adjunction between completion lattices, not an inverse between graph descriptions. Below: state transfer is a commuting square in $\TS_E$. The closure $q^*q_!$ adds source occurrences whose values are available from already completed shared nodes.}\label{fig:sharing-adjunction}
\end{figure}

\begin{proof}
The proof of Proposition~\ref{prop:graph-state-transfer} shows that $q(U)$ is predecessor-closed and contains all target inputs. For inverse image, all source inputs map to target inputs. If an operation node $v$ maps into $V$, the image of its defining record is the unique defining record of $q(v)$. Every operand image is consequently a predecessor of $q(v)$, hence belongs to $V$. Thus $q^{-1}(V)$ is predecessor-closed. The set-theoretic equivalence~\eqref{eq:completed-adjunction} now gives the poset adjunction.

For every $U$, $c_q(U)=q^{-1}(q(U))$ is the union of all $q$-fibres meeting $U$. This description proves extensivity, monotonicity and idempotence, and the preceding closure argument keeps it in $\mathcal L_G$. Surjectivity gives $q(q^{-1}(V))=V$. The mutually inverse restrictions of $q_!$ and $q^*$ then yield~\eqref{eq:saturated-completions}. Finally, $v_{G,a}=v_{G',a}\circ q$ makes the right-hand side of~\eqref{eq:saturation-values} independent of the choice of $u$ and equal to $v_{G,a}(v)$. Its restriction is the unique partial valuation on the completed set $c_q(U)$.
\end{proof}

Direct image preserves unions of completed sets. Inverse image preserves both unions and intersections here; it also sends the bottom to the bottom, since an operation node cannot map to an input node in a valid graph. Direct image need not preserve intersections. For a concrete example, let $G$ have input $x$, two nodes $u=\operatorname{NOT}(x)$ and $v=\operatorname{NOT}(x)$, and output ports $(u,v)$. Let $G'$ share the two operations as one node $w$ and have output ports $(w,w)$. With $U=\{x,u\}$ and $V=\{x,v\}$,
\[
q_!(U\cap V)=\{x\},\qquad
q_!U\cap q_!V=\{x,w\}.
\]
Sharing can therefore identify two independent source events. It preserves the stated valuation and adjunction, but must not be claimed to preserve every execution pullback or every dimension of a ready cube.

The saturation adds \emph{available values}, not a historical claim that all their source operations physically ran. The transition from $U$ to $c_q(U)$ is a certified reuse step; it need not be one ordinary ready batch. It depends on the graph map and its operation records, not on accidental numerical equality between unrelated nodes. The adjunction is the classical direct-image/inverse-image adjunction restricted to the paper's legal completion sets. The content to verify is that these sets are closed under both maps and that the resulting closure supports correct state reuse; it is not a new general theorem about adjoints.

For the sharing map of Proposition~\ref{prop:sharing}, one completed occurrence supplies the value of its shared subterm even if other occurrences have not yet executed. This gives a structural certificate for reusing completed work, rather than merely restarting an equivalent program. The runtime still needs a consistent snapshot and ownership of any in-flight tasks before identifying their result locations. Associative reorganization need not admit a graph Token map; for such changes the residual-interface rule below applies.

\begin{theorem}[Replacement at a completed boundary]\label{thm:residual-replacement}
Let $(U,\sigma)$ be a partial execution state for $G$ on $a$. Let $r'$ be a valid description on the same ordered interface as $G/U$, with a certificate that
\begin{equation}\label{eq:residual-rewrite}
\delta(r')=\delta(G/U):D^k\to D^m.
\end{equation}
Starting $r'$ on the retained boundary tuple $\beta_U$ returns $\delta(G)(a)$. If $r'$ is a graph description, its evaluation may use any legal progressing schedule. Repeating such certified boundary replacements a finite number of times preserves the final result, provided the last remaining evaluation completes.
\end{theorem}
\begin{proof}
Equation~\eqref{eq:residual-rewrite} and Proposition~\ref{prop:residual} give
$\operatorname{eval}(r',\beta_U)=\delta(G)(a)$.
For a graph evaluator, Theorem~\ref{thm:schedule-independence} applies. Each later replacement preserves the result of its then-current residual invocation by the same argument; induction over the finite sequence proves the last assertion.
\end{proof}

For a fixed completed set $U$, define the external-to-frontier function
\[
b_U:D^n\to C_U,\qquad
b_U(a)=(v_{G,a}(b_1),\ldots,v_{G,a}(b_k)),
\]
using its fixed boundary order. This is the mathematical function supplying the stored tuple $\beta_U$. The residual-correctness theorem for every external input gives $\delta(G)=f_U\circ b_U$. An adaptation certificate $\delta(r')=f_U$ therefore gives Figure~\ref{fig:online-replacement}.

\begin{figure}[htbp]
\centering
\begin{tikzcd}[column sep=huge,row sep=large]
D^n \arrow[r,"\delta(G)"] \arrow[d,"b_U"'] & D^m \arrow[d,equal] \\
C_U \arrow[r,"\delta(r')=f_U"'] & D^m
\end{tikzcd}
\caption{Replacing unfinished work at a completed boundary, in $\mathbf{FinSet}$. The left arrow summarizes values already obtained; the bottom arrow can use a new implementation of the same residual mapping. At runtime the boundary tuple is retained, so this square does not prescribe recomputation from the original input.}\label{fig:online-replacement}
\end{figure}

For $U\subseteq V$, Theorem~\ref{thm:frontier-functor} also gives $b_V=b_{U,V}b_U$. Thus boundaries within one graph compose coherently. A replacement can change the remaining graph and hence its future completion lattice; subsequent adaptations use the new graph's frontier construction. The proof of repeated replacement uses the common remaining result, and does not assume that all adapted states belong to the original $\mathcal L_G$.

The certificate is deliberately for every boundary input, not merely for the actual tuple currently stored. It is therefore reusable when the same residual structure recurs. A weaker certificate valid on a subset of boundary states is possible only together with an invariant proving that the current boundary belongs to that subset. Such specializations are not supplied by equality on the original external interface alone.

In a concrete concurrent implementation, the boundary must be a consistent snapshot. One simple protocol stops new dispatch, lets issued operations finish, records the completed set and its frontier, constructs and checks the replacement, then resumes. For a subgraph replacement, independent outside work can continue if its interface and ownership conditions are separately respected. The theorem proves the completed-boundary case; it does not establish arbitrary lock-free graph rewriting, distributed snapshot correctness or safe reuse of every old internal value.

For example, consider a four-input parity chain computing $u_2=x_1\oplus x_2$, $u_3=u_2\oplus x_3$, and $y=u_3\oplus x_4$. After $u_2$ completes, a frontier for the single final output is $(u_2,x_3,x_4)$. The residual mapping is
\[
(b,c,d)\longmapsto(b\oplus c)\oplus d.
\]
It may be replaced, using associativity, by $(b,c,d)\mapsto b\oplus(c\oplus d)$. No value for a new internal node is obtained merely by renaming $u_3$; it is computed from the retained boundary. This small example distinguishes a valid mid-execution transformation from an unsupported transfer of internal state. For several independent invocations, choosing an implementation at the start of each invocation is an even simpler safe policy.

\subsection{Exact criteria for changing a retained boundary}\label{v13:mig:criteria}

The replacement theorem keeps the ordered boundary interface fixed. A more flexible machine may change that interface as well: sharing can merge retained values, expansion can expose additional ones, and another representation can organize the remaining work around different intermediate results. Such a change requires a map between \emph{boundary assignments}. A graph Token map acts on node occurrences, whereas a boundary migration acts on assignments of values to a whole frontier. The latter must be constructed and checked; it is not obtained merely by treating the node map as a state map.

Retain the pure, total graph assumptions and the frontier notation of Sections~\ref{sec:residual-mapping}--\ref{sec:online-replacement}. Write $A=D^n$ for the external input space, $Y=D^m$ for the output space, and $F:A\to Y$ for the common denotation of the representations under comparison. For a completed set $U$, distinguish the full Cartesian boundary space from its reachable part:
\begin{equation}\label{v13:mig:reachable}
 C_U=D^{|\partial_G U|},\qquad
 R_U=\im(b_U)\subseteq C_U.
\end{equation}
The tuple currently retained by a correctly executed invocation belongs to $R_U$. Values chosen independently in $C_U$ need not arise from any external input. By Theorem~\ref{thm:frontier-functor}, advancement restricts to a surjection $R_U\to R_V$ whenever $U\subseteq V$: every $b_V(a)$ is $b_{U,V}(b_U(a))$. These restrictions form a subfunctor of $\mathcal B_G$. Increasing the number of stored coordinates can therefore accompany a loss of information about the original input; frontier dimension alone is not an information order.

\begin{definition}[Output-correct boundary migration]\label{v13:mig:output-definition}
For residual mappings $f:C\to Y$ and $f':C'\to Y$, an output-correct migration is a function $m:C\to C'$ satisfying
\begin{equation}\label{v13:mig:output-equation}
 f=f'\circ m.
\end{equation}
One may instead declare invariant subsets $R\subseteq C$ and $R'\subseteq C'$ and require $m:R\to R'$ with $f|_R=(f'|_{R'})m$. The chosen domains are part of the certificate.
\end{definition}

Output-correct migrations compose: if $f=f'm$ and $f'=f''n$, then $f=f''(nm)$. With fixed $Y$, they are precisely the arrows of the familiar slice category $\mathbf{FinSet}/Y$. They also have a direct Token encoding. Fix an evaluation label and distinct core labels $\bar y$ for $y\in Y$, and equip $C$ with
\[
 H_f=\{(\mathrm{eval},c,\overline{f(c)}):c\in C\}.
\]
A function $m:C\to C'$ is a Token map between these evaluation classes exactly when~\eqref{v13:mig:output-equation} holds. This is a compositional specification of migration, while the represented machine retains an executable description of $m$ and its cost. For example, the arrow $f:(C,f)\to(Y,\id_Y)$ satisfies the equation by performing the entire remaining computation. Existence of a semantic arrow therefore does not establish that adaptation is inexpensive.

Output correctness need not reproduce the internal state that the target representation would have reached on the same original input. This stronger requirement has an exact information criterion. For any function $g$, write
\[
 \ker(g)=\{(a,a'):g(a)=g(a')\}
\]
for its equal-value equivalence relation.

\begin{theorem}[Exact reconstruction of a reachable boundary]\label{v13:mig:kernel}
Suppose two boundaries of representations of $F$ satisfy
\[
 b:A\to C,\qquad b':A\to C',\qquad F=fb=f'b'.
\]
Put $R=\im(b)$ and $R'=\im(b')$. A map $m:R\to R'$ satisfying
\begin{equation}\label{v13:mig:exact-state}
 m\circ b=b'
\end{equation}
exists if and only if
\begin{equation}\label{v13:mig:kernel-condition}
 \ker(b)\subseteq\ker(b').
\end{equation}
When it exists, $m$ is unique and surjective, and it is output-correct on the reachable states. Exact reconstruction in both directions is possible if and only if the kernels are equal; in that case the two reconstruction maps are inverse bijections.
\end{theorem}

\begin{proof}
If~\eqref{v13:mig:exact-state} holds and $b(a)=b(a')$, applying $m$ gives $b'(a)=b'(a')$. Thus the kernel inclusion is necessary. Conversely, that inclusion makes
\[
 m(b(a))=b'(a)
\]
well defined: different external inputs naming the same retained source state name the same target state. Every element of $R$ has the form $b(a)$, proving uniqueness. Every element of $R'$ has the form $b'(a)$, proving surjectivity. Finally, for $c=b(a)\in R$,
\[
 f(c)=F(a)=f'(b'(a))=f'(m(c)),
\]
which proves output correctness. Equal kernels give constructions in both directions, and their composites fix every element represented by an external input. Conversely, maps in both directions imply both kernel inclusions.
\end{proof}

\begin{figure}[htbp]
\centering
\begin{tikzcd}[column sep=large,row sep=large]
 & A \arrow[dl,"b"'] \arrow[dr,"b'"] & \\
 R \arrow[rr,dashed,"m"] \arrow[dr,"f|_R"'] && R' \arrow[dl,"f'|_{R'}"] \\
 & Y &
\end{tikzcd}
\caption{Two obligations for a change of retained boundary, in $\mathbf{FinSet}$. The upper triangle reconstructs exactly the target state belonging to each external input; it commutes precisely under the kernel inclusion. The lower triangle preserves the remaining output. Under the common-denotation assumptions, the upper obligation implies the lower one.}
\label{v13:mig:state-output-diagram}
\end{figure}

The criterion describes information after reclamation. If the source boundary has identified two external inputs but the requested target boundary still distinguishes them, no deterministic reconstruction from that boundary alone is possible. Keeping additional data can change the source map $b$ and remove the obstruction. Increasing the number of available workers cannot reconstruct a distinction that the retained state has lost.

\begin{table}[htbp]
\centering\small\renewcommand{\arraystretch}{1.14}
\begin{tabularx}{\linewidth}{@{}>{\raggedright\arraybackslash}p{0.22\linewidth}>{\raggedright\arraybackslash}p{0.32\linewidth}X@{}}
\toprule
Required preservation & Exact condition & What the condition establishes \\
\midrule
Remaining output & $f|_R=(f'|_{R'})m$ on the declared state domains. & Finishing after migration returns the same answer. The target state may differ from its state on the original input. \\
Exact reachable boundary & $mb=b'$, equivalently $\ker(b)\subseteq\ker(b')$. & The retained source state determines the designated target state for every original input. With $F=fb=f'b'$, this also preserves the reachable output. \\
Every source ready cube & For a valid graph map $q$, injectivity on nodes. & Each local cube of independent source completions maps isomorphically onto its image. The target may have additional work, and exact boundary reconstruction remains a further question. \\
\bottomrule
\end{tabularx}
\caption{Three contracts for representation changes. The second entails the reachable version of the first under their common-denotation assumptions. The third is an occurrence-level condition, proved in Corollary~\ref{v13:mig:all-cubes}, rather than a statement about recoverable values.}
\label{v13:mig:contracts}
\end{table}

\subsection{Reachable frontier migrations induced by sharing}\label{v13:mig:sharing}

The state-transfer rule of Proposition~\ref{prop:graph-state-transfer} uses the full partial valuation $\sigma:U\to D$. After reclaiming $U\setminus\partial_G U$, only $\beta_U$ remains. A target can require values absent from that smaller store. For surjective graph maps, the source residual contains enough information to recover the target frontier, although recovery can require evaluating some unfinished source work.

\begin{theorem}[Natural migration of reachable frontiers]\label{v13:mig:graph-natural}
Let $q:G\to G'$ be a graph Token map between valid graphs on the same interpreted primitive library and external interfaces, and suppose $q$ is surjective on nodes. Put $\phi(U)=q(U)$. There is a unique family
\[
 m_U:R_U\longrightarrow R'_{\phi(U)},\qquad
 m_Ub_U=b'_{\phi(U)}.
\]
Each component is surjective and preserves residual outputs. For $U\subseteq V$ it satisfies
\begin{equation}\label{v13:mig:naturality}
 m_V\circ(b_{U,V}|_{R_U})
 =\bigl(b'_{\phi(U),\phi(V)}|_{R'_{\phi(U)}}\bigr)\circ m_U.
\end{equation}
Thus these components are a natural transformation between the reachable frontier functors, with the target functor reindexed by $\phi$. The transformations for composable surjective graph maps compose componentwise.
\end{theorem}

\begin{proof}
Theorem~\ref{thm:sharing-adjunction} makes $q(U)$ completed. Fix $U$ and write $U'=q(U)$. For each target frontier node $t\in\partial_{G'}U'$, we first find a source node $s$ in
\[
 (N_G\setminus U)\cup\partial_G U
\]
with $q(s)=t$. If $t$ is named by a target output port, choose the source node at the same output port. It is either unfinished or a completed source output, hence belongs to the displayed carrier. Otherwise, $t$ is an operand of an unfinished target node $w'\notin U'$. Node surjectivity supplies $w$ with $q(w)=w'$, necessarily $w\notin U$. Preservation of the ordered operation record supplies an operand $s$ of $w$ with $q(s)=t$. If $s$ is completed, its unfinished consumer makes it a source frontier node; otherwise it is itself unfinished. This proves the required existence in both cases.

Now evaluate $G/U$ from a boundary tuple and read the selected $s$-values. For the reachable tuple $b_U(a)$, residual correctness identifies these values with $v_{G,a}(s)$. Theorem~\ref{thm:graph-evaluation} gives
\[
 v_{G,a}(s)=v_{G',a}(q(s))=v_{G',a}(t).
\]
Consequently $b'_{U'}(a)$ is determined by $b_U(a)$, proving the kernel inclusion. Theorem~\ref{v13:mig:kernel} now gives the unique exact, surjective, output-correct map $m_U$. Its restriction to reachable states is independent of the chosen source nodes.

For naturality, precompose both sides of~\eqref{v13:mig:naturality} with the surjection $b_U:A\to R_U$. By $b_{U,V}b_U=b_V$ and the corresponding target identity, both resulting maps are $b'_{\phi(V)}$. Surjectivity gives the asserted equality. If $q':G'\to G''$ is also surjective, composing the two exact state equations yields the exact state equation for $q'q$. Uniqueness identifies that composite with its migration component; the identity case is immediate.
\end{proof}

The proof supplies a possible implementation, not a bound on its cost. A selected source node $s$ can be outside $U$, so its value may need recomputation from the retained frontier. The required residual subgraph can be evaluated just far enough to obtain the selected values; in some cases this is substantial work. A runtime that still has a consistent full partial valuation can instead use the direct transfer of Proposition~\ref{prop:graph-state-transfer}. These are different storage situations.

The theorem is deliberately stated for the reachable subfunctors. Chosen source lifts can define functions on full Cartesian boundary spaces, but the theorem does not claim that such extensions are canonical, output-correct, or natural on arbitrary tuples. The following sharing example shows that a full-space output-correct extension may not exist at all.

\begin{example}[Sharing excludes inconsistent boundary assignments]\label{v13:mig:diagonal}
Over $D=\{0,1\}$, let $G$ have input $x$, operations $u=\operatorname{NOT}(x)$ and $v=\operatorname{NOT}(x)$, and outputs $(u,v)$. Let $G'$ share them as $w=\operatorname{NOT}(x)$ with outputs $(w,w)$. The node map fixes $x$ and identifies $u,v$ with $w$.

At full completion, the source residual is $\id:D^2\to D^2$, while the target residual is $\mathrm{diag}:D\to D^2$. No function $m:D^2\to D$ satisfies $\id=\mathrm{diag}\circ m$, since $(0,1)$ is not in the image of $\mathrm{diag}$. The reachable source frontier is only $\{(0,0),(1,1)\}$, and there the exact migration is $m(z,z)=z$. The disagreement of an arbitrary pair of boundary coordinates is an impossible state of this source graph, but it is a legitimate element of its Cartesian boundary interface.
\end{example}

\begin{example}[An injective graph change can require reclaimed information]\label{v13:mig:reclaimed}
Let $G$ have Boolean input $x$, a nullary constant node $c=0$, and the sole output $c$. At its completed set $U=\{x,c\}$, the live frontier is $(c)$ and $b_U(x)=0$. Form $G'$ by adding an unused operation $w=\operatorname{NOT}(x)$ and retaining the same output. Inclusion defines an injective graph Token map $q:G\to G'$. At $q(U)$, the target frontier is $(x,c)$ because $w$ is unfinished, so
\[
 b'_{q(U)}(x)=(x,0).
\]
The source frontier has identified the two inputs, whereas the target frontier distinguishes them. The kernel inclusion fails, and there is no exact reconstruction from the reclaimed source state. An output-correct migration can still choose either value for the missing coordinate, since the reported output remains $0$. Retaining $x$ in addition to the source frontier would remove this obstruction. The inclusion also preserves every source ready cube, as shown below; that structural property does not recover discarded values.
\end{example}

This example identifies the role of node surjectivity in Theorem~\ref{v13:mig:graph-natural}: an additional target consumer may require a source value that reclamation correctly discarded for the original representation. For a nonsurjective change, the runtime can check the kernel criterion directly, retain a richer boundary, or remove the extra target work using an independently justified output-demand restriction. The existence of the graph map alone is insufficient after reclamation.

\subsection{The exact change in a ready cube}\label{v13:mig:concurrency}

The migration of values and the organization of independent events supply different data to an elastic runtime. Sharing equal subcomputations preserves the answer and can reduce work, while identifying the occurrences that offered separate scheduling choices. The completion map $q_!$ makes that change measurable without identifying equal numerical values with equal occurrences.

\begin{theorem}[Image of a ready cube]\label{v13:mig:cube-image}
Let $q:G\to G'$ be any graph Token map between valid graphs with the fixed library and interfaces. For $U\in\mathcal L_G$ and a ready batch $R\subseteq\mathcal R_G(U)$, put
\[
 U'=q(U),\qquad S=q(R)\setminus q(U).
\]
Then $S\subseteq\mathcal R_{G'}(U')$, and direct image gives a surjective monotone map of Boolean intervals
\[
 [U,U\cup R]\longrightarrow[U',U'\cup S].
\]
Under the identifications of Proposition~\ref{prop:execution-lattice}, its formula is
\begin{equation}\label{v13:mig:cube-formula}
 \mathcal P(R)\longrightarrow\mathcal P(S),
 \qquad A_0\longmapsto q(A_0)\setminus q(U).
\end{equation}
The image therefore has dimension $|q(R)\setminus q(U)|$. This restriction is an isomorphism precisely when $q|_R$ is injective and $q(R)\cap q(U)=\varnothing$.
\end{theorem}

\begin{proof}
Each operand of a ready $r\in R$ belongs to $U$. Its operation record maps to the unique defining record of $q(r)$, so every target operand of $q(r)$ belongs to $q(U)$. Unless already completed, $q(r)$ is ready, proving the first assertion. The equation $q(U\cup A_0)=q(U)\cup q(A_0)$ gives~\eqref{v13:mig:cube-formula}. To obtain any subset of $S$, choose one source representative in $R$ of each desired target node. This proves surjectivity.

If some $r\in R$ maps into $q(U)$, the empty subset and $\{r\}$ have the same image. If two distinct members of $R$ have the same target, their singleton subsets have the same image. Either event prevents injectivity. Conversely, if neither occurs, $q$ is a bijection from $R$ onto $S$, and its powerset map is a Boolean-interval isomorphism.
\end{proof}

\begin{corollary}[Faithful source cubes characterize node injectivity]\label{v13:mig:all-cubes}
A valid graph Token map carries every source ready cube isomorphically onto its image if and only if it is injective on nodes.
\end{corollary}

\begin{proof}
Node injectivity gives both conditions of Theorem~\ref{v13:mig:cube-image} for every disjoint completed set and ready batch. Conversely, suppose $q$ is not injective. Distinct fixed input labels prevent two source inputs from being identified, and unique input or operation definitions prevent an input from being identified with an operation. There are therefore distinct operation nodes $u,v$ with $q(u)=q(v)=t$.

Let $V'$ contain every strict causal predecessor of $t$ and all target input nodes. It is completed and does not contain the operation node $t$. By Theorem~\ref{thm:sharing-adjunction}, $U=q^{-1}(V')$ is completed. Every operand of $u$ and $v$ maps to an operand of $t$ and hence belongs to $U$. Thus both $u,v$ are ready at $U$, whereas $t\notin q(U)$. Their two-dimensional cube has a one-dimensional image by Theorem~\ref{v13:mig:cube-image}, contradicting the claimed preservation of all source cubes.
\end{proof}

In Example~\ref{v13:mig:diagonal}, the two independent NOT operations give a two-dimensional cube at the initial cut, and sharing maps it to the one-dimensional interval of completing $w$. After only $u$ has completed, firing $v$ adds no new target completion at all, because its image $w$ is already in $q(U)$. These are the two mechanisms in~\eqref{v13:mig:cube-formula}: distinct ready occurrences can merge, and a ready occurrence can already be represented by completed shared work.

A noninjective graph map therefore changes some potential concurrency, while it can still improve execution time by removing redundant work. Conversely, an injective graph map faithfully represents the source's local choices inside the target without asserting that those are all the target's choices. Example~\ref{v13:mig:reclaimed} simultaneously has this fidelity and a failure of exact state reconstruction. The contracts of Table~\ref{v13:mig:contracts} make these distinctions available to the machine instead of hiding them behind denotational equality.

\paragraph{Using these certificates in an adaptive machine.}
For T5, a proposed expansion, fusion, or sharing step should specify its source and target representations, the retained-state invariant, the implemented migration, and the contract being preserved. A reachable certificate is usable only while the runtime maintains its stated reachability invariant. A family of migrations that commutes with advancement can be composed across successive representation choices; finite sequences preserve the answer when each component preserves it. Checking all finite state pairs can establish the kernel criterion, but may be expensive. The structural theorem for surjective graph maps supplies one alternative certificate; other representation languages need their own sufficient conditions.

Correct migration alone does not guarantee progress. An invocation can alternate equivalent representations indefinitely without completing its requested observation. An adaptation budget can exclude an infinite block of replacements without an execution step; overall completion still needs a progress argument. One sufficient form uses a well-founded rank decreased by productive execution and adaptation and not increased by administrative transitions, with no nonterminal deadlock and no indefinitely administrative suffix. In a continuing population of cores, progress should follow the finite request and its dependencies: repeatedly changing task identifiers must not defeat the scheduling obligation. Section~\ref{sec:unbounded-cores} states the finite-request conditions for the unbounded population.

The runtime must also account for the adaptation it selects. If migration costs $K_m(p)$ on $p$ workers and the replacement residual costs $T_{r'}(p)$, the relevant comparison includes $K_m(p)+T_{r'}(p)$, synchronization, and preparation. Peak storage includes any overlap between old and new descriptions, retained data, and temporary migration values. The intended machine can choose among compressed or expanded representations with different granularity and visible independence; these results specify what must remain coherent when that choice occurs during execution. They provide correctness tests and structural tradeoffs, while measured or proved resource bounds determine whether a particular choice is beneficial.

\subsection{A uniform family with selectable block granularity}\label{sec:parallel-prefix}

The prefix-XOR mapping $P_n$ of Section~\ref{sec:parity} illustrates the difference between storage compression, work reduction and parallel depth. Its occurrence graph has $n(n-1)/2$ XOR operations and span $n-1$. Its shared chain has $n-1$ operations and the same span. These facts show exactly what syntactic sharing improves. Associativity permits a further reorganization with much smaller span. Parallel prefix computation is classical \cite{blelloch1993}; its role here is to supply a verified family of alternative descriptions of the same mapping.

First consider any associative binary operation $\star:D^2\to D$ on a finite set. A recursive prefix graph on $n$ inputs can be constructed as follows. For $n=1$, return the input. For $n\ge2$, put $h=\lfloor n/2\rfloor$, compute the pairs $y_i=x_{2i-1}\star x_{2i}$ for $1\le i\le h$, and recursively compute their prefixes $z_i$. Output $x_1$ in position $1$, output $z_i$ in position $2i$, and output $z_i\star x_{2i+1}$ in each existing odd position $2i+1\ge3$. All pair nodes are independent, and so are the final odd-position nodes once their respective operands are ready. Distinct operation occurrences use distinct nodes; input and output routing costs are excluded from the operation count.

\begin{proposition}[Parallel prefix construction]\label{prop:parallel-prefix}
The recursive graph is valid and computes all inclusive prefixes under $\star$. If $w_n,s_n$ are its operation count and span, then
\begin{equation}\label{eq:prefix-recurrence}
\begin{split}
w_1=s_1&=0,\\
w_n&=w_{\lfloor n/2\rfloor}+n-1\le2n-2,\\
s_n&\le s_{\lfloor n/2\rfloor}+2
\le2\lceil\log_2 n\rceil\qquad(n\ge2).
\end{split}
\end{equation}
For XOR it consequently gives another element of $\mathcal I(P_n)$ with linear work and logarithmic span.
\end{proposition}

For a selectable granularity, split the input into $b$ contiguous nonempty blocks, $1\le b\le n$, each of length at most $\ell=\lceil n/b\rceil$. Compute all local prefixes by a chain within each block, retaining the block totals $z_1,\ldots,z_b$. Compute the inclusive prefixes of these totals by the recursive construction. For every element of block $j>1$, prepend the preceding-block total $z_1\star\cdots\star z_{j-1}$ to its local prefix, applying $\star$ with that total as the left operand. The first block's local prefixes are already final outputs. Call the resulting graph $R_{n,b}$, and let $\ell_1$ be the first block length.

\begin{proposition}[Block elasticity]\label{prop:block-prefix}
All $R_{n,b}$ compute the same prefix mapping. With the described construction and no subsequent removal of unused intermediate results, for $b\ge2$,
\begin{equation}\label{eq:block-prefix-cost}
\begin{split}
W(R_{n,b})&=(n-b)+w_b+(n-\ell_1)\le3n-3,\\
S(R_{n,b})&\le\lceil n/b\rceil+2\lceil\log_2 b\rceil.
\end{split}
\end{equation}
For $b=1$, $R_{n,1}$ is the shared chain with work and span $n-1$. For XOR, every $R_{n,b}$ preserves the causal heap and the output behavior of $P_n$.
\end{proposition}

The proofs are routine inductions and are given in Appendix~\ref{app:scheduling}.

\begin{table}[htbp]
\centering\small\renewcommand{\arraystretch}{1.13}
\begin{tabularx}{\linewidth}{@{}>{\raggedright\arraybackslash}p{0.27\linewidth}>{\raggedright\arraybackslash}p{0.23\linewidth}>{\raggedright\arraybackslash}p{0.21\linewidth}>{\raggedright\arraybackslash}X@{}}
\toprule
Description of $P_n$ & XOR work & Span & Correctness certificate \\
\midrule
Separate left-associated prefixes & $n(n-1)/2$ & $n-1$ & Direct term evaluation. \\
Shared prefix chain & $n-1$ & $n-1$ & Syntactic sharing Token map. \\
Recursive parallel prefixes & $w_n\le2n-2$ & $\le2\lceil\log_2n\rceil$ & Associativity and recursive construction. \\
Block family $R_{n,b}$, $b\ge2$ & $(n-b)+w_b+n-\ell_1$ & $\le\lceil n/b\rceil+2\lceil\log_2b\rceil$ & Local prefixes, totals and associative recombination. \\
\bottomrule
\end{tabularx}
\caption{One behavior with several dependency structures. Work counts XOR operations; input delivery, output materialization, scheduling and communication are excluded. The bounds concern these constructions, not optimality over every implementation.}\label{tab:parallel-prefix-cost}
\end{table}

Under the fixed-resource model, choose $b=\min\{p,n\}$. Equations~\eqref{eq:work-span} and~\eqref{eq:block-prefix-cost} give $O(n/p+\log p)$ slots for $2\le p\le n$, the single chain gives $O(n)$ for $p=1$, and $b=n$ gives $O(\log n)$ for $p>n$. This is a resource-parameterized family of exact descriptions of one mapping---elasticity of representation and granularity in the sense of Table~\ref{tab:elasticity-levels}---rather than a speed-up of the unchanged chain. The prefix mapping remains causal even when its full known input is processed in parallel: causality constrains dependence of an output on future input coordinates, not the wall-clock order of independent operations; tasks whose inputs have not arrived are simply not ready. The associativity certificate applies to XOR and other specified associative operations; reassociating floating-point additions can change their bit-level result, so it does not transfer to numerical neural computation without further assumptions. The block parameter is selected before an invocation or before a residual invocation; changing $b$ while operations are in flight is handled by Theorem~\ref{thm:residual-replacement}, not by Proposition~\ref{prop:block-prefix}.

\subsection{A first execution protocol and its invariant}\label{sec:elastic-protocol}

The following protocol makes the scope of an initial machine prototype concrete. A graph invocation retains its current description, completed-node set, live values, ready tasks and invocation identifier; a resource manager supplies workers. The complete mathematical valuation can be used as a proof object even when the implementation stores only live values. The protocol has four phases or transitions:
\begin{enumerate}
\item \emph{Dispatch and completion.} Dispatch ready tasks to available workers. Accept each task result once, mark its node completed, update readiness, and retain values needed by unfinished consumers or outputs.
\item \emph{Allocation change.} Adjust the number or placement of workers without changing the logical graph or losing live values. In-flight tasks finish or use a separately verified migration mechanism.
\item \emph{Certified adaptation.} At a consistent completed boundary, construct the residual graph. Choose a candidate implementation with the same residual interface, validate its certificate, retain the ordered boundary values, and replace the remaining invocation. If preparation fails or exceeds its budget, resume the original residual graph.
\item \emph{Completion.} Return the ordered output tuple and reclaim invocation-specific data once no caller requires it.
\end{enumerate}

Each adaptation preserves the remaining result by Theorem~\ref{thm:residual-replacement}. Dispatch preserves the partial valuation invariant. Allocation changes preserve values and dependencies. Consequently any terminating execution with finitely many such certified adaptations returns the original specified output. Termination follows, for example, if only finitely many adaptations are allowed, their preparation terminates, and after the last one a fixed finite residual graph receives unbounded cumulative bounded worker capacity as in Proposition~\ref{prop:variable-capacity}. A policy making infinitely many individually correct expansions may never return an answer; correctness of each rewrite alone is not a termination argument.

This protocol can be implemented centrally before investigating decentralized variants. Nothing in schedule independence requires a centralized program counter, but efficient readiness tracking, memory reclamation, backpressure and load balancing remain real algorithmic tasks. Unrestricted eagerness can exceed memory even when much work is independent. A memory-constrained scheduler may delay ready nodes or choose a different representation; its progress and performance then need analysis beyond the greedy bounds above.

\subsection{An unbounded population of Token computing cores}\label{sec:unbounded-cores}

The machine envisaged here has potentially unboundedly many Token computing cores. This is an architectural conviction: computational organization should be able to grow with the task, while interfaces, dependencies and permitted changes remain explicit. Thesis~\ref{thesis:cores} brings that conviction into the foreground. The finite graph results above supply its local execution theory; this subsection specifies a first operational foundation for the population itself. The word \emph{unbounded} means that the machine specification imposes no common finite bound across all admitted workloads and runs. It does not assert that one physical machine supplies infinitely many simultaneous operations.

\subsubsection{What a computing core contains}

\begin{definition}[Token computing core]\label{def:computing-core}
Fix an effective typed representation system as in Definition~\ref{def:representation-system}. A \emph{Token computing core} is a logical invocation record
\[
\kappa=(\iota,T_A,T_B,r,\sigma),\qquad r\in R_{A,B},
\]
where $\iota$ is a fresh invocation identifier, $T_A,T_B$ are effectively presented finite interfaces, $r$ is an executable description, and $\sigma$ contains its finite input bindings, execution status and retained local values. Input bindings may refer to outputs of other invocations. A core denotes a particular invocation of $\delta(r):T_A\to T_B$; it is ready when its required inputs are available. In the pure fragment, different ready invocations have no additional shared mutable effects.
\end{definition}

Here \emph{computing core} names an execution unit. The \emph{signature core} $E$ still names the fixed symbols used in Token records. Neither is a tokenizer vocabulary item or a physical processor. A structural Token is a record; the operational interpretation of selected records supplies the computing behavior. Two invocations of the same description on the same values may have distinct identifiers, retaining the distinction emphasized by Thesis~\ref{thesis:occurrence}. Identifying their results by sharing is a separate certified transformation.

Core descriptions, interfaces, bindings and statuses can be carried by a finite Token class. For example, with disjoint carrier tags for identifiers, ports, codes and values, records may have forms
\[
\begin{gathered}
(\mathrm{code},\iota,\ulcorner r\urcorner),\quad
(\mathrm{interface},\iota,\ulcorner T_A\urcorner,\ulcorner T_B\urcorner),\\
(\mathrm{bind},\iota,j,\iota',k),\quad
(\mathrm{value},\iota,j,a),\quad (\mathrm{status},\iota,s).
\end{gathered}
\]
Only the record heads need be fixed signature symbols; varying indices, status values and finite codes are carrier data. Multiport interfaces use declared component types, with one port as the default. Validity checks enforce uniqueness of the appropriate fields, correct code and value types, and compatible bindings. A general map of such record classes preserves their records, but is not automatically a legal runtime transition. In particular, changing a status or deleting a completed invocation is specified by an operational rule, rather than asserted to be a heap-preserving inclusion.

The execution boundary is relative. A description may be evaluated as one core, expanded into a finite network of smaller cores, or replaced by a certified fused description. Theorems~\ref{thm:module-expansion} and~\ref{thm:residual-replacement} justify specified instances of these changes. Tables, circuits, symbolic programs and numerical models are possible descriptions. A Transformer is one possible evaluator within this population, and its own computation can be exposed at a finer granularity when an appropriate description and correctness argument are available.

\subsubsection{Finite configurations and effective growth}

\begin{definition}[A locally finite core regime]\label{def:core-regime}
A configuration has a finite set $C$ of active invocation identifiers, finite encoded core descriptions and dependency records, a finite store of retained values and completed results, and a finite encoded control state. A fresh-name counter ranges over $\mathbb N$ without a fixed maximum. Active cores are waiting or executing; completed results can remain in the store after their cores leave $C$. The initial configuration is finite. Each transition is effective, terminates after finite work, and constructs only finitely many records and fresh identifiers. The following transitions are admitted with explicit validity checks:
\begin{enumerate}[label=(\roman*),leftmargin=2.2em]
\item \emph{Generate or admit work.} Create finitely many fresh cores and typed bindings under a specified workload generator or request protocol. New requests may change the overall workload; their admission is not a claim of equivalence to a previous fixed mapping.
\item \emph{Dispatch and complete.} Execute a ready core with its declared evaluator and publish its result once. A small-step realization may retain finite intermediate evaluator state; completion agrees with $\delta(r)$.
\item \emph{Refine or fuse.} Change an existing computation using a certified finite description transformation. During execution, retain a consistent boundary and satisfy the state conditions of Section~\ref{sec:online-replacement}.
\item \emph{Reclaim or reallocate.} Remove completed core records when their live results and references remain available, and change worker assignments without altering required values or dependencies.
\end{enumerate}
The pure acyclic fragment requires the dependency graph of each current invocation to be valid. A resource policy may delay admission, expansion or dispatch when storage or workers are unavailable. Infinite storage availability and progress are not assumed by a record format.
\end{definition}

\begin{proposition}[Finite stages of an unbounded regime]\label{prop:finite-core-stages}
Every configuration reached after a finite number of transitions in a locally finite core regime is finite. If the initial active population has size $c_0$ and transition $t$ introduces $b_t$ fresh cores, then
\[
|C_t|\le c_0+\sum_{s<t}b_s<\infty.
\]
An infinite sequence of transitions introduces at most countably many invocation identities. These conclusions do not impose a uniform finite bound on $|C_t|$ across all admitted runs.
\end{proposition}
\begin{proof}
Induct on transition count. A finite configuration and a terminating transition that constructs finitely many records give another finite configuration. Removing records cannot increase the bound; introducing $b_t$ fresh identities increases it by at most $b_t$. Over an infinite sequence, the identities lie in a countable union of finite sets, or directly in the chosen name set $\mathbb N$. The bound depends on the run and the number of transitions, so it is not a common architecture constant.
\end{proof}

For a run $\rho$, write $C_t^\rho$ for its active cores and $Q_t^\rho\subseteq C_t^\rho$ for its waiting ready cores. Two distinct requirements can now be stated:
\begin{align}
\text{unbounded active population:}\quad
&\forall N\ \exists\rho,t\quad |C_t^\rho|>N,\label{eq:unbounded-population}\\
\text{unbounded available concurrency:}\quad
&\forall N\ \exists\rho,t\quad |Q_t^\rho|>N.\label{eq:unbounded-ready}
\end{align}
The second condition implies the first. Neither follows from indefinitely many successive invocations: a service that creates a new core only after the previous one is reclaimed can run forever with $|C_t|\le1$. Even a large active population may form a serial dependency chain. Fresh identities, active width and ready width are consequently different observables.

If $p_t$ workers are available, a dispatch step starts at most $\min\{p_t,|Q_t|\}$ new core evaluations. The earlier work--span estimates apply to a fixed exposed graph with their declared primitive costs; treating a complex core as one task does not make its evaluator unit-cost. Proposition~\ref{prop:variable-capacity} assumes a finite bound on worker capacity during the invocation. Equations~\eqref{eq:unbounded-population}--\eqref{eq:unbounded-ready} instead range over workloads and configurations. They do not extend that scheduling theorem to arbitrary graph growth or infinite instantaneous capacity.

\subsubsection{A uniform witness and the role of representation}

\begin{proposition}[Arbitrarily wide finite core populations]\label{prop:unbounded-core-width}
Let $r\in R_{A,B}$ denote $f:T_A\to T_B$, with nonempty finite carrier $A$. Suppose the regime admits a uniform generator that, for each finite input $(n,a_1,\ldots,a_n)$, creates $n$ distinct ready invocations of $r$ before dispatch, retaining their output order. Then it satisfies~\eqref{eq:unbounded-ready}. Whenever all generated cores complete, their ordered outputs realize
\begin{equation}\label{eq:core-product-family}
F_n=f^{\otimes n}:T_A^{\otimes n}\longrightarrow T_B^{\otimes n},
\qquad (a_1,\ldots,a_n)\longmapsto(f(a_1),\ldots,f(a_n)).
\end{equation}
The constructor is effective and uses only finite configurations at each finite stage. Here $\otimes$ denotes the Cartesian product of Token classes, not a tensor imposing a physical execution schedule.
\end{proposition}
\begin{proof}
One fixed loop reads the finite input and allocates one invocation with a fresh identifier and bound input at each iteration. No invocation depends on another, and each description and binding is finite. After $n$ allocations all $n$ are ready. Taking $n=N+1$ proves~\eqref{eq:unbounded-ready}. By evaluator correctness the $i$th result is $f(a_i)$, independently of the order of completion. Keeping the output indices yields~\eqref{eq:core-product-family}, which is a Token map by the product construction. For $n=0$ the empty population realizes the identity of the terminal object. If available resources allow every allocated invocation to finish, this gives the asserted completed output.
\end{proof}

The proposition assumes no fixed admission cap across the family; an implementation with a fixed cap may admit the requests in smaller batches. Constructing the cores, supplying their inputs, storing results and collecting outputs all have costs. The example proves a realizable structural width, not an arbitrary speedup. It also makes a distinction essential to the machine proposal: changing $n$ changes the external interface and the mapping $F_n$. For each fixed $n$, however, one can evaluate $F_n$ as a single composite module or expose its $n$ independent component calls. Correct product constructors or module certificates connect these descriptions. This is representation elasticity at a fixed interface, while varying $n$ is growth of a uniform interface family.

For a concrete finite mapping, take Boolean negation and obtain $F_n(x_1,\ldots,x_n)=(1-x_1,\ldots,1-x_n)$. A generator can create one core per coordinate, a graph can expose all the calls, and a packaged evaluator can return the same vector through one invocation. Equal input values do not erase coordinate identities. This small witness isolates the organization principle without attributing novelty to coordinatewise evaluation. More substantial examples must show how representation choice, communication, storage and adaptation affect total cost.

\begin{proposition}[Finite requests in a continuing population]\label{prop:finite-request-progress}
Consider an execution that may continue admitting other work indefinitely. At a finite stage, suppose a requested finite output tuple has a finite, fully generated acyclic dependency closure $G$, with all external inputs available. Assume that its core descriptions are pure and total, its defining records remain fixed, and required values are retained. Already executing cores of $G$ eventually complete; every waiting core of $G$ that remains ready is eventually dispatched and completes. Then the requested output tuple is available at a finite later stage and equals the evaluation of $G$ on those inputs, independently of the completion order of its cores.
\end{proposition}
\begin{proof}
Order the finitely many core nodes topologically. Inputs and retained completed values are available by hypothesis, and already executing cores eventually finish. Once the predecessors of a waiting node have completed, their retained values make it persistently ready until dispatch. The stated progress condition gives its completion at a finite later stage, and evaluator correctness determines its value. Induction over the finite ordering gives a finite completion stage for every required node; the maximum of these finitely many stages supplies the requested output. The value induction is the one in Theorem~\ref{thm:graph-evaluation}, so the resulting tuple is independent of completion order. Work outside the dependency closure cannot change its values under the purity and retention hypotheses.
\end{proof}

This result concerns completion of an observation inside a continuing computation, not termination of the whole population. The progress condition includes protection against starvation by newly admitted work. A stream of correct allocations can otherwise postpone an old ready request forever. The proposition supplies no latency bound and does not justify arbitrary rewrites of the closure during execution; certified replacement still needs the state and progress conditions discussed above. A useful finite-resource implementation must enforce the stated progress condition through admission and scheduling policies.

\subsubsection{A categorical view of the unbounded horizon}

Finite configurations also admit a precise relationship to an infinite structural description. We state it for compatible inclusions; an elastic runtime can additionally delete and replace records, so its successive live states need not have this form.

\begin{proposition}[Union of a growing sequence of Token classes]\label{prop:core-union}
Fix a signature core $E$. Let $T_n=\langle B_n,E,H_n\rangle$ be a sequence of finite Token classes, presented so that $B_n\subseteq B_{n+1}$ and $H_n\subseteq H_{n+1}$. Then in $\TS_E$ the inclusion diagram has colimit
\[
T_\infty=\left\langle\bigcup_n B_n,E,\bigcup_n H_n\right\rangle.
\]
For any $\mathcal V\in\TS_E$, compatible Token maps $v_n:T_n\to\mathcal V$ extend to a unique Token map $v_\infty:T_\infty\to\mathcal V$. Every finite subset of the union carrier and heap is already present at some finite stage.
\end{proposition}
\begin{proof}
Compatibility defines a function $v_\infty(b)=v_n(b)$ for any stage containing $b$, independent of the chosen stage. Each record of $H_\infty$ belongs to an $H_n$, so its image belongs to the heap of $\mathcal V$. This proves the Token-map condition. Uniqueness holds on every carrier element because it appears at a finite stage. The same argument for every compatible cocone is the colimit universal property. A finite collection of elements and records appears in finitely many stages; their largest index contains them all.
\end{proof}

\begin{figure}[htbp]
\centering
\begin{tikzcd}[column sep=large,row sep=large]
T_n \arrow[r,hook,"j_n"] \arrow[dr,"v_n"'] &
T_{n+1} \arrow[r,hook,"j_{n+1,\infty}"] \arrow[d,"v_{n+1}"] &
T_\infty \arrow[dl,dashed,"v_\infty"] \\
& \mathcal V &
\end{tikzcd}
\caption{Finite structural stages and their unbounded union. The displayed triangle is part of the full inclusion chain. Compatible finite-stage values determine the unique dashed Token map; no step evaluates the infinite union all at once. This diagram describes monotone structural growth or an append-only history, not arbitrary replacement of live configurations.}\label{fig:core-union}
\end{figure}

This is an elementary colimit calculation in the record category, not a new general colimit theorem. It gives the phrase ``theoretically unboundedly many cores'' a structural interpretation: every finite collection of recorded cores, bindings and events can lie in a larger finite description, while their compatible union gives a mathematical horizon. That union belongs to $\TS_E$ and need not belong to its finite bounded variant. A history can use immutable event identifiers to record allocation, completion and reclamation; changing a live status need not overwrite a past event. Garbage-collected live storage is a separate object and need not grow by inclusions.

The compatible maps in Proposition~\ref{prop:core-union} are an assumption to establish for the chosen execution records. When they are valuations, the proposition assembles their already compatible values; it does not create values for unfinished work. Nor does finite arity alone guarantee that an infinite dependency structure is executable: it can contain an infinite predecessor chain. The operational result proved here covers a requested observation when it has a finite generated dependency closure, terminating local evaluators and the progress conditions of Proposition~\ref{prop:finite-request-progress}. These requirements connect the unbounded horizon to the finite execution theorems without asserting infinite work in finite time.

T5 therefore adds a population and growth principle to the earlier execution theory. T1 exposes the cores' internal records; T2 makes their descriptions available to computation; T3 fixes the interfaces and preservation obligations of replacement; T4 retains invocation identity; T5 asks the machine to generate and reorganize these units without a preset global scale. Its further development concerns effective generation, finite-resource progress and useful representation policies. The finite-stage, width, finite-request progress and union propositions establish a foundation for that programme, while the choice of a generally advantageous dynamic architecture remains open.

\subsection{Research position and an evaluation programme}\label{sec:elastic-research}

Concurrent logical units and their creation already have a foundational precedent in the Actor model \cite{hewitt1973}. Category theory has also long supplied semantics and comparisons for parallel computation \cite{winskel1986}. T5 therefore locates our architectural ambition in an established research area; an unbounded population alone is not a novelty claim. Our proposed direction combines typed finite-mapping representations, inspectable structure, invocation identity and certified changes of execution granularity within Token Space. Dependency-directed execution, work--span scheduling and parallel prefix networks are established ideas \cite{blelloch1999,blelloch1993}. Dynamically scheduled elastic circuits also use local readiness and communication protocols to accommodate execution timing \cite{josipovic2017,josipovic2018}. The word ``elastic'' in that setting does not by itself mean changing the number of physical processors. Our four distinctions in Table~\ref{tab:elasticity-levels} are intended to avoid that ambiguity. Certified compilation and equality-based optimization supply further precedents, as discussed in Section~\ref{sec:finite-representations}.

The present contribution is the explicit connection among the Token-level constructions, unique graph valuation, the lattice of partial execution states, and implementation changes. Theorem~\ref{thm:execution-gluing} realizes execution by universal gluing in Token Space and by a fully faithful diagram over the fixed graph. Theorem~\ref{thm:sharing-adjunction} restricts image and inverse image to legal completions and gives a sound saturation for state reuse. Theorem~\ref{thm:frontier-functor} proves coherent advancement between changing live interfaces, with residual outputs as a natural transformation. These are concrete results of the present formalization. Their ingredients---configuration orders, universal properties and adjunctions---are established mathematics \cite{maclane1998,winskel1987}; the results are not presented as a new general theory of concurrency or as a priority claim over all related formulations. Ranked module expansion and term compilation supply proved bridges between descriptions, callable components and primitive graphs; they do not claim priority for inlining or term-graph evaluation. The scheduling estimates are standard consequences of the graph model, included with proofs rather than claimed as new parallel algorithms. The residual construction supplies a specific sufficient rule for adapting an unfinished invocation. The prefix family provides a complete exact example in which the same mapping admits selectable granularity and different parallel depths. Together they make the proposed machine's concurrency and elasticity mathematically testable. A generally efficient choice of representation, distributed adaptation mechanism and practical advantage over existing runtimes remain research questions.

The population results of Section~\ref{sec:unbounded-cores} add an effective growth model and a structural horizon to this account. A useful prototype should retain the representation and its certificates as data, expose the scheduler's state, and record the following measurements under a declared primitive library and hardware model:
\begin{itemize}
\item Compare the occurrence graph, shared chain, recursive prefixes and several $R_{n,b}$ descriptions on identical inputs. Check exact outputs and record work, critical path, total latency and peak live storage.
\item Vary worker supply during a fixed-graph invocation, including temporary zero-capacity intervals. Compare the observed idealized schedule with Proposition~\ref{prop:variable-capacity}, then report actual scheduling and communication overhead separately.
\item Replace residual graphs at chosen completed boundaries. Record retained boundary size, discarded storage, certificate cost and migration cost; compare with uninterrupted execution and restarting from the original input.
\item Include a low-parallelism chain, a wide independent graph, graphs with heavy sharing, and a small irregular mapping represented by both a table and a graph. These cases expose different tradeoffs instead of assuming one representation wins everywhere.
\item Implement the uniform core generator and interleave generation, completion and reclamation. Record cumulative identities, peak active population, ready width, physical worker count and retained storage separately. Include a serial stream with many successive identities but width one, and the independent family of Proposition~\ref{prop:unbounded-core-width}. Add admission caps to measure backpressure and include the cost of generating descriptions.
\item In subsequent neural experiments, specify numerical precision, deterministic kernels, and whether the required relation is exact bit-level equality or a quantified approximation. Measure representation preparation and repeated-use amortization as well as forward latency.
\end{itemize}

These are proposed experiments, not reported hardware results. The mathematical examples establish behavior and ideal operation bounds. Empirical performance must also include graph construction, certificate checking, allocation, data movement and the control protocol. The motivating objective is a machine whose computational behavior is stable while its representation, active structure and execution scale can change under explicit correctness conditions.

\section{Causal Token Structures and Their Computational Role}
\label{v13:causal:section}

The execution theory has two structural questions behind it: when can all completed fragments reconstruct their induced union, and when can ordered record structures themselves support function objects? The mathematical companion \cite{pancausal} develops these questions independently, with complete definitions, proofs and counterexamples. This section states the results needed for the machine programme and explains their scope. The graph evaluation, execution gluing, frontier and migration proofs in Sections~\ref{sec:finite-representations}--\ref{sec:elastic-execution} remain self-contained.

\subsection{A support criterion for universal gluing}\label{v13:causal:relative}

An \emph{ordered Token structure} is $T=(B,\leq,I,H)$: a poset, a lower set $I$ of initially available elements, and a heap on the carrier $B$ with fixed core $E$. Its maps are monotone Token maps preserving initial availability. For $U\supseteq I$ lower, write $T_U=\langle U,E,H_U\rangle$ for the induced class. The support $\supp(\tau)$ of a record is the finite set of its carrier entries. Core entries have no support.

A record is \emph{relatively rooted} if its support is contained in $I$, or the remaining support $\supp(\tau)\setminus I$ has a greatest member. The root must occur in the record; an upper bound elsewhere does not suffice. The qualification relative to $I$ admits arbitrary relations between initial inputs, since every completed fragment already contains them.

\begin{theorem}[Universal gluing criterion; companion paper]\label{v13:causal:gluing}
The following conditions on an ordered Token structure are equivalent:
\begin{enumerate}[label=(\roman*),leftmargin=2em]
\item every record is relatively rooted;
\item $H_{U\cup V}=H_U\cup H_V$ for all lower $U,V\supseteq I$;
\item $H_{\bigcup_\alpha U_\alpha}=\bigcup_\alpha H_{U_\alpha}$ for every nonempty family of such sets;
\item every induced-fragment square in Figure~\ref{v13:causal:gluing-figure} is both a pullback and a pushout in $\TS_E$.
\end{enumerate}
\end{theorem}

For orientation, if a remaining root belongs to one piece, lower closure places the whole record there. Conversely, covering a finite support by the sets $I\cup\downarrow b$ forces one support member to dominate all the noninitial members whenever all such covers reconstruct the heap. The companion proves both implications and the universal properties, including empty supports. The empty family is excluded: pure-core and initial records already belong to $T_I$.

\begin{figure}[htbp]
\centering
\begin{tikzcd}[column sep=large,row sep=large]
T_{U\cap V}\arrow[r,hook]\arrow[d,hook]&T_U\arrow[d,hook]\\
T_V\arrow[r,hook]&T_{U\cup V}
\end{tikzcd}
\caption{Completed fragments reconstruct the induced union exactly under relative rootedness. The square glues their records as well as their carriers. This is a property of the specified family of restrictions, not an assertion of ambient adhesivity.}
\label{v13:causal:gluing-figure}
\end{figure}

This criterion explains why the operation records of a valid graph glue: each record includes its result, and all its operands precede that result. If an additional heap record relates two independent noninitial occurrences without containing a greatest one, it can cross two completed fragments and defeat the property. Such a record may be appropriate for another interface, but the completed-fragment argument no longer applies to that enlarged heap. For relatively rooted structures the records decompose as the disjoint union of initial records and records indexed by their remaining root. This suggests a local organisation of records by completion event; it does not itself specify a scheduler or a communication cost.

\subsection{Function objects with designated result roles}\label{v13:causal:designated}

Root existence alone does not ensure cartesian closure. The category of rooted ordered structures with empty initial set already fails to be cartesian closed: a product can delete a paired record whose two component roots occur in opposite positions, preventing the required exponential evaluation. The companion gives an explicit two-chain counterexample and proves the product construction.

A sufficient alternative is to fix, for every positive-arity record shape $R$, one carrier position $j_R$, and require
\[
R(a_1,\ldots,a_n)\ \Longrightarrow\ a_i\leq a_{j_R}\quad\text{for all }i.
\]
This compares initial entries too. Denote by $\mathsf{DRTS}_{E,j}$ the resulting category, with the same monotone, initial-preserving Token maps.

\begin{theorem}[Designated-root function objects; companion paper]\label{v13:causal:designated-ccc}
The category $\mathsf{DRTS}_{E,j}$ has finite limits, finite coproducts and exponentials. The carrier of $B^A$ is $\operatorname{Mon}(|A|,|B|)$ with pointwise order and initial set $\{f:f(I_A)\subseteq I_B\}$. Its positive-arity records are
\begin{equation}\label{v13:causal:exponential-records}
\begin{split}
R^{B^A}(f_1,\ldots,f_n)\ \Longleftrightarrow\ &f_i\leq f_{j_R}\ \text{for every }i,\quad\text{and}\\
&R^A(a_1,\ldots,a_n)\Rightarrow R^B(f_1(a_1),\ldots,f_n(a_n))\ \text{for every tuple}.
\end{split}
\end{equation}
For nullary shapes use $R^{B^A}=(R^A\Rightarrow R^B)$.
\end{theorem}

\begin{figure}[htbp]
\centering
\begin{tikzcd}[column sep=large,row sep=large]
X\times A\arrow[rr,"h"]\arrow[dr,"\widehat h\times\id_A"']&&B\\
&B^A\times A\arrow[ur,"\operatorname{ev}"']&
\end{tikzcd}
\caption{Currying and evaluation preserve the dependency order, initial availability and the declared records. The function object supplies structural semantics; choosing a compact executable representation requires further work.}
\label{v13:causal:exponential-figure}
\end{figure}

The carrier consists of all monotone functions; the morphisms $A\to B$ correspond to the \emph{global points} of this structured object. These sets should not be identified. The designated-position condition is sufficient, not necessary: the full subcategory assigning to each poset its symmetric comparability relation, empty initial set and no other records is isomorphic to $\mathbf{Pos}$ and is cartesian closed, although no single position roots both orientations of a two-element comparable pair. No local cartesian closure or quasitopos assertion for $\mathsf{DRTS}_{E,j}$ is needed here.

\subsection{Occurrence structures, execution and the finite horizon}

\begin{proposition}[Graph and tree presentations; companion paper]\label{v13:causal:graph-bridge}
For a fixed finite primitive library and input/output interface, valid computation graphs with graph Token maps lift fully faithfully to $\mathsf{DRTS}_{E,j}$ by taking dependency order, input nodes as the initial set, and each operation's result slot as designated root. Supported tree reification likewise lifts fully faithfully to the designated-root category over its extended node/root core, using the subtree order, the tuple identifier as designated root of each node record, and empty initial set.
\end{proposition}

The companion independently defines and proves both lifts, including recovery of substitutions from reified maps. For trees the empty initial set permits a base letter to be substituted by a compound term; marking every base letter initial would impose an additional restriction. For graphs the carrier records \emph{occurrences}. Boolean NOT, for example, can validly relate an input occurrence to a later result occurrence, although imposing its two value records on the value carrier would force both $0\leq1$ and $1\leq0$.

A full embedding preserves the chosen structures and arrows, but does not make executable graphs closed under ambient exponentials. With one input, one output and a NOT primitive, the exponential of the one-node identity graph with itself acquires a NOT self-loop; it violates strict acyclicity and unique definition. More generally, graph maps preserve denotation, so a fixed-interface category containing both identity and NOT graphs need not even have a terminal object. Function objects in the structure category consequently organise higher-order preservation obligations; the represented-mapping machine still needs its own valid descriptions and evaluators. Certified implementation rewrites also need not be graph Token maps.

For the unbounded core horizon of Thesis~\ref{thesis:cores}, one further result isolates the structural finiteness assumptions. If $I$ is finite and every record is relatively rooted, induced completed stages with both finite carrier and finite heap form a nonempty cover exactly when each principal down-set $\downarrow b$ is finite, the initial heap $H_I$ is finite, and the records with any fixed remaining root are finite in number. Their filtered union recovers $T$ as a colimit in $\TS_E$. If, in addition, the nodes in $I$ are input nodes and each noninitial node has a unique operation definition with strictly earlier finite operands and total primitives on a nonempty value set, every fixed input assignment determines a unique compatible valuation. Every finite observation is determined in its finite predecessor closure together with $I$. The valuation diagram concerns the input--operation reduct; any further relation requires its own interpretation.

These statements clarify T5 without replacing the execution theory. Local finiteness permits finite observations of an unbounded structural horizon; a serial chain still has ready width one. Effective generation, sufficient resources and progress require the operational conditions of Section~\ref{sec:unbounded-cores}. The main paper retains its migration criteria because changing a running representation also requires retained information and concurrency to be analysed. The companion's gluing and function-object theorems supply structural criteria; they do not alone guarantee economical or terminating adaptation.

\section{From Tokenizer Outputs to Transformer Computations}\label{attention-in-token-space}

\subsection{The computational pipeline and its carriers}\label{sec:pipeline}

Fix trained parameters $\theta$ and a deterministic inference protocol. Dropout is disabled, normalization includes a positive numerical stabilizer, and positional conventions and attention masks are specified. We first treat a fixed length $n$ within the supported context window, without external memory or retrieval. Define
\[
\mathcal H_{\ell,n}=(\mathbb R^{d_\ell})^n,
\qquad
\Delta(V)=\{p\in[0,1]^V:\textstyle\sum_{v\in V}p(v)=1\}.
\]
The following diagram separates tokenizer IDs, embeddings, block states, and output distributions:
\begin{figure}[htbp]
\centering
\begin{tikzcd}[column sep=huge,row sep=large]
\mathcal X_n \arrow[r,"\tau"] & V^n \arrow[r,"e_\theta^n"] & \mathcal H_{0,n} \arrow[d,"F_{0,\theta}"] \\
\Delta(V)^n & \mathcal H_{L,n} \arrow[l,"r_\theta^n"'] & \mathcal H_{1,n} \arrow[l,"\mathcal F_\theta"']
\end{tikzcd}
\caption{A deterministic forward pass, with $\mathcal F_\theta=F_{L-1,\theta}\circ\cdots\circ F_{1,\theta}$. $\mathcal X_n$ is the set of texts tokenized to length $n$. Position information is supplied at the input or within the blocks according to the architecture. Each arrow after tokenization is a function on a whole sequence. Selected heaps turn particular arrows into Token maps; being a function alone does not specify which structure is preserved.}\label{fig:pipeline}
\end{figure}

For an embedding table $W_{\mathrm{emb}}\in\mathbb R^{|V|\times d}$,
\[
e_\theta(v)=W_{\mathrm{emb}}^\top\operatorname{onehot}(v),
\qquad H_i^{(0)}=e_\theta(v_i).
\]
With additive position embeddings one instead starts from $H_i^{(0)}=e_\theta(v_i)+p_i$. Architectures using rotary position information can leave this lookup unchanged and transform queries and keys inside attention. For a column state $h$, the output readout is
\begin{equation}\label{eq:readout}
r_\theta(h)=\Softmax\bigl(W_{\mathrm{out}}N_f(h)+b\bigr)\in\Delta(V),
\end{equation}
where $W_{\mathrm{out}}\in\mathbb R^{|V|\times d}$ and $N_f$ is the final normalization. The distribution at position $i$ predicts the token following the prefix through $i$, rather than recovering an intrinsic word assigned to $h_i^{(L)}$. The last row supplies the next-token distribution during generation. Embedding and output readout are not assumed to be inverse maps \cite{vaswani2017}.

\begin{example}[Why the carrier cannot generally be one token vector]\label{ex:context}
Consider scalar, two-position causal attention with zero queries and keys and identity values. At the second position it returns $(a+b)/2$. Inputs $(a,0)$ and $(a',0)$ have the same vector at that position but produce different outputs when $a\ne a'$. There is no function of that vector alone which represents the attention update for all contexts. A sequence carrier $\mathcal H_{\ell,n}$ contains the needed context. Alternatively, a carrier of located occurrences must include the computational context on which the update depends.
\end{example}

The categorical entrywise action on a structural Token must not be confused with entrywise neural computation. If the carrier elements are whole matrices $H$, a Token map acts on each matrix occurring in a record. Inside the carrier function $F(H)$, attention can mix information from many positions. This is how contextual dependence and the entrywise definition of $F_E(f)$ coexist.

\paragraph{The same vector does not determine the same Token occurrence.}
Thesis~\ref{thesis:occurrence} distinguishes numerical coordinates from the identity of a state in a computation. At a fixed deterministic protocol and parameter choice $\theta$, let $\Omega_\theta$ be the set of located states
\[
\omega=(\theta,\ell,s,i,h_i^{(\ell)}(s)),
\qquad 0\le\ell\le L,\quad 1\le i\le |s|,
\]
over the admitted inputs. The occurrence identity includes the layer, input and position. If repeated invocations of the same input must also be distinguished, add an invocation identifier, as in Section~\ref{sec:elastic-execution}. A structural Token can record a fact about one or several such occurrences; it is not identical to the occurrence it records.

The coordinate projection $\nu$ forgets location. For a common width it is $\nu:\Omega_\theta\to\mathbb R^d$, $\omega\mapsto h_i^{(\ell)}(s)$; with varying widths use $\coprod_d\mathbb R^d$ and retain the dimension tag. Two occurrences can be different even when their coordinate projections are equal. The set-theoretic tuple supplies this identity data. The categorical description then makes the forgetting map, the structural relations and the permitted comparisons explicit. Merely putting different subscripts on equal sets does not do this: a diagram can use the same object at several stages.

\begin{proposition}[When a vector suffices for an observation]\label{prop:observation-factorization}
For an observation $q:\Omega_\theta\to Q$, there is a unique function $\overline q:\nu(\Omega_\theta)\to Q$ with $q=\overline q\circ\nu$ iff
\[
\nu(\omega)=\nu(\omega')\ \Longrightarrow\ q(\omega)=q(\omega')
\quad\text{for all }\omega,\omega'\in\Omega_\theta.
\]
Here the codomain of $\nu$ is restricted to its image. In particular, any observation that distinguishes a pair of equal-coordinate occurrences cannot be recovered from their vectors alone.
\end{proposition}
\begin{proof}
A factorization gives the displayed implication. Conversely define $\overline q(h)$ to be the common value of $q$ on the nonempty fiber $\nu^{-1}(h)$. The implication makes this definition unambiguous, and surjectivity onto the image proves uniqueness.
\end{proof}

\begin{figure}[htbp]
\centering
\begin{tikzcd}[column sep=large,row sep=large]
\Omega_\theta \arrow[r,two heads,"\nu"] \arrow[dr,"q"'] &
\nu(\Omega_\theta) \arrow[d,dashed,"\overline q"]\\
&Q
\end{tikzcd}
\caption{Forgetting location, in $\Set$. The dashed factorization exists exactly when the observation is constant on every coordinate fiber. Equal vectors do not by themselves establish that condition. These are observation functions; Token-map status additionally requires specified heaps.}\label{fig:occurrence-forgetting}
\end{figure}

For example, equal initial embeddings at two positions of a repeated vocabulary item do not recover the position label. A layer label similarly cannot be recovered from the vector when the same vector occurs at two layers. An observation deliberately defined as a fixed readout of the vector does factor. The proposition therefore asks a concrete question of each interpretation, rather than denying that vectors can support useful observations. The contextual attention example above supplies a further reason to retain context: a row's update may depend on other rows. A row $h_i$ is obtained by projection from $H$; it is not itself an element of the whole-sequence carrier $\mathcal H_{\ell,n}$.

\paragraph{Located layer objects and comparison arrows.}
At fixed $\theta,n$, distinguish stages by the tagged carriers
\[
\widehat{\mathcal H}_{\ell,n}=\{\ell\}\times\mathcal H_{\ell,n},
\qquad \kappa_\ell(H)=(\ell,H).
\]
Transport the chosen sequence heap along $\kappa_\ell$. A block then acts by
\[
\widehat F_{\ell,\theta}(\ell,H)=(\ell+1,F_{\ell,\theta}(H)).
\]
It is a Token map between the transported classes exactly when $F_{\ell,\theta}$ is a Token map between the original classes. Thus tagging retains stage identity without changing the previously proved causal or equivariance conditions. The untagged pipeline is recovered through the carrier bijections $\kappa_\ell^{-1}$.

When adjacent residual-stream widths agree, their coordinate identification gives a specific comparison function
\[
\iota_\ell:\widehat{\mathcal H}_{\ell,n}\to\widehat{\mathcal H}_{\ell+1,n},
\qquad(\ell,H)\mapsto(\ell+1,H).
\]
For classes with untagged heaps $\Upsilon_\ell,\Upsilon_{\ell+1}$ over the same carrier and core, this function is a Token map iff $\Upsilon_\ell\subseteq\Upsilon_{\ell+1}$; it is a Token isomorphism iff the heaps are equal. This follows by applying the identity carrier function to their records. In particular, using the same causal heap at each equal-width stage makes $\iota_\ell$ an isomorphism. Different heaps or widths require a separately specified comparison.

\begin{figure}[htbp]
\centering
\begin{tikzcd}[column sep=large,row sep=large]
\widehat{\mathcal H}_{\ell,n} \arrow[r,"\iota_\ell"] \arrow[d,"\kappa_\ell^{-1}"'] &
\widehat{\mathcal H}_{\ell+1,n} \arrow[d,"\kappa_{\ell+1}^{-1}"]\\
(\mathbb R^d)^n \arrow[r,equal] & (\mathbb R^d)^n
\end{tikzcd}
\caption{A coordinate comparison between located stages, in $\Set$. With the same causal heaps transported to both stages, it is also a square of Token maps. The top arrow transports the same coordinates; it is distinct from the actual block update $\widehat F_{\ell,\theta}$.}\label{fig:layer-comparison}
\end{figure}

The residual equation $F_{\ell,\theta}(H)=H+\Delta_{\ell,\theta}(H)$ uses this common numerical coordinate system. It does not assert that a general Token category is additive, or that the comparison arrow equals the computational update. Reading an intermediate state with the final output map~\eqref{eq:readout}, as in a logit lens, uses the chosen coordinate comparison to make a readout well typed. Whether that readout is a faithful interpretation of the intermediate computation is a further empirical and relational question. The tuned lens, for example, learns layer-specific affine probes to improve the reliability of direct intermediate readout \cite{belrose2023}. Section~\ref{sec:interpretation} records the readout procedure and its location explicitly, so that a comparison protocol can be evaluated without assigning a vector a privileged intrinsic meaning.

\subsection{Permutation heaps and their precise scope}\label{sec:equivariance}

Let $S_n$ act on $X^n$ by $(\sigma s)_i=s_{\sigma^{-1}(i)}$, and let $E_n^{\mathrm{sym}}=\{\bar\sigma:\sigma\in S_n\}$ be core symbols naming the permutations.

\begin{definition}[Sequence classes]\label{def:seq}
The \emph{sequence class} and the \emph{positional sequence class} over $X$ are
\begin{align}
\Seq_n(X)&=\left\langle X^n,E_n^{\mathrm{sym}},
\{(\bar\sigma,s,\sigma s):\sigma\in S_n,\ s\in X^n\}\right\rangle,\label{eq:seq}\\
\Seq_n^{\mathrm{pos}}(X)&=\left\langle X^n,E_n^{\mathrm{sym}},\{(\bar1,s,s):s\in X^n\}\right\rangle .\nonumber
\end{align}
\end{definition}

Both carriers are sets of ordered sequences. The heap of $\Seq_n(X)$ records the action of each named permutation, making reindexing explicit at the Token level. The heap of $\Seq_n^{\mathrm{pos}}(X)$ records only self-identity; it neither pins individual positions nor expresses a successor or absolute-position relation. We retain the name ``positional sequence class'' to indicate a background in which position-sensitive functions are allowed without an equivariance requirement. This is a convention about admissible maps, not a claim that the diagonal heap itself encodes positional structure.

\begin{proposition}\label{prop:equivariance}
A function $F:X^n\to Y^n$ is a Token map $\Seq_n(X)\to\Seq_n(Y)$ iff it is permutation equivariant, $F(\sigma s)=\sigma F(s)$ for all $\sigma,s$. Every function $X^n\to Y^n$ is a Token map $\Seq_n^{\mathrm{pos}}(X)\to\Seq_n^{\mathrm{pos}}(Y)$.
\end{proposition}
\begin{proof}
The image of $(\bar\sigma,s,\sigma s)$ is $(\bar\sigma,F(s),F(\sigma s))$. It lies in the target heap iff $F(\sigma s)=\sigma F(s)$. For the positional class only $\sigma=1$ occurs and the condition is void.
\end{proof}

For rows collected in a matrix $H$, unmasked attention without position-dependent terms has
\[
\Attn_\theta(H)=\Softmax\!\left(\frac{(HW_Q)(HW_K)^\top}{\sqrt{d_k}}\right)HW_V,
\]
with softmax applied separately to each row. A permutation matrix $\Pi$ replaces the score matrix $S$ by $\Pi S\Pi^\top$, and rowwise softmax obeys
\[
\Softmax(\Pi S\Pi^\top)=\Pi\Softmax(S)\Pi^\top.
\]
Multiplication by $\Pi HW_V$ then gives $\Attn_\theta(\Pi H)=\Pi\Attn_\theta(H)$. Shared positionwise feed-forward maps and normalization, residual addition, and feature concatenation also preserve this equivariance. Thus unmasked blocks without positional asymmetries are Token maps between the corresponding sequence classes \cite{zaheer2017,bronstein2021}.

The choice of heap decides which functions count as admissible layers. For a fixed $X$, the inclusion of heaps has the direction
\begin{equation}\label{eq:sequence-inclusion}
H_n^{\mathrm{pos}}(X)\subseteq H_n^{\mathrm{sym}}(X),\qquad
\id_{X^n}:\Seq_n^{\mathrm{pos}}(X)\longrightarrow\Seq_n(X).
\end{equation}
The displayed identity-carrier Token map goes from the smaller heap to the larger heap. Replacing permutation heaps by diagonal heaps at both source and target enlarges the admitted hom sets from equivariant functions to all functions. This relaxation of a family of specifications must not be confused with reversing the arrow in~\eqref{eq:sequence-inclusion}.

Even position reversal is a Token map between positional classes. With a fixed additive positional matrix $P$, the map $H\mapsto\Attn_\theta(H+P)$ is generally not permutation equivariant, but it is a Token map between positional classes. With a fixed mask $M$ and no other positional asymmetry, the calculation yields equivariance under permutations satisfying $\Pi M\Pi^\top=M$; the ordinary finite causal mask admits no nonidentity such permutation. The absence of an equivariance requirement says nothing about prefix causality. The latter is an independent relation on sequence pairs, not an intermediate heap in a linear scale from full permutation symmetry to diagonal records. The next subsection defines it and gives counterexamples in both directions.

\subsection{Prefix-agreement heaps characterize causal computation}\label{sec:causal}

Let $E_n^{\mathrm{pre}}=\{c_1,\ldots,c_n\}$ be fresh labels. Define the \emph{causal sequence class}
\begin{equation}\label{eq:causal-class}
\Cau_n(X)=\left\langle X^n,E_n^{\mathrm{pre}},
\{(c_k,s,t):s_{1:k}=t_{1:k},\ 1\le k\le n\}\right\rangle.
\end{equation}
A Token here compares two complete possible sequence states. The relation is equality of their prefixes, not a claim that the model stores or explicitly enumerates pairs of states.

\begin{theorem}[Causal maps]\label{thm:causal}
For nonempty $X$, a function $F:X^n\to Y^n$ is a Token map $\Cau_n(X)\to\Cau_n(Y)$ iff, for every $i$, there is a function $f_i:X^i\to Y$ such that
\begin{equation}\label{eq:causal-factor}
F(s)_i=f_i(s_{1:i}).
\end{equation}
Equivalently, agreeing on an input prefix forces agreement on the corresponding output prefix.
\end{theorem}
\begin{proof}
Preservation of $(c_k,s,t)$ is exactly the implication
\[
s_{1:k}=t_{1:k}\ \Longrightarrow\ F(s)_{1:k}=F(t)_{1:k}.
\]
For $k=i$, this makes $F(s)_i$ independent of the chosen extension of $s_{1:i}$. Every prefix extends because $X$ is nonempty, so it defines $f_i$. Conversely, \eqref{eq:causal-factor} gives equality of each output coordinate up to $k$ whenever those input prefixes agree.
\end{proof}

\begin{figure}[htbp]
\centering
\begin{tikzcd}[column sep=large,row sep=large]
X^n \arrow[r,"F"] \arrow[d,"\pi_{n,k}"'] & Y^n \arrow[d,"\pi_{n,k}"] \\
X^k \arrow[r,"\overline F_k^{(n)}"'] & Y^k
\end{tikzcd}
\caption{Prefix causality as factorization in $\Set$. The projection $\pi_{n,k}$ keeps the first $k$ coordinates. Theorem~\ref{thm:causal} is equivalent to existence of the lower map for every $1\le k\le n$. The superscript records the fixed original length $n$; equality with a separately specified shorter-length model is the additional consistency condition of Section~\ref{sec:lengths}.}\label{fig:causal-factorization}
\end{figure}

Explicitly, $\overline F_k^{(n)}(x_1,\ldots,x_k)=(f_1(x_1),\ldots,f_k(x_1,\ldots,x_k))$. The square concerns which inputs can affect outputs. It places no additional execution order on independent operation nodes inside a valid implementation of $F$. This separates the causal square from the execution diamond of Figure~\ref{fig:ready-diamond}.

This is a substantive restriction on functions. The full self-attention example with uniform weights returns $((a+b)/2,(a+b)/2)$ and fails it: inputs $(0,0)$ and $(0,2)$ agree at the first position but their first outputs differ. Causal attention with uniform permitted weights returns
\begin{equation}\label{eq:two-position}
F(a,b)=\left(a,\frac{a+b}{2}\right).
\end{equation}
It preserves prefix agreement but is not permutation equivariant: for the swap $\sigma$, $F(2,1)=(2,3/2)$ differs from $\sigma F(1,2)=(3/2,1)$. These calculations separate the two heaps without appealing to an empirical benchmark.

The classes in \eqref{eq:causal-class} encode non-anticipation on tokenized sequences. They do not claim that tokenizing a text prefix always gives a prefix of the final tokenization: that is a distinct property of $\tau$. A beginning-of-sequence convention supplies an initial tokenized context where generation needs one.

\subsection{A decoder block with position information}\label{sec:decoder}

We now connect Theorem~\ref{thm:causal} to a concrete architecture. Consider a pre-normalized decoder block of the kind used in LLaMA, with RMSNorm, rotary query/key processing, and a gated feed-forward sublayer \cite{touvron2023}. This is a representative architecture, not a claim that all deployed LLMs use the same block.

Let $H\in\mathbb R^{n\times d}$ have rows $h_i^\top$, put $U=N_1(H)$, and for head $a$ set
\begin{align*}
Q^{a}&=UW_Q^{a},&K^{a}&=UW_K^{a},&Z^{a}&=UW_V^{a},\\
\widetilde q_i^{a}&=R_i^{a}q_i^{a},&
\widetilde k_i^{a}&=R_i^{a}k_i^{a}.&&
\end{align*}
Here $q_i^{a},k_i^{a}$ are the column versions of the matrix rows, $W_Q^{a},W_K^{a}\in\mathbb R^{d\times d_k}$, $W_V^{a}\in\mathbb R^{d\times d_v}$, and $R_i^{a}$ is a prescribed rotary position transformation \cite{su2021}. For the causal mask,
\[
M_{ij}=\begin{cases}0,&j\le i,\\-\infty,&j>i,\end{cases}
\]
define
\begin{equation}\label{eq:masked-attention}
\alpha_{ij}^{a}(H)=
\begin{cases}
\displaystyle\frac{\exp(\langle\widetilde q_i^{a},\widetilde k_j^{a}\rangle/\sqrt{d_k})}
{\sum_{r\le i}\exp(\langle\widetilde q_i^{a},\widetilde k_r^{a}\rangle/\sqrt{d_k})},&j\le i,\\[8pt]
0,&j>i,
\end{cases}
\qquad
o_i^{a}=\sum_{j\le i}\alpha_{ij}^{a}(H)z_j^{a}.
\end{equation}
The multihead output and residual updates are
\begin{align}
\MHA_\theta(U)&=\operatorname{Concat}(O^1,\ldots,O^m)W_O,\nonumber\\
\widetilde H&=H+\MHA_\theta(N_1(H)),\label{eq:block}\\
F_{\ell,\theta}(H)&=\widetilde H+\FFN_\theta(N_2(\widetilde H)).\nonumber
\end{align}
In \eqref{eq:block}, $\MHA$ includes the rotary and mask operations specified above; $W_O\in\mathbb R^{md_v\times d}$. Normalization is applied independently to each position. For example, with learned scale $g$ and $\epsilon>0$,
\[
\RMSNorm(x)=g\odot\frac{x}{\sqrt{d^{-1}\sum_{r=1}^d x_r^2+\epsilon}}
\]
is a total positionwise function \cite{zhang2019}. A rowwise gated feed-forward map can take the form
\[
\FFN(x)=\bigl(\operatorname{SiLU}(xW_g)\odot(xW_u)\bigr)W_d,
\]
with conformable matrices \cite{shazeer2020}. The causal argument requires their positionwise dependence, not this particular activation.

\begin{proposition}[Decoder blocks preserve the causal heap]\label{prop:decoder-causal}
For fixed parameters and fixed positional conventions, the block \eqref{eq:block} is a Token map
\[
F_{\ell,\theta}:\Cau_n(\mathbb R^d)\longrightarrow\Cau_n(\mathbb R^d).
\]
The same statement holds with additive position embeddings and for post-normalized variants when their normalizations remain positionwise.
\end{proposition}
\begin{proof}
If $H$ and $H'$ agree through position $k$, their normalized rows and projected queries, keys, and values agree through $k$. The prescribed $R_i^a$ are the same in both executions. For $i\le k$, every allowed key/value index in \eqref{eq:masked-attention} satisfies $j\le i\le k$, so its numerator, denominator, weights, and weighted sum agree. Concatenation and output projection act within positions. Residual addition, positionwise normalization, and the feed-forward maps preserve equality of each corresponding row. Thus the output prefixes agree, and Theorem~\ref{thm:causal} applies. Adding the same $p_i$ at each corresponding input position also preserves prefix equality.
\end{proof}

The proof covers fixed causal sliding-window masks as well. A mask or positional rule chosen using future token content could violate the hypothesis. Cross-token batch statistics, global routing capacity decisions, retrieval, or external memory require separate dependence analysis. Stochastic training operations can be considered with randomness fixed as additional data, but the proposition concerns deterministic inference.

\subsection{Composition, multiple heads, and output probabilities}

Theorems~\ref{thm:tse} and \ref{thm:causal} make the dependence argument compositional. If $F$ and $G$ preserve prefix agreement, so does $G\circ F$. Pairing heads with the categorical product yields
\[
\Cau_n(Y_1)\otimes\Cau_n(Y_2)\cong\Cau_n(Y_1\times Y_2),
\]
because a product Token compares prefixes with the same label $c_k$ in both factors. The bijection pairs the two feature vectors at each position. Feature concatenation implements the ordinary identification of a product of Euclidean feature spaces with a larger Euclidean space. Residual addition is the Token map
\[
\Cau_n(\mathbb R^d)\otimes\Cau_n(\mathbb R^d)
\longrightarrow\Cau_n(\mathbb R^d),\qquad(H,G)\mapsto H+G.
\]
Consequently, sharing the input through a diagonal, applying parallel heads, projecting their concatenation, and adding a residual are explicit Token maps. The analogous product statement holds for permutation heaps, with a matched label $\bar\sigma$ in each factor.

Embedding lookup and the positionwise output map \eqref{eq:readout} preserve prefix agreement. The deterministic model therefore gives a composite
\begin{equation}\label{eq:causal-pipeline}
\Cau_n(V)\xrightarrow{e_\theta^n}\Cau_n(\mathbb R^d)
\xrightarrow{F_{L-1,\theta}\circ\cdots\circ F_{0,\theta}}
\Cau_n(\mathbb R^d)\xrightarrow{r_\theta^n}\Cau_n(\Delta(V)).
\end{equation}
The preserved relation is prefix agreement. This does not assert preservation of each numerical activation, attention weight, or output probability across layers. Their values are expected to change.

Sampling a token from $\Delta(V)$ is stochastic and is not a deterministic Token map on that carrier. A fixed decoding rule gives a function, while random sampling can be modeled by an additional randomness input or by a stochastic extension of the category. Those choices should not be hidden inside the word ``map.'' A specified training update can also be represented in the record language. Its interpretation as gradient-based learning additionally uses differential and optimization structure, as explained in Section~\ref{sec:learning-dynamics}.

\subsection{Position-indexed presentations and prefix naturality}
\label{llm:prefix}\label{sec:lengths}

Computational dependency and sequence causality are different structures. A
bidirectional Transformer has an acyclic forward computation for each fixed
input length, even though an early output can depend on later input
positions. Designated roots therefore do not by themselves distinguish a
causal decoder from a bidirectional encoder. We add positions and a coherent
family of interpreted presentations. The construction applies to fixed
parameters in deterministic inference; it is not a statement about the
training procedure or about sampled generation. Section~\ref{sec:learning-dynamics} distinguishes this sequence causality from training-time causality and the dependencies of a training computation.

Let $N\in\mathbb N\cup\{\infty\}$ be a context bound. When $N=\infty$, all length and position indices in what follows are still finite; $\infty$ is not an additional object or stage. Use one common core
$E$ containing all operation, input, parameter and output tags needed up to
that bound. In particular the tags for position $i$ have the same identity
at all lengths $n\ge i$. This permits a varying-interface family inside
$\mathsf{DRTS}_{E,j}$; it is not a family of objects of one fixed-interface graph
category from Section~\ref{sec:finite-representations}.

\begin{definition}[Prefix-compatible occurrence presentation]
\label{llm:def-prefix-presentation}
A prefix-compatible occurrence presentation consists of a finite set $P$
of shared parameter occurrences, input occurrences $x_i$ and output
occurrences $y_i$ for $1\le i\le N$, a carrier $B$, and a position map
$p:B\to\{0,\ldots,N\}$ (to $\mathbb N$ when $N=\infty$), subject to the
following conditions.
\begin{enumerate}
\item $p^{-1}(0)=P$, $p(x_i)=p(y_i)=i$, and every
      $B_n=\{v:p(v)\le n\}$ is finite.
\item Every occurrence is defined exactly once: as a parameter, as one of
      the distinct inputs, or by one finite operation record
      $(\omega,u_1,\ldots,u_k,v)$ with designated result $v$.
      The operand-to-result graph is acyclic and
      $p(u_r)\le p(v)$ for every operand. Repeated operands are permitted.
\item The heap consists of those defining records and the unary-carrier
      parameter, input and output records. Its restriction to $B_n$ is
      denoted $H_n$; no defining record or label changes with $n$.
\end{enumerate}
Equip $B_n$ with the reflexive transitive closure of its operand edges and
with initial set $I_n=P\cup\{x_1,\ldots,x_n\}$. Write $G_n$ for the
resulting ordered Token structure. An interpretation assigns a fixed set
of values to each occurrence sort, a total function of the indicated
operand sorts to each operation, and fixed values to $P$. These assignments
are the same at every length. Inputs have a common value set $X$ and
outputs a common value set $Y$.
\end{definition}

The values need not be finite sets. They may be tagged vector spaces,
finite vocabulary symbols or tensor values. The finite condition concerns
occurrences and records. The order is on occurrences, not on numerical
values. Totality can be replaced by a specified common domain on which
all the operations are defined; the theorem below then requires that
domain to be stable under the restrictions being used.

\begin{proposition}[Structural prefix embeddings]
\label{llm:prefix-embedding}
Every $G_n$ lies in $\mathsf{DRTS}_{E,j}$. For $m\le n$, $B_m$ is lower in $B_n$,
and its induced heap, order and initial set are exactly those of $G_m$.
The inclusions $\iota_{m,n}:G_m\hookrightarrow G_n$ are morphisms and
satisfy $\iota_{n,r}\iota_{m,n}=\iota_{m,r}$.
\end{proposition}
\begin{proof}
Each operand is below its result in the generated order, so the operation
records satisfy the designated-root comparisons; unary-carrier records
satisfy them automatically. Inputs and parameters have no incoming
operand edges, hence $I_n$ is lower. Position cannot decrease along an
operand edge. Thus every ancestor of a member of $B_m$ is in $B_m$.
Every path between members of $B_m$ lies there, proving equality of the
restricted orders. The heaps are induced by construction and
$I_n\cap B_m=I_m$. These facts prove record, order and initial-set
preservation. The last identity is an identity of inclusions.
\end{proof}

There is a small but essential qualification when invoking the completed
fragment results of Section~\ref{v13:causal:section}. In general $B_m$ does
not contain $I_n$: the inputs at positions $m+1,\ldots,n$ are missing.
Thus $G_m$ itself is not a completed fragment relative to $I_n$. At fixed
length $n$ the appropriate completed prefix is
\[
 U_m^{(n)}=I_n\cup B_m.
\]
It is lower, contains $I_n$, and its induced structure is $G_m$ with the
additional unused future input occurrences and all induced unary records on them.
The completed-fragment gluing criterion of Theorem~\ref{v13:causal:gluing} applies to these padded fragments and their completed subfragments. No extra future-input
dependence is thereby introduced into the outputs of $G_m$.

Topological evaluation of $G_n$ gives a function
$F_n:X^n\to Y^n$; the fixed parameter assignment is suppressed in the
notation. Write $\rho^X_{n,m}:X^n\to X^m$ and
$\rho^Y_{n,m}:Y^n\to Y^m$ for prefix projections, including the unique
maps to the zero-length products.

\begin{theorem}[Naturality of interpreted causal prefixes]
\label{llm:prefix-naturality}
For every interpreted prefix-compatible occurrence presentation,
\begin{equation}\label{eq:prefix-natural}
 \rho^Y_{n,m} F_n = F_m\rho^X_{n,m}\qquad(0\le m\le n\le N).
\end{equation}
Let $\mathcal L_N$ have objects $0,\ldots,N$ and one arrow $n\to m$
exactly when $m\le n$. Then $n\mapsto X^n$ and $n\mapsto Y^n$, with
prefix projections, are functors $\mathcal L_N\to\Set$, and the family
$F=(F_n)_n$ is a natural transformation. The occurrence presentations,
with their inclusions, instead form a functor
$\mathcal L_N^{\mathrm{op}}\to\mathsf{DRTS}_{E,j}$.
\end{theorem}
\begin{proof}
Fix $a=(a_1,\ldots,a_n)$. Evaluate $G_n$ at $a$ and restrict its valuation
to $B_m$. The parameter values remain fixed and the retained input values
are exactly $a_1,\ldots,a_m$. Every operation defining a retained
occurrence has all its operands in $B_m$, the same operation label and
the same interpretation as in $G_m$. The restricted valuation therefore
satisfies all defining equations of $G_m$. Uniqueness of topological
evaluation identifies it with the evaluation of $G_m$ on the prefix of
$a$. In particular the values at $y_1,\ldots,y_m$ agree, proving
\eqref{eq:prefix-natural}. Prefix projections and inclusions compose as
stated, which proves the functorial and naturality assertions.
\end{proof}

\begin{figure}[htbp]
\centering
\begin{tikzcd}[column sep=large,row sep=large]
 X^n\arrow[r,"F_n"]\arrow[d,"\rho^X_{n,m}"'] &
 Y^n\arrow[d,"\rho^Y_{n,m}"]\\
 X^m\arrow[r,"F_m"'] & Y^m
\end{tikzcd}
\caption{Prefix naturality in numerical semantics. Restricting the input before evaluation agrees with restricting the output after evaluation. The occurrence presentations have inclusion arrows in the opposite length direction.}
\label{llm:prefix-naturality-diagram}
\end{figure}

\begin{proposition}[The extensional meaning of prefix naturality]
\label{llm:prefix-kernels}
A family $F_n:X^n\to Y^n$ is natural under prefix projections if and only
if there are functions $f_i:X^i\to Y$ such that
\[
 F_n(a_1,\ldots,a_n)_i=f_i(a_1,\ldots,a_i)
 \quad(1\le i\le n).
\]
When they exist, the functions are uniquely determined by
$f_i=\pi_i F_i$.
\end{proposition}
\begin{proof}
Apply naturality with $m=i$ and take its final component. Conversely,
the displayed component formula gives equality of all $m$ coordinates
in~\eqref{eq:prefix-natural}.
\end{proof}

This proposition characterizes a behavior, whereas
Definition~\ref{llm:def-prefix-presentation} supplies a sufficient
structural realization of that behavior. The converse from behavior to
the structure of a given presentation is false: a presentation may read
a future input and cancel its contribution. For example, over $\mathbb R$
one can compute an identity family by setting
\[
 y_i=x_i+\sum_{j=i+1}^{n}(x_j-x_j).
\]
The denotation is prefix-natural, but this displayed implementation reads
future positions. Likewise a bidirectional architecture can acquire
prefix-natural behavior for special parameters, such as parameters
giving constant outputs. A structural architecture guarantee means that
the property holds for every admissible parameter assignment under the
declared execution assumptions, not merely for one trained model.

For the decoder block of Proposition~\ref{prop:decoder-causal}, give
each input, hidden state, projection and attention result its layer,
position and operation occurrence identity. Fixed weights may be shared
initial occurrences. A causal row reads only positions at most its own,
and the local normalization, residual and feed-forward operations use
the same position. The resulting finite operation graph satisfies the
definition, provided each retained operation and positional convention
is unchanged by extending the prefix. The definition also permits fixed
external encoder memory. Every finite prefix is finite even when the
chosen family is unbounded; an actual context limit restricts its index
set. The order is on occurrences, not on untagged numerical vectors.

Two explicit failures distinguish the assumptions. Bidirectional
averaging $F_n(x)_i=n^{-1}\sum_jx_j$ has acyclic operation graphs, but
$F_1(0)=0$ and $F_2(0,2)_1=1$. The pointwise family $F_n(x)_i=x_i+n$
is causal at every fixed length, but its first output changes from $1$
to $2$ under the same extension. Thus an acyclic graph, a causal
fixed-length map and a coherent family are different requirements.
The naturality square is standard in categorical accounts of causal
streams \cite[Section~2.1]{spiwack2018}; the construction identifies
the occurrence and interpretation assumptions that establish it here.
Appendix~\ref{llm:generic-prefix} gives a free-algebra criterion when
interpretations are allowed to vary universally.

Equation~\eqref{eq:prefix-natural} compares outputs at already present
positions. It does not keep the next token's distribution unchanged
after a new input: that distribution uses the new last row. Collecting
length classes in a coproduct does not itself impose coherence. The
following invariant makes precise the additional relation to cached
execution, an established implementation of incremental attention
\cite[Section~2.4]{shazeer2019}.

\subsection{A deterministic causal decoder and its cache invariant}
\label{llm:llm-cache}

The next result isolates the numerical assumptions needed to compare two
execution organizations. It is a theorem about deterministic maps of sets;
it does not assert that arbitrary numerical operations are morphisms of
$\mathsf{OTS}_E$. Dependency order is carried by the occurrences of
Section~\ref{llm:prefix}; the typed numerical valuation is separate.

Fix an input set $X$, an output set $Y$, a number of layers $L$, and sets
$H_0,\ldots,H_L$. At layer $\ell$ let $Q_\ell,K_\ell,V_\ell,A_\ell$
be query, key, value, and attention-result sets. For every position $i\geq1$
fix total deterministic maps
\begin{align*}
 e_i&:X\longrightarrow H_0,\\
 q_{\ell,i}&:H_{\ell-1}\longrightarrow Q_\ell,
 &k_{\ell,i}&:H_{\ell-1}\longrightarrow K_\ell,
 &v_{\ell,i}&:H_{\ell-1}\longrightarrow V_\ell,\\
 a_{\ell,i}&:Q_\ell\times(K_\ell\times V_\ell)^i
                  \longrightarrow A_\ell,
 &r_{\ell,i}&:H_{\ell-1}\times A_\ell\longrightarrow H_\ell,
 &o_i&:H_L\longrightarrow Y.
\end{align*}
All these maps are fixed independently of the total prefix length.
For $x=(x_1,\ldots,x_n)$ define, for $1\leq i\leq n$,
\begin{align}
 h_i^{0,n}&=e_i(x_i),\notag\\
 (q_i^{\ell,n},k_i^{\ell,n},v_i^{\ell,n})
  &=\bigl(q_{\ell,i}(h_i^{\ell-1,n}),
          k_{\ell,i}(h_i^{\ell-1,n}),
          v_{\ell,i}(h_i^{\ell-1,n})\bigr),\notag\\
 h_i^{\ell,n}
  &=r_{\ell,i}\!\left(h_i^{\ell-1,n},
       a_{\ell,i}\!\left(q_i^{\ell,n},
                   (k_j^{\ell,n},v_j^{\ell,n})_{j=1}^{i}\right)\right),
       \qquad 1\leq\ell\leq L,\label{llm:decoder-recurrence}\\
 y_i^n&=o_i(h_i^{L,n}).\notag
\end{align}
Write $F_n(x)=(y_1^n,\ldots,y_n^n)$ and let $F_0$ be the unique map
$X^0\to Y^0$.

The query, key and value sets may bundle several heads. Feature-wise
normalization, residual operations and position-wise feed-forward maps
can be absorbed into the displayed local maps; both usual pre-normalized
and post-normalized blocks have this form. Absolute position information,
rotary transforms and relative attention biases are allowed when their
rules for fixed positions are unchanged by extending the prefix. A causal
mask may select any fixed subset of $\{1,\ldots,i\}$ inside $a_{\ell,i}$.
The hypothesis excludes a positional rescaling depending on the current
total length, normalization using future positions, and other operations
that alter old keys or values after an extension. Parameters are fixed,
execution is deterministic, and this first invariant uses a full cache without truncation or lossy
conversion. Safe deletion for a declared access family is established
separately in Theorem~\ref{llm:eviction-theorem}. Equality is interpreted in exact arithmetic;
floating-point implementations require a separate error or operation-order
analysis.

\begin{lemma}[Layer-wise prefix consistency]\label{llm:layer-prefix}
If $m\leq n$, $i\leq m$ and $x\in X^n$, every state and every query,
key and value at $(\ell,i)$ computed on $x$ equals the corresponding
quantity computed on $x_{1:m}$. Consequently
\[
 \rho^Y_{n,m}F_n=F_m\rho^X_{n,m},
\]
where $\rho$ denotes prefix restriction.
\end{lemma}
\begin{proof}
At layer zero both quantities are $e_i(x_i)$. Suppose the claim holds at
layer $\ell-1$ for all $i\leq m$. The fixed local maps give equal queries,
keys and values at those positions. For a position $i\leq m$, the ordered
attention operand list contains exactly the positions $1,\ldots,i$ in
both computations. Applying the same $a_{\ell,i}$ and $r_{\ell,i}$ gives
equal layer-$\ell$ states. Induction over the layers and then the fixed
output maps prove the assertions.
\end{proof}

For $n\geq0$ define the cache space and the full-prefix cache map by
\begin{equation}\label{llm:cache-space}
 C_n=\prod_{\ell=1}^{L}(K_\ell\times V_\ell)^n,
 \qquad
 c_n(x)=\bigl((k_j^{\ell,n},v_j^{\ell,n})_{j=1}^{n}\bigr)_{\ell=1}^{L}.
\end{equation}
The space $C_0$ is a singleton. An incremental transition
\[
 u_n:C_n\times X\longrightarrow Y\times C_{n+1}
\]
is defined as follows. Given an old cache and a new input $z$, put
$i=n+1$ and $h=e_i(z)$. At layer $\ell$, compute
$q=q_{\ell,i}(h)$, $k=k_{\ell,i}(h)$ and $v=v_{\ell,i}(h)$, append
$(k,v)$ to the old layer-$\ell$ cache, and replace $h$ by
$r_{\ell,i}(h,a_{\ell,i}(q,\text{appended cache}))$. After the last layer,
return $o_i(h)$ and all the appended cache lists. This definition is total
on $C_n\times X$, including cache arrays that are not reachable from any
input prefix.

\begin{theorem}[Full-prefix and cache simulation]\label{llm:cache-simulation}
For $x\in X^n$ and $z\in X$, write $xz=(x_1,\ldots,x_n,z)$. Then
\begin{equation}\label{llm:cache-invariant}
 u_n(c_n(x),z)=\bigl(y_{n+1}^{n+1}(xz),c_{n+1}(xz)\bigr).
\end{equation}
In particular, repeated incremental execution from $C_0$ yields exactly
$F_n(x)$ as its output sequence and $c_n(x)$ as its final cache. If
\[
 r_n:X^n\times X\longrightarrow Y\times X^{n+1},\qquad
 r_n(x,z)=\bigl(y_{n+1}^{n+1}(xz),xz\bigr),
\]
then the square in Figure~\ref{llm:cache-square} commutes in $\Set$.
\end{theorem}

\begin{figure}[htbp]
\centering
\begin{tikzcd}[column sep=large]
 X^n\times X \arrow[r,"r_n"] \arrow[d,"c_n\times\id_X"']
   &Y\times X^{n+1} \arrow[d,"\id_Y\times c_{n+1}"]\\
 C_n\times X \arrow[r,"u_n"'] &Y\times C_{n+1}
\end{tikzcd}
\caption{An exact simulation in $\Set$. The upper path extends the full input history; the lower path updates only the cache. The cache map preserves both the emitted output and the state required for the next step.}
\label{llm:cache-square}
\end{figure}
\begin{proof}
By Lemma~\ref{llm:layer-prefix}, every old entry of $c_n(x)$ equals the
entry at the same layer and position in $c_{n+1}(xz)$. At the new position,
the incremental layer-zero state is $e_{n+1}(z)=h_{n+1}^{0,n+1}$.
Suppose the current incremental state is $h_{n+1}^{\ell-1,n+1}$.
Its query, key and value therefore agree with their full-prefix
counterparts. Together with the old-entry equality, this identifies the
entire appended layer-$\ell$ cache with the corresponding list in
$c_{n+1}(xz)$. The attention and local block maps then produce
$h_{n+1}^{\ell,n+1}$. Layer induction proves equality of every appended
cache and of the final output. This proves~\eqref{llm:cache-invariant}
and the displayed square. Induction on $n$, starting with the unique empty
cache, proves the assertion about repeated execution.
\end{proof}

\begin{remark}[State semantics and the meaning of the square]
\label{llm:cache-state-semantics}
For the unbounded idealized recurrence just specified, let
$\mathsf{Hist}=\coprod_{n\geq0}X^n$ and $S=\coprod_{n\geq0}C_n$. The length tag
in $S$ supplies the next position. The maps $u_n$ assemble into
$u:S\times X\to Y\times S$; equivalently, by symmetry and currying,
$X\to(Y\times S)^S$. This is a Kleisli arrow for the ordinary state monad
\cite{moggi1991}
$T_S(B)=(B\times S)^S$ on $\Set$. Similarly the $r_n$ assemble into a
full-history transition on $\mathsf{Hist}$, and $c=\coprod_n c_n:\mathsf{Hist}\to S$ is a
state-machine simulation preserving outputs and transitions. It need not
be injective, invertible, or a minimal-state representation. No state-monad
structure specific to Transformers is claimed. A finite context limit
instead gives the indexed transitions only for the allowed lengths, unless
an overflow policy is explicitly added. Sampling requires a further
random-state or probabilistic interpretation; the deterministic output
$Y$ may already consist of logits or probability distributions.
\end{remark}

\begin{corollary}[Coalgebraic simulation and future sufficiency]
\label{llm:cache-coalgebra}
Under that unbounded assumption, let $\mathcal H(Z)=(Y\times Z)^X$. Currying the assembled transitions
in their input argument gives $\mathcal H$-coalgebras on $\mathsf{Hist}$ and $S$,
and $c:\mathsf{Hist}\to S$ is a coalgebra morphism. If $c(p)=c(p')$, every finite
common continuation produces the same future output sequence from $p$
and $p'$.
\end{corollary}
\begin{proof}
For $f:Z\to Z'$, the functor action is
$\mathcal H(f)(g)=(\id_Y\times f)\circ g$. The coalgebra-morphism equation
is therefore precisely the assembled simulation square, after currying.
For the last assertion, equal cache states and equal next inputs give
equal outputs and equal next cache states. Induction on continuation
length completes the proof.
\end{proof}

\begin{figure}[htbp]
\centering
\begin{tikzcd}[column sep=huge,row sep=large]
 \mathsf{Hist}\arrow[r,"c"]\arrow[d,"\bar r"'] & S\arrow[d,"\bar u"]\\
 (Y\times\mathsf{Hist})^X\arrow[r,"(\id_Y\times c)^X"'] & (Y\times S)^X
\end{tikzcd}
\caption{The cache simulation after currying in the input: a morphism of Mealy coalgebras. Here $\bar r$ and $\bar u$ are the curried history and cache transitions. Currying the cache transition in its state argument instead gives a state-monad Kleisli arrow.}
\label{llm:coalgebra-square}
\end{figure}

This is the ordinary deterministic Mealy-machine construction
\cite[Section~3.8]{rutten2000}. The two
currying choices of $S\times X\to Y\times S$ exhibit, respectively,
a state-monad computation and a coalgebra of input/output behavior.
The sufficiency assertion concerns future outputs; previously emitted
outputs need not be recoverable from the cache. It neither identifies
all behaviorally equivalent states nor supplies a cache-minimization
algorithm.

\paragraph{From occurrences to retained values.}
Expanding the recurrence gives a typed call DAG with separate identities
for keys, values and hidden activations. A fixed finite prefix admits
the graph-to-designated-root construction of
Proposition~\ref{v13:causal:graph-bridge}; the same record construction
works for total maps of nonempty value sets. Its valuation forgets
the occurrence order and sends operation records to their true numerical
graphs. This connects the structural and numerical accounts without
declaring all numerical functions to be ordered Token maps.

At a completed position cut, old keys and values are the operands of
later rows in the displayed recurrence. They therefore form the
corresponding retained-value interface of
Proposition~\ref{prop:residual}, with fixed parameters and future inputs
supplied separately. The cache array is a store of values, not an
induced lower subgraph containing all their past dependencies. The
residual presentation removes the old defining records of boundary
nodes and makes them inputs. Universal gluing reconstructs historical
records; the residual factorization determines sufficient values for
continuation. The latter is the result needed for deletion.

\subsection{Safe eviction from a prescribed attention interface}
\label{llm:safe-eviction}

The full cache in the preceding simulation is sufficient, but a declared
access policy can permit exact deletion. We now specialize the live
frontier of Definition~\ref{def:residual} to that policy. The result is a
statement about an explicit operand interface. It is not a lower bound on
all representations of the same input--output mapping.

Fix a horizon $N\in\mathbb N\cup\{\infty\}$ and, for every layer $\ell$
and position $1\leq i\leq N$, a fixed access set
\[
 \mathcal A_{\ell,i}\subseteq\{1,\ldots,i\}.
\]
All finite subsets below are ordered increasingly. Replace the attention
map in the decoder recurrence by the total map
\[
 \widehat a_{\ell,i}:
 Q_\ell\times\prod_{j\in\mathcal A_{\ell,i}}(K_\ell\times V_\ell)
 \longrightarrow A_\ell.
\]
The other maps and prefix-coherence assumptions are unchanged. An empty
access set is allowed when a total empty-context operation is specified;
ordinary normalized softmax requires a nonempty context. The record of
this operation lists exactly these operands. A full-list operation that
internally masks entries still displays full-list dependencies; replacing
it by the selected-operand record requires the stated factorization
through coordinate restriction. If a mask depends on data, an exact
policy-independent interface requires the union of its possible accesses,
or a separately proved policy-specific invariant.

At the quiescent boundary after position $n$, let
\begin{equation}\label{llm:future-read-set}
 B_{\ell,n}
 =\{j\leq n:\text{there exists }i\text{ with }n<i\leq N
                       \text{ and }j\in\mathcal A_{\ell,i}\}.
\end{equation}
For $N=\infty$, the upper inequality means any later finite position.
These are sets of occurrence indices, not sets of distinct numerical
values. A key--value pair is treated as one addressable operand; a finer
presentation can separate its components or heads.

\begin{theorem}[Exact future-read interface and safe eviction]
\label{llm:eviction-theorem}
For the fixed access family, the following statements hold.
\begin{enumerate}
\item For $0\leq n<N$, the old operands for step $n+1$ and the next retained indices satisfy
\begin{equation}\label{llm:boundary-update}
 \mathcal A_{\ell,n+1}\cap\{1,\ldots,n\}\subseteq B_{\ell,n},
 \qquad
 B_{\ell,n+1}\subseteq B_{\ell,n}\cup\{n+1\}.
\end{equation}
Consequently all prescribed reads are supported by retaining exactly
$B_{\ell,n}$ after each step, computing the new pair, and restricting to
$B_{\ell,n+1}$ after its last current use.
\item Suppose an implementation stores old pairs by their occurrence
indices, must perform every read in the displayed access family, and
cannot regenerate, reconstruct, or obtain an omitted old pair from
another store. A retained set $M_{\ell,n}\subseteq\{1,\ldots,n\}$ can
support all those future reads if and only if
$B_{\ell,n}\subseteq M_{\ell,n}$. Thus $B_{\ell,n}$ is the unique least
set of retained old indices for this representation and horizon.
\item Let $C_n$ be the full cache space, $u_n^{\mathcal A}$ its incremental
transition for these access maps, and define
\begin{equation}\label{llm:pruned-cache}
 C_n^B=\prod_{\ell=1}^L\prod_{j\in B_{\ell,n}}(K_\ell\times V_\ell),
 \qquad \pi_n:C_n\longrightarrow C_n^B
\end{equation}
by coordinate restriction. There is a total transition
$u_n^B:C_n^B\times X\to Y\times C_{n+1}^B$ such that, for every $n<N$,
\begin{equation}\label{llm:eviction-factorization}
 u_n^B\circ(\pi_n\times\mathrm{id}_X)
 = (\mathrm{id}_Y\times\pi_{n+1})\circ u_n^{\mathcal A}.
\end{equation}
The identity holds on all full cache arrays. On reachable arrays it gives
exactly the full-prefix output and the retained full-prefix pairs.
\end{enumerate}
\end{theorem}
\begin{proof}
If an old index occurs at position $n+1$, that position witnesses its
membership in $B_{\ell,n}$. If $j\in B_{\ell,n+1}$ and $j\leq n$, its
witness is also later than $n$, so $j\in B_{\ell,n}$. The only remaining
possibility is $j=n+1$. This proves~\eqref{llm:boundary-update}.

For necessity in the second assertion, an omitted $j\in B_{\ell,n}$ has
a future displayed read. By the stipulated storage discipline it can
never again become available. That read therefore cannot be performed.
For sufficiency, retain its pair until its last required read, or retain
it throughout if there is no last read. At each step the old operands
are available by the first inclusion, the current pair is computed, and
the second inclusion permits the next retained state to be formed. This
inductively supports every finite continuation, including every finite
prefix when the horizon is unbounded.

Construct $u_n^B$ layer by layer. Starting with the embedding of the new
input at position $n+1$, compute its query, key, and value. Supply each
requested old pair from the layer's $B_{\ell,n}$ coordinates and the
current pair when requested; apply $\widehat a_{\ell,n+1}$ and the fixed
local update. Output only the old and new pairs indexed by
$B_{\ell,n+1}$. Equation~\eqref{llm:boundary-update} makes this a total
operation even on unreachable cache arrays.

Compare this construction with $u_n^{\mathcal A}$ on any full array.
At layer zero their new hidden states agree. If the new hidden states
agree before a layer, its local projections agree, each requested old
pair is supplied by the same coordinate, and its attention and local
update therefore agree. The retained new pairs and all retained old
pairs also agree. Induction over the layers proves the output equality
and every coordinate equality in~\eqref{llm:eviction-factorization}.
The full-prefix/cache invariant then gives the reachable-state assertion.
\end{proof}

\begin{figure}[htbp]
\centering
\begin{tikzcd}[column sep=large,row sep=large]
 C_n\times X \arrow[r,"u_n^{\mathcal A}"]
   \arrow[d,"\pi_n\times\mathrm{id}_X"']
 &Y\times C_{n+1}\arrow[d,"\mathrm{id}_Y\times\pi_{n+1}"]\\
 C_n^B\times X\arrow[r,"u_n^B"']&Y\times C_{n+1}^B
\end{tikzcd}
\caption{Coordinate deletion is an exact simulation for the prescribed
access family. The square commutes on arbitrary full cache arrays; the
full-prefix invariant identifies its behavior on reachable states.}
\label{llm:eviction-square}
\end{figure}

\paragraph{Relation to the general frontier.}
Fix a finite future continuation and display its numerical operations as
an occurrence graph. At a completed prefix cut, the old key--value pairs
read by its remaining operations are exactly the corresponding part of
the live frontier. Proposition~\ref{prop:residual} proves their sufficiency
for residual evaluation, and Theorem~\ref{thm:frontier-functor} explains
why successive changes of frontier compose. The proof above additionally
exhibits a uniform transition as future inputs arrive. The future inputs
are supplied externally, and parameters are fixed in the local maps; if
parameters are displayed as value occurrences, their required values must
also be retained or supplied.

The frontier in Definition~\ref{def:residual} also retains an already
computed output that must still be returned. Formula~\eqref{llm:future-read-set}
counts only the state needed for future computation. Outputs already
emitted to an external observer need not remain in that state, but an
implementation promising the entire final output list must store or
otherwise preserve that list separately. Deletion is performed at a
quiescent boundary: outstanding reads and independent observers impose
additional lifetime requirements.

\begin{corollary}[Persistent window cache]
\label{llm:window-eviction}
Suppose the layer's fixed attention window includes the current position
and has width $w_\ell\geq1$:
\[
 \mathcal A_{\ell,i}
 =\{\max(1,i-w_\ell+1),\ldots,i\}.
\]
For $0\leq n<N$,
\begin{equation}\label{llm:window-boundary}
 B_{\ell,n}
 =\{j:1\leq j\leq n,\ j\geq n-w_\ell+2\},
 \qquad |B_{\ell,n}|=\min(n,w_\ell-1).
\end{equation}
For a finite horizon, $B_{\ell,N}=\varnothing$.
\end{corollary}
\begin{proof}
The first future position $n+1$ reads exactly the old indices in the
displayed interval. The lower window endpoint is nondecreasing at all
later positions, so no other old index will be read. At the final horizon
there is no future read. These observations also cover $w_\ell=1$, when
the interval of old indices is empty.
\end{proof}

Thus a persistent window cache holds at most $w_\ell-1$ old pairs per
layer between steps; the current pair and temporary attention values are
additional during a step. This count uses a window that includes the
current position. A convention counting only past positions changes the
formula. It concerns local access in each layer: retained contextual
pairs may encode effects of much earlier inputs, so it does not justify
recomputing from only the last $w_\ell$ raw input tokens.

\paragraph{Effective construction and the limit of minimality.}
For a finite, explicitly given horizon, initialize each old position's
last-use value as $d_\ell(j)=0$ and scan all access lists to update it to
$\max\{i:j\in\mathcal A_{\ell,i}\}$. Then
\[
 B_{\ell,n}=\{j\leq n:d_\ell(j)>n\}.
\]
The scan takes time linear in the total number of listed accesses, apart
from allocation and indexing. For an unbounded horizon, membership in
$B_{\ell,n}$ contains an existential quantifier over all future positions.
An arbitrary computable access family does not supply a decision
procedure for that question; a window rule, a finite-state policy with
an appropriate analysis, or another proved future-access summary is
additional input. For example, a computable policy can read old position
$1$ at a future step precisely when a fixed program has halted by that
step. Deciding whether $1$ ever belongs to its future-read set would then
decide halting. A conservative overapproximation supports safe retention
but may retain more pairs than necessary.

The unique least index set in Theorem~\ref{llm:eviction-theorem} is not
the least state representation of the numerical behavior. If all output
maps are constant, the displayed graph can still prescribe reads of
every past pair, although the future output behavior needs none of them.
An alternative implementation of prefix sums can retain a running sum
instead of all summands; that implementation changes the retained state
and the computation. Numerical equality between two pairs likewise does
not erase their occurrence identities from a direct-address interface.
Semantic identification and exact re-encoding require the behavioral
quotient developed below, together with a separate resource model for
storage or execution costs.

\subsection{Parameterised state for growing and pruned interfaces}
\label{sec:parameterised-cache-state}

The cache interface changes as a computation advances.  A growing prefix
may require more retained values, while a windowed or otherwise pruned
execution may release old values at the same step.  Both are instances of
typed state change.  We use the standard parameterised state construction
of Atkey~\cite{atkey2009}, making its indexing and its relation to the
retained interfaces explicit.  This construction organizes a previously
proved update; its monad laws do not themselves prove that an eviction or
a replacement preserves the model's outputs.

Let $J$ be a set of stage labels and let $C_j$ be a set of admissible state
values at stage $j$.  For a decoder one may take $J=\mathbb N$, or
$J=\{0,\ldots,N\}$ for a bounded context, with $j$ denoting the number of
processed inputs.  Equip $J$ with the \emph{discrete} category structure.
There is therefore no implicit map $C_j\to C_k$ merely because $j\le k$.
The state sets may be the full well-typed arrays, provided that the update
is total on them, or reachable-state sets closed under all admitted next
inputs.  A genuinely partial update requires a partiality or error
semantics in addition to the construction below.

\begin{proposition}[Parameterised state structure]
\label{prop:parameterised-cache-state}
For each $j,k\in J$, define an endofunctor on $\Set$ by
\[
  T_{j,k}B=\Set(C_j,B\times C_k),\qquad
  T_{j,k}f(p)=(f\times\id_{C_k})\circ p.
\]
The family is a parameterised strong monad, with unit,
multiplication and strength
\begin{align*}
 \eta_{j,B}&:B\longrightarrow T_{j,j}B,
 &\eta_{j,B}(b)(s)&=(b,s),\\
 \mu_{j,k,l,B}&:T_{j,k}(T_{k,l}B)\longrightarrow T_{j,l}B,
 &\mu_{j,k,l,B}(P)(s)&=q(t)
       \quad\text{where }P(s)=(q,t),\\
 \operatorname{st}_{A,j,k,B}&:A\times T_{j,k}B
                    \longrightarrow T_{j,k}(A\times B),
 &\operatorname{st}_{A,j,k,B}(a,p)(s)&=((a,b),t)
       \quad\text{where }p(s)=(b,t).
\end{align*}
Equivalently, if $p\in T_{j,k}B$ and $g:B\to T_{k,l}D$,
their sequential composite is
\[
 (p\mathbin{\star}g)(s)=g(b)(t)
       \quad\text{where }p(s)=(b,t).
\]
Here $\star$ denotes this state-passing operation.
\end{proposition}
\begin{proof}
The identity and composition laws for $T_{j,k}$ follow by
postcomposition with $f\times\id$. For $p\in T_{j,k}B$ and
$g:B\to T_{k,l}D$, the correctly indexed unit equations are
\[
 \eta_{k,B}(b)\mathbin{\star}g=g(b),
 \qquad
 p\mathbin{\star}\eta_{k,B}=p.
\]
Both follow by evaluation at a state of their respective domains.
For $g:B\to T_{k,l}D$ and $h:D\to T_{l,r}E$, both
\[
 ((p\mathbin{\star}g)\mathbin{\star}h)(s)
 \quad\text{and}\quad
 (p\mathbin{\star}(b\mapsto g(b)\mathbin{\star}h))(s)
\]
evaluate first $p(s)=(b,t)$, then $g(b)(t)=(d,u)$, and finally
$h(d)(u)$.  This proves the two unit laws and associativity, equivalently
the parameterised multiplication laws.  Naturality in the value sets
holds because a value map is applied to the returned value without
altering the state.  The strength retains an extra value $a$ throughout
this same evaluation; pairing two retained values, or retaining a
terminal value, gives the product coherence laws, and retaining $a$
before or after a unit or a composite gives the compatibility laws.
Naturality and dinaturality in the discrete stage category have only
identity cases.  These are precisely the parameterised strong-monad
axioms.
\end{proof}

\noindent This structure has a pre-state and a post-state index; it is not being identified with
an ordinary monoid-graded monad.  Such an identification would require
additional data relating stages to grades and their composition.

A total incremental update $u_n:C_n\times X\to Y\times C_{n+1}$
now has the curried interface
\[
  \widehat u_n:X\longrightarrow T_{n,n+1}Y,
  \qquad \widehat u_n(x)(s)=u_n(s,x).
\]
Two calls can be composed only when the first post-state type matches
the second pre-state type.  If their inputs are supplied in advance,
the outputs can both be retained:
\begin{align*}
 \operatorname{run}_{n,2}(x_1,x_2)
  ={}&\widehat u_n(x_1)\mathbin{\star}
       \bigl(y_1\mapsto
       \widehat u_{n+1}(x_2)\mathbin{\star}
       (y_2\mapsto\eta_{n+2,Y^2}(y_1,y_2))\bigr)\\
 &\in T_{n,n+2}(Y^2).
\end{align*}
More generally, ordinary products retain a fixed-length output trace:
\[
 \operatorname{run}_{n,r}:X^r\longrightarrow
                  T_{n,n+r}(Y^r),
 \qquad \operatorname{run}_{n,0}(*)=\eta_{n,1}(*).
\]
Induction on $r$ gives exactly the repeated update, and the associativity
law makes regrouping consecutive blocks immaterial.  For a bounded
context only $n+r\le N$ is admitted.  A variable-length supplied trace
has the dependent output family $Y^r\times C_{n+r}$: it can be collected
as a tagged coproduct over the admitted $r$, without pretending that
all final states have the same type.  Autoregressive generation further
supplies each next input by a declared deterministic policy, or by an
explicit random-state or probabilistic effect.

The proposition defines $T_{j,k}$ for \emph{all} stage pairs, including
when $J$ is bounded.  Actual forward execution may select only arrows
with $n\le m$; this selection is closed under sequential composition and
contains every identity.  It should not be confused with replacing the
discrete parameter category by the ordered category: the latter would
require additional reindexing maps and coherence proofs.  Retained-state
size may decrease during a forward step.  The stage index records
elapsed input, not the number of cache entries, and equal cache
cardinalities need not mean equal coordinate identities or equal
positional protocols.

For example, if $B_{\ell,n}$ is the certified retained set of old
positions at layer $\ell$, a pruned value interface has the form
\[
 C_n^{\mathrm{pr}}=
 \prod_{\ell=1}^{L}\prod_{j\in B_{\ell,n}}
      (K_\ell\times V_\ell),
\]
with any further required metadata included.  Its coordinates retain
their position labels even when the sets $B_{\ell,n}$ shift.  Theorem~\ref{llm:eviction-theorem} supplies an example. Suppose
projection $q_n:C_n^{\mathrm{full}}\to C_n^{\mathrm{pr}}$ and the two
update families satisfy the safe-eviction simulation
\[
 u_n^{\mathrm{pr}}\circ(q_n\times\id_X)
   =(\id_Y\times q_{n+1})\circ u_n^{\mathrm{full}}
\]
on the declared states.  Repeated state-passing composition then
preserves this equation by induction on the number of steps, and hence
preserves the output trace.  When this equation is established only on
reachable states, the induction has that same scope.  Such a simulation
does not generally induce a transformation between the \emph{entire}
parameterised state monads: an arbitrary full-state computation may
inspect exactly the information discarded by $q_n$.  The simulation
compares the certified update families.

\begin{remark}[Tagged state and structural interfaces]
\label{rem:parameterised-cache-boundary}
An ordinary state presentation with $S=\coprod_{n\in J}C_n$ retains the
stage label in its coproduct tag.  It is a valid alternative when the
global transition and, for a bounded run, terminal behavior are supplied.
The parameterised presentation strengthens the interface discipline:
the required pre-state and resulting post-state are visible in the
arrow's type rather than solely in an invariant on $S$.  Neither
presentation proves a cache-minimality or memory bound.

The construction here is in $\Set$, where its arguments concern stored
values.  Replacing $C_j$ by Token classes and function sets by internal
function objects is legitimate in a suitable cartesian closed category,
but then each concrete update and state change must preserve its
declared heaps.  Numerical update equations do not establish these
extra obligations automatically.  For the type-oriented companion,
the interface $A\to T_{j,k}B$ therefore gives a precise point of contact:
it types a state change, while selected heaps specify additional
structural obligations on that change.
\end{remark}

\subsection{Concurrency, numerical execution and elastic resources}\label{sec:llm-elasticity}

The decoder equations also identify where a represented-mapping machine can expose concurrent tasks. Given a layer input, distinct heads can compute their projections and attention outputs independently before head combination. Many computations at different known positions also admit parallel evaluation: a causal mask restricts which input rows a position can read, but does not require that the output at the previous position be computed first. Layer composition, head combination and residual additions introduce their own dependencies. An exact primitive-level graph records these obligations more precisely than a single arrow for the whole layer. Whenever this finite-precision description satisfies Definition~\ref{def:computation-graph}, Proposition~\ref{v13:causal:graph-bridge} places it in the designated-root category. This occurrence-level order concerns execution dependencies, whereas the prefix-agreement heap constrains behavior across alternative input sequences.

During autoregressive generation, the next sampled or selected vocabulary item becomes part of the subsequent input. That value is not available before the corresponding decision, so a readiness-based evaluator cannot remove this dependence by relabelling the graph. Independent requests can expose additional concurrency, while cached execution reuses state under the prefix-consistency assumptions of Section~\ref{sec:lengths}. Different request lengths change activity and storage requirements; they are not automatically changes in the model's parameters or architecture.

At a fixed finite-precision protocol, legal rescheduling of a graph with fixed pure primitive functions preserves its result by Theorem~\ref{thm:schedule-independence}. Replacing that graph by a numerically different kernel is a stronger change. A parallel reduction with different parenthesization, a lower-precision representation, or altered overflow behavior can change outputs even when the corresponding real-valued formula is equal. Exact certificates must refer to the actual numerical primitives; approximate replacements require the error framework described in Section~\ref{sec:reencoding}. The parallel-prefix associativity certificate is therefore not a general permission to reassociate neural arithmetic.

\begin{table}[htbp]
\centering\small\renewcommand{\arraystretch}{1.13}
\begin{tabularx}{\linewidth}{@{}>{\raggedright\arraybackslash}p{0.24\linewidth}>{\raggedright\arraybackslash}p{0.32\linewidth}>{\raggedright\arraybackslash}X@{}}
\toprule
LLM component or protocol & Representation in the framework & What the execution account permits \\
\midrule
Heads and known-position operations & Independent subgraphs joined at explicit consumers & Concurrent scheduling subject to actual dependencies and resource limits. \\
Layer stack and generation loop & Composed maps and explicitly supplied next inputs & Waiting on real dependencies; no general elimination of layer or generated-token dependence. \\
Key--value cache & Invocation state tied to parameters, prefix and positions & Reuse by exact simulation; deletion by future reads; other encodings require an output and update invariant. \\
Kernel or representation choice & Alternative descriptions at a fixed interface & Certified exact substitution, or a separately quantified approximation. \\
Request length and worker availability & Interface family, active tasks and resource protocol & Distinct forms of elasticity, with their respective consistency and progress conditions. \\
\bottomrule
\end{tabularx}
\caption{How the concurrent-machine proposal relates to the decoder model. These are formal correspondences and correctness conditions, not measured claims of faster LLM inference.}\label{tab:llm-elasticity}
\end{table}

In the terminology of T5, a layer call, a head computation or a kernel invocation can be assigned a logical computing-core identity at a declared granularity. This is not a one-to-one identification between tokenizer outputs and execution units. A fixed bounded inference graph has finite width; unbounded organization concerns an effective family of workloads, requests or generated descriptions. A trained Transformer does not acquire unlimited context length or eliminate autoregressive dependencies merely by entering that family.

Thus a Transformer is one structured evaluator among the machine's possible implementation forms. It can be invoked as a module, expanded into a graph at a chosen numerical granularity, and combined with exact tables, symbolic components or other represented mappings. Efficiently selecting those boundaries and implementations is part of the research programme, rather than an automatic consequence of Cartesian closedness. A learned model that only approximates a target must carry that approximation obligation throughout any larger composition.

\paragraph{What these categorical structures classify.}
The distinctions give a sequence of tests, rather than architecture
names. Rooted occurrence structure tests which descriptions can be
assembled and compared; prefix naturality tests information access
across lengths; the cache square tests a proposed state realization;
the safe-eviction square tests one deletion protocol. The parameterised
state interface records the types before and after a step. These
properties do not classify trained weights by quality, and different
architectures can pass the same tests. Fixed parameters or declared
parameter quantifiers are part of each claim. A stronger comparison of
retained states is supplied by the future-observation construction in
Section~\ref{obs:future-observation}, which identifies precisely the
distinctions needed for the chosen future outputs.

Taken together, these tests provide comparison coordinates for a model
description and its inference protocol. A causal map at each fixed
length and a coherent family across lengths express different
properties; two exact state realizations may retain different amounts
of information. Specifying which property is under comparison helps
locate a failed simulation or an unavailable transformation. The
classification programme is to organize such properties and their
implications across implementations.

\subsection{Parameters as a family of admissible layers}\label{sec:parameters}

The general parameter construction of Section~\ref{sec:parameter-objects} applies to neural architectures as one representation family. For fixed sequence length and positional conventions, it gives
\[
\Phi_\ell:P_E(\Theta)\otimes\Cau_n(\mathbb R^d)\to\Cau_n(\mathbb R^d),
\qquad \Phi_\ell(\theta,H)=F_{\ell,\theta}(H).
\]
For causal decoder blocks, take $E=E_n^{\mathrm{pre}}$ and $T_A=T_B=\Cau_n(\mathbb R^d)$. Proposition~\ref{prop:decoder-causal} then supplies the joint map for every parameter value in the specified architecture. One can also use the smaller parameter heap $\{(c_k,\theta,\theta):k\le n,\theta\in\Theta\}$ for these particular classes; it expresses that the parameter is held fixed when comparing two prefixes. Definition~\ref{def:parameter} has the advantage of working for arbitrary Token classes.

Different parameter choices may have the same denotation. The map into the exponential need not be injective, just as different finite descriptions may compute the same mapping in Section~\ref{sec:representation-system}. Stacking layers is the internal composition Token map $\Cau_n(Z)^{\Cau_n(Y)}\otimes\Cau_n(Y)^{\Cau_n(X)}\to\Cau_n(Z)^{\Cau_n(X)}$ of the Cartesian closed category, and applying a layer to data is evaluation. In this way the \emph{architecture} of a network is expressed by Token maps built from products, exponentials, evaluation and composition, while the \emph{parameters} are points of bases; this is the perspective of \cite{fong2019,cruttwell2022,gavranovic2024} enriched by the Token-level structure carried by heaps. Admissibility of the forward family does not itself provide a training algorithm.

\subsection{Learning Dynamics, Training Protocols, and Causal Structure}
\label{sec:learning-dynamics}

The parameter family above describes what a model can compute at each
parameter value. Learning adds a trajectory through that family and a
protocol for producing the trajectory. Token Space can describe all three:
the maps evaluated by the model, the state transitions that change its
parameters, and the occurrence graphs that implement those transitions.
This extends the connection between represented programs and data: a
parameter that is fixed during one invocation becomes data consumed and
produced by another. The useful question is which structure is being
represented or preserved at each level.

\paragraph{Language modelling supplies a structured learning problem.}
Consider autoregressive training on a tokenized sequence
$x=(x_1,\ldots,x_n)$. With a fixed start convention, a basic objective is
\begin{equation}\label{eq:training-loss}
 L(\theta;x,w)=-\sum_{i=1}^{n}w_i
       \log p_\theta(x_i\mid x_{<i}),\qquad w_i\geq0.
\end{equation}
Here the model assigns positive probabilities, and the weights specify
which predictions are scored and their normalization. Tokenization,
document boundaries, attention masks and selected target positions are
parts of this learning specification. A text occurrence can be a target
for one prediction and an input to later predictions. Thus a corpus
supplies a family of related tasks through its sequence structure; the
learning problem is richer than the choice of an optimizer.

The same parameters are used across positions and examples. When the
objective weights $w_i$ are fixed and the summands are differentiable, their gradient
contributions add to a shared parameter gradient. One update can therefore
change predictions in many contexts. In teacher-forced training the input
prefixes are supplied from the training sequence. Many position-wise
computations can proceed in parallel subject to the layer dependencies
and the causal mask \cite{vaswani2017}. Autoregressive generation instead
waits for each newly selected token before the next prefix is available.
These are different execution organizations of related model evaluations.

\paragraph{A training protocol has a state interface.}
For one stage with fixed interfaces, write a sufficiently complete state as
$s_k=(\theta_k,\mu_k,k,\pi_k)\in S$. The component $\mu_k$ contains optimizer
memory, such as momentum or moment estimates; $\pi_k$ contains any further
protocol state needed to determine the next transition. A specified total
update has the form
\begin{equation}\label{eq:training-update}
 U:S\times\mathcal B\times\Xi\longrightarrow S,
 \qquad s_{k+1}=U(s_k,b_k,\xi_k).
\end{equation}
The input $b_k$ carries a batch and any external feedback used by the
update. Random choices are supplied explicitly by $\xi_k$, or generated
from a random-generator state included in $S$. The data-selection rule,
learning-rate schedule and numerical conventions are part of the
specified protocol. If a stage change alters the state interface, the
parameterised state construction of
Section~\ref{sec:parameterised-cache-state} supplies the appropriate
organization. Partial updates additionally require an explicit error or
partiality semantics.

This interface accommodates LLM training beyond a single objective.
Pretraining, supervised instruction tuning, preference optimization and
reinforcement learning can use different data, scoring rules and feedback
loops. For example, InstructGPT combines demonstrations with subsequent
preference-based reinforcement learning \cite{ouyang2022}; Llama~3 uses
staged pretraining and iterative post-training involving supervised
fine-tuning, rejection sampling and direct preference optimization
\cite{llama2024}. Generated candidates, evaluations and selected training
examples can be recorded as distinct occurrences with their provenance.
An adaptive data-collection process becomes another stateful component,
with human or environmental responses exposed as inputs. The interface
then shows how a model helps produce the evidence used to update itself.

\paragraph{From a state update to Token records.}
At a coarse level, a record
$(\mathrm{trainstep},o_{k+1},o_k,o_b,o_\xi)$ designates its first occurrence
as the result and the remaining occurrences as operands. Its numerical
interpretation is equation~\eqref{eq:training-update}. Expanding the
operation exposes a finer occurrence graph. For a gradient-based instance
of $U$, Figure~\ref{fig:training-update} shows one such decomposition.
Here $W_0=S\times\mathcal B\times\Xi$, and the later $W_j$ can be taken
to be the successive images of the specified stage maps. The state $W_1$
also carries predictions and the forward information needed by the
derivative computation; $W_2$ adds the loss and its derivative seed;
$W_3$ adds the accumulated parameter gradient. Each stage retains the
inputs needed by later stages. The seed is computed from predictions
and targets by the specified loss rule, not inferred from the scalar
loss value alone.

\begin{figure}[htbp]
\centering
\begin{tikzcd}[column sep=large]
 W_0 \arrow[r,"\mathrm{forward}"]
     \arrow[rrrr,bend left=28,"U"]
 & W_1 \arrow[r,"\substack{\mathrm{loss}\\\mathrm{and\ seed}}"']
 & W_2 \arrow[r,"\mathrm{backward}"']
 & W_3 \arrow[r,"\mathrm{update}"']
 & S
\end{tikzcd}
\caption{A training update as a composition of functions with explicit
intermediate state. The outer arrow is the composite. Operation records
expose the dependencies inside these stages; structural claims about
the functions additionally specify heaps that they preserve.}
\label{fig:training-update}
\end{figure}

These arrows are first functions in $\Set$. A heap on their numerical
interfaces determines which arrows are also Token maps, just as in the
forward-layer construction. Alternatively, the operation-graph encoding
records a chosen update algorithm and its execution. This distinction
separates representing an update from asserting that it preserves a
previously selected relation. Core symbols name operation roles;
trainable values belong to the carrier or to versioned value occurrences.
An initial parameter occurrence may be a designated input of one step,
while the resulting occurrence is an input to the next. Designation as
an input does not make its numerical value immutable across steps.

To interpret the backward stage mathematically, one additionally supplies
derivative rules and a loss; a gradient also uses a chosen identification
of covectors with vectors, conventionally the Euclidean inner product.
Nondifferentiable primitives require specified treatment. For a numerical
implementation, automatic differentiation executes declared rules under
chosen arithmetic; it need not compute a derivative of the rounded
floating-point map. Cartesian closedness organizes evaluation and
composition, while this differential structure explains the update.
The accounts of backpropagation and gradient-based learning in
\cite{fong2019,cruttwell2022} provide complementary categorical machinery.

\paragraph{Three causal questions.}
The word ``causal'' refers to different orders in this discussion.
Table~\ref{tab:training-causality} separates them. During a training
forward pass, the target $x_i$ may enter the loss although the prediction
of $x_i$ uses only $x_{<i}$. Backpropagation carries derivative information
against forward dataflow, but each backward operation waits for its own
operands. Its gradient occurrences are new nodes. Unrolling a finite
training run with distinct $\theta_k$ and $\theta_{k+1}$ can therefore give
an acyclic computation presentation even though a schematic learning
loop is drawn with a feedback arrow.

\begin{table}[htbp]
\centering\small\renewcommand{\arraystretch}{1.12}
\begin{tabularx}{\linewidth}{@{}>{\raggedright\arraybackslash}p{0.22\linewidth}>{\raggedright\arraybackslash}X>{\raggedright\arraybackslash}X@{}}
\toprule
Order & What is held fixed & Causal requirement \\
\midrule
Sequence positions $i$
& Parameters and the forward protocol
& An output prefix depends only on the corresponding input prefix. \\
Learning steps $k$
& Initial training state and the update protocol
& State after $k$ steps depends only on the first $k$ supplied training events. \\
Operation dependencies
& A specified implementation and its inputs
& Each operation uses values already produced or supplied at its interface. \\
\bottomrule
\end{tabularx}
\caption{Distinct causal structures for a language model and its training.
A training event can contain an entire sequence or batch.}
\label{tab:training-causality}
\end{table}

For equation~\eqref{eq:training-update}, two event streams with the same
first $k$ pairs $(b_j,\xi_j)$ produce the same states through step $k$
from the same initial state, by induction. This is causality in training
time. It permits each event to contain a complete sequence. It does not
assert prefix causality with respect to the positions inside that
sequence, as the following small example makes explicit.

\begin{example}[Updating on a sequence and then predicting it]
\label{ex:training-prefix-boundary}
For every fixed real $\theta$, the family
$F_{\theta,n}(x)_i=\theta x_i$ on $\mathbb R^n$ is prefix-causal and
coherent across lengths. Consider the illustrative squared loss
\[
 L_n(\theta;x)=\tfrac12\sum_{j=1}^n(\theta x_j-1)^2.
\]
From $\theta_0=0$, one gradient step of size $\eta>0$ gives
\[
 \theta'_n(x)=\theta_0-\eta\,\partial_\theta
     L_n(\theta_0;x)=\eta\sum_{j=1}^n x_j.
\]
Now update using the complete sequence and then evaluate
$H_n(x)=F_{\theta'_n(x),n}(x)$. The inputs $(1,0)$ and $(1,1)$ have the
same first entry, but the first outputs of $H_2$ are $\eta$ and $2\eta$.
Thus $H_2$ is not prefix-causal. Each frozen-parameter map remains
prefix-causal. The example concerns the composite of learning and
subsequent prediction, not information access in the original forward
pass; the squared loss is chosen only to make that distinction explicit.
\end{example}

\paragraph{What the execution framework contributes to training.}
Once a finite training computation is presented by pure total operations
with explicit state and randomness, the typed extension in
Remark~\ref{rem:typed-frontier} makes the evaluation and frontier
constructions available. Independent forward or backward subcomputations
expose concurrency; shared parameter-gradient accumulation exposes their
joins. Schedule independence concerns legal schedules of that fixed
graph. Changing a floating-point reduction tree changes the numerical
operation and requires its own comparison. Microbatch organization,
communication and gradient accumulation must likewise be represented
when they affect the update.

The retained-state question also becomes precise. For a training graph,
future reads include backward operations and requested diagnostics.
Consequently an activation that is no longer needed by forward
evaluation may still belong to the live frontier. Checkpointing replaces
some retention by a recomputation graph, whose correctness includes
reproducing the needed forward state and random choices. The inference
KV eviction formula is therefore specific to its stated access protocol;
training uses the frontier of the complete training implementation.
Likewise, changing parameters generally invalidates a cache justified
for the previous parameters, unless it is recomputed or a separate
migration invariant is proved.

A learning step deliberately moves from $F_{\theta_k}$ to
$F_{\theta_{k+1}}$ and may change its predictions. The relevant execution
comparison is instead between implementations of the same declared
training update $U$. With a specified transition and observation
interface, the future-observation construction of
Section~\ref{obs:future-observation} can also ask which retained training
states support the same future observations under the admitted event
streams. This connects learning to the paper's questions about
representation, concurrency and exact retained information. Convergence,
generalization and empirical model quality remain questions about the
chosen learning problem and protocol, rather than consequences of its
record representation.

\subsection{Occurrence-level records and attention snapshots}\label{sec:attention-records}

Sequence classes describe transformations of all possible sequence states. A complementary description records one execution. For a fixed input $s$, layer $\ell$, and head $a$, let $O_s=\{o_1,\ldots,o_n\}$ be distinct occurrence identifiers. The class
\begin{equation}\label{eq:attention-heap}
A_{s,\ell,a}=\left\langle O_s,E_w,
\{(\mathrm{weight},o_i,o_j,\overline{\alpha_{ij}^{\ell,a}(s)}):1\le i,j\le n\}\right\rangle
\end{equation}
uses tagged core labels for weights in $[0,1]$, in the manner of $\Rep_{\mathbf{Vct}}$ (Remark~\ref{rem:signature}). A map to a second such class preserves these weights exactly: it is a comparison of weighted relations, not automatically the neural layer that computed the matrix. The operators $\Cup,\Cap,\ominus$ and the renaming $\oslash$ compare the attention patterns of different heads or layers; if multiple inputs, layers, or heads are collected, keep those indices in occurrence identifiers or in separate tagged record fields.

One can add $(\mathrm{next},o_i,o_{i+1})$, $(\mathrm{pos}_i,o_i)$, or $(\mathrm{type},o_i,\bar v_i)$ to record adjacency, a fixed position, or a vocabulary identity. Position labels force position correspondence under a Token map, but do not forbid a layer from using other positions' values. That dependency is instead controlled by the causal sequence heap. Adding more kinds of records strengthens the morphism condition and may exclude a computation that was admissible for a smaller heap.

Row sums equal to one and the presence of a unique weight for each ordered pair are properties imposed on the class of attention snapshots. They are not satisfied by arbitrary heaps over $E_w$. Likewise, exact equality of real-valued labels does not express numerical closeness. Comparing attention patterns approximately requires a discrepancy criterion or an enriched framework; heap intersection with exact real labels is generally too rigid for that purpose.

\section{Distillation as Structure-Preserving Compression}\label{distillation}

\subsection{Teacher-induced heaps and faithful students}

Knowledge distillation trains a small \emph{student} network to imitate a large \emph{teacher} \cite{hinton2015,gou2021}. Different distillation methods differ mainly in \emph{what} of the teacher has to be imitated: its output probabilities, the ranking of its logits, the similarities it assigns to pairs of samples, or its attention patterns \cite{park2019,tian2020,jiao2020}. In the Token Space all of these choices are choices of a heap, and ``imitating'' becomes ``being a Token map''. The framework describes the exact relational specification associated with such a choice; it does not identify a small empirical loss with exact preservation.

Let $f_T:X\to Y_T$ and $f_S:X\to Y_S$ be teacher and student functions. Equip $Y_T$ and $Y_S$ with classes $T_T$ and $T_S$ over the same core $E$. Their carriers may have different dimensions; the core specifies which relations are being compared. The teacher induces the \emph{teacher-induced class}
\[
f_T^*T_T=\langle X,E,F_E(f_T)^{-1}(\heap(T_T))\rangle
\]
on the common input domain (the pullback \eqref{eq:pullback}); its Tokens are exactly the structural facts about inputs that the teacher asserts through its outputs, and by \eqref{eq:heap-adjunction} it is the largest class on $X$ through which $f_T$ is a Token map.

\begin{definition}[Structural distillation]\label{def:faithful}
Let $K\subseteq f_T^*\heap(T_T)$ be a \emph{knowledge heap}, the part of the teacher's structural facts that is to be transferred. The student is \emph{structurally faithful} to the teacher on $K$ if
\begin{equation}\label{eq:distill}
K\subseteq f_S^*\heap(T_S),
\end{equation}
i.e.\ iff $f_S:\langle X,E,K\rangle\to T_S$ is a Token map. For a specified candidate domain $\mathcal D\subseteq F_E(X)$, it is \emph{exactly faithful on $\mathcal D$} when $\mathcal D\cap f_S^*\heap(T_S)=\mathcal D\cap f_T^*\heap(T_T)$. This compares both presence and absence within the chosen domain; selecting a positive heap $K$ alone does not specify negative targets.
\end{definition}

These conditions define positive transfer and exact structural agreement on an explicit test domain; they are not consequences of ordinary distillation training.

For example, logit rankings can be encoded by $(\mathrm{greater},y,\bar i,\bar j)$ when $y_i>y_j$, using fixed labels for output classes; a structurally faithful student then preserves every ranking of the teacher, the ``dark knowledge'' \cite{hinton2015} about which wrong classes are more plausible. Pairwise proximity can be encoded by $(\mathrm{close},y,y')$ when a specified distance is below a threshold, the qualitative core of relational distillation \cite{park2019,tian2020}. For attention distillation \cite{jiao2020}, the output carrier should be a set of attention matrices, not the position carrier of the single-snapshot class \eqref{eq:attention-heap}. A well-typed exact target uses
\[
(\mathrm{weight},A,\bar i,\bar j,\bar a)\in\heap(T_Y)
\quad\Longleftrightarrow\quad A_{ij}=a.
\]
Pulling this heap back along a teacher $f_T:X\to Y$ records its matrix entries for each input. Different numbers of heads, layers, or positions require an explicit alignment or aggregation before this comparison. Full probability matching can similarly use $(\mathrm{prob},p,\bar v,\bar a)$ iff $p(v)=a$, at the cost of demanding exact probabilities.

To distinguish positive preservation from agreement on tested absences, specify a domain $\mathcal D\subseteq F_E(X)$ of candidate records. Define
\[
K_T=\mathcal D\cap f_T^*\heap(T_T),\qquad
K_S=\mathcal D\cap f_S^*\heap(T_S).
\]
Exact agreement on $\mathcal D$ is $K_T=K_S$; this is the exact-faithfulness clause of Definition~\ref{def:faithful}. If $\mathcal D$ is chosen as all records of selected shapes, $K_T$ must include every teacher-positive record of those shapes, not merely a selected subset. Selecting only some positive teacher records, by contrast, does not make every omitted record a negative claim.

\begin{proposition}[Distillation transfers symmetries]\label{prop:distill-equivariance}
Suppose $f_T:X^n\to Y_T^n$ and $f_S:X^n\to Y_S^n$ use permutation classes on their outputs, with the same permutation labels. If the student is structurally faithful on the full teacher-induced heap, then
\[
f_T(t)=\sigma f_T(s)\quad\Longrightarrow\quad
f_S(t)=\sigma f_S(s).
\]
In particular, an equivariant teacher has an equivariant student, and more generally every symmetry of the teacher's outputs is inherited by the student. If $f_S=g\circ f_T$ for a Token map $g:T_T\to T_S$, then $f_S$ is structurally faithful on the whole of $f_T^*\heap(T_T)$.
\end{proposition}
\begin{proof}
The premise is exactly membership of $(\bar\sigma,s,t)$ in the teacher-induced heap. Its preservation gives the conclusion. If the teacher is equivariant, apply this implication with $t=\sigma s$. The last statement is the composite of the Token maps $f_T:f_T^*T_T\to T_T$ and $g$.
\end{proof}

The full induced heap can be a much stronger requirement than the selected training targets. Conversely, symmetry preservation alone can be weak evidence of predictive quality: a constant student whose output is fixed by all permutations preserves every induced permutation record. The relation one elects to preserve must therefore match the intended learning objective. The operators of Section~\ref{sec:operators} describe the common variants of the setting: $\ominus$ removes facts that should not be transferred (selective distillation), the meet $\Cap$ of several teacher-induced classes is the knowledge on which an ensemble of teachers agrees, $\varnothing$ removes every structural requirement; label-only training instead requires retaining explicit label-target records, and the tree Tokens of Section~\ref{interior-structure-mapping-and-tree-Token-classes} together with reification describe the transfer of nested reasoning traces to a student that reads flat sequences.

\subsection{Compression as quotient}\label{sec:quotient-compression}

State aggregation is one possible form of distillation. A smaller implementation need not merge semantic states: Section~\ref{sec:parity} preserves every external input while sharing its implementation nodes. Here we study the distinct case where a specified map identifies elements of a semantic carrier. Let $T=\langle B,E,H\rangle$ describe the teacher's internal representations $B$ and the structural facts $H$ among them. By Proposition~\ref{prop:mono} a surjective state aggregation of $T$ is an epimorphism $q:T\to T'$, a Token map with surjective base, and every epimorphism factors as $T\to T/\ker q\to T'$ with the second map bijective on bases. For an equivalence relation $\eta$ on $B$ with quotient map $q:B\to B/\eta$ define
\[
T/\eta=\langle B/\eta,E,F_E(q)(H)\rangle .
\]
Up to a carrier bijection and adding target Tokens, these state aggregations are the quotients $T/\eta$. The quotient map always preserves Tokens: merging never loses a fact of the teacher, but it may \emph{create} facts. It \emph{reflects} Tokens when $F_E(q)(r)\in\heap(T/\eta)$ implies $r\in H$; reflection rules out tuples becoming indistinguishable from a recorded tuple solely through merging their entries.

\begin{definition}[Congruence]\label{def:congruence}
An equivalence $\eta$ is a \emph{congruence} of $T$ (is \emph{relation-saturating} for $H$) if, whenever $r,r'$ have the same shape and core entries and corresponding carrier entries are $\eta$-related,
\[
r\in H\quad\Longleftrightarrow\quad r'\in H.
\]
This is stronger than the operation-compatible congruence of universal algebra when operations are represented by their graphs.
\end{definition}

\begin{proposition}[Faithful compressions]\label{prop:congruence}
The quotient $q:T\to T/\eta$ reflects Tokens iff $\eta$ is a congruence of $T$. The congruences of $T$ form a complete lattice under inclusion; its greatest element is the \emph{indiscernibility relation} $\approx$, where $b\approx b'$ iff replacing $b$ by $b'$ at any one carrier position of any word leaves membership in $H$ unchanged. Consequently $T/\!\approx$ is the smallest quotient of $T$ that neither loses nor invents a structural fact, and $|B/\!\approx|$ bounds below the target carrier size of any map that preserves and reflects membership of all these records. It does not bound parameter count or description length.
\end{proposition}
\begin{proof}
A word has its image in the quotient heap iff some recorded word has the same image, i.e.\ is componentwise $\eta$-related to it with the same core entries; so reflection is exactly saturation. The replacement condition defining $\approx$ is reflexive, symmetric, and transitive by successive replacements. It is a congruence because a word has finitely many positions, so componentwise replacement can be performed one position at a time; and every congruence permits each single replacement, hence is contained in $\approx$. Saturation is preserved under nonempty intersections and under the join (transitive closure of the union) of congruences, because a componentwise chain can be traversed one position and one step at a time; the empty meet is $\approx$ and the empty join is equality, giving a complete lattice. For any preserving-and-reflecting map $e:B\to C$, equal $e$-images can be substituted without changing membership, so $\ker e$ is a congruence contained in $\approx$. Consequently $B/\ker e\cong e(B)$ surjects onto $B/\!\approx$, proving the carrier bound, and the quotient by $\approx$ attains it.
\end{proof}

The relation $\approx$ expresses indiscernibility under single-position replacement in the chosen atomic relations. This invites comparison with Leibniz congruences \cite{burris1981,blok1989} and behavioral equivalences \cite{sangiorgi2011}. For an arbitrary heap no such identification follows: bisimulation has different matching conditions, and ordinary algebraic congruences need not saturate operation graphs. Section~\ref{obs:future-observation} constructs a particular observation heap for which indiscernibility does coincide with future behavioral equivalence.

\begin{example}[Scope of the bound]\label{ex:no-merge}
If $H$ records equality through $(\mathrm{eq},b,b)$ for every $b$, replacing one occurrence of $b$ by $b'\ne b$ produces an unrecorded tuple, so $\approx$ is equality. Both the permutation sequence heap (through the identity permutation) and the causal sequence heap (through $c_n$) contain such a diagonal relation; their exact preserving-and-reflecting quotients therefore cannot merge distinct complete sequences. The lower bound of Proposition~\ref{prop:congruence} is informative for heaps of relational facts---rankings, proximities, attention graphs---and trivial for classes that contain a diagonal.
\end{example}

This example separates structural quotients from the representation compression of Section~\ref{sec:finite-representations}. A smaller circuit or network can compute the same function on the same state domain without being a quotient of the teacher's set of representations; cardinality of a carrier is not parameter count, storage cost, or arithmetic cost. Heaps are sets and Token maps are all-or-nothing, so the Token Space specifies \emph{which} structure a distillation method transfers and proves what is preserved, but does not measure how closely a student approximates its teacher; that requires heaps weighted by real numbers, or Token classes enriched over metric spaces, which we leave to future work (Section~\ref{conclusion}).

\subsection{Future observation heaps, behavioral quotients, and exact compression}
\label{obs:future-observation}

The quotient in Proposition~\ref{prop:congruence} is determined by the
\emph{chosen heap}. The cache simulation supplies a different datum: a
state suffices to reproduce every future output. These notions agree for
a particular observational heap, and only for the observations it records.
The following construction connects the structural quotient, deterministic
state minimization, and the migration criterion of
Theorem~\ref{v13:mig:kernel}. The Mealy-machine and behavioral-quotient
ingredients are standard coalgebraic semantics
\cite[Sections~3.8 and~8--9]{rutten2000}; the purpose here is to identify
their exact realization by Token records and their additional update
obligation.

Fix an input set $A$, an output set $Y$, and a deterministic, total history
response
\[
 \lambda:A^*\times A\longrightarrow Y.
\]
Write $\mathsf H=A^*$, with empty history $\epsilon$, and append an input
$a$ as $ha$. The history machine is
$r(h,a)=(\lambda(h,a),ha)$. Define its \emph{future observation}
recursively by
\begin{equation}\label{obs:future-recursion}
 \mathsf O(h,\epsilon)=\epsilon,
 \qquad
 \mathsf O(h,au)=\lambda(h,a)\cdot\mathsf O(ha,u),
\end{equation}
where the dot on the right prepends one output to an output word.
Thus $\mathsf O(h,u)\in Y^{|u|}$ contains only outputs produced
\emph{after} history $h$, during the continuation $u$. It does not include
outputs already delivered. For a deterministic decoder, $Y$ can consist
of logits or probability vectors; sampling is not included unless its
additional state and inputs are specified. A compression boundary is
assumed to have delivered the output of its last step. If a current
logit vector or another result is still owed to a caller or sampler,
the interface must retain or reproduce that pending observation;
omitting it from the chosen behavior would change the correctness
question.

Use a fixed, disjointly tagged core
\[
 E_{\mathrm{obs}}=\{\mathrm{obs}\}
       \sqcup\{\bar u:u\in A^*\}
       \sqcup\{\bar v:v\in Y^*\}.
\]
The input and output tags are distinct even when their underlying words
coincide. The \emph{future observation class} is
\begin{equation}\label{obs:observation-heap}
 T_{\mathrm{obs}}=
 \langle\mathsf H,E_{\mathrm{obs}},H_{\mathrm{obs}}\rangle,
 \qquad
 H_{\mathrm{obs}}=
 \{(\mathrm{obs},\bar u,h,\overline{\mathsf O(h,u)}):
                                      h,u\in A^*\}.
\end{equation}
Every record has exactly one carrier position. Continuations and output
words are fixed observations in the core, not carrier elements available
for identification. The heap is generally infinite; its definition is a
semantic specification, not an assertion that it can be materialized or
decided effectively.

Define \emph{future behavioral equivalence} by
\begin{equation}\label{obs:behavior-equivalence}
 h\sim h'
 \quad\Longleftrightarrow\quad
 \mathsf O(h,u)=\mathsf O(h',u)\quad\hbox{for every }u\in A^*.
\end{equation}
This is the residual, or Nerode-style, equivalence appropriate to the
specified deterministic input/output behavior. A right congruence on
$A^*$ here means an equivalence $\eta$ satisfying
$h\mathrel\eta h'\Rightarrow ha\mathrel\eta h'a$ for every $a\in A$;
it need not be a two-sided monoid congruence.

\begin{theorem}[Observation, behavior, and exact incremental compression]
\label{obs:observation-compression}
For the total history machine above, the following statements hold.
\begin{enumerate}[label=\textup{(\roman*)},leftmargin=2em]
\item The indiscernibility relation $\approx$ of
Proposition~\ref{prop:congruence}, computed on $T_{\mathrm{obs}}$, is
exactly $\sim$. Consequently, for any equivalence $\eta$ on $\mathsf H$,
the quotient $T_{\mathrm{obs}}\to T_{\mathrm{obs}}/\eta$ reflects Tokens
if and only if $\eta\subseteq\sim$.

\item The relation $\sim$ is a right congruence. With
$Q=\mathsf H/{\sim}$ and $q(h)=[h]$, the formula
\begin{equation}\label{obs:behavior-transition}
 \tau_Q([h],a)=(\lambda(h,a),[ha])
\end{equation}
defines a total Mealy transition $Q\times A\to Y\times Q$ and makes $q$
an output-preserving state-machine map. Starting at $[\epsilon]$, every
state of $Q$ is reachable.

\item Let $c:\mathsf H\to C$ be any specified encoding, and put
$C_{\mathrm r}=\im(c)$. There is a transition
$\tau_C:C_{\mathrm r}\times A\to Y\times C_{\mathrm r}$ satisfying
\begin{equation}\label{obs:compression-invariant}
 \tau_C(c(h),a)=(\lambda(h,a),c(ha))
 \qquad(h\in\mathsf H,\ a\in A)
\end{equation}
if and only if both conditions hold:
\begin{align}
 c(h)=c(h')&\ \Longrightarrow\
           \lambda(h,a)=\lambda(h',a)
                 &&\text{for every }a\in A,
                 \label{obs:output-fiber}\\
 c(h)=c(h')&\ \Longrightarrow\ c(ha)=c(h'a)
                 &&\text{for every }a\in A.
                 \label{obs:update-fiber}
\end{align}
Equivalently, $\ker(c)$ is a right congruence contained in $\sim$.
When it exists, $\tau_C$ is unique on the reachable states, and
$C_{\mathrm r}$ equipped with the image heap
$F_{E_{\mathrm{obs}}}(c)(H_{\mathrm{obs}})$ is an exact
preserving-and-reflecting quotient of the observation class.

\item For every encoding satisfying \textup{(iii)}, there is a unique
surjective Mealy homomorphism
\begin{equation}\label{obs:minimal-factor}
 p:C_{\mathrm r}\longrightarrow Q,
 \qquad p(c(h))=[h],\qquad pc=q.
\end{equation}
It preserves and reflects the future observation heaps. More generally,
the pointed machine $(Q,[\epsilon],\tau_Q)$ is terminal in the category
of reachable, pointed deterministic Mealy machines with the same
input/output behavior as the history machine at $\epsilon$, whose
arrows preserve the initial state, transitions and outputs.
In particular, $Q$ is a least-state exact realization, unique up to
isomorphism among reachable realizations whose distinct states have
distinct future behaviors.
\end{enumerate}
\end{theorem}

\begin{proof}
For \textup{(i)}, replacing $h$ by $h'$ in a record with observation
labels $\bar u,\bar v$ preserves its membership exactly when
\[
 v=\mathsf O(h,u)\quad\Longleftrightarrow\quad
 v=\mathsf O(h',u).
\]
For every $u$, exactly one output word is recorded for each history.
Thus replacement preserves membership for all $u,v$ if and only if the
two future observation functions agree. All other record shapes are
absent and remain absent after a carrier substitution. This proves
$\approx\,=\,\sim$. Since the only recorded shapes have one carrier
position, every equivalence contained in $\sim$ saturates this heap;
the converse follows by the same test. Proposition~\ref{prop:congruence}
then proves the quotient assertion.

For \textup{(ii)}, suppose $h\sim h'$. Applying the definition to a
one-letter continuation $a$ yields
$\lambda(h,a)=\lambda(h',a)$. Applying it to $au$ and cancelling the
equal first output in~\eqref{obs:future-recursion} yields
$\mathsf O(ha,u)=\mathsf O(h'a,u)$ for every $u$, hence $ha\sim h'a$.
Both coordinates of~\eqref{obs:behavior-transition} are therefore
independent of the chosen representative. Its defining formula is the
homomorphism equation for $q$. Reading $h$ from $[\epsilon]$ reaches
$[h]$, by induction on $|h|$.

For \textup{(iii)}, equation~\eqref{obs:compression-invariant}
immediately implies both fiber conditions, because equal arguments
of $\tau_C$ have equal output pairs. Conversely, the two conditions
make that equation a well-defined definition of $\tau_C$ on every
element of $C_{\mathrm r}\times A$. Its second coordinate belongs to
$C_{\mathrm r}$, and surjectivity of $c$ onto this image gives
uniqueness. Induction on continuation length, using the first condition
at the first step and the second condition at the next state, shows
that $c(h)=c(h')$ implies $h\sim h'$. The second condition is precisely
right stability of $\ker(c)$. Conversely, right stability and
$\ker(c)\subseteq\sim$ imply the second condition and, by one-letter
observations, the first. Part~\textup{(i)} gives preservation and
reflection for the image heap.

For \textup{(iv)}, part~\textup{(iii)} gives
$\ker(c)\subseteq\sim=\ker(q)$. Hence $p(c(h))=[h]$ is well defined,
unique, and surjective. For every $h,a$,
\[
 (\id_Y\times p)\tau_C(c(h),a)
     =(\lambda(h,a),[ha])
     =\tau_Q(p(c(h)),a),
\]
so $p$ is a homomorphism. Each reachable state $c(h)$ has future
observations $\mathsf O(h,-)$, proving preservation and reflection
of its observation records.

For the categorical assertion, take any reachable pointed machine
$(D,d_0,\tau_D)$ with the specified behavior from $d_0$. Let $d(h)$ be
the state reached after reading $h$. Reachability makes
$d:\mathsf H\to D$ surjective. Equality of the output traces on all
words $hu$, followed by removal of their already emitted $|h|$ outputs,
shows that the observations from $d(h)$ are exactly
$\mathsf O(h,u)$. In particular its next output is $\lambda(h,a)$,
and its next state is $d(ha)$. Thus $d$ satisfies \textup{(iii)}, and
the preceding construction yields a homomorphism $D\to Q$. Any
initial-state-preserving homomorphism must send $d(h)$ to $[h]$,
by induction on $h$; it is therefore unique. If $D$ has no two distinct
states with the same future behavior, this map is also injective,
and hence an isomorphism. The surjection from every reachable exact
realization proves the stated least-state property.
\end{proof}

\begin{figure}[htbp]
\centering
\begin{tikzcd}[column sep=large,row sep=large]
 \mathsf H\times A \arrow[r,"r"]
       \arrow[d,"c\times\id_A"']
   & Y\times\mathsf H \arrow[d,"\id_Y\times c"]\\
 C_{\mathrm r}\times A \arrow[r,"\tau_C"]
       \arrow[d,"p\times\id_A"']
   & Y\times C_{\mathrm r} \arrow[d,"\id_Y\times p"]\\
 Q\times A \arrow[r,"\tau_Q"']
   & Y\times Q
\end{tikzcd}
\caption{Exact incremental compression in $\Set$. The upper square
requires both preservation of the next output and preservation of the
chosen next-state encoding. The lower square maps every exact reachable
realization onto the behavioral quotient. Currying in $A$ gives
homomorphisms of coalgebras for $(Y\times -)^A$. The arrows also preserve
and reflect the selected future observation records.}
\label{obs:observation-square}
\end{figure}

For an exact encoding in part~\textup{(iii)}, its factor
$p:C_{\mathrm r}\to Q$ is an isomorphism exactly when
$\ker(c)={\sim}$. Indeed $p$ is already surjective, and it is injective
exactly when $h\sim h'$ implies $c(h)=c(h')$. This is the criterion
for the encoding to retain exactly the selected behavioral distinctions,
as opposed to a sufficient but redundant state. It does not identify
an inexpensive representation of that state.

\begin{corollary}[Migration between exact realizations]
\label{obs:behavior-migration}
Suppose $c_i:\mathsf H\to C_i$ for $i=1,2$ satisfy
Theorem~\ref{obs:observation-compression}\textup{(iii)}, with reachable
parts $C_{i,\mathrm r}$. A map $m:C_{1,\mathrm r}\to C_{2,\mathrm r}$
satisfying $mc_1=c_2$ exists exactly when
\[
 \ker(c_1)\subseteq\ker(c_2).
\]
It is then unique, surjective, and automatically a Mealy homomorphism
preserving and reflecting future observation records. Writing
$p_i:C_{i,\mathrm r}\to Q$ for the maps to the behavioral quotient,
one has $p_2m=p_1$.
\end{corollary}
\begin{proof}
The set-theoretic argument of Theorem~\ref{v13:mig:kernel} applies to
the common history domain, whether finite or infinite, and supplies
the unique map $m(c_1(h))=c_2(h)$. The two compression invariants give
\[
 (\id_Y\times m)\tau_1(c_1(h),a)
     =(\lambda(h,a),c_2(ha))
     =\tau_2(mc_1(h),a).
\]
Every reachable state is represented by a history, so this is the
homomorphism equation on the whole stated domain. Future observations
and the identity $p_2m=p_1$ follow by evaluation on $c_1(h)$.
\end{proof}

This corollary separates two exact operations. Sending any correct
cache to $Q$ forgets distinctions irrelevant to future outputs.
Migrating to a particular other cache must also reproduce that cache's
chosen internal distinctions; it can fail even though both caches
realize the same input/output behavior. Their common map to $Q$ alone
does not provide a migration between them. Conversely, the existence
of $m$ is an information criterion; as in
Section~\ref{v13:mig:criteria}, it supplies neither an inexpensive
algorithm for $m$ nor a compact representation of $Q$.

\begin{example}[Observational sufficiency without an update invariant]
\label{obs:insufficient-kernel}
Let $A=\{a\}$, $Y=\{0\}$, and let every step output $0$. All histories
are behaviorally equivalent. Define
\[
 c(a^n)=
 \begin{cases}
  0,&n=0,1,\\
  1,&n\ge2.
 \end{cases}
\]
Then $\ker(c)\subseteq\sim$, so the observation heap descends exactly
and every future output can be obtained from $c(h)$. But
$c(\epsilon)=c(a)=0$, while $c(a)=0\ne1=c(aa)$. A deterministic update
of state $0$ on input $a$ cannot simultaneously produce both required
next states. Thus no transition maintains this \emph{specified}
encoding. A different one-state machine correctly emits only zeros;
the failure concerns the proposed representation invariant, not the
existence of any compressed implementation.
\end{example}

\begin{example}[A finite semantic state for unbounded histories]
\label{obs:parity-history}
For $A=Y=\{0,1\}$, let
$\lambda(h,a)$ be the parity of all bits in $ha$.
Two histories have equal future observations exactly when their
parities agree: the same parity gives the same output sequence under
every continuation, and the one-letter continuation $0$ distinguishes
different parities. Thus $Q=\{0,1\}$ with transition
$(s,a)\mapsto(s\mathbin\oplus a,s\mathbin\oplus a)$ is the exact
behavioral quotient of infinitely many histories. This is an example
of semantic summarization. It differs from retaining selected entries
of a historical cache and illustrates why a liveness argument for
particular entries is not a general state-minimality theorem.
\end{example}

\paragraph{Heap dependence and computational scope.}
The identity $\approx\,=\,\sim$ is proved for
\eqref{obs:observation-heap}; it does not change
Example~\ref{ex:no-merge}. If exact equality records
$(\mathrm{eq},h,h)$ are added, the greatest relation-saturating
equivalence becomes equality, even for the constant-output machine
of Example~\ref{obs:insufficient-kernel}. Likewise, saturating a graph
of a state-update operation is stronger than descending that operation
to a behavioral quotient. Consequently, exact compression of future
behavior need not preserve and reflect all occurrence identities,
intermediate activations, or records of the original computation.
Selecting $H_{\mathrm{obs}}$ states explicitly which knowledge is to
survive distillation or migration.

If the behavioral quotient has finite size $k$, every exact reachable
realization has at least $k$ distinct states. Under an additional
fixed-length binary encoding model, this entails at least
$\lceil\log_2 k\rceil$ state bits. Neither assertion alone gives a
lower bound on parameter count, arithmetic operations, KV entries,
or real-valued coordinates. Arbitrary set-theoretic encodings impose
no continuity, computability, finite-precision, or decoding-cost
condition. In particular, no impossibility of constant-dimensional
softmax state follows from this theorem. Such a result would require
a specified representation and computational model.

The full heap also quantifies over all finite continuations. Agreement
on a finite test sample establishes only the selected observations;
it need not establish $\sim$ or right stability. The theorem gives
a semantic universal property, not a general effective minimization
algorithm for neural states. Finite-state instances can be checked by
the standard refinement of behavioral partitions, whereas an arbitrary
numerical decoder requires additional effective hypotheses.

\paragraph{Finite horizons and admissible continuations.}
The theorem assumes total transitions on all of $A^*$. A model with
maximum length $N$ must retain this domain distinction. One precise
formulation uses histories $A^n$ and, for each $n\le N$, compares all
continuations in $A^{\le N-n}$. Define $h\sim_n h'$ only for
$h,h'\in A^n$ by equality of these future observations. The same
one-step argument gives
\[
 h\sim_n h'\ \Longrightarrow\ ha\sim_{n+1}h'a
 \qquad(n<N),
\]
and yields indexed quotients $Q_n=A^n/{\sim_n}$ and transitions
$Q_n\times A\to Y\times Q_{n+1}$. There is no transition required
beyond $N$. If no undelivered terminal observation remains, $Q_N$ is a
singleton when $A^N$ is nonempty, and is empty otherwise. An alternative
is to totalize the machine with
an explicit observable termination or rejection protocol. Comparing
unequal-length histories only on the intersection of their admissible
continuations, while omitting availability or termination, does not
justify the total-machine theorem. Other restricted input protocols
likewise require the protocol state or an explicit observation of
admissibility. These choices connect the behavioral quotient to the
length-indexed state interfaces rather than silently discarding their
typing information.

\subsection{Structural records and interpretation}\label{sec:interpretation}

Three claims should be kept apart: that a relation is included in a formal description; that a measurement or intervention supports that relation; and that the relation determines a semantic content. The first is specified by a heap, the second requires an empirical account, and the third requires further argument. Successful readout, reconstruction, or causal intervention do not by themselves establish one uniquely privileged natural-language content for an individual hidden-state occurrence, and the framework equally does not prohibit well-supported local content attributions \cite{cao2022}.

For interpretability, take as carrier a collection of located occurrences $\omega=(\theta,\ell,s,i,h_i^{(\ell)}(s))$ rather than bare vectors when layer and context affect the claim. A heap can then record a probe readout $(\mathrm{readout},\bar P,\omega,\bar c)$, meaning that the specified procedure $P$ returns $c$ for $\omega$, a reconstruction result, or an intervention outcome whose site, replacement rule, baseline and observable are specified in the record or in the protocol named by a label. Different query classes and protocols induce different heaps on the same carrier, and the operators of Section~\ref{sec:operators} compare them explicitly: the meet $\Cap$ of the heaps supported by several protocols is the structure they agree on. Core labels have fixed identity in the category; the signature can name a measurement protocol or a formal relation while remaining neutral about a theory of representation. We therefore keep the term \emph{structural Token} and do not define it as a ``Token without a referent.''

\section{Explicit Knowledge Structures and Hidden-State Computation}\label{sec:explicit-hidden}

The record language of Token Space applies both to explicit knowledge structures and to internal states obtained through learning. A knowledge graph presents named entities and relations; an ontology supplies axioms and a semantics for their consequences; a language model computes through parameterized transformations of contextual states. We want these structures and computations to enter a common framework without identifying their sources of structure or their correctness conditions. In particular, the interior made available as data in Thesis~\ref{thesis:interior} need not be restricted to relations that have already received familiar concept names. This section develops that connection and states what would make it informative beyond the possibility of encoding each system.

\subsection{Explicit records and ontology computation}\label{sec:explicit-records}

Fix a vocabulary of relation labels in a core $E$. Let $B$ contain the entities and other nodes allowed to vary under the chosen maps. A fact $R(a,b)$ has the record
\[
(R,a,b)\in\Upsilon\subseteq(B\sqcup E)^*.
\]
Unary records can specify sorts, and relations of other finite arities use the same construction. A knowledge graph is thus represented first as a \emph{Token class}. A map of its carrier is a Token map precisely when it preserves the chosen records with the core labels fixed. If named nodes must also be fixed, that requirement must be expressed in the encoding, for example by distinguished constants as in algebraic tokenization. This is the fixed-signature relational case of Theorem~\ref{thm:tokenization}, rather than an additional representation theorem for every feature of a knowledge representation language.

An ontology adds axioms and a consequence relation. The syntax of its axioms can be recorded, including nested expressions through tree Tokens, but recording an axiom does not give it a logical interpretation automatically. OWL~2, for example, specifies entailment through its models \cite{owl2012}. A proposed reasoning computation therefore requires both a representation of its states and an account of the semantics that makes its steps sound. Neither soundness nor completeness follows from the Token map condition alone.

For a concrete positive fragment, consider the rule
\begin{equation}\label{eq:knowledge-rule}
(\mathrm{type},x,c),\ (\mathrm{sub},c,d)
\quad\Longrightarrow\quad
(\mathrm{type},x,d),
\end{equation}
where $\mathrm{type},\mathrm{sub}\in E$ are fixed labels and $x,c,d\in B$. This records the familiar propagation of membership along a subclass relation, one of the RDFS entailment patterns \cite{rdf2014}. Given a heap $\Upsilon$, let $\operatorname{Cl}(\Upsilon)$ be the least heap containing it and closed under this declared rule. With the carrier fixed and no fresh nodes generated, the identity function defines
\begin{equation}\label{eq:knowledge-closure-map}
\langle B,E,\Upsilon\rangle
\xrightarrow{\ \id_B\ }
\langle B,E,\operatorname{Cl}(\Upsilon)\rangle,
\end{equation}
because $\Upsilon\subseteq\operatorname{Cl}(\Upsilon)$. For example, the two records $(\mathrm{type},a,\mathrm{Cat})$ and $(\mathrm{sub},\mathrm{Cat},\mathrm{Animal})$, with the three named nodes in $B$, generate $(\mathrm{type},a,\mathrm{Animal})$. For a finite initial fact set over finite $B$ and the displayed vocabulary, exhaustive addition of these facts terminates. This example concerns the declared positive rule, not a complete OWL reasoner.

The arrow in~\eqref{eq:knowledge-closure-map} expresses retention of existing facts as the object changes. Its carrier function is an identity, so it is not itself the algorithm that constructs the enlarged heap. An executable account makes the heap representation part of a machine state and specifies the function or sequence of transitions that computes the closure. This distinction connects the example to Thesis~\ref{thesis:executable}: structural records, descriptions of their transformations, and executions of those descriptions have related but different roles.

\subsection{Learned internal structure}\label{sec:learned-interior}

In a conventional symbolic knowledge system, computation is organized around declared entities, concepts, relations and rules. A language model also computes through internal representations whose organization is obtained through training and whose values depend on the input context. Those states need not correspond one by one to previously named concepts. This is the useful distinction behind the contrast between explicit knowledge and hidden-state computation. The words \emph{explicit} and \emph{hidden} describe how a representation is made available; they do not establish a hierarchy of semantic depth.

These are not disjoint classes of systems. Knowledge graph embeddings already learn internal representations of entities and relations. TransE, for example, models relations by translations between learned entity embeddings \cite{bordes2013}. Whether a relation has an explicit name, whether its representation is learned, and whether a state varies with context are separate modeling choices. Token Space can expose these choices without forcing every system onto one side of a symbolic--neural divide.

For a language model with fixed parameters $\theta$, inference evaluates a specified family of state transformations. Training selects or modifies those parameters; computing a new contextual state during inference is not, by itself, relearning the parameterized transformation. The state carriers must also contain the information required by the computation. A complete layer-state sequence under a fixed inference protocol can provide such a carrier. An isolated vector at one position generally does not determine the next state without the rest of the context. Located occurrences are appropriate for observations that depend on layer, position or input, as in Section~\ref{sec:interpretation} and Thesis~\ref{thesis:occurrence}.

Having a state function $F_\theta:H\to H'$ is only the start. To regard it as a Token map, choose classes $\mathcal H=\langle H,E,\Upsilon_H\rangle$ and $\mathcal H'=\langle H',E,\Upsilon_{H'}\rangle$ and establish
\[
F_E(F_\theta)(\Upsilon_H)\subseteq\Upsilon_{H'}.
\]
Prefix-agreement heaps make this condition a statement about information dependence. Heaps of specified readouts or interventions instead concern the corresponding experimental protocols. Preservation of the first kind of relation does not imply preservation of ontology consequences, and a successful readout of the second kind does not establish a uniquely privileged content. The heap specifies which claim is being made.

\subsection{Correspondences between explicit and internal computations}\label{sec:knowledge-correspondence}

Suppose a task specifies explicit state sets $S,S'$ and a computation $g:S\to S'$. For example, a state can encode a finite fact set and a query, and $g$ can compute a declared inference step or a query answer. Let $F_\theta:H\to H'$ describe the model computation on a domain containing the required context. Fix total readouts $\pi:H\to S$ and $\pi':H'\to S'$ on these declared domains. A candidate correspondence is the equation
\begin{equation}\label{eq:knowledge-correspondence}
\pi'\circ F_\theta=g\circ\pi.
\end{equation}

\begin{figure}[htbp]
\centering
\begin{tikzcd}[column sep=large,row sep=large]
H \arrow[r,"F_\theta"] \arrow[d,"\pi"'] & H' \arrow[d,"\pi'"] \\
S \arrow[r,"g"'] & S'
\end{tikzcd}
\caption{A specified correspondence between hidden-state computation and explicit computation. Commutation is a condition to prove or test on the declared domains. The diagram is in $\Set$ initially; it is a diagram in $\TS_E$ only after all four state sets are equipped with classes over the same core and each arrow preserves the specified heap.}\label{fig:explicit-hidden-square}
\end{figure}

Equation~\eqref{eq:knowledge-correspondence} says that computing internally and then reading out agrees with reading out first and performing the specified explicit computation. It is not assumed for arbitrary models or tasks. If the four functions are also Token maps between specified classes over a common core, the square is a commuting diagram in $\TS_E$. Partial readouts require a restricted domain on which all compositions are defined, and stochastic computations require a different account of arrows; neither is silently identified with the total deterministic square.

There is a concrete obstruction to this correspondence. If
\[
\pi(h_1)=\pi(h_2)
\quad\text{but}\quad
\pi'(F_\theta(h_1))\ne\pi'(F_\theta(h_2)),
\]
then no deterministic $g$ on those explicit states can make the square commute. Conversely, constancy of $\pi'\circ F_\theta$ on every fiber of $\pi$ defines a unique set function $g$ when $S=\pi(H)$. This is the factorization argument of Proposition~\ref{prop:observation-factorization} applied to another observation. It does not show that the induced function equals an independently prescribed reasoner or that it preserves the heaps chosen for $S$ and $S'$. Those remain additional conditions.

Declared correspondences also compose. For internal transitions $F_1:H_0\to H_1$ and $F_2:H_1\to H_2$, explicit transitions $g_1:S_0\to S_1$ and $g_2:S_1\to S_2$, and compatible readouts $\pi_i:H_i\to S_i$, the two local equations imply
\begin{equation}\label{eq:knowledge-correspondence-composition}
\pi_2F_2F_1=g_2\pi_1F_1=g_2g_1\pi_0.
\end{equation}
When all arrows are admissible, this is ordinary pasting of commuting squares in Token Space. Its role is to make the interface agreement needed for a composite claim explicit, not to assert that every neural layer performs a named logical inference.

\subsection{What makes a correspondence informative?}\label{sec:knowledge-evidence}

Every function is a Token map from an empty-heap class, and induced heaps can make preservation true by construction. An application therefore has to explain why its relations and observations matter. For comparisons of explicit and learned computation, the task, domains, readouts and preserved relations should be stated independently of the particular outcomes used to test the claim. In an empirical study, learned readouts can be fixed using designated training data and evaluated on separate cases. A correspondence obtained only by choosing a convenient description after observing every tested outcome has a different evidential status from a predeclared comparison.

A concrete study could use finite class hierarchies and membership queries. The explicit computation is fixed by rules such as~\eqref{eq:knowledge-rule}. The internal computation is a specified model segment or a complete inference procedure, with the necessary context retained. The study would compare readout results with the rule-based reference, examine failures of fiber constancy, and test whether conclusions persist across designated changes of names, contexts or hierarchy depth. Intervention protocols can additionally investigate whether the recorded information participates in producing the output. Agreement on a finite sample does not prove the universal equation, and decodability alone does not establish causal use.

Approximate agreement needs a metric or loss, an evaluation domain or distribution, and conditions for propagating error through composition. The exact Token map condition provides none of these automatically. Similarly, the commuting square establishes a particular computational correspondence, not uniqueness of interpretation or identity of internal algorithms. These distinctions retain the separation between records, evidence and content in Section~\ref{sec:interpretation}.

\subsection{Consequences for represented finite-mapping machines}\label{sec:knowledge-machines}

This connection supports the broader machine envisioned in Theses~\ref{thesis:interfaces} and~\ref{thesis:cores}. On declared finite interfaces, rules, tables, circuits and parameterized models can provide different executable representations of a mapping. An implementation might expose a named relation at one boundary and perform distributed internal computation at another. The relevant question is which interface obligations it satisfies and which transformations preserve its specified behavior. Neither finiteness alone nor the presence of learned parameters makes these representations equally compact or efficient.

Let $Q:X\to Y$ be an implementation, let $\alpha:X\to S$ and $\beta:Y\to S'$ be interface maps, and suppose $\beta Q=g\alpha$. If a second implementation $Q'$ on the same interfaces satisfies $Q'=Q$ on the declared input domain, then $\beta Q'=g\alpha$ as well. Equality of implementations is a sufficient condition here; it is stronger than merely agreeing after the output observation $\beta$. More general replacements can use the re-encoding conditions of Section~\ref{sec:finite-representations}. Output agreement alone does not guarantee preservation of independent events, resource cost or online execution state.

Under the hypotheses already established in Section~\ref{sec:elastic-execution}, valid pure acyclic descriptions permit independent ready work, and certified residual replacement preserves the required result. The symbolic or neural origin of a module does not remove those hypotheses. This gives a route for combining explicit knowledge procedures and learned modules while making their dependencies, representation choices and adaptation conditions available to the machine.

We do not require all hidden-state computation to be reconstructed as a knowledge graph. Explicit knowledge is one possible observation and interface among others. The framework leaves room for internal relations that have not been named or are not usefully expressed as entity--relation graphs. Its intended contribution is to make those selected interiors, the correspondences between computations, and changes of executable representation available in one mathematical setting. Their practical value must be demonstrated by relations that explain a success or failure, replacements whose correctness can be established, or execution choices that improve a measured resource objective.

\section{Scope}\label{sec:limits}

Several boundaries are part of the proposed account.
\begin{enumerate}
\item \emph{Exact relational conditions versus numerical performance.} Token maps preserve membership in heaps exactly. A numerical loss, a threshold error rate, or an approximate commuting diagram needs the quantitative extension discussed below; recording real numbers as labels still gives exact equality.
\item \emph{Modeling choices versus discovered structure.} Every function is a map from an empty-heap class or from a suitably induced class. Applications must say why the heap matters and what observations support it; causal and permutation heaps encode architectural constraints, algebraic tokenization encodes stipulated structure, and neither establishes semantic content or predictive ability.
\item \emph{Set-based description versus effective implementation.} Carriers and heaps may be infinite, and the heap of all pairs of equal-prefix sequences is a specification, not a data structure enumerated during inference. Exponential carriers include functions that no finite architecture realizes; parameter objects select realized families but do not solve realizability or optimization.
\item \emph{Concurrency versus resource guarantees.} Schedule independence is proved for valid acyclic graphs of pure total operations. Work--span and variable-capacity bounds assume unit-cost operations, free data access and adequate memory. The completed-boundary theorem justifies one form of online replacement; arbitrary asynchronous graph mutation, failures, recursion and resource-optimal adaptation require further protocols and proofs.
\item \emph{An unbounded population versus infinite simultaneous execution.} T5 admits arbitrarily large finite core populations through effective rules. Its finite-stage invariant, independent product family and union colimit establish specified parts of that position. They do not establish arbitrary asynchronous mutation, progress of endlessly growing workloads, infinite available hardware, or computation of a result with infinitely many unresolved prerequisites.
\item \emph{A fixed inference protocol versus every LLM system.} The decoder results assume specified masks, positional rules, and positionwise operations. Retrieval, persistent memory, cross-example operations, conditional routing, tool use, and stochastic policies need enlarged state spaces and different relations. Section~\ref{sec:learning-dynamics} supplies a state and execution description for specified training protocols, distinguishing changes to model predictions from correctness of implementing a training update.
\item \emph{Elementary classes versus tree Token classes.} The elementary $\TS_E$ is a quasitopos; the finite bicartesian closed constructions also hold in the bounded variant of Proposition~\ref{prop:bounded}. Shared-subterm reification is fully faithful on the supported substitution category into its explicitly extended-core Token Space. No such claim is made for arbitrary unused generators, or for all occurrence-graph maps. Limits and exponentials of the tree categories and nested extensions require separate analysis.
\end{enumerate}

The ambient quasitopos classification applies to the full elementary category with a fixed set of core symbols. It does not transfer merely by full faithfulness to valid computation graphs, tree categories or nested extensions. Designated-root exponentials represent structured higher-order families but may introduce records incompatible with executable graph axioms. An exact reachable-state migration, together with a common denotation, implies reachable output correctness; concurrency fidelity is an additional property. None of these certificates alone bounds recomputation cost or ensures progress under endless reorganization.

\section{Conclusion}\label{conclusion}

The framework developed here supplies a mathematical language for
studying internal structure, execution state and permitted changes of
representation. Its construction brings algebraic encoding, program
evaluation, concurrent execution and observations of neural computation
into a common setting. The distinctions it makes are part of its
contribution: sharing a description, deleting an unused address, and
merging behaviorally indistinguishable states ask different questions,
even when all three are called compression. The proved connections
show how these questions can be related without losing their stated
assumptions. This gives the language a role in formulating and connecting
research problems as well as in establishing particular results.

The Token Space is built by making structural facts explicit as finite records. Its fixed-core elementary category is complete, cocomplete and locally cartesian closed, with a strong-subobject classifier: it is a quasitopos. This classification applies general relational theory, with nullary facts retained explicitly. Its bijective noninvertible maps still show that it is not a topos. Exponential carriers contain all functions; global elements select the admissible computations and parameter objects organize their families. A bounded version keeps these mathematical objects finite, although their representations may still be large.

Algebraic tokenization remains the central device. For fixed finitary signatures and all homomorphisms, it gives a full embedding; componentwise products and induced equalizers are preserved under the stated closure hypotheses. A common core and the operator calculus make the choices of carrier, signature and heap visible and permit comparisons across encodings. The tree construction supplies nested syntax, and shared-subterm reification gives a full embedding under an explicit support condition. These constructions retain our motivating question: which internal relationships should be visible in a computation, and how can the machine manipulate them?

The proposal for represented finite-mapping machines distinguishes behaviors, implementations and certified changes of implementation. A valid graph evaluates through a unique Token map. Its completed sets form a finite distributive lattice realized fully faithfully over the graph. Completed pieces glue by bicartesian squares, ready batches determine Boolean intervals, and partial valuations agree along subclass inclusions. Sharing induces a completed-set adjunction and a value-reuse closure. Frontier interfaces and their advancements form a functor, while residual output maps give its compatible observation. Independent ready firings therefore commute. Explicit coding makes descriptions available as data, and a ranked module library supports denotation-preserving expansion. Primitive term compilation connects supported reification to this executable fragment without treating arbitrary tree records as programs. This gives intrinsic concurrency a precise scope: dependencies constrain execution while admissible schedules may vary. Worker supply can change under the explicit capacity model. A live frontier identifies a residual mapping and supports a sufficient rule for replacing unfinished computation without changing its result. These results connect structural descriptions to an execution invariant rather than leaving the machine at the level of input--output functions.

Computational elasticity includes worker allocation, active work, representation granularity and changing interfaces. The prefix-XOR family separates sharing from parallel reorganization, and its block descriptions supply a uniform choice of granularity with explicit work and span bounds. Neither common-subexpression sharing nor parallel prefix computation is claimed as a new algorithm. Their role is to make the framework's construction and correctness conditions concrete. The broader objective is to choose and transform implementations in response to resources and demand while preserving the specified behavior and guaranteeing progress. Token computing cores make the units of that organization explicit. Effective finite transitions preserve finite configurations, a uniform product construction supplies unbounded available concurrency across workloads, and compatible growing record classes assemble by a union colimit. These results connect the locally finite execution theory to the unbounded machine envisioned by T5; they do not assume that a finite physical installation realizes all of this width at once.

Transformers remain one implementation family among others. Their
decoder computations preserve prefix-agreement heaps under the stated
architectural assumptions. Coherent occurrence families connect this
fixed-length fact to prefix naturality. The full-prefix/cache invariant
and the prescribed-access deletion square then establish correctness
of two changes in retained state. Parameterised state exposes the
changing interfaces. The future-observation heap completes the semantic
connection: its structural quotient is the behavioral quotient, and an
exact incremental encoding must both preserve future observations and
remain stable under input extension. Its reachable state maps uniquely
onto the behavioral quotient. Direct-address liveness and behavioral
minimality answer different questions; neither by itself proves a
memory advantage for a chosen numerical representation. Structural
distillation can choose these observations explicitly, while approximate
and cost-sensitive compression require further hypotheses.

Explicit knowledge structures and learned internal states provide another use of the selected-interior perspective. Knowledge graph records, ontology computation and contextual state transformations can be compared through declared interfaces. The correspondence square of Section~\ref{sec:explicit-hidden} distinguishes agreement with a specified computation from merely encoding both systems. It also connects certified implementation changes to that agreement. Establishing useful correspondences in actual models, including their failures and costs, remains an empirical and formal task.

The causal results summarized from the mathematical companion explain which parts of a record presentation make a structural difference. Relative rootedness characterizes reconstruction from completed fragments. Fixing the result role by record shape supplies one cartesian closed category of ordered structures, with an explicit bridge from graph occurrences and supported reification. Exponentials need not remain executable graphs. Exact state migration is controlled by retained information, and the ready-cube formula measures how a graph map changes exposed independence. These distinctions turn interface and representation choices into obligations that a machine can check or certify, while compact evaluation and economical adaptation remain further tasks.

\paragraph{The five theses revisited.}
The five theses specify a single research direction. Thesis~\ref{thesis:interior} makes the selected interior available as records, with full faithfulness for the specified finitary structures and a calculus that acts on the records. Thesis~\ref{thesis:executable} makes categorical objects data in that same language: small categories and functors have full encodings, natural transformations have endpoint-constrained encodings, and finite products, slices, arrow categories and functor categories have terminating constructions. Thesis~\ref{thesis:interfaces} puts the relations to be preserved into computational interfaces, so different implementations can realize the same admissible mapping. Thesis~\ref{thesis:occurrence} retains the location of a language-model state and makes the forgetting of that location an explicit map, with a precise factorization test for observations. Thesis~\ref{thesis:cores} organizes these structures into a potentially unbounded population of logical execution units, with finite operational stages and explicit distinctions among identity, activity, readiness and hardware allocation. These are the ways in which the framework maintains an internal, operational, situated and extensible view of computation.

The broader ambition remains a machine that can compute with its own structured representations through a dynamically organized population of Token computing cores. Its internal descriptions may be compressed, expanded, shared, scheduled or replaced while specified behavior is preserved. The finite categorical constructions and the concurrent graph constructions establish parts of this ambition at different levels. Their combination suggests that categorical specifications, executable descriptions and runtime state should be available to one computational process, with each transformation carrying its own correctness conditions.

Several avenues remain open.
\begin{itemize}
\item Develop the categorical self-description programme beyond finite tables: specify effective presentations, algorithms for the desired constructions and size bounds at each level. Study natural transformations and their horizontal and vertical composition as structured data, and investigate a relation-lifting semantics that connects structural interface refinement to an abstraction theorem.
\item Extend the explicit dependent-product semantics of the elementary quasitopos to effective interface languages, specifying which families admit finite or compressed representations and effective evaluation. Determine which additional axioms and nested constructions preserve the required closure properties. The elementary classification is established in Theorem~\ref{v13:thm-quasitopos}; these effective and constrained extensions remain open.
\item Develop quantitative specifications for approximate re-encoding and distillation, with conditions for composing error bounds. Weighted records alone do not provide such a composition theorem.
\item Instantiate the explicit--hidden correspondence of Section~\ref{sec:knowledge-correspondence} on a finite ontology task with independently specified rules, domains and readouts. Compare rule-based and model computations, identify failures of fiber constancy, and distinguish readout agreement from intervention evidence. Investigate whether certified hybrid implementations offer a measured benefit in representation size, work or adaptation cost.
\item Characterize the concrete functors for which nested $\SX^{\&}$ constructions preserve Cartesian closedness, and investigate limits and exponentials of tree Token categories with their specified substitution classes.
\item Let $\mathcal C_{Cartesian}$ and $\mathcal C_{Topos}$ be the \emph{super operators} assigning to a subcategory $\mathbf P$ of a Token category $\mathbf{TC}$ the least bicartesian closed subcategory, resp.\ the least topos subcategory, of $\mathbf{TC}$ containing $\mathbf P$, where they exist. Given a category $\mathbf C$ of structured objects, how does $\Rep_{\mathbf C}(\mathbf C)$ relate to $\mathcal C_{Cartesian}(\Rep_{\mathbf C}(\mathbf C))$ and $\mathcal C_{Topos}(\Rep_{\mathbf C}(\mathbf C))$? Existence is part of the question: a subcategory containing a bijective noninvertible arrow of $\TS$ cannot be a topos, so $\mathcal C_{Topos}$ is defined at most for images that avoid such arrows, and external completions are the alternative.
\item Extend the ranked module construction to product-valued and multi-output bodies, and the representation library to algebraically justified rewrites, decision diagrams, tables and hybrid evaluators. Study partial orders of resource profiles in $\mathsf{Impl}(f)$ and effective selection policies that include transformation and certificate costs.
\item Implement the concurrent protocol of Section~\ref{sec:elastic-protocol} and the experiments of Section~\ref{sec:elastic-research}. Investigate memory-constrained progress, non-unit primitives, communication and locality, and conditions under which adaptation yields an amortized benefit.
\item Extend the proved reachable-state migration and ready-cube criteria beyond the pure, total acyclic graph language, including certified rewrites that are not graph Token maps, local or distributed adaptation and failures. Develop policies that account for information recovery, recomputation, synchronization and progress when allowing expansion and compression. Use the finite causal-stage criterion to study dynamically generated dependencies, admission control and distributed generation; compare concrete protocols with Actor and dataflow models.
\end{itemize}
The Token-level approach aims to connect a machine's structural specification to the representations it can actually execute and transform, including the unbounded horizon of its computational organization. The definitions, proofs and finite examples above provide a starting point for that investigation. The mathematical companion develops the universal constructions for dependency presentations, while the type-theoretic companion develops a separate programming interpretation; this paper retains the machine and AI programme that motivates both.

\appendix
\section{Proofs of the Scheduling and Prefix Bounds}\label{app:scheduling}

\begin{proof}[Proof of Proposition~\ref{prop:work-span}]
At most $p$ nodes complete in a slot, and the nodes of a dependency path require distinct slots, giving the lower bound. Call an unfinished slot underfull if it executes fewer than $p$ nodes; it then executes every ready node. In the remaining graph a longest path starts at a ready node, otherwise a predecessor would extend it, so removing all ready nodes decreases the remaining span by at least one. Hence there are at most $S$ underfull slots. Each executes at least one node and wastes at most $p-1$ worker-slots; full slots waste none. Thus $pT_p-W\le(p-1)S$.
\end{proof}

\begin{proof}[Proof of Proposition~\ref{prop:variable-capacity}]
While work remains, a positive-capacity underfull slot executes every ready node and decreases the remaining span, so there are at most $S$ such slots, each wasting at most $P-1$ worker-slots; zero-capacity slots waste nothing. If work remained after horizon $h$, at most $W-1$ operations would have completed, while the capacity accounting gives at least $C(h)-(P-1)S\ge W$, a contradiction. The zero-work case is immediate, and the output assertion is Theorem~\ref{thm:schedule-independence}.
\end{proof}

\begin{proof}[Proof of Proposition~\ref{prop:parallel-prefix}]
The recursive graph depends only on pair results, and final odd-position nodes depend only on recursive prefixes and original inputs, giving an acyclic graph with unique definitions. Associativity identifies $z_i$ with the prefix through $x_{2i}$; appending $x_{2i+1}$ yields the next odd prefix. There are $\lfloor n/2\rfloor$ pair operations and $\lfloor(n-1)/2\rfloor$ final operations, whose sum is $n-1$; the work bound follows by induction. At most one layer is added before and one after the recursive subgraph, proving the span recurrence and bound. No identity element or padding is required, so odd lengths are included.
\end{proof}

\begin{proof}[Proof of Proposition~\ref{prop:block-prefix}]
Associativity combines each preceding-block total with a local prefix to produce the corresponding full prefix, with the original operand order unchanged. Local chains require $\sum_j(\ell_j-1)=n-b$ operations and span at most $\ell-1$. The totals network uses $w_b\le2b-2$ operations and span at most $2\lceil\log_2 b\rceil$. There are $n-\ell_1$ final combinations, executable in one further layer. Summing yields the work formula and span bound; since $b\le n$ and $\ell_1\ge1$, work is at most $3n-3$. Denotational equality to $P_n$ implies preservation of its specified relations.
\end{proof}

\section{Reproducible Finite Checks}\label{app:finite-checks}

The validation script, its recorded results and a README with reproduction instructions are supplied in the accompanying ancillary directory. The files are named \texttt{validate\_V16.py}, \texttt{validation\_results\_V16.json} and \texttt{README.txt}. The script supports the complete source bundle and an explicitly specified manuscript path when files are downloaded separately.

The accompanying standard-library Python script \texttt{validate\_V16.py} implements the recursive and block prefix graph constructors, ready-node schedules and residual graphs used in Section~\ref{sec:elastic-execution}. It also constructs the syntactic sharing map for separate prefix expressions. These checks supplement the proofs by testing their finite constructions and boundary cases; they are not machine-checked proofs, complexity lower bounds or hardware performance measurements. The following scope makes the reported checks reproducible.

\begin{itemize}
\item For every binary input of lengths $1$ through $10$, evaluate the recursive prefix graph and every block choice $1\le b\le n$, comparing each result with sequential prefix XOR. Verify their operation counts and span bounds. This gives $20{,}480$ graph--input comparisons.
\item Repeat the prefix comparisons for lengths $1$ through $6$ over the four-element monoid of functions from $\{0,1\}$ to itself, with function composition as the operation. This gives $36{,}408$ graph--input comparisons. Composition is associative but not commutative, so this test detects an incorrect reversal of operands that XOR would conceal.
\item Generate $32$ acyclic graphs with three input nodes, six operation nodes and three output ports, using seed $5001$. The primitive library consists of Boolean XOR, AND, NOT, and the two nullary Boolean constants. Enumerate predecessor-closed completed sets and all eight inputs. Check partial valuations, readiness, commutation of each ready pair and residual output equality. Output aliases, completed outputs, nullary operations and unused work are permitted.
\item For the same graphs and all block-prefix graphs with $2\le n\le14$, check eight fixed worker counts and sixteen seeded variable-capacity schedules per graph. The checks cover $1{,}088$ fixed schedules and $2{,}176$ variable schedules, including zero-capacity slots, against the corresponding inequalities. These are idealized schedules under the paper's cost model.
\item Enumerate completed sets and inputs for the occurrence-prefix sharing maps at lengths $2$ through $5$, checking the image completed set, consistency on identified nodes, and transferred valuations. Also check all eight residual inputs and sixteen original inputs in the four-input mid-execution parity example of Section~\ref{sec:online-replacement}.
\end{itemize}

Further checks for the module and lattice constructions are described below. The script also verifies the direction of~\eqref{eq:sequence-inclusion} and distinguishes causality from equivariance on two-position Boolean sequences.

\begin{itemize}
\item Enumerate all $512$ directed acyclic graphs on five topologically ordered nodes with the first two designated as inputs. Check closure under union and intersection for their input-containing down-sets, their cover relation, and the ready-set Boolean intervals. These are structural order checks, independent of a choice of numerical primitives.
\item Build a ranked Boolean module library containing argument projection with an unused formal input, a nullary constant, balanced parity, nested calls and repeated arguments. Use seed $7001$ to generate $16$ graphs at each input arity $0$ through $4$, with eight operations and four output ports (including a repeated output); add one direct invocation of each of the four modules. The resulting $84$ graphs give $521$ full-expansion input comparisons and $1{,}824$ selected-call comparisons. Preserve output aliases and reject both a forward reference and a recursive reference.
\item Check packing and expansion of primitive graphs, including bodies that return an input; and compile primitive terms by shared-subterm reification, including variable roots, nullary roots, repeated subterms and unused declared inputs. Packing uses $35$ primitive graphs and $277$ input comparisons. Term compilation uses $123$ terms and $761$ input comparisons, with all inputs enumerated for each term. Compare the resulting graph values with direct evaluation.
\item In the same $32$ primitive graphs used for the partial-state checks, enumerate all nonempty operation-node subsets and candidate exits. The fusion conditions admit $241$ regions. Check all eight original inputs for each region, giving $1{,}928$ comparisons of direct, fused and re-expanded outputs. Also reject the unique-exit region with a cycle after contraction described in Section~\ref{sec:module-expansion}.
\end{itemize}

The diagrammatic results have additional finite checks. For the seeded primitive graphs, the script verifies the carrier and heap union/intersection formulas on completed sets and constructs advancement graphs. It checks the iterated residual law, identity and composition of frontier maps, and residual naturality on all boundary tuples within the stated finite sample. A separate small span enumerates arbitrary two-element target heaps and compatible maps to check unique gluing. For syntactic sharing of separate prefix expressions at lengths $2$ through $5$, it checks both directions of the completed-set adjunction, saturation and reused values, including a nonsurjective graph map and the example where direct image fails to preserve an intersection. The Boolean two-position causal square is checked over all $256$ functions. The gluing identities cover $38{,}703$ pairs of completed sets in all $32$ seed-$5001$ graphs; the arbitrary-target test covers $64$ heaps and $96$ compatible map pairs. The frontier sample contains the first eight of those graphs and four boundary cases: an empty graph, a nullary constant with aliased outputs, an input-only graph with unused inputs, and a graph with no outputs but unfinished dead work. It gives $26{,}274$ boundary-input checks of residual naturality, $123{,}283$ of transition composition, and $2{,}747$ structural iterated-residual checks. The six graph maps in the sharing sample give $742$ adjunction comparisons and $4{,}296$ input checks of saturation and value reuse. The companion script and result file retain these counts together with the earlier checks.

The categorical self-description and located-state constructions have separate checks. Seven finite categories are used: the empty category, the terminal category, a discrete two-object category, the walking arrow, the one-object cyclic group of order two, the one-object idempotent monoid of order two, and a category with two parallel nonidentity arrows. Over all $49$ ordered pairs, enumerate $83{,}737$ functions on the tagged carriers and check that the $85$ record-preserving maps are exactly the functors obtained by independently checking the category tables. Construct and check the category laws for $49$ products, $49$ functor categories, seven arrow categories and nine slices. Independently enumerate functors into the arrow categories and compare them with the $195$ natural-transformation lifts. The $49$ functor-category evaluation maps are checked as both functors and Token maps. Empty-source and empty-target cases are included, as is the two-versus-four carrier counterexample for ambient exponentials.

For located observations, use four distinct occurrences with two coordinate fibers and enumerate all $81$ observations into a three-element set. Exactly nine factor through the coordinate projection, matching an independent enumeration of the vector-only functions. Over all $4{,}096$ pairs of heaps on a two-element carrier with one unary and one binary relation symbol, check that the coordinate identity is a Token map exactly for heap inclusion, and an isomorphism exactly for heap equality. These checks supplement Propositions~\ref{prop:cat-tokenization}, \ref{prop:finite-cat-constructions} and~\ref{prop:observation-factorization}. The manuscript check reads \texttt{TokenSpace\_V16.tex}, requires it to exist, and checks label uniqueness, references and bibliography keys.

The T5 extension checks all Boolean inputs of lengths $0$ through $8$ for the independent negation-core family, including empty input and repeated values. It also runs finite interleavings of generation, completion and reclamation, checking fresh identities, readiness, preservation of values still needed by waiting consumers, and ordered output correctness against a separately evaluated dependency graph. A serial stream is included to distinguish cumulative identities from available width. These finite checks exercise the operational distinctions; the unbounded assertion comes from the uniform constructor and proof of Proposition~\ref{prop:unbounded-core-width}.

All listed checks pass. The companion \texttt{validation\_results\_V16.json} records their counts. Tests over these finite ranges do not establish claims for arbitrary graph sizes; the propositions and their stated assumptions provide those arguments. The experimental programme in Section~\ref{sec:elastic-research}, including measured adaptation overhead and distributed execution, remains future work.

\subsection*{Migration checks and companion structural checks}
The additional standard-library script \texttt{check\_residual\_migrations.py} checks $961$ finite output-fiber cases and $11{,}141$ kernel-factorization cases. Its $128$ syntactic-sharing samples, generated with seed $20261002$, include $126$ noninjective maps and check $42{,}604$ ready cubes with $196{,}697$ vertices, together with two information/concurrency obstruction examples. The graph sample checks completion and concurrency structure; the general natural migration of reachable frontier diagrams is established by the proof. The script and its recorded results accompany this paper.

The checks of general causal gluing, designated-root exponentials and dependent products accompany the mathematical companion \cite{pancausal}, where their full proofs now appear. They are not counted as additional runtime or machine experiments in this manuscript. Both bundles specify finite test ranges and seeds. The tiny decoder checks below test declared mathematical models; no trained-model benchmark or hardware measurement is claimed.

\subsection*{Prefix and cache semantics}\label{llm:checks}

The accompanying script \texttt{check\_llm\_prefix\_cache.py} uses two
independent implementations of the same scalar decoder: one computes all
positions together layer by layer, while the other computes one new
position using stored keys and values. They do not share an attention
or layer evaluator. This checks the distinction relevant to
Theorem~\ref{llm:cache-simulation}, rather than comparing one routine
against itself. It is not a formal proof or a test of a trained model.

There are two layers, one head and scalar hidden states. At position
$i\geq1$ the embedding is $x_i+i/10$. With layer index $\ell=0,1$, the
projections are
\[
 q=((\ell+1)h+1)/3,\qquad k=(2h-\ell)/5,\qquad v=(h+\ell+1)/4.
\]
If $a$ is the attention result, the update is $r=h+a$ followed by
$h'=r+r^2/13$. The exact rational version uses positive weights
$w(q,k)=1+(qk)^2$, normalized over allowed positions. It is an instance
of the general recurrence, explicitly not softmax. The floating version
uses stable softmax with score $qk/2$.

\begin{table}[htbp]
\centering\small
\begin{tabularx}{\textwidth}{@{}Xr@{}}
\toprule
Finite check & Recorded result\\
\midrule
Exact rational input sequences, lengths $0$--$4$ & $121$\\
Prefix squares & $547$\\
Exact prefix-output equalities & $1{,}002$\\
Exact incremental-output equalities & $1{,}002$\\
Exact cached key/value equalities & $4{,}008$\\
Floating softmax input sequences & $126$\\
Maximum absolute incremental-output discrepancy & $8.882\times10^{-16}$\\
Maximum absolute cache discrepancy & $2.220\times10^{-16}$\\
Floating absolute tolerance & $10^{-12}$\\
Explicit QKV graphs / checked prefix cuts & $4\,/\,14$\\
\bottomrule
\end{tabularx}
\caption{Finite checks for the declared scalar model. The exact tests
use rational arithmetic; the softmax tests use floating arithmetic.}
\label{llm:checks-table}
\end{table}

The five additional floating sequences beyond the rational-domain
catalogue are listed in the script. Every cached entry is compared with
the corresponding full-prefix layer value. On the explicit QKV graphs,
each position prefix is lower and every old occurrence read directly
by a later position is a key or value. This is checked for the declared
graph granularity; it does not automatically cover every possible
primitive expansion or compiler representation.

Three boundary witnesses are evaluated. Bidirectional averaging changes
the first output from $0$ on $(0)$ to $1$ on $(0,2)$. The length-dependent
pointwise rule $x_i+n$ changes the first output from $1$ to $2$ on those
same inputs. Finally, deleting all but two cached positions in the
full-attention floating model changes its output on $(1,0,1,0)$, first at
position three, by approximately $0.102321$. This last test compares
cache deletion against the original full-attention semantics, not
against a model consistently defined to use a local attention window.

The JSON output records all counts, errors, assumptions and graph cuts.
These finite checks can detect implementation mistakes and make the
examples reproducible. They do not establish an error bound for arbitrary
sequence lengths, certify the general proofs, or measure model quality
or resource performance.

\subsection{Checks for eviction, indexed state, and behavioral compression}
\label{obs:checks-state-compression}

Three additional dependency-free Python scripts accompany the new state
results. Each writes a JSON report with its own filename stem:
\nolinkurl{check_safe_eviction.py},
\nolinkurl{check_parameterised_state.py}, and
\nolinkurl{check_observation_quotient.py}.
All reported cases pass. Table~\ref{obs:state-check-counts} records cases,
not claims about the number of independent proofs or test assertions.

\begin{table}[htbp]
\centering\small
\begin{tabularx}{\linewidth}{@{}Xr@{}}
\toprule
Check family & Cases\\
\midrule
All causal access families through horizon $5$ & $33{,}867$\\
Future-read availability, over every retained subset and cut
 & $2{,}097{,}151$\\
Boundary update inclusions & $168{,}146$\\
Sliding-window cuts & $5{,}544$\\
Arbitrary-array eviction factorizations (seeded samples) & $2{,}500$\\
Exact two-layer inputs / evaluated prefix cuts
 & $32{,}768\,/\,98{,}304$\\
\midrule
Indexed function-space cardinalities & $27$\\
Indexed value-functor laws / unit laws & $607\,/\,644$\\
Indexed composition associativity & $2{,}149$\\
Strength naturality and coherence laws & $3{,}723$\\
\midrule
All binary Mealy machines with two / three labelled states
 & $256\,/\,46{,}656$\\
State partitions of these machines & $233{,}792$\\
Exact update quotients and their reachable maps to $Q$ & $60{,}692$\\
Observationally exact but update-unstable partitions & $6{,}048$\\
\bottomrule
\end{tabularx}
\caption{Finite validation of the state results. Most families are
exhaustive over the stated finite domains. The arbitrary-array
factorization family consists of reproducible random samples.}
\label{obs:state-check-counts}
\end{table}

The eviction script compares the required-set formula with an operational
oracle that processes future reads without using that formula. It
enumerates every causal access family for horizons $0$ through $5$ and
every retained subset at every cut. Window widths $1$ through $24$ are
checked at all cuts of horizons $0$ through $20$, including terminal
cuts and gaps before a later read. A separate integer-arithmetic decoder
uses total finite sums, including empty access rows: a layer-major
full-prefix evaluator is compared with a token-major pruned-cache
evaluator for all pairs of layer access families at horizon $3$ and
all inputs in $\{-1,1\}^3$. These are exact tests of the declared
recurrence, not softmax experiments. The $2{,}500$ additional cases use
seed $20261002$, horizons $1$ through $7$, one through three layers,
and independently assigned integer cache entries; no reachability
assumption is imposed on these arrays. This sampled check is
distinguished from the exhaustive access-set and decoder families.

The indexed-state script enumerates total function tables. State and
value cardinalities $0$, $1$, and $2$ test cardinality formulas, value
functor laws, units, and the specified strength identities.
Associativity is exhaustively checked over nine explicitly listed
patterns of four state cardinalities and three value cardinalities,
including different pre-, intermediate-, and post-state sizes; the
JSON report lists the patterns. Its oracle expands the three
uncurried transitions directly. Empty-domain functions and cases
with no admissible function triples are reported separately, so a
vacuous equality is not presented as evidence for an inhabited
execution. This is exhaustive enumeration within the chosen patterns,
not enumeration of arbitrary state sizes.

The observation script exhausts every binary-input, binary-output
Mealy transition table on two and three labelled states and every
partition of each carrier. A greatest-fixed-point computation of
behavioral equivalence is compared with explicit observations of all
continuation words through length equal to the number of states.
It checks membership reflection for the unary observation records,
descent of the prescribed transition to each proposed quotient,
and the reachable factorization to the behavioral quotient.
The $6{,}048$ unstable cases concretely separate observational
sufficiency from the update condition of
Theorem~\ref{obs:observation-compression}; these are expected
counterexamples to the weaker claim, not test failures. Additional
examples check the effect of adding equality records and the
history-encoding counterexample of
Example~\ref{obs:insufficient-kernel}.

These finite checks supplement the proofs. They establish neither
formal proof certification nor validation on a trained LLM, and do
not yield a floating-point error guarantee, an effective algorithm
for arbitrary neural-state minimization, or a lower bound on
semantic memory outside the stated computational model.

\section{A Generic Interpretation Test for Prefix Naturality}
\label{llm:generic-prefix}

Prefix naturality is a behavioral property of interpreted diagrams. Its
relationship to syntactic dependencies can be stated exactly when the
interpretation is allowed to vary freely. The following is a standard
free-algebra argument, used here to delimit what a structural
classification can recover.

Fix a finitary, single-sorted operation signature $\Omega$ and a finite
set $P$ of parameter names. For each $n\le N$ let $Q_n$ be a valid finite
DAG presentation with input names $x_1,\ldots,x_n$, parameter names $P$
and $n$ ordered output ports. The operation signature and parameter names
are common to all lengths, but no embeddings or compatible node names
between $Q_m$ and $Q_n$ are assumed. Unfold the DAG at its $i$th output to
obtain a finite term
\[
 t_{n,i}\in\operatorname{Term}_{\Omega}
                 (P\sqcup\{x_1,\ldots,x_n\}).
\]
Sharing affects the DAG but not this unfolded term. An $\Omega$-algebra
$A$, together with $\theta:P\to|A|$, interprets $Q_n$ as
$F_n^{A,\theta}:|A|^n\to|A|^n$.

\begin{proposition}[Universal-interpretation criterion]
\label{llm:free-prefix-test}
The following are equivalent.
\begin{enumerate}
\item For every $\Omega$-algebra $A$ and parameter assignment $\theta$,
      the family $(F_n^{A,\theta})_n$ is prefix-natural.
\item For every $m\le n$ and $i\le m$, $t_{n,i}=t_{m,i}$ as formal
      terms, regarding the latter term as a term over the larger
      variable set.
\end{enumerate}
If these conditions hold, no input $x_j$ with $j>i$ has a directed
operand-to-result path to the $i$th output of $Q_n$.
\end{proposition}
\begin{proof}
Equality of formal terms is preserved by every interpretation and
substitution, proving the second condition implies the first. Conversely,
fix $m\le n$ and use the free $\Omega$-algebra on
$P\sqcup\{x_1,\ldots,x_n\}$. Assign each parameter and input its own
generator. The $i$th output values are then exactly $t_{n,i}$ and
$t_{m,i}$; naturality forces their equality. In a freely unfolded term,
an input variable occurs if and only if there is an operand path from
its input occurrence to the output: induction along the DAG proves this
by taking the union of operand variable sets at each constructor.
Equality with $t_{i,i}$ therefore excludes all later variables and their
paths.
\end{proof}

This criterion is stronger than testing all numerical parameter values
of a fixed neural architecture: the latter still fixes the meanings of
addition, attention and other primitives, whose equations may hide a
syntactic dependency. Nor does the criterion recover a unique occurrence
presentation. A shared subexpression and two separately recomputed copies
can unfold to the same term. Thus the three questions---equality of
numerical behavior, generic algebraic dependency, and identity and sharing
of occurrences---remain distinct classification questions.


\begin{thebibliography}{99}

\bibitem{lawvere1997}
F. W. Lawvere and S. H. Schanuel, \emph{Conceptual Mathematics: A First Introduction to Categories}, Cambridge University Press, 1997.

\bibitem{awodey2006}
S. Awodey, \emph{Category Theory}, Oxford University Press, 2006.

\bibitem{maclane1998}
S. Mac Lane, \emph{Categories for the Working Mathematician}, 2nd ed., Springer, 1998.

\bibitem{johnstone2002}
P. T. Johnstone, \emph{Sketches of an Elephant: A Topos Theory Compendium}, Vols. 1--2, Oxford University Press, 2002.

\bibitem{goldblatt1984}
R. Goldblatt, \emph{Topoi: The Categorial Analysis of Logic}, revised ed., North-Holland, 1984.

\bibitem{lambek1986}
J. Lambek and P. J. Scott, \emph{Introduction to Higher Order Categorical Logic}, Cambridge University Press, 1986.

\bibitem{hodges1993}
W. Hodges, \emph{Model Theory}, Cambridge University Press, 1993.

\bibitem{lawvere1963}
F. W. Lawvere, Functorial semantics of algebraic theories, \emph{Proceedings of the National Academy of Sciences} 50 (1963), 869--872; reprinted in \emph{Reprints in Theory and Applications of Categories} 5, 2004.

\bibitem{adamek1990}
J. Ad\'amek, H. Herrlich, and G. E. Strecker, \emph{Abstract and Concrete Categories: The Joy of Cats}, Wiley, 1990; reprinted in \emph{Reprints in Theory and Applications of Categories} 17, 2006.

\bibitem{burris1981}
S. Burris and H. P. Sankappanavar, \emph{A Course in Universal Algebra}, Springer, 1981.

\bibitem{blok1989}
W. J. Blok and D. Pigozzi, \emph{Algebraizable Logics}, Memoirs of the American Mathematical Society 396, 1989.

\bibitem{sangiorgi2011}
D. Sangiorgi, \emph{Introduction to Bisimulation and Coinduction}, Cambridge University Press, 2011.

\bibitem{comon2007}
H. Comon, M. Dauchet, R. Gilleron, F. Jacquemard, D. Lugiez, C. L\"oding, S. Tison, and M. Tommasi, \emph{Tree Automata Techniques and Applications}, 2007. \url{https://inria.hal.science/hal-03367725}.

\bibitem{riley2018}
M. Riley, Categories of optics, arXiv:1809.00738, 2018.

\bibitem{fong2019}
B. Fong, D. Spivak, and R. Tuy\'eras, Backprop as functor: A compositional perspective on supervised learning, \emph{LICS}, 2019.

\bibitem{cruttwell2022}
G. S. H. Cruttwell, B. Gavranovi\'c, N. Ghani, P. Wilson, and F. Zanasi, Categorical foundations of gradient-based learning, \emph{ESOP}, LNCS 13240, 2022.

\bibitem{shiebler2021}
D. Shiebler, B. Gavranovi\'c, and P. Wilson, Category theory in machine learning, arXiv:2106.07032, 2021.

\bibitem{gavranovic2024}
B. Gavranovi\'c, P. Lessard, A. Dudzik, T. von Glehn, J. G. M. Ara\'ujo, and P. Veli\v{c}kovi\'c, Position: Categorical deep learning is an algebraic theory of all architectures, \emph{ICML}, PMLR 235, 2024.

\bibitem{sennrich2016}
R. Sennrich, B. Haddow, and A. Birch, Neural machine translation of rare words with subword units, \emph{ACL}, 2016, 1715--1725. \url{https://aclanthology.org/P16-1162/}.

\bibitem{kudo2018}
T. Kudo and J. Richardson, SentencePiece: A simple and language independent subword tokenizer and detokenizer for neural text processing, \emph{EMNLP: System Demonstrations}, 2018, 66--71. \url{https://aclanthology.org/D18-2012/}.

\bibitem{vaswani2017}
A. Vaswani et al., Attention is all you need, \emph{NeurIPS} 30, 2017. \url{https://arxiv.org/abs/1706.03762}.

\bibitem{ouyang2022}
L. Ouyang et al., Training language models to follow instructions with human feedback, arXiv:2203.02155, 2022. \url{https://arxiv.org/abs/2203.02155}.

\bibitem{llama2024}
A. Grattafiori et al., The Llama 3 herd of models, arXiv:2407.21783, 2024. \url{https://arxiv.org/abs/2407.21783}.

\bibitem{zaheer2017}
M. Zaheer et al., Deep sets, \emph{NeurIPS} 30, 2017.

\bibitem{bronstein2021}
M. M. Bronstein, J. Bruna, T. Cohen, and P. Veli\v{c}kovi\'c, Geometric deep learning: Grids, groups, graphs, geodesics, and gauges, arXiv:2104.13478, 2021.

\bibitem{touvron2023}
H. Touvron et al., LLaMA: Open and efficient foundation language models, arXiv:2302.13971, 2023. \url{https://arxiv.org/abs/2302.13971}.

\bibitem{su2021}
J. Su, Y. Lu, S. Pan, A. Murtadha, B. Wen, and Y. Liu, RoFormer: Enhanced Transformer with rotary position embedding, arXiv:2104.09864, 2021. \url{https://arxiv.org/abs/2104.09864}.

\bibitem{zhang2019}
B. Zhang and R. Sennrich, Root mean square layer normalization, \emph{NeurIPS} 32, 2019. \url{https://arxiv.org/abs/1910.07467}.

\bibitem{shazeer2020}
N. Shazeer, GLU variants improve Transformer, arXiv:2002.05202, 2020. \url{https://arxiv.org/abs/2002.05202}.

\bibitem{hinton2015}
G. Hinton, O. Vinyals, and J. Dean, Distilling the knowledge in a neural network, arXiv:1503.02531, 2015.

\bibitem{gou2021}
J. Gou, B. Yu, S. J. Maybank, and D. Tao, Knowledge distillation: A survey, \emph{International Journal of Computer Vision} 129, 2021, 1789--1819.

\bibitem{park2019}
W. Park, D. Kim, Y. Lu, and M. Cho, Relational knowledge distillation, \emph{CVPR}, 2019, 3967--3976.

\bibitem{tian2020}
Y. Tian, D. Krishnan, and P. Isola, Contrastive representation distillation, \emph{ICLR}, 2020.

\bibitem{jiao2020}
X. Jiao et al., TinyBERT: Distilling BERT for natural language understanding, \emph{Findings of EMNLP}, 2020, 4163--4174.

\bibitem{cao2022}
R. Cao, Putting representations to use, \emph{Synthese} 200, article 151, 2022. \url{https://doi.org/10.1007/s11229-022-03522-3}.

\bibitem{bryant1986}
R. E. Bryant, Graph-based algorithms for Boolean function manipulation, \emph{IEEE Transactions on Computers} C-35(8), 1986, 677--691. \url{https://doi.org/10.1109/TC.1986.1676819}.


\bibitem{darwiche2002}
A. Darwiche and P. Marquis, A knowledge compilation map, \emph{Journal of Artificial Intelligence Research} 17 (2002), 229--264. \url{https://arxiv.org/abs/1106.1819}.

\bibitem{pavlovic2013}
D. Pavlovic, Monoidal computer I: Basic computability by string diagrams, \emph{Information and Computation} 226 (2013), 94--116. \url{https://arxiv.org/abs/1208.5205}.

\bibitem{leroy2009}
X. Leroy, A formally verified compiler back-end, \emph{Journal of Automated Reasoning} 43 (2009), 363--446. \url{https://arxiv.org/abs/0902.2137}.

\bibitem{willsey2021}
M. Willsey, C. Nandi, Y. R. Wang, O. Flatt, Z. Tatlock, and P. Panchekha, egg: Fast and extensible equality saturation, \emph{Proceedings of the ACM on Programming Languages} 5 (POPL), article 23, 2021. \url{https://arxiv.org/abs/2004.03082}.


\bibitem{blelloch1993}
G. E. Blelloch, Prefix sums and their applications, in J. H. Reif (ed.), \emph{Synthesis of Parallel Algorithms}, Morgan Kaufmann, 1993, 35--60. \url{https://www.cs.cmu.edu/~guyb/papers/Ble93.pdf}.

\bibitem{blelloch1999}
G. E. Blelloch, P. B. Gibbons, and Y. Matias, Provably efficient scheduling for languages with fine-grained parallelism, \emph{Journal of the ACM} 46(2), 1999, 281--321. \url{https://doi.org/10.1145/301970.301974}.

\bibitem{josipovic2017}
L. Josipovi\'c, P. Brisk, and P. Ienne, From C to elastic circuits, \emph{51st Asilomar Conference on Signals, Systems, and Computers}, 2017, 121--125.

\bibitem{josipovic2018}
L. Josipovi\'c, R. Ghosal, and P. Ienne, Dynamically scheduled high-level synthesis, \emph{ACM/SIGDA International Symposium on Field-Programmable Gate Arrays}, 2018, 127--136.

\bibitem{hewitt1973}
C. Hewitt, P. Bishop, and R. Steiger, A universal modular ACTOR formalism for artificial intelligence, \emph{Proceedings of the Third International Joint Conference on Artificial Intelligence}, 1973, 235--245. \url{https://www.ijcai.org/Proceedings/73/Papers/027B.pdf}.

\bibitem{winskel1986}
G. Winskel, \emph{Category Theory and Models for Parallel Computation}, University of Cambridge Computer Laboratory Technical Report UCAM-CL-TR-85, 1986. \url{https://doi.org/10.48456/tr-85}.

\bibitem{winskel1987}
G. Winskel, Event structures, in W. Brauer, W. Reisig, and G. Rozenberg (eds.), \emph{Petri Nets: Applications and Relationships to Other Models of Concurrency}, LNCS 255, Springer, 1987. \url{https://doi.org/10.1007/3-540-17906-2_31}.

\bibitem{lack2003}
S. Lack and P. Soboci\'nski, Adhesive categories, \emph{BRICS Report Series} 10(31), 2003. \url{https://doi.org/10.7146/brics.v10i31.21799}.

\bibitem{rydeheard1988}
D. E. Rydeheard and R. M. Burstall, \emph{Computational Category Theory}, Prentice Hall, 1988.

\bibitem{spivak2012}
D. I. Spivak, Functorial data migration, \emph{Information and Computation} 217 (2012), 31--51. \url{https://arxiv.org/abs/1009.1166}.

\bibitem{reynolds1983}
J. C. Reynolds, Types, abstraction and parametric polymorphism, in R. E. A. Mason (ed.), \emph{Information Processing 83}, North-Holland, 1983, 513--523.

\bibitem{freeman1991}
T. Freeman and F. Pfenning, Refinement types for ML, \emph{PLDI}, 1991, 268--277. \url{https://doi.org/10.1145/113445.113468}.


\bibitem{belrose2023}
N. Belrose et al., Eliciting latent predictions from Transformers with the tuned lens, 2023. \url{https://arxiv.org/abs/2303.08112}.


\bibitem{owl2012}
B. Motik, P. F. Patel-Schneider, and B. Cuenca Grau (eds.), \emph{OWL 2 Web Ontology Language: Direct Semantics}, 2nd ed., W3C Recommendation, 11 December 2012. \url{https://www.w3.org/TR/owl2-direct-semantics/}.

\bibitem{rdf2014}
W3C, \emph{RDF 1.1 Semantics}, W3C Recommendation, 25 February 2014. \url{https://www.w3.org/TR/rdf11-mt/}.

\bibitem{bordes2013}
A. Bordes, N. Usunier, A. Garcia-Dur\'an, J. Weston, and O. Yakhnenko, Translating embeddings for modeling multi-relational data, \emph{Advances in Neural Information Processing Systems} 26 (2013), 2787--2795. \url{https://proceedings.neurips.cc/paper_files/paper/2013/hash/1cecc7a77928ca8133fa24680a88d2f9-Abstract.html}.

\bibitem{parker2022}
J. Parker, Exponentiability in categories of relational structures, \emph{Theory and Applications of Categories} 39(16) (2023), 493--518. Preprint 2022. \url{https://arxiv.org/abs/2208.07350}.

\bibitem{garner2012}
R. Garner and S. Lack, On the axioms for adhesive and quasiadhesive categories, \emph{Theory and Applications of Categories} 27(3) (2012), 27--46. \url{https://arxiv.org/abs/1108.2934}.

\bibitem{pancausal}
W. Pan, \emph{Causal Token Structures: Universal Gluing, Designated Roots, and Function Objects for Dependency Presentations}, working manuscript, version 3, structural revision, October 2026. A self-contained mathematical development of the dependency-presentation results. Local working versions are distinct from public arXiv version numbers.

\bibitem{pantypes}
W. Pan, \emph{Types as Token Classes: Data, Codata and Constraint Satisfaction in Token Space}, working manuscript, version 3, October 2026.

\bibitem{spiwack2018}
A.~Spiwack, \emph{Circuits via topoi},
arXiv:1807.07159 (2018).
\url{https://arxiv.org/abs/1807.07159}.

\bibitem{shazeer2019}
N.~Shazeer, \emph{Fast Transformer Decoding: One Write-Head is All You Need},
arXiv:1911.02150 (2019).
\url{https://arxiv.org/abs/1911.02150}.

\bibitem{moggi1991}
E.~Moggi, \emph{Notions of computation and monads},
Information and Computation \textbf{93}(1) (1991), 55--92.
\url{https://doi.org/10.1016/0890-5401(91)90052-4}.

\bibitem{rutten2000}
J. J. M. M. Rutten.
Universal coalgebra: a theory of systems.
\emph{Theoretical Computer Science}, 249(1):3--80, 2000.
\url{https://ir.cwi.nl/pub/48/}.

\bibitem{atkey2009}
R. Atkey, Parameterised notions of computation,
\emph{Journal of Functional Programming} 19(3--4) (2009), 335--376.
\url{https://doi.org/10.1017/S095679680900728X}.

\end{thebibliography}
\end{document}